%% file: 21-1091.tex
\documentclass[twoside,11pt]{article}

\usepackage{jmlr2e}
\usepackage{amsmath,bbm}
 \usepackage{algorithm,amssymb}
\usepackage{algpseudocode}
\usepackage[dvipsnames]{xcolor}
 \usepackage{multirow}
 \usepackage{cleveref}
\usepackage{graphicx}
\usepackage{subcaption}
\usepackage{float}

\algnewcommand\algorithmicinput{\textbf{INPUT:}}
\algnewcommand\INPUT{\item[\algorithmicinput]}
\algnewcommand\algorithmicoutput{\textbf{OUTPUT:}}
\algnewcommand\OUTPUT{\item[\algorithmicoutput]}

\newcommand{\mk}{\mathbb K}

\newcommand{\lt}{{\mathcal L^2}}
\newcommand{\me}{\mathbf e}
\newcommand{\my}{\mathbf y}

\newcommand{\dint}{\, \mathrm{d}}

\newcommand{\mb}{\mathbf b}
\newcommand{\mh}{\mathbf h}
\newcommand{\wsi}{{ w_s(i) }}
\newcommand{\wrj}{{ w_r(j) }}
\newcommand{\wsj}{{ w_s(j) }}
\newcommand{\wri}{{ w_r(i) }}
\DeclareMathOperator*{\argmin}{arg\,min}

\allowdisplaybreaks

\newcommand{\lb}{ \langle}
\newcommand{\rb}{ \rangle}
\newcommand{\h}{\mathcal H}
\newcommand{\E}{\mathbb E}
\newcommand{\hk}{{\mathcal H   ( \mathbb K ) }}
\newcommand{\hb}{{\mathcal H (\mathbb{K}_\beta) }}
\newcommand{\kmh}{{ \mathbb K^{-1/2}}}
\newcommand{\kh}{{ \mathbb K^{1/2}}}
\newcommand{\wot}{{W^{\alpha ,2}}}
\newcommand{\fk}{ {\mathrm{F} (\mathbb K )}}

\newtheorem{assume}{Assumption}

\usepackage{lastpage}
\jmlrheading{23}{2022}{1-\pageref{LastPage}}{9/21; Revised
7/22}{7/22}{21-1091}{Daren Wang, Zifeng Zhao, Yi Yu  and  Rebecca Willett}
\ShortHeadings{Functional Linear Regression with Mixed Predictors}{Wang, Zhao, Yu  and Willett}

\begin{document}

\title{Functional Linear Regression with Mixed Predictors}

\author{\name Daren Wang \email dwang24@nd.edu \\
       \addr Department of ACMS\\
       University of Notre Dame
\\Indiana, USA
       \AND
       \name Zifeng Zhao \email zifeng.zhao@nd.edu \\
       \addr Mendoza College of Business\\
        University of Notre Dame
\\Indiana, USA     
       \AND
	   \name Yi Yu \email Yi.Yu.2@warwick.ac.uk  \\	
	   \addr Department of Statistics   \\
	    University of Warwick\\
	   Coventry, UK
	   \AND Rebecca Willett \email willett@uchicago.edu \\
	   \addr Department of Statistics  \\
	   University of Chicago\\
	     Illinois, USA}
\editor{Pradeep Ravikumar} 
\maketitle 

 \begin{abstract}

We study a functional linear regression model that deals with functional responses and allows for both functional covariates and high-dimensional vector covariates. The proposed model is flexible and nests several functional regression models in the literature as special cases. Based on the theory of reproducing kernel Hilbert spaces (RKHS), we propose a penalized least squares estimator that can accommodate functional variables observed on discrete sample points. Besides a conventional smoothness penalty, a group Lasso-type penalty is further imposed to induce sparsity in the high-dimensional vector predictors. We derive finite sample theoretical guarantees and show that the excess prediction risk of our estimator is minimax optimal. Furthermore, our analysis reveals an interesting phase transition phenomenon that the optimal excess risk is determined jointly by the smoothness and the sparsity of the functional regression coefficients. A novel efficient optimization algorithm based on iterative coordinate descent is devised to handle the smoothness and group penalties simultaneously. Simulation studies and real data applications illustrate the promising performance of the proposed approach compared to the state-of-the-art methods in the literature.

\end{abstract} 
\begin{keywords}  Reproducing kernel Hilbert space, functional data analysis, high-dimensional regression, minimax optimality.
\end{keywords}
 
\section{Introduction}
\input{introduction}

\section{Background}\label{sec-background}

\input{background}

\section{Main results}\label{sec-main-results}
\input{main}

\section{Optimization}\label{sec-optimization}

\input{optimization}

\section{Numerical results}\label{sec-simulation}
\input{numerical}

\section{Conclusion}\label{sec-discussion}
 \input{conclusion}

 \bibliography{citations}

 \clearpage

\appendix

\section*{Appendices}

\Cref{sec-app-disc-assump} provides more discussions on the connection between bivariate functions and compact linear operators, Assumptions~\ref{assume:mixed} and \ref{assume:joint smooth}.  \Cref{sec-app-proof-sec3} collects proofs of results in \Cref{sec-main-results}.  \Cref{sec-app-func}  contains a series of deviation bounds, which are interesting \emph{per se}. \Cref{sec-app-lb} contains the proof of \Cref{lem-lb111}.  \Cref{sec-app-additional} includes additional details on technical proofs.  \Cref{sec-app-F} collects additional details on optimization.  Additional numerical details are exhibited in Appendices~\ref{subsec:additional_simu}, \ref{subsec:realdata_cv}, \ref{subsec:simu_me}, \ref{subsec:withLargeBasis}, \ref{subsec:withFPCA} and \ref{subsec:phase_transition}.

\input{appendix}

\end{document}

%% file: introduction.tex
Functional data analysis is the collection of statistical tools and results revolving around the analysis of data in the form of (possibly discrete samples of) functions, images and more general objects. Recent technological advancement in various application areas, including neuroscience \citep[e.g.][]{petersen2019frechet, dai2019age}, medicine \citep[e.g.][]{chen2017modelling, ratcliffe2002functional}, linguistic \citep[e.g.][]{hadjipantelis2015unifying, tavakoli2019spatial}, finance \citep[e.g.][]{fan2014functional, benko2007functional}, economics \citep[e.g.][]{ming2007economic, ramsay2002functional}, transportation \citep[e.g.][]{chiou2014functional, wagner2017functional}, climatology \citep[e.g.][]{fraiman2014detecting, bonner2014modeling}, and others, has spurred an increase in the popularity of functional data analysis.

The statistical research in functional data analysis has covered a wide range of topics and areas.  We refer the readers to a recent comprehensive review \citep{wang2016functional}. In this paper, we are concerned with a general functional linear regression model that deals with functional responses and accommodates both functional and vector covariates.

By rescaling if necessary, without loss of generality, we assume the domain of the functional variables is $[0,1]$. Let $A^*(\cdot, \cdot): \, [0, 1] \times [0, 1] \to \mathbb{R}$ be a bivariate coefficient function and $\{\beta^*_j(\cdot): [0, 1] \to \mathbb{R}\}_{j = 1}^p$ be a collection of $p$ univariate coefficient functions. The functional linear regression model concerned in this paper is as follows:
\begin{align}\label{eq-model-mixed-cov}
	Y_t(r)  = \int_{[0, 1]} A^*(r, s)X_t(s)\dint s + \sum_{j = 1}^p \beta^*_j(r)Z_{tj} + \epsilon_t(r), \quad r\in [0,1],
\end{align}
where $Y_t(r ): [0, 1] \to \mathbb{R}$ is the functional response, $X_t(s ): [0, 1] \to \mathbb{R}$ is the functional covariate, $Z_t = (Z_{tj})_{j = 1}^p \in \mathbb{R}^p$ is the vector covariate and $\epsilon_t(r):[0,1]\to \mathbb{R}$ is the functional noise   such that $\mathbb{E}(\epsilon_t(r))=0 $  and $\mathrm{Var} (\epsilon_t(r)   )< \infty  $ for all $ r\in[0,1]$. 
The index $t \in \{1, 2, \ldots, T\}$ denotes $T$ independent and identically distributed samples. The vector dimension $p$ is allowed to diverge as the sample size $T$ grows unbounded.  Furthermore, instead of assuming the functional variables are fully observed, we consider a more realistic setting where the functional covariates $\{X_t\}_{t=1}^T$ and responses $\{Y_t\}_{t=1}^T$ are only observed on discrete sample points $\{s_i\}_{i = 1}^{n_1}$ and $\{r_j\}_{j = 1}^{n_2}$, respectively. The two collections of sample points do not need to coincide. 

As an important real-world example, consider a dataset collected from the popular crowdfunding platform kickstarter.com~(see more details in \Cref{subsec:realdata}). The website provides a platform for start-ups to create fundraising campaigns and charges a 5\% service fee from the final fund raised by each campaign over its 30-day campaign duration. Denote $N_t(r)$ as the pledged fund for the campaign indexed by $t$ at time $r$.  Note that $\{N_t(r), \, r\in[0,30]\}$ forms a fundraising curve. For both the platform and the campaign creators, it is of vital interest to generate an accurate prediction of the future fundraising path $\{N_t(r), r\in(s,30]\}$ at an early time $s\in(0,30)$, as the knowledge of $\{N_t(r), \, r\in(s,30]\}$ helps the platform better assess its future revenue and further suggests timing along $(s, 30]$ for potential intervention by the creators and the platform to boost the fundraising campaign and achieve better outcome.

\begin{figure}[h]
\begin{center}
	 \hspace*{-0.1cm}                                                           
	\includegraphics[angle=270, width=\textwidth]{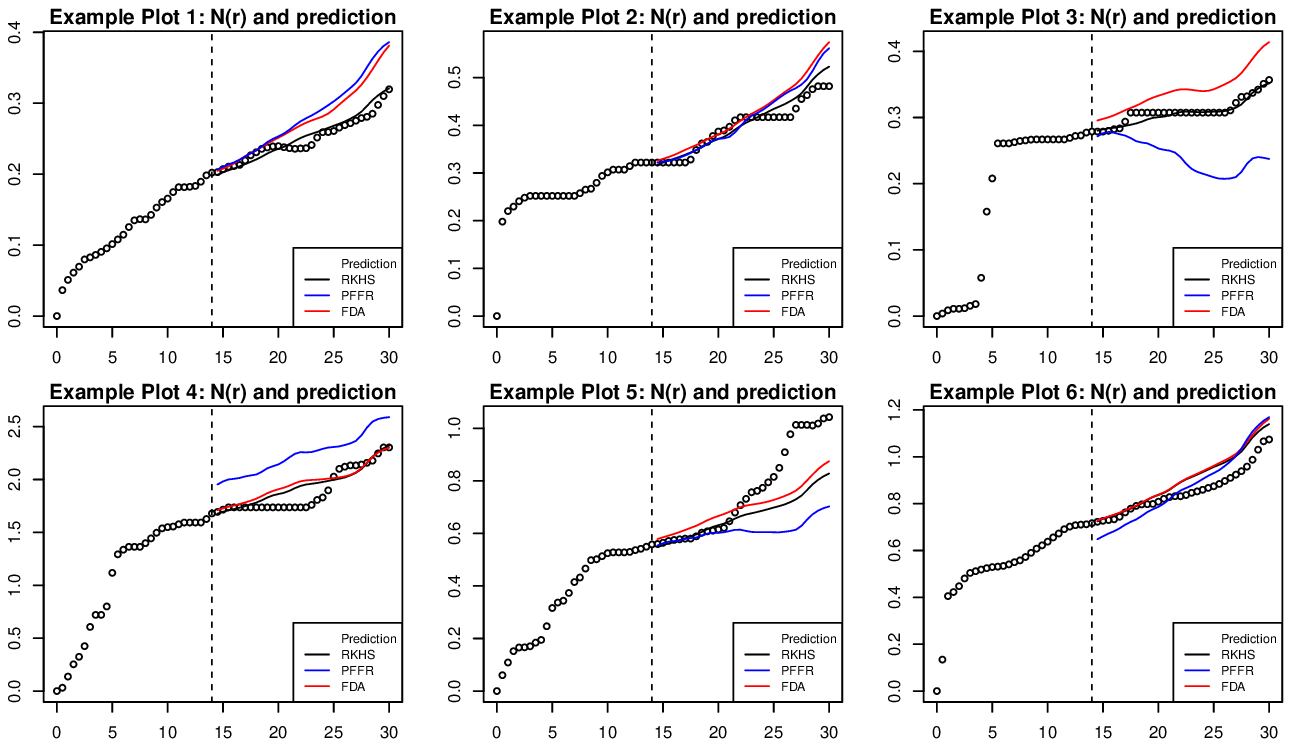}
	\vspace{-0.2cm}
	\caption{Observed (normalized) fundraising curves $\{N_t(r)\}$~(dots) of six representative campaigns and functional predictions~(solid lines) given by RKHS, FDA and PFFR at $s=14$th day~(dashed vertical line) based on the functional covariate $\{N_t(r),\, r\in[0,s]\}$.}
	\label{fig:representative_plot}
\end{center}
	\vspace{-0.3cm}
\end{figure}

At time $s$, to predict the functional response $Y_t$, i.e.~$\{N_t(r), \, r\in(s,30]\}$, a functional regression as proposed in model \eqref{eq-model-mixed-cov} can be built based on the functional covariate $X_t$, i.e.~$\{N_t(r), \, r\in[0,s]\}$, and vector covariates $Z_t$ such as the number of creators and product updates of campaign $t$. \Cref{fig:representative_plot} plots the normalized fundraising curves of six representative campaigns and the functional predictions given by our proposed method~(RKHS) and two competitors, FDA in \cite{Ramsay2005} and PFFR in \cite{ivanescu2015penalized}. In general, the proposed method achieves more favorable performance. More detailed real data analysis is presented in \Cref{subsec:realdata}.

\subsection{Literature review}	
As mentioned in recent reviews, e.g.~\cite{wang2016functional} and \cite{morris2015functional}, there are in general three types of functional linear regression models: 1) functional covariates and  scalar or vector responses \citep[e.g.][]{cardot2003spline, cai2006prediction, hall2007methodology, yuan2010reproducing, raskutti2012minimax, cai2012minimax,  jiang2014inverse,  fan2015functional}; 2) functional covariates and functional responses \citep[e.g.][]{wu1998asymptotic, liang2003relationship, yao2005functional, fan2014functional, ivanescu2015penalized, sun2018optimal} and 3) vector covariates and functional responses \citep[e.g.][]{laird1982random, li2010uniform, faraway1997regression, wu2000kernel}. Another closely related area is functional time series. For instance, functional autoregressive models also preserve the regression format. The literature on functional time series is also abundant, including \cite{van2017nonparametric}, \cite{aue2015prediction}, \cite{bathia2010identifying}, \cite{wang2020functional} and many others.

Regarding the sampling scheme, to establish theoretical guarantees, it is often assumed in the literature that functions are fully observed, which is generally untrue in practice. To study the properties of functional regression under the  more realistic  scenario where functions are observed only on discrete sample points, one usually imposes certain smoothness conditions on the underlying functions. The errors introduced by the discretization of functions can therefore be controlled as a function of the smoothness level and $n$, the number of discretized observations available for each function.  Some fundamental tools regarding this aspect were established in \cite{mendelson2002geometric} and \cite{bartlett2005local}, and were further applied to various regression analysis problems \citep[e.g.][]{raskutti2012minimax, koltchinskii2010sparsity, cai2012minimax}.  

Regarding the statistical inference task, in the regression context, estimation and prediction are two indispensable pillars. In the functional regression literature, \cite{valencia2013radial}, \cite{cai2011optimal}, \cite{yuan2010reproducing}, \cite{lin2006convergence}, \cite{park2018singular}, \cite{fan2014functional},  \cite{fan2015functional}, and \cite{reimherr2019optimal}, among many others, have studied different aspects of the estimation problem.  As for prediction, the existing literature includes \cite{cai2012minimax}, \cite{cai2006prediction}, \cite{ferraty2009additive},   \cite{sun2018optimal}, \cite{reimherr2018optimal} to name but a few.

In contrast to the aforementioned three paradigms of functional regression, in this paper, we study the functional linear regression problem in model \eqref{eq-model-mixed-cov} with functional responses and mixed predictors that consist of both functional and vector covariates. We assume the functions are only observed on discrete sample points and aim to derive optimal upper bounds on the excess prediction risk~(see \Cref{def-excess-risk-mixed}).  More detailed comparisons with existing literature are deferred till we present the main results in \Cref{sec-main-results}.

\subsection{Main contributions}
The main contributions of this paper are summarized as follows.

Firstly, to the best of our knowledge, the model we study in this paper, as defined in \eqref{eq-model-mixed-cov}, is among the most flexible ones in the literature of functional linear regression. In terms of predictors, we allow for both functional and vector covariates. In terms of model dimensionality, we allow for the dimension~$p$ of vector covariates to grow exponentially as the sample size $T$ diverges. In terms of function spaces, which will be elaborated later, we allow for the coefficients of the functional and vector covariates to be from different reproducing kernel Hilbert spaces~(RKHS).  In terms of dependence between the functional and vector covariates, we allow the correlation to be of order up to $O(1)$, which is a rather weak condition, especially considering that the vector covariates are of high dimensions. In terms of the sampling scheme, we allow for the functional covariate and response to be observed on (possibly different) discretized sample points. This general and unified framework imposes new challenges on both theory and optimization, as we elaborate in detail later.

Secondly,  we develop new peeling techniques in our theoretical analysis, which is crucial and fundamental in dealing with the potentially exponentially growing vector covariate dimension. Existing asymptotic analysis techniques and results in the functional regression literature are insufficient to provide Lasso-type guarantees with the presence of high-dimensional covariates. See \Cref{remark:analysis} for more details.

Thirdly, we demonstrate an interesting phase transition phenomenon in terms of the smoothness of the functional covariate and the dimensionality of the vector covariate. To be specific, let~$ \mathfrak {s}$ be the sparsity of $\{\beta^*_j\}_{j  = 1}^p$, i.e.~the number of non-zero univariate coefficient functions. Let $\delta_T$ be a quantity jointly determined by the complexities and the alignment of the two RKHS's where the functional covariates $\{ X_t\}_{t=1}^T$ and the bivariate coefficient function $A^*$ reside; see \Cref{thm-main-mixed} for the detailed definition of $\delta_T$. We show that the excess prediction risk of our proposed estimator is upper bounded by 
\[
	O_{\mathrm{p}}\left\{  \mathfrak {s} \log(p  \vee T ) T ^{-1}  +  \delta_T  \right\} + \mbox{discretization error}.
\] 
\begin{itemize}
	\item When $\delta_T \lesssim    \mathfrak {s} \log(p  \vee T ) T ^{-1}$,  the excess risk is dominated by $  \mathfrak {s} \log(p \vee T ) T ^{-1}$, which is the standard excess risk rate in the high-dimensional parametric statistics literature \citep[e.g.][]{buhlmann2011statistics}. Therefore, the  difficulty is dominated by estimating $\{\beta^*_ j \}_{j  = 1}^p$, the $p$ univariate coefficient functions. 
	\item When $\delta_T  \gtrsim    \mathfrak {s} \log(p \vee T ) T ^{-1}$, the excess risk is dominated by $ \delta_T  $, which suggests that the difficulty is dominated by estimating the bivariate coefficient function $ A^*$. We show that in this regime, the optimal excess risk of model \eqref{eq-model-mixed-cov} is of order $ \delta_T  $. As a result, in this regime the  difficulty is dominated by estimating the bivariate coefficient function $ A^*$. 
	\item We further develop matching lower bounds to justify that this phase transition  between  high-dimensional parametric rate and the non-parametric rate is indeed minimax optimal; see \Cref{sec-lower-bounds}. To the best of our knowledge,  this is the first time such phenomenon is observed in the functional regression literature. 
\end{itemize} 
Note that, this phase transition and its associated optimality are not with respect to the discretization error, i.e.~$n$ - the minimum number of observations obtained for each function.

Lastly, we derive a representer theorem and further propose a novel optimization algorithm via iterative coordinate descent, which efficiently solves a sophisticated penalized least squares problem with the presence of both a ridge-type smoothing penalty and a group Lasso-type sparsity penalty.

The rest of the paper is organized as follows. In \Cref{sec-background}, we introduce relevant definitions and quantities. The main theoretical results are presented in \Cref{sec-main-results}. The optimization procedure is discussed in \Cref{sec-optimization}. Numerical experiments and real data applications are conducted in \Cref{sec-simulation} to illustrate the favorable performance of the proposed estimator over existing approaches in the literature. \Cref{sec-discussion} concludes with a discussion.
 
\subsection*{Notation}

Let $\mathbb{N}_*$ denote the collection of positive integers. Let $a(n)$ and $b(n)$ be two quantities depending on $n$. We say that $a(n) \lesssim b(n)$ if there exists $n_0 \in \mathbb{N}_*$ and an absolute constant $C > 0$ such that for any $n \geq n_0$, $a(n) \leq C b(n)$.  We denote~$a(n) \asymp b(n)$, if $a(n) \lesssim b(n)$ and $b(n) \lesssim a(n)$. For any function $g: [0,1]\to \mathbb R$, denote  
\[
   \|g\|_{\mathcal{L}^2} = \sqrt{\int_{[0, 1]} g^2(r) \dint r} \quad \text{and} \quad \|g\|_\infty = \sup_{r\in [0,1]} |g(r)|.
\]
 For any $p \in \mathbb{N}_*$ and $a, b \in \mathbb R^p$, let $\lb a, b\rb_p = \sum_{i=1}^p a_ib_i$. 

Let $\mathcal{L}^2 = \mathcal{L}^2([0, 1])$ be the space of all square integrable functions with respect to the uniform distribution on $[0, 1]$, i.e.~$\mathcal{L}^2 = \{f: [0, 1] \to \mathbb{R}, \,  \|f\|_\lt < \infty\}$. For any $\alpha >0$, let $W^{\alpha, 2}$ be the Sobolev space of order $\alpha$. For a  positive integer $\alpha$, we have $$W^{\alpha, 2} = \{f \in \lt [0,1], \, f^{(\alpha - 1)} \mbox{ is absolutely continuous} \mbox{ and }  \|f^{(\alpha)} \|_{\lt } < \infty \},$$
where $f^{(k)}$ is the $k$-th weak derivative of $f$, $k \in \mathbb{N}$. 
In this case the Sobolev norm of $f$ is defined as 
$$\|f\|_{\wot }^2 =  \|f\|_\lt^2 + \|f^{(\alpha)}\|_\lt^2 . $$
  For non-integer valued $\alpha > 0$, the Sobolev space $W^{\alpha, 2}$ and the norm $\| \cdot   \|_{\wot }$ are also well-defined \citep[e.g.~see][]{brezis2010functional}.  Note that if $0 < \alpha_1 < \alpha_2$, then $W^{\alpha_2, 2} \subset W^{\alpha_1, 2}$.

%% file: background.tex
In this section, we provide some background on the fundamental tools used in this paper, with RKHS and compact linear operators studied in Section \ref{sec-rkhs} and \Cref{sec-bivfunc-linmap}, respectively.

\subsection{Reproducing kernel Hilbert spaces (RKHS)}	\label{sec-rkhs}

Consider a Hilbert space $\mathcal{H} \subset \mathcal{L}^2$ and its associated inner product $\langle \cdot, \cdot \rangle_{\mathcal{H}}$, under which $\mathcal{H}$ is complete.  We assume that there exists a continuous symmetric nonnegative-definite kernel function $\mathbb{K}: [0, 1] \times [0, 1] \to \mathbb{R}_+$ such that the space  $\mathcal{H}$ is an RKHS, in the sense that for each $r \in [0, 1]$, the function $\mathbb{K}(\cdot, r) \in \mathcal{H}$ and $g(r) = \langle g(\cdot), \mathbb{K}(\cdot, r)\rangle_{\mathcal{H}}$, for all $g \in \mathcal{H}$.  To emphasize this relationship, we write~$\mathcal{H}(\mathbb{K})$ as an RKHS with $\mathbb{K}$ as the associated kernel.  For any $\mathbb{K}$, we say it is a bounded kernel if $\sup_{r \in [0, 1]} \mathbb{K}(r, r) < \infty$.  

It follows from Mercer's theorem \citep[][]{mercer1909xvi} that there exists an orthonormal basis of~$\mathcal{L}^2$, $\{\phi_k\}_{k = 1}^{\infty} \subset \mathcal{L}^2$, such that a non-negative definite kernel function $\mathbb{K}(\cdot, \cdot)$ has the representation $\mathbb{K}(s, r ) = \sum_{k = 1}^{\infty} \mu_k \phi_k(s) \phi_k(r ),~ s, r  \in [0, 1],$ where $\mu_1 \geq \mu_2 \geq \cdots \geq 0$ are the eigenvalues of $\mathbb{K}$ and $\{\phi_k\}_{k = 1}^{\infty}$ are the corresponding eigenfunctions.  In the rest of this paper, when there is no ambiguity, we drop the dependence of $\mu_k$'s and $\phi_k$'s on their associated kernel function $\mathbb{K}$ for ease of notation.
	
Note that any function $f \in \mathcal{H}(\mathbb{K})$ can be written as 
	\[
		f(r) = \sum_{k = 1}^{\infty} \left\{\int_{[0, 1]} f(r)\phi_k(r)\dint r\right\} \phi_k(r) = \sum_{k = 1}^{\infty} a_k\phi_k(r), \quad r \in [0, 1],
	\]
and its RKHS norm is defined as $ \|f\|_{\hk } = \sqrt{\sum_{k = 1}^{\infty} {a_k^2}/{\mu_k}}.$ Thus, for the eigenfunctions, we have $\|\phi_k\|_\hk ^2 = \mu_k^{-1}$.  Throughout this section, we further denote $\psi_k = \sqrt{\mu_k}\phi_k$ and note that $\|\psi_k\|_\hk=1$ for $k\in \mathbb{N}_*$.  
   
Define the linear map $L_\mk : \lt \to \lt$, associated with $\mk$, as $L_\mk  (f) (\cdot)  = \int _{[0,1] }\mk(\cdot , r )f (r ) \dint r$. It holds that $L_\mk (\phi_k ) = \mu_k  \phi_k$, for $k \in \mathbb{N}_*$. Furthermore, define $L_{\mk^{1/2}}: \lt\to\hk$ such that $L_{\mk^{1/2}} (\phi_k )  = \sqrt {\mu_k } \phi_k$, 
and define $L_{\mk ^{ - 1/2 }}: \hk  \to \lt $ such that $L_{\mk ^{ - 1/2 }} (\phi_k )  =  \mu_k^{-1/2} \phi_k$, for $k \in \mathbb{N}_*$.
For any two bivariate functions $R_1(\cdot, \cdot), R_2(\cdot, \cdot): [0, 1] \times [0, 1] \to \mathbb{R}$, define $R_1R_2(r, s) =\int _{[0,1] }R_1(r , u)R_2(u,s )\dint u$, $r, s \in [0, 1]$. It holds that $L_{R_1R_2} = L_{R_1} \circ L_{R_2}$, where $ L_{R_1} \circ L_{R_2}$ represents the  composition of  $R_1$ and $R_2$, which is a linear map from $\lt$ to $\lt$ given that both $L_{R_1}$ and $L_{R_2}$ are linear maps from $\lt$ to $\lt$.

\subsection{Bivariate functions and compact linear operators}\label{sec-bivfunc-linmap}

To regulate the bivariate coefficient function $A^*$ in model \eqref{eq-model-mixed-cov}, we consider a class of compact linear operators in $\hk$.  Specifically, for any compact linear operator $A_2: \hk \to \hk$, denote 
\begin{align*}
	A_2[f,  g] = \langle A_2[g], f\rangle_\hk, \quad f, g \in \hk . 
\end{align*}
Note that $A_2[f, g]$ is well defined for any $f, g \in \hk $ due to the compactness of $A_2$.  Let $\{\psi_i\}_{i=1}^\infty $ be the eigenbasis of  $\hk $ and   $a_{ij} = A_2[\psi_i, \psi_j] = \lb A_2[\psi_j], \psi_i\rb_{\mathcal{H}(\mathbb{K})}$, $i, j \in \mathbb{N}_*$.  We thus have for any $f, g \in \hk $, it follows from \Cref{remark:bivariate functions} in the appendix that  
\begin{align}
	A_2[f, g] &  = \sum_{i, j = 1}^{\infty} a_{ij} \langle f, \psi_i\rangle_{\hk } \langle g, \psi_j\rangle_{\hk }.\label{eq-comp}
\end{align}

Note that via \eqref{eq-comp}, we can define a bivariate function $A_1(r,s),~r, s \in [0, 1]$, affiliated with the compact operator $A_2$. Specifically, plugging $f=\mathbb{K}(r,\cdot)$ and $g=\mathbb{K}(s,\cdot)$ into \eqref{eq-comp}, we have that
\begin{align*}
	A_2[\mathbb K (r,\cdot), \mathbb K(s, \cdot) ]  = \sum_{i, j=1}^\infty a_{ij} \langle \mathbb{K}(r, \cdot), \psi_i\rangle_{\hk } \langle\mathbb{K}(s, \cdot), \psi_j\rangle_{\hk } = \sum_{i, j = 1}^{\infty} a_{ij}\psi_i(r) \psi_j(s) =A_1(r,s).
\end{align*}
From the above set up, it holds that for all  $v \in \hk $, we have
\begin{align*}
	A_2[v](r) = \langle A_1(r, \cdot), v(\cdot) \rangle_\hk, \quad r \in [0, 1].
\end{align*}
We have established an equivalence between a compact linear operator $A_2$ and its corresponding bivariate function $A_1(r,s)$, therefore any compact linear operator $A_2: \hk  \to \hk $ can be viewed as a bivariate function $A_1: [0, 1] \times [0, 1] \to \mathbb{R}$.

To facilitate the later formulation of penalized convex optimization and the derivation of the representer theorem, we further focus on the Hilbert--Schmidt operator, which is an important subclass of compact operators.  A compact operator $A_2$ is Hilbert--Schmidt if 
	\begin{equation}\label{eq-FK-def}
		\|A_2\|_\fk ^2 = \sum_{i, j = 1}^\infty \lb A_2[\psi_j], \psi_i\rb_\hk ^2 = \sum_{i, j = 1}^\infty a_{ij}^2< \infty.
	\end{equation}

In the rest of this paper, for ease of presentation, we adopt some abuse of notation and refer to compact linear operators and the corresponding bivariate functions by the same notation $A$.

%% file: main.tex
\subsection{The constrained/penalized estimator and the representer theorem}\label{sec-main-setup}
Recall that in model \eqref{eq-model-mixed-cov}, the response is the function $Y$, the covariates include the function $X$ and vector $Z$, and the unknown parameters are the bivariate coefficient function $A^* $ and univariate coefficient functions $\beta^*=\{\beta_ l ^*\}_{ l = 1}^p$. Given the observations $\{X_t(s_i), Z_t, Y_t(r_j)\}_{t = 1, i = 1, j = 1}^{T, n_1, n_2}$, our main task is 
to estimate $A^*$ and $\beta^*$.

Define the weight functions $  w_s(i) = (n_1+1) (s_i -s_{i-1}) $  for  $1\le i \le n_1$   and  $w_r(j) = (n_2+1)(   r_j-r_{j-1}  )$  for $ 1\le j \le n_2 $, where by convention we set $s_0=r_0=0$. In addition, for any $1\le l \le p$, define 
$$  \| \beta_l \|_{n_2}   =  \sqrt{ \frac{1}{n_2 } \sum_{j=1}^{n_2 }w_r(j)  \beta_l^2 ( r_j )  }  .$$
We propose the following constrained/penalized  least squares estimator
\begin{align}\label{eq:approx A beta}
& (\widehat A , \widehat \beta ) = \argmin_{\substack{A \in \mathcal{C}_A, \\  \beta = \{\beta_ l \}_{l  = 1}^p \in \mathcal{C}_\beta}} \Bigg[  \frac{1}{Tn_2} \sum_{t= 1}^T \sum_{j=1}^{n_2}  w_r(j)  \left \{Y_t(r_j) -  \frac{1}{n_1}\sum_{i=1}^{n_1} w_s(i)  A  (r_j, s_i)  X_{t  } (s_i)   - \lb \beta (r_j),  Z_t  \rb_p \right\}^2\nonumber\\
& \hspace{10cm} + \lambda \sum_{ l =1}^p \| \beta_ l \|_{n_2} \Bigg],
\end{align}
where $\lambda > 0$ is a tuning parameter that controls the group Lasso penalty,  $\mathcal{C}_A$ and $\mathcal{C}_\beta$ characterize the spaces of coefficient functions $A$ and $\beta$ respectively such that
	\begin{align}\label{eq-ca}
		\mathcal C_A = \{A: \hk \to \hk, \, \|A\|_{\fk}^2     \leq C_A\}, \, \mathcal{C}_{\beta} = \Bigg\{ \{\beta_l   \}_{l=1}^p \subset \hb :\, \sum_{ l =1}^p \| \beta_ l\|  _{\hb} \le C_\beta \Bigg\}.
	\end{align}
Here  $C_A, C_{\beta} > 0$ are two absolute constants and $\|\cdot\|_{\mathrm{F}(\mathbb{K})}$ is defined in \eqref{eq-FK-def}. Note that we allow the RKHS's $\hk$ and $\hb$ to be generated by two different kernels $\mk$ and $\mk_\beta$. The detailed optimization scheme is deferred to \Cref{sec-optimization}.
  
The above optimization problem essentially constrains three tuning parameters. The tuning parameters $C_A$ and $C_{\beta}$ control the smoothness   of the regression coefficient functions. These constrains are standard and ubiquitous in the functional data analysis and non-parametric statistical literature. The tuning parameter $\lambda$ controls the group Lasso penalty, which is used to encourage sparsity when the dimensionality $p$ diverges faster than the sample size $T$. 

The optimization problem in \eqref{eq:approx A beta} makes use of the weight functions $\wsi $ and $\wrj $, which are determined by the discrete sample points $\{ s_i\}_{i=1}^{n_1} $ and $\{ r_j\}_{j=1}^{n_2} $ respectively. This is designed to handle the scenario where the functional variables are observed on unevenly spaced sample points. Indeed, for evenly spaced sample points $ \{ s_i\}_{i=1}^{n_1}$  and $ \{ r_j\}_{j=1}^{n_2}$, we have that $\wsi=\wrj=1 $ for all $i,j$. We remark that it is well known that the integral of  a regular function can be well approximated by its weighted  sum evaluated at discrete sample points.  

Note that \eqref{eq:approx A beta} is an infinite dimensional optimization problem. Fortunately, \Cref{lemma:representer 3} states that the estimator $(\widehat A , \widehat \beta)$ can in fact be written as linear combinations of their corresponding kernel functions evaluated at the discrete sample points $\{s_i\}_{i=1}^{n_1}$ and $\{r_j\}_{j=1}^{n_2}$.
\begin{proposition} \label{lemma:representer 3}
Denote $\mk$ and $\mk_{\beta}$ as the RKHS kernels of $\hk$ and $\hb$ respectively. There always exists a minimizer $(\widehat{A}, \widehat{\beta})$ of \eqref{eq:approx A beta} such that,
	\[
		\widehat{A}(r, s) = \sum_{i = 1}^{n_1} \sum_{j = 1}^{n_2}\widehat{a}_{ij} \mk (r, r_j) \mk(s, s_i), \quad (r, s) \in [0, 1] \times [0, 1], \, \{\widehat{a}_{ij}\}_{i, j = 1}^{n_1, n_2} \subset \mathbb{R}, 
	\]
	and
	\[
		\widehat \beta_ l (r ) = \sum_{j =1}^{n_2} \widehat  b_{lj }\mk_ \beta  (r_j , r), \quad r \in [0, 1], \, l \in \{1, \ldots, p\}, \, \{\widehat{b}_{lj}\}_{l = 1, j = 1}^{p, n_2} \subset \mathbb{R}. 
	\]
\end{proposition}
\Cref{lemma:representer 3} is a generalization of the well-known representer theorem  for RKHS \citep{wahba1990spline}. Various versions of representer theorems are derived and used in the functional data analysis literature \citep[e.g.][]{yuan2010reproducing}.

\subsection{Model assumptions}\label{sec-assumption}
To establish the optimal theoretical guarantees for the estimator $(\widehat{A}, \widehat{\beta})$, we impose some mild model assumptions, on the coefficient functions (\Cref{assume:regularity}), functional and vector covariates (\Cref{assume:mixed}) and sampling scheme (\Cref{assume:joint smooth}).

\begin{assume}[Coefficient functions] \label{assume:regularity}
\
	\begin{itemize}
		\item [(a)] The bivariate coefficient function $A^* $ belongs to $\mathcal{C}_A$ defined in \eqref{eq-ca}, where $\hk$ is the RKHS associated with a bounded kernel $\mk$ and $C_A>0$ is an absolute constant.  
		\item [(b)] The univariate coefficient functions $\{\beta^*_ l \}_{l = 1}^p$ belong to $\mathcal{C}_\beta$  defined in \eqref{eq-ca}, where $\hb$ is the RKHS associated with a bounded kernel $\mk_\beta$ and     $C_\beta >0$ is an absolute constant. 
		
		In addition,     there exists a set $S \subset \{1, \ldots, p\}$ such that $\beta_j^* = 0$ for all~$j \in  \{1, \ldots, p\}\setminus S$ and there exists a sufficiently   large    absolute constant $C_{\mathrm{snr}}> 0$ such that
			\begin{equation}\label{eq-spar-lev-con}
				 T \ge C_{\mathrm{snr}}  \mathfrak {s} \log(p \vee T ),
			\end{equation}
		where $ \mathfrak {s}  = |S|$ denotes the cardinality of the set $S$.		 
	\end{itemize}
\end{assume}

\Cref{assume:regularity}(a) requires the bivariate coefficient function $A^*$ to be a Hilbert--Schmidt operator mapping from $\hk$ to $\hk$.  \Cref{assume:regularity}(b) imposes smoothness and  sparsity conditions on the univariate coefficient functions $\{\beta_l ^*\}_{l =1}^p$ to handle the potential high-dimensionality of the vector covariate. Note that we allow $\{\beta_l ^*\}_{l =1}^p$ to be from a possibly different RKHS $\hb$ than $\hk$ and thus would allow  users to choose   kernels based on their practical  needs.

\begin{assume} [Covariates]\label{assume:mixed}
\
\begin{itemize}
	\item[(a)] The functional noise $\{\epsilon_ t  \}_{t =1}^ T   $ is  a collection of independent and identically distributed   Gaussian processes such that $\mathbb{E}(\epsilon_1(r))=0 $  and $\mathrm{Var} (\epsilon_1(r)   )< C_\epsilon $ for all $ r\in[0,1]$. In addition, $\{\epsilon_t \}_{ t  =1}^ T $ are independent of  $\{X_t, Z_t \}_{t=1}^T  $.
		
	\item [(b)] The functional covariate $\{X_ t  \}_{ t  =1}^T     $  is a collection of independent and identically distributed centered     Gaussian processes with covariance operator  $\Sigma_X$ and $\mathbb{E}\big( \|X_1\|_\lt^2 \big) =C_X<\infty.  $

	\item [(c)]  The vector covariate $\{Z_t\}_{t = 1}^T \subset \mathbb{R}^p$ is a collection of independent and identically distributed Gaussian random vectors from $\mathcal{N}(0, \Sigma_Z)$, where $\Sigma_Z \in \mathbb R^{p\times p}$ is a positive definite matrix such that 
		\[
			c_z  \| v\|^2_2 \le  v^\top \Sigma_Zv \le C_z  \|v\|_2 ^2, \quad v \in \mathbb{R}^p,
		\] 
		and $c_z, C_Z > 0$ are absolute constants.
	
	\item [(d)] For any deterministic $f \in \hk $ and deterministic $v \in \mathbb R^p$, it holds that
		\[
			\mathbb{E}( \lb X_1, f \rb_\lt Z^\top _1 v)  \le \frac{3}{4} \sqrt { \mathbb{E}\{\langle X_1, f\rangle_{\mathcal{L}^2}^2\} (v^\top  \Sigma_Zv)} = \frac{3}{4} \sqrt { \Sigma_X[f,f] (v^\top  \Sigma_Zv)}.
		\]
\end{itemize} 
\end{assume}

\Cref{assume:mixed}(a) and (b) state that both the functional noise $ \{\epsilon_t\}_{t=1}^T $ and functional covariates $ \{X_t\}_{t=1}^T $ are  Gaussian processes.  In fact, one could further relax them to be  sub-Gaussian processes.  
Such assumptions  are frequently used in high-dimensional functional literature such as \cite{kneip2016functional} and  \cite{wang2020low}.
\Cref{assume:mixed}(d) allows that the functional and vector covariates to be correlated up to $3/4$, which means that the correlation, despite the functional and high-dimensional nature of the problem, can be of order $O(1)$.  We do not claim the optimality of the constant $3/4$ but emphasize that this correlation cannot be equal to one, as detailed in \Cref{sec-appendix-A.2}.

Assumptions~\ref{assume:regularity} and \ref{assume:mixed} are sufficient for establishing theoretical guarantees if we require all functions~(i.e.\ the functional responses and functional covariates) to be fully observed, which is typically not realistic in practice. To allow for discretized observations, we further introduce assumptions on the sampling scheme and on the smoothness of the functional covariates.

\begin{assume}[Sampling scheme]\label{assume:joint smooth}  
\
\begin{itemize}
\item [(a)] The discrete sample points $\{s_i\}_{i = 1}^{n_1}$ and $\{r_j\}_{j = 1}^{n_2}$ are collections of   points with $0\le s_1 < s_2\ldots< s_{n_1}=  1   $ and $0\le r_1 < r_2\ldots< r_{n_2}= 1   $, and there exists  an absolute constant    $C_{d} $ such that 
$$   s_ i -s_{i-1} \le \frac{C_d}{n_1} \quad\text{for all $1\le i \le n_{1}+1  $ }  \quad \text{and} \quad 
    r_ j -r_{j -1} \le \frac{C_d}{n_2} \quad \text{for all 
$1\le j \le n_{2}+1 $,}$$
where by convention we set  $s_0= r_0=0$.
\item [(b)]Suppose that $ \hk, \hb \subset W^{\alpha, 2}$ for some $\alpha>1/2$.  In addition, suppose that 
	\begin{align} \label{eq:X almost surely}
		 \mathbb E (\| X_1\|_\wot ^2 )<\infty.
	\end{align}

\end{itemize} 
\end{assume}

\Cref{assume:joint smooth}(a) allows the functional variables to be partially observed on discrete sample points. Importantly, \Cref{assume:joint smooth}(a) can accommodate both fixed and random sampling schemes. In particular, suppose the sample points  are randomly generated from an unknown distribution on $[0,1]$ with a density function $\mu: \, [0,1] \to \mathbb R $ such that $\inf_{r \in[0,1]}\mu(r)>0$, we have \Cref{assume:joint smooth}(a) holds with high probability. We refer to Theorem 1 of \cite{wang2014falling} for more details.

To handle the partially observed functional variables, \Cref{assume:joint smooth}(b) imposes the smoothness assumption that $X$, $A^*$ and $\{\beta^*_j\}_{j=1}^p$ can be enclosed by a common superset, the Sobolev space $W^{\alpha,2} $. In particular,  \eqref{eq:X almost surely}     is a commonly used assumption for   bivariate function estimation  in the functional data analysis  literature. See  for example, \cite{cai2010nonparametric}  and \cite{wang2020low} and references therein.  We emphasize that \Cref{assume:joint smooth}(b) indeed allows different smoothness levels for $X$, $A^*$ and $\{\beta^*_j\}_{j=1}^p$, and only requires that the least smooth space among $X$, $A^*$ and $\{\beta^*_j\}_{j=1}^p$ is covered by $W^{\alpha,2}$. 
 
The condition \eqref{eq:X almost surely} requires that the second moment of $\| X\|_{ W^{\alpha, 2}} $ is finite, which implies $X \in W^{\alpha,2}$ almost surely. Due to the fact that the functional covariates $\{X_t\}_{t=1}^T$ are partially observed, to derive finite-sample guarantees, we need to establish uniform control over the approximation error
$$ \left|\int_{[0,1]}A(r, s)X_t(s)\dint  s -   \frac{1}{n_1}\sum_{i=1}^{n_1} w_s(i) A(r, s_i) X_t(s_i )\right|\quad \text{for all } t\in \{  1 ,\ldots,T\}.$$
This requires the realized sample paths $\{X_t\}_{t=1}^T$ to be regular.   In particular, the second moment condition in  \eqref{eq:X almost surely} can be used to show that $\{\|X_t\|_{\wot}\}_{t=1}^T$ are bounded with high probability, which implies $\{X_t\}_{t=1}^T$ are H\"{o}lder smooth with H\"{o}lder parameter $ \alpha -1/2>0$  by the Morrey inequality~(see \Cref{theorem:Morrey} in Appendix \ref{sec:Holder_smooth} for more details). In \Cref{example:joint smooth} in \Cref{sec-app-disc-assump}, we further provide a concrete example to illustrate a sufficient condition for   \eqref{eq:X almost surely}.

\subsection{Theoretical  guarantees}\label{sec-opt-func-res}

With the assumptions in hand, we establish theoretical guarantees and investigate the minimax optimality of the proposed estimator \eqref{eq:approx A beta}, via the lens of excess risk defined in \Cref{def-excess-risk-mixed}.  The derivation of \eqref{eq-excess-risk111111} is collected in \Cref{lem-derivation-of-eq8}.

\begin{definition}[Excess prediction risk]\label{def-excess-risk-mixed}
	Let $(\widehat{A}, \widehat{\beta})$ be any estimator of $(A^*, \beta^*)$ in model \eqref{eq-model-mixed-cov}. The excess risk of $(\widehat{A}, \widehat{\beta})$ is defined as
	\begin{align}
		 \mathcal{E}^*(\widehat{A}, \widehat{\beta}) =  &\mathbb{E}_{X^*, Z^*, Y^*}\left\{\int_{[0, 1]} \left(Y^*(r) - \int_{[0, 1]} \widehat{A}(r, s) X^*(s)\dint s - \langle Z^*, \widehat{\beta}(r) \rangle_p\right)^2\dint r\right\} \nonumber\\
		& \hspace{0mm} - \mathbb{E}_{X^*, Z^*, Y^*}\left\{\int_{[0, 1]} \left(Y^*(r) - \int_{[0, 1]} A^*(r, s) X^*(s)\dint s - \langle Z^*, \beta^*(r) \rangle_p \right)^2\dint r\right\} \nonumber \\
		= & \mathbb{E}_{X^*, Z^*, Y^*} \left\{  \int_{[0, 1]}  \left (\int_{[0, 1]}  \Delta_A  (r,s)   X^ *(s) \dint s   +    \langle Z^*,  \Delta_\beta (r)\rangle_p     \right)^2          \dint r\right\},  \label{eq-excess-risk111111}
	\end{align}
	where  
	\begin{itemize}
		\item $\Delta_A(r, s) = \widehat{A}(r, s) - A^*(r, s)$ and $\Delta_{\beta}(r)=\widehat{\beta}(r)-\beta^*(r)$ for $s, r \in [0, 1]$;
		\item the random objects $(X^*, Z^*, Y^*)$ are independent from and identically distributed as the observed data generated from model \eqref{eq-model-mixed-cov}.
	\end{itemize}
\end{definition}
  
  \begin{remark}
Definition \ref{def-excess-risk-mixed} is ubiquitously used to measure prediction accuracy in regression settings.  Excess risks are usually considered to  provide   better quantification of    prediction  accuracy   for  the estimators  of interest than the  $\ell_2$ error bounds. 
Consequently, throughout our paper, we evaluate our  proposed  estimators using excess risks. We remark  that as pointed out by
\cite{cai2006prediction} and  \cite{cai2012minimax}, much of the practical interest in the  functional regression parameters is centered around  applications for the purpose of  prediction.   
\end{remark}

 We are now ready to present our main results, which provide upper and lower bounds on the excess risk $\mathcal{E}^*(\widehat{A}, \widehat{\beta})$ of the proposed estimator~$(\widehat{A}, \widehat{\beta})$ defined in \eqref{eq:approx A beta}. For notational simplicity, in the following, we assume without loss of generality that $n_1 = n_2 = n$. For $n_1\neq n_2$, all the statistical guarantees continue to hold by setting $n=\min\{ n_1, n_2\}$.

\subsubsection*{Upper bounds}\label{sec-upper-bound}

\begin{theorem}\label{thm-main-mixed}
Under Assumptions~\ref{assume:regularity}, \ref{assume:mixed} and \ref{assume:joint smooth}, suppose that the eigenvalues $\{ \xi_k\}_{k=1}^\infty$ of the linear operator  $L_{ \mk ^{1/2} \Sigma_X  \mk ^{1/2}}$ satisfy
	\begin{align} \label{eq:eigen decay}
		\xi_k \asymp   k^{-2r },
	\end{align} 
	for some $r>1/2$.
Let $(\widehat{A}, \widehat{\beta})$ be any solution to \eqref{eq:approx A beta} with the tuning parameter $\lambda =C_\lambda \sqrt {{\log(p\vee T)}{T}^{-1}}$ for some sufficiently large constant $C_\lambda$. For any $T\gtrsim \log(n)$, there exists an absolute constant $C > 0$ such that with probability at least $1 - 8T^{-4}   $,  it holds that
	\begin{align}\label{eq-thm-main-mixed-result}
		\mathcal{E}^*(\widehat{A}, \widehat{\beta}) \le  C \log(T) \big \{  \delta_T +  \mathfrak {s} \log(p\vee T ) T^{-1}   +   \zeta_n    \big\},    
	\end{align}
	where $\zeta_n =  n^{-\alpha+1/2}$, $ \mathfrak {s} $ is the sparsity parameter defined in \Cref{assume:regularity}({\bf b}) and  $\delta_T = T^{-2r/(2r+1)}$. 
\end{theorem}


\Cref{thm-main-mixed} provides a high-probability upper bound on the excess risk $\mathcal{E}^*(\widehat{A}, \widehat{\beta})$. The result is stated in a general way for any RKHS's $\hb $ and $\hk$, provided that the spectrum of $L_{ \mk ^{1/2} \Sigma_X  \mk ^{1/2}}$ satisfies the polynomial decay rate in \eqref{eq:eigen decay}. There are three components in the upper bound in \eqref{eq-thm-main-mixed-result}.
\begin{itemize}
	\item The term $\delta_T$ is an upper bound on the error associated with the estimation of the bivariate coefficient function $A^*$ in the presence of functional noise. As formally stated in~\eqref{eq:eigen decay}, $\delta_T$ is determined by the alignment between the kernel $\mk$ of the RKHS that $A^*$ resides in and the covariance operator of $X$. This is the well-known nonparametric rate frequently seen in the functional regression literature, see e.g.\ \cite{cai2012minimax}.  
		
	\item The term $\mathfrak {s} \log(p)T^{-1} $ is an upper bound on the error associated with the estimation of the high-dimensional sparse univariate coefficient functions $ \beta^* $ and is a parametric rate frequently seen in the high-dimensional linear regression settings \citep[e.g.][]{buhlmann2011statistics}.  
	
	\item The term  $\zeta_n$ is an  upper bound  on the error   due to the fact that the functional variables $\{ X_t, Y_t\}_{t=1}^T$ are only observed on discrete sample points. Recall model \eqref{eq-model-mixed-cov} consists of two components $\int_{[0, 1]}A^*(\cdot, s)X_t(s)\dint s$ and $\langle Z_{t}, \beta^*(\cdot)\rangle$.  The discretization errors  are captured through $\zeta_n =n^{-\alpha +1/2}$, which  only depends on the smoothness of $\wot $ with $ \alpha>1/2$.
\end{itemize}
We show later in \Cref{lem-lb111} that, provided the sample points are dense enough, e.g.~$n \gg   T$, up to a logarithmic factor of $T$, this upper bound achieves minimax optimality. In \Cref{subsec:phase_transition}, we further provide numerical illustration for the nonparametric rate $\delta_T$ and the high-dimensional parametric rate $\mathfrak{s}\log(p\vee T)T^{-1}$ in Theorem \ref{thm-main-mixed}.

\begin{remark}[New peeling techniques]\label{remark:analysis} 
	To prove \Cref{thm-main-mixed}, we develop new peeling techniques to obtain new exponential tail bounds, which are crucial in dealing with the potentially exponentially growing dimension $p$.  
\end{remark} 
 
\begin{remark}[Phase transition]
For sufficiently many samples, i.e.~large $n$, the upper bound in \eqref{eq-thm-main-mixed-result} implies that, there exists an absolute constant  $C  > 0$ such that with large probability 
 	\[
		\mathcal{E}^*(\widehat{A}, \widehat{\beta}) \le  C \log(T) \big \{  \delta_T +  \mathfrak {s} \log(p \vee T) T^{-1}\big\}.
 	\]
	This unveils a phase transition between the nonparametric regime and the high-dimensional parametric regime, governed by the eigen-decay of the linear operator  $L_{ \mk ^{1/2} \Sigma_X  \mk ^{1/2} }  $ and the sparsity of the univariate coefficient functions $\beta^*$. To be specific, if $\delta_T  \gtrsim  \mathfrak {s}   \log(p \vee T )T^{-1}$, the high-probability upper bound on the excess risk is determined by the nonparametric rate~$\delta_T $; otherwise, the high-dimensional parametric rate dominates. 
\end{remark}

\Cref{cor-main} in \Cref{subsection:proof of main} further presents a formal theoretical guarantee  which quantifies a discretized version of the excess risk defined in \Cref{def-excess-risk-mixed} for the proposed estimators $(\widehat{A}, \widehat{\beta})$. Due to the fact that the functional variables are only partially observed on discrete sample points, the  discretized  version of the excess risk can be more relevant in certain practical applications.

\subsubsection*{Lower bounds}\label{sec-lower-bounds}  

In this section, we derive a matching lower bound on the excess risk $\mathcal{E}^*(\widehat{A} ,\widehat{\beta})$ and thus show that the upper bound provided in \Cref{thm-main-mixed} is nearly minimax optimal in terms of the sample size $T$, the dimension $p$ and the sparsity parameter $\mathfrak{s}$, saving for a logarithmic factor.  We establish the lower bound under the assumption that the functional variables $\{ X_t, Y_t \}_{t=1}^T$ are fully observed, which is equivalent to setting   $n=\infty$.  Thus, we do not claim optimality in terms of the number of sample points  $n$ for the result in \Cref{thm-main-mixed}.

\begin{proposition} \label{lem-lb111}
Under Assumptions~\ref{assume:regularity} and \ref{assume:mixed}, suppose that the functional variables $\{ X_t, Y_t \}_{t=1}^T$ are fully observed and the eigenvalues $\{ \xi_k\}_{k=1}^\infty$ of the linear operator  $L_{ \mk ^{1/2} \Sigma_X  \mk ^{1/2}}$ satisfy $\xi_k \asymp   k^{-2r }$, for some $r>1/2$. There exists a sufficiently small constant $c > 0$  such that 
	\[
	 \inf_{ \widehat A , \widehat \beta}    \sup_{A^*\in \mathcal C_A, \beta ^* \in \mathcal C_\beta } \mathbb   E\{ \mathcal{E}^*(\widehat{A}, \widehat{\beta})\}    \geq  c\{   T^{-\frac{2r }{2r  + 1}}   +  \mathfrak {s}  \log(p) T^{-1}\},
	\]
	where the infimum is taken over all  possible estimators of $A^* $ and $\beta ^*$ based on the observations $\{X_t, Z_t, Y_t\}_{t=1}^T $.   
\end{proposition}

\subsection{Extension to functional responses with measurement errors}\label{subsec:extension_me}
In this subsection,  we consider the setting where  the functional responses are corrupted with measurement errors.  More precisely, let the functional response $Y_t(r)$ be generated as in \eqref{eq-model-mixed-cov}.  Assume we observe
	\[
		y_{t, j } = Y_t (r_j) +\mathfrak   E   _{t,j} ,   \quad t \in \{1, \ldots, T\}, \, j \in \{1, \ldots, n_2\},
	\]
where $\{ \mathfrak   E  _{t,j}\}_{t=1, j=1}^{T, n_2} $ is a collection of independent and identically distributed sub-Gaussian measurement errors with mean zero and $\mathrm{Var}(\mathfrak   E   _{t,j}) \le C_{\mathfrak   E}  $.
\\
\\
Given the observations $\{X_t(s_i), Z_t, y_{t,j}\}_{t = 1, i = 1, j = 1}^{T, n_1, n_2}$,
consider 
\begin{align}\label{eq:approx A beta noise}
(\widehat A , \widehat \beta ) = \argmin_{\substack{A \in \mathcal{C}_A, \\  \beta = \{\beta_ l \}_{l  = 1}^p \in \mathcal{C}_\beta}} \Bigg[  \frac{1}{Tn_2} \sum_{t= 1}^T \sum_{j=1}^{n_2}  w_r(j)  \left \{ y_{t,j} -  \frac{1}{n_1}\sum_{i=1}^{n_1} w_s(i)  A  (r_j, s_i)  X_{t  } (s_i)   - \lb \beta (r_j),  Z_t  \rb_p \right\}^2\nonumber\\
 + \lambda \sum_{ l =1}^p \| \beta_ l \|_{n_2} \Bigg],
\end{align}
where  $\mathcal{C}_A$ and $\mathcal{C}_\beta$  are defined in \eqref{eq-ca}.  Note that without the presence of measurement  errors, i.e.\ $\mathfrak   E  _{t,j} = 0$ for all $t$ and $j$, the  estimator  \eqref{eq:approx A beta noise} is identical to \eqref{eq:approx A beta} proposed for the setting with only functional noise.  In what follows, we  show that the excess risk of  the  estimator   \eqref{eq:approx A beta noise} also achieves the  same  convergence rate   as that in \Cref{thm-main-mixed}.
	 
\begin{theorem}\label{thm-main-mixed noise}
Suppose all  the  assumptions in  \Cref{thm-main-mixed} hold. 
Let $(\widehat{A}, \widehat{\beta})$ be any solution to \eqref{eq:approx A beta noise} with the tuning parameter $\lambda =C_\lambda \sqrt {{\log(p\vee T)}{T}^{-1}}$ for some sufficiently large constant $C_\lambda$. For any $T\gtrsim \log(n)$, there exists an absolute constant $C > 0$ such that with probability at least $1 - 8T^{-4}   $,  it holds that
	\begin{align}\label{eq-thm-main-mixed-result noise}
		\mathcal{E}^*(\widehat{A}, \widehat{\beta}) \le  C \log(T) \big \{  \delta_T +  \mathfrak {s} \log(p\vee T ) T^{-1}   +   \zeta_n    \big\},    
	\end{align}
	where $\zeta_n =  n^{-\alpha+1/2}$, $ \mathfrak {s} $ is the sparsity parameter defined in \Cref{assume:regularity}({\bf b}) and  $\delta_T = T^{-2r/(2r+1)}$. 
\end{theorem}

%% file: optimization.tex
In this section, we propose an efficient convex optimization algorithm for solving \eqref{eq:approx A beta} given the observations $\{X_t(s_i), Z_t, Y_t(r_j)\}_{t = 1, i = 1, j = 1}^{T, n_1, n_2}$.  We remark that identical optimization can be applied to \eqref{eq:approx A beta noise}.   To ease presentation, in the following we assume without loss of generality that $\h(\mk)=\h(\mk_\beta)$, and thus $\mathbb{K} = \mathbb{K}_\beta$. The general case where $\mathbb{K} \neq \mathbb{K}_\beta$ can be handled in exactly the same way with more tedious notation. Section \ref{subsec:formulation_convex} formulates \eqref{eq:approx A beta} as a convex optimization problem and Section \ref{subsec:coordinateD} further proposes a novel iterative coordinate descent algorithm to efficiently solve the formulated convex optimization.


\subsection{Formulation of the convex optimization}\label{subsec:formulation_convex}


By the equivalence between constrained and penalized optimization \citep[see e.g.][]{Hastie2009}, we can reformulate \eqref{eq:approx A beta} into a penalized optimization such that
\begin{align}\label{eq:approx A beta_penalized}
	(\widehat A , \widehat \beta ) = \argmin_{A, \beta} \Bigg[& \sum_{t= 1}^T \sum_{j=1}^{n_2} w_r(j)\left \{Y_t(r_j) -  \frac{1}{n_1}\sum_{i=1}^{n_1} w_s(i) A  (r_j, s_i)  X_{t  } (s_i)   - \lb \beta (r_j),  Z_t  \rb_p \right\}^2 \nonumber\\
	& + \lambda_1\|A\|_{\mathrm F(\mk)}^2 + \lambda_2 \sum_{l=1}^{p}\|\beta_l\|_{\h(\mk)}+\lambda_3\sum_{l=1}^p\|\beta_l\|_{n_2}   \Bigg],
\end{align}
where $\lambda_1,\lambda_2,\lambda_3> 0$ denote tuning parameters. Note that for notational simplicity, we drop the factor ${1}/{(Tn_2)}$ of the squared loss. In the following, we show that \eqref{eq:approx A beta_penalized} can be solved via convex optimization.

We first define some necessary notation. For $t=1,\ldots, T$, denote the functional curves observed on discrete grids as
\[
Y_t=[Y_t(r_1),Y_t(r_2),\cdots,Y_t(r_{n_2})]^\top \in \mathbb{R}^{n_2} \quad \mbox{and} \quad X_t=[X_t(s_1),X_t(s_2),\cdots,X_t(s_{n_1})]^\top \in \mathbb{R}^{n_1}.
\] 
Define $Y=[Y_1,\cdots,Y_T] \in \mathbb{R}^{n_2 \times T}, \, X=[X_1,\cdots,X_T] \in \mathbb{R}^{n_1 \times T}$ and $Z=[Z_1,\cdots,Z_T]  \in \mathbb{R}^{p \times T}.$

For any $r, s \in [0, 1]$, denote the RKHS kernels as
\begin{align*}
	k_1(r)=[\mk(r,r_1),\mk(r,r_2),\cdots,\mk(r,r_{n_2})]^\top \in \mathbb{R}^{n_2} \text{ and } k_2(s)=[\mk(s,s_1),\mk(s,s_2),\cdots,\mk(s,s_{n_1})]^\top \in \mathbb{R}^{n_1}.
\end{align*}
Denote $K_1=[k_1(r_1),k_1(r_2),\cdots,k_1(r_{n_2})] \in \mathbb{R}^{n_2 \times n_2}$ and $K_2=[k_2(s_1),k_2(s_2),\cdots,k_2(s_{n_1})]  \in \mathbb{R}^{n_1 \times n_1}$.  Note that $K_1=\lb k_1(r), k_1(r)  \rb_{\h(\mk)}$ and $K_2=\lb k_2(s), k_2(s)  \rb_{\h(\mk)}$, thus both are symmetric and positive definite matrices. Furthermore, denote $W_S=\text{diag}(w_s(1),w_s(2),\cdots, w_s(n_1))\in \mathbb{R}^{n_1\times n_1}$ and $W_R=\text{diag}(\sqrt{w_r(1)},\sqrt{w_r(2)},\cdots,\sqrt{w_r(n_2)})\in \mathbb{R}^{n_2\times n_2}$ as the diagonal weight matrices. Define $Y^*=W_RY$, $K_1^*=W_RK_1$ and $K_2^*=K_2W_S$.
	
By the representer theorem~(\Cref{lemma:representer 3}), we have the minimizer of \eqref{eq:approx A beta_penalized} taking the form
$${A}(r,s)=k_1(r)^\top Rk_2(s) \quad \text{and} \quad  {\beta}_l(r)=k_1(r)^\top \mb_l, \quad l=1,\cdots,p,$$
where $R \in \mathbb{R}^{n_2 \times n_1}$ is an $n_2 \times n_1$ matrix and $\mb_l=[b_{l1},b_{l2},\cdots,b_{ln_2}]^\top \in \mathbb{R}^{n_2}$ is an $n_2$-dimensional vector for $l=1,\cdots,p$. Denote $\beta(r)=[{\beta}_1(r), {\beta}_2(r),\cdots, {\beta}_p(r)]^\top=[\mb_1,\mb_2,\cdots,\mb_p]^\top k_1(r)=B^\top k_1(r)$, where $B=[\mb_1,\mb_2,\cdots,\mb_p]\in \mathbb{R}^{n_2\times p}$.

By straightforward algebra, we can rewrite the optimization problem in \eqref{eq:approx A beta_penalized} as
\begin{align}\label{eq:optimization}
	\left\|Y^*-\frac{1}{n_1}K_1^{*}RK_2^*X-K_1^*B Z\right\|_{\mathrm{F}}^2+\lambda_1 \mathrm{tr}(R^\top K_1 RK_2) + \lambda_2 \sum_{l=1}^p\sqrt{\mb_l^\top K_1 \mb_l} + \lambda_3 \sum_{l=1}^p\sqrt{\frac{1}{n_2}\mb_l^\top K_1^{*\top}K_1^{*}\mb_l},
\end{align}
where $\|\cdot\|_{\mathrm{F}}$ and $\text{tr}(\cdot)$ are the Frobenius norm and trace of a matrix. We refer to \Cref{subsec:convex_formulation} for the detailed derivation.

It is easy to see that \eqref{eq:optimization} is a convex function of $R$ and $B$. Note that the first two terms of \eqref{eq:optimization} are quadratic functions and can be handled easily, while the main difficulty of the optimization lies in the group Lasso-type penalty $\lambda_2 \sum_{l=1}^p\sqrt{\mb_l^\top K_1 \mb_l}+\lambda_3 \sum_{l=1}^p\sqrt{\mb_l^\top K_1^{*\top}K_1^{*}\mb_l/n_2}$.

\subsection{Iterative coordinate descent}\label{subsec:coordinateD}
In this section, we propose an efficient iterative coordinate descent algorithm which solves \eqref{eq:optimization} by iterating between the optimization of $R$ and $B$.

\textbf{Optimization of $R$}~(i.e.\ the bivariate coefficient function $A$): Given a fixed $B=\left[\mb_1,\mb_2,\cdots,\mb_p\right]$, denote $\widetilde{Y}=Y^*-K_1^*BZ \in \mathbb{R}^{n_2\times T}$. We have that \eqref{eq:optimization} reduces to a function of $R$   that 
\begin{align}\label{eq:FR_onlyA}
	\|\widetilde{Y}-n_1^{-1}K_1^*RK_2^* X\|_{\mathrm{F}}^2+\lambda_1 \mathrm{tr}(R^\top K_1 RK_2).
\end{align}
Define $E=K_1^{1/2}RK_2^{1/2} \in \mathbb{R}^{n_2 \times n_1}$ and $\me=\mathrm{vec}(E).$ Denote $\widetilde{\my}=\mathrm{vec}(\widetilde{Y})=[\widetilde{Y}_1^\top,\widetilde{Y}_2^\top,\cdots,\widetilde{Y}_T^\top]^\top$ and denote $S_1=n_1^{-1}(X^\top W_S K_2^{1/2})\otimes (W_R K_1^{1/2})\in \mathbb{R}^{Tn_1 \times n_1n_2}$. We can rewrite \eqref{eq:FR_onlyA} as
\begin{align}\label{eq:optimization_nolasso1}
	\|\widetilde{\my}-S_1{\me}\|_2^2 + \lambda_1\|{\me}\|_2^2,
\end{align}
which can be seen as a classical ridge regression with a structured design matrix $S_1$. This ridge regression has a closed-form solution and can be solved efficiently by exploiting the Kronecker structure of $S_1$, see \Cref{subsec:opt_alg} for more details. Thus, given $B$, $R$ can be updated efficiently.

\textbf{Optimization of $B$}~(i.e.\ the univariate coefficient functions $\beta$): Given a fixed $R$, with some abuse of notation, denote $\widetilde{Y}=Y^*- n_1^{-1}K_1^*RK_2^* X \in \mathbb{R}^{n_2 \times T}$. For $l=1, \ldots, p$, define $\mh_l=K_1^*\mb_l \in \mathbb{R}^{n_2}$ and $H=K_1^*B=[\mh_1,\mh_2,\cdots,\mh_p] \in \mathbb{R}^{n_2 \times p}$. Further define $K_3=(K_1^*)^{-1}K_1(K_1^*)^{-1}$. We have that \eqref{eq:optimization} reduces to a function of $H$~(and thus $B$), i.e.
\begin{align}\label{eq:cd_nest1}
	&\left\|\widetilde{Y}- H Z\right\|_{\mathrm{F}}^2+ \lambda_2 \sum_{l=1}^p\sqrt{\mh_l^\top K_3 \mh_l}+ \lambda_3 \sum_{l=1}^p\sqrt{\frac{1}{n_2}\mh_l^\top \mh_l}\nonumber \\
	=&\sum_{t=1}^{T}\left\|\widetilde{Y}_t-\sum_{l=1}^pZ_{tl}\mh_l\right\|^2_2+\lambda_2 \sum_{l=1}^p\sqrt{\mh_l^\top K_3 \mh_l}+ \frac{\lambda_3}{\sqrt{n_2}} \sum_{l=1}^p \|\mh_l\|_2,
\end{align}
which can be seen as a linear regression with two group penalties: a weighted group Lasso penalty and a standard group Lasso penalty. We solve the optimization of \eqref{eq:cd_nest1} by performing coordinate descent on $\mh_l, l=1,2\cdots,p$. See \cite{Friedman2010} for a similar strategy used to solve the sparse group Lasso problem for a linear regression with both a Lasso and a group Lasso penalty.

Specifically, for each $l=1,2,\cdots,p$, the Karush--Kuhn--Tucker condition for $\mh_l$ is
\begin{align*}
	-2\sum_{t=1}^{T}Z_{tl}\left(\widetilde{Y}_t-\sum_{l=1}^pZ_{tl}\mh_l\right) + \lambda_2 K_4\mathbf s_l^{(1)} + \frac{\lambda_3}{\sqrt{n_2}}\mathbf s_l^{(2)}=0.
\end{align*}
Define $K_4$ as the root of $K_3,$ i.e. $K_4K_4=K_3.$ We have that $K_4\mathbf{s}_l^{(1)}$ and $\mathbf{s}_l^{(2)}$ are the subgradients of $\sqrt{\mh_l^\top K_3 \mh_l}$ and $\|\mh_l\|_2$ respectively, such that
\begin{align*}
	\mathbf s_l^{(1)}=\begin{cases}
		{K_4\mh_l}/{\|K_4\mh_l\|_2}, & \mh_l\neq 0,\\
		\text{a vector with } \|\mathbf s_l^{(1)}\|_2 \leq 1, & \mh_l=0,
	\end{cases}	
    ~\text{ and }~
	\mathbf s_l^{(2)}=\begin{cases}
		{\mh_l}/{\|\mh_l\|_2}, & \mh_l\neq 0,\\
		\text{a vector with } \|\mathbf s_l^{(2)}\|_2 \leq 1, & \mh_l=0. 
	\end{cases}
\end{align*}
Define $\widetilde{Y}_{t}^{l}=\widetilde{Y}_t-\sum_{j\neq l}Z_{tj}\mh_j$. Thus, given $\{\mh_j, j\neq l\}$, we have that $\mh_l=0$ is the optimal solution if there exist two vectors $\mathbf s_l^{(1)}$ and $\mathbf s_l^{(2)}$ with $\|\mathbf s_l^{(1)}\|_2 \leq 1$ and $\|\mathbf s_l^{(2)}\|_2 \leq 1$, and
\begin{align*}
	\mathbf s_l^{(1)}=\frac{1}{\lambda_2}K_4^{-1}\left(2\sum_{t=1}^{T}Z_{tl}\widetilde{Y}_{t}^l-\frac{\lambda_3}{\sqrt{n}}\mathbf s_l^{(2)}\right).
\end{align*}
This is equivalent to checking
\begin{align*}
	\min_{\|\mathbf s\|_2\leq 1} \left\|2\sum_{t=1}^{T}Z_{tl}\widetilde{Y}_{t}^l-\frac{\lambda_3}{\sqrt{n}}\mathbf s\right\|_2 \leq 1,
\end{align*}
which is a standard constrained optimization problem and can be solved efficiently.

Otherwise, $\mh_l\neq0$ and to update $\mh_l$, we need to optimize
\begin{align}\label{eq:cd_nest2}
	\sum_{t=1}^{T}\|\widetilde{Y}_t^l-Z_{tl}\mh_l\|^2_2+\lambda_2 \sqrt{\mh_l^\top K_3 \mh_l} + \frac{\lambda_3}{\sqrt{n_2}} \|\mh_l\|_2.
\end{align}
We again solve this optimization by performing coordinate descent on $\mh_{lk}, k=1,2,\cdots, n_2$, where for each $k$, given $\{\mh_{lj},j\neq k\}$, we can update $\mh_{lk}$ by solving a simple one-dimensional optimization
\begin{align}\label{eq:cd_nest3}
	\min\left\{\sum_{t=1}^{T}(\widetilde{Y}_{tk}^l-Z_{tl}\mh_{lk})^2+\lambda_2 \sqrt{K_{3,kk}\mh_{lk}^2+2\mh_{lk}\sum_{j\neq k} K_{3,kj}\mh_{lj}+\sum_{i,j\neq k}\mh_{li} K_{3,ij}\mh_{lj}} + \frac{\lambda_3}{\sqrt{n_2}} \sqrt{\mh_{lk}^2+\sum_{j\neq k}\mh_{lj}^2}\right\}.
\end{align}
Thus, given $R$, $B$ can also be updated efficiently.

\textbf{The iterative coordinate descent algorithm}: Algorithm \ref{alg:cd} formalizes the above discussion and outlines the proposed iterative coordinate descent algorithm. The convergence of coordinate descent for convex optimization is guaranteed under mild conditions, see e.g.\ \cite{Wright2015}. Empirically, the proposed algorithm is found to be efficient and stable, and typically reaches a reasonable convergence tolerance within a few iterations.

\begin{algorithm}
	\caption{Iterative coordinate descent}\label{alg:cd}
	\begin{algorithmic}[1]
		\State \textbf{input}: {Observations $\{X_t(s_i), Z_t, Y_t(r_j)\}_{t = 1, i = 1, j = 1}^{T, n_1, n_2}$, tuning parameters $(\lambda_1,\lambda_2,\lambda_3)$, the maximum iteration $L_{max}$ and tolerance $\epsilon$.}
		\State \textbf{initialization}: $L=1, B_0=R_0=0.$
		\Repeat \Comment{First level block coordinate descent}
			\State Given $B=B_{L-1}$, update $R_L$ via the ridge regression formulation \eqref{eq:optimization_nolasso1}.
			\State Given $R=R_L$, set $\widetilde{Y}=Y^*-{1}/{n_1}K_1^*RK_2^* X$ and initialize $H=K_1^*B_{L-1}$. 
			\Repeat \Comment{Second level coordinate descent}
				\For{$l=1,2,\cdots,p$}
					\State Given $\{\mh_j, j\neq l\}$, set $\widetilde{Y}_{t}^{l}=\widetilde{Y}_t-\sum_{j\neq l}Z_{tj}\mh_j,$ for $t=1,\cdots,T$.
					\If{$\min_{\|\mathbf s\|_2\leq 1} \left\|2\sum_{t=1}^{T}Z_{tl}\widetilde{Y}_{t}^l-\frac{\lambda_3}{\sqrt{n}}\mathbf s\right\|_2 \leq 1$}
						\State Update $\mh_l=0$.
					\Else
						\Repeat \Comment{Third level coordinate descent}
							\For{$k=1,2,\cdots,n_2$}
							    \State Given $\{\mh_{lj},j\neq k\}$, update $\mh_{lk}$ via the one-dimensional optimization \eqref{eq:cd_nest3}.
							\EndFor
						\Until{Decrease of function value \eqref{eq:cd_nest2}  $<\epsilon$.}
					\EndIf
				\EndFor
			\Until{Decrease of function value \eqref{eq:cd_nest1} $<\epsilon$.}
			\State Update $B_L=K_1^{*-1}H$ and set $L \leftarrow L+1$.
		\Until{Decrease of function value \eqref{eq:optimization} $<\epsilon$ or $L\geq L_{max}$.}
		\State \textbf{output}: $\widehat{R}=R_L$ and $\widehat{B}=B_L$.
	\end{algorithmic}
\end{algorithm}

%% file: numerical.tex
 In this section, we conduct extensive numerical experiments to investigate the performance of the proposed RKHS-based penalized estimator~(hereafter RKHS) for the functional linear regression with mixed predictors. Sections \ref{subsec:simu_setting}-\ref{subsec:mfr} compare RKHS with popular methods in the literature via simulation studies. Section \ref{subsec:realdata} presents a real data application on crowdfunding prediction to further illustrate the potential utility of the proposed method.
The implementations of our numerical experiments can be found at \url{https://github.com/darenwang/functional_regression}.
\subsection{Simulation settings}\label{subsec:simu_setting}
\textbf{Data generating process}: We simulate data from the functional linear regression model
\begin{align}\label{eq:FR_mixed}
	Y_t(r)=\int_{[0,1]} A^*(r,s)X_t(s)\dint s + \sum_{j=1}^{p} \beta^*_j(r)Z_{tj} + \epsilon_t(r), \quad r \in [0,1],
\end{align}
for $t=1,2,\ldots, T.$ Note that for $p=0$, \eqref{eq:FR_mixed} reduces to the classical function-on-function regression.

We generate the functional covariate $\{X_t\}$ and the functional noise $\{\epsilon_t\}$ from a $q$-dimensional subspace spanned by basis functions $\{u_i(s)\}_{i=1}^{q}$, where $\{u_i(s)\}_{i=1}^{q}$ consists of orthonormal basis of $\mathcal L^2[0,1]$. Following \cite{yuan2010reproducing}, we set $u_i(s)=1$ if $i=1$ and $u_i(s)=\sqrt{2}\cos((i-1)\pi s)$ for $i=2, \ldots, q$. Thus, we have $X_t(r)=\sum_{i=1}^qx_{ti}u_i(r)$ and $\epsilon_t(r)=\sum_{i=1}^{q} e_{ti}u_i(r)$. For $x_t=(x_{t1},\ldots,x_{tq})^\top$, we simulate $x_{ti} \stackrel{\mbox{i.i.d.}}{\sim} \text{Unif}[-1/i,1/i]$, $i=1,\ldots, q$. For $e_t=(e_{t1},\ldots,e_{tq})^\top$, we simulate $e_{ti} \stackrel{\mbox{i.i.d.}}{\sim} \text{Unif}[-0.2/i,0.2/i]$,  $i=1,\ldots, q$. For the vector covariate $\{Z_{tj}\}$, we simulate $Z_{tj} \stackrel{\mbox{i.i.d.}}{\sim} \text{Unif}[-1/\sqrt{3}, 1/\sqrt{3}]$, $j=1,\ldots,p$.

For the coefficient functions $A^*(r,s)$ and $\{\beta^*_j(r)\}_{j=1}^p$, we consider Scenarios A and B.
\begin{itemize} 
	\item Scenario A (Exponential): We set $A^*(r,s)=\kappa \sqrt{3}e^{-(r+s)}$, $\beta^*_1(r)=\kappa \sqrt{3}e^{-r}$ and $\beta_j^*(r)\equiv 0$, $j=2,\ldots, p$.
	\item Scenario B (Random): We set $A^*(r,s)=\kappa\sum_{i,j=1}^{q} \lambda_{ij} u_i(r)u_j(s)$, $\beta^*_1(r)=\kappa\sum_{i=1}^{q} b_{i}u_i(r)$ and $\beta^*_j(r)\equiv 0$, $j=2,\ldots,p$. Define matrix $\Lambda = (\lambda_{ij})$ and $b=(b_{1},\ldots,b_{q})^\top.$ We simulate $\lambda_{ij} \stackrel{\mbox{i.i.d.}}{\sim} \mathcal{N}(0,1)$ and rescale $\Lambda$ such that its spectral norm $\|\Lambda\|_{\mathrm{op}}=1.$ For $b$, we simulate $b_{i} \stackrel{\mbox{i.i.d.}}{\sim} \mathcal{N}(0,1)$ and rescale $b$ such that $\|b\|_2=1.$

\end{itemize}

Note that $(A^*,\beta^*)$ in Scenario B is more complex than that in Scenario A, especially when $q$ is large, while for Scenario A, its complexity is insensitive to $q$. The parameter $\kappa$ is later used to control the signal-to-noise ratio~(SNR). Here, we define the SNR for $A^*(r,s)$ as 
\begin{align*}
	\sqrt{\mathbb{E}\int_{[0,1]}\bigg[\int_{[0,1]} A^*(r,s)X_t(s)\dint s\bigg]^2dr} \bigg/ \sqrt{\mathbb{E}\int_{[0,1]}\epsilon_t(r) ^2\dint r},
\end{align*}
which roughly equals to 1, 2 and 4 as we vary $\kappa=0.5,1,2.$ Similarly, we define SNR for $\beta^*$ (note that only $\beta_1^*$ is a non-zero function) as
\begin{align*}
	\sqrt{\mathbb{E}\int_{[0,1]} \bigg[\beta_1^*(r)Z_{t1}\bigg]^2\dint r} \bigg/ \sqrt{\mathbb{E}\int_{[0,1]}\epsilon_t(r) ^2\dint r},
\end{align*}
which also roughly equals to 1, 2 and 4 as we vary $\kappa=0.5,1,2.$

For simplicity, we set the discrete sample points $\{r_j\}_{j=1}^{n_2}$ for $Y_t$ and $\{s_i\}_{i=1}^{n_1}$ for $X_t$ to be evenly spaced grids on $[0,1]$ with the same number of grids $n=n_1=n_2$. The simulation result for random sample points where $\{r_j\}_{j=1}^{n_2}$ and $\{s_i\}_{i=1}^{n_1}$ are generated independently via the uniform distribution on [0,1] is similar and thus omitted.

\textbf{Evaluation criteria}: We evaluate the performance of the estimator by its excess risk. Specifically, given the sample size $(n, T)$, we simulate observations $\{X_t(s_i), Z_t, Y_t(r_j)\}_{t=1, i=1, j=1}^{T+0.5T, n, n}$, which are then split into the training data $\{X_t(s_i), Z_t, Y_t(r_j)\}_{t=1, i=1, j=1}^{T, n, n}$ for constructing the estimator $(\widehat{A},\widehat{\beta})$ and the test data $\{X_t(s_i), Z_t, Y_t(r_j)\}_{t=T+1, i=1, j=1}^{T+0.5T, n, n}$ for the evaluation of the excess risk.  Based on $(\widehat{A},\widehat{\beta})$ and the predictors $\{X_t(s_i), Z_t\}_{t=T+1, i=1}^{T+0.5T, n}$, we generate the prediction $\{\widehat{Y}_t(r)\}_{t=T+1}^{T+0.5T}$ and define
\begin{align}
	&\text{RMISE}(\widehat{A}, \widehat{\beta})= \sqrt{ \frac{1}{0.5T }\sum_{t=T+1}^{T+0.5T}\int_{[0,1]}\left(Y_t^{\text{oracle}}(r)-\widehat{Y}_t(r) \right)^2\dint r},\label{eq:RMISE}\\
	&\text{nRMISE}(\widehat{A}, \widehat{\beta})= \sqrt{ \frac{\frac{1}{0.5T }\sum_{t=T+1}^{T+0.5T}\int_{[0,1]}\left(Y_t^{\text{oracle}}(r)-\widehat{Y}_t(r) \right)^2\dint r}{\frac{1}{0.5T}\sum_{t=T+1}^{T+0.5T}\int_{[0,1]}\left(Y_t^{\text{oracle}}(r)-0 \right)^2\dint r  }},
\end{align} 
where $Y_t^{\text{oracle}}(r)=\int_{[0,1]}A(r,s)X_t(s)\dint s +\beta_1(r)Z_{t1}$ is the oracle prediction of $Y_t(r)$. Note that nRMISE is a normalized RMISE, which can be viewed as a percentage error and thus is easy to assess and interpret. Smaller RMISE and nRMISE indicate a better recovery of the signal.

\textbf{Simulation settings}: We consider three simulation settings for $(n,T,q)$ where $(n,T,q)\in \{(5,50,5), (20,100,20),(40,200,50)\}.$ For each setting, we vary $\kappa\in \{0.5,1,2\}$, which roughly corresponds to $\text{SNR}\in\{1,2,4\}.$ As for the number of scalar predictors $p$, Section \ref{subsec:ffr} considers the classical function-on-function regression, which is a special case of \eqref{eq:FR_mixed} with $p=0$ and Section \ref{subsec:mfr} considers functional regression with mixed predictors and sets $p=3,10,50,100,200$. For each setting, we conduct 500 experiments.

\textbf{Implementation details of the RKHS estimator}: We set $\mathbb{K}=\mathbb{K}_\beta$ and use the rescaled Bernoulli polynomial as the reproducing kernel such that

\centerline{$\mk(x,y)=1+k_1(x)k_1(y)+k_2(x)k_2(y)-k_4(x-y),$}
\noindent where $k_1(x)=x-0.5$, $k_2(x)=2^{-1}\{k_1^2(x)-1/12\}$, $k_4(x)=1/24\{k_1^4(x)- k_1^2(x)/2 + 7/240\}$, $x\in[0,1]$, and $k_4(x-y)=k_4(|x-y|)$, $x,y\in [0,1]$. Such $\mk$ is the reproducing kernel for $W^{2,2}$. See Chapter 2.3.3 of \cite{Gu2013} for more details. In \Cref{alg:cd}, we set the tolerance parameter $\epsilon=10^{-8}$ and the maximum iterations $L_{\max}=10^4.$ A standard 5-fold cross-validation (CV) on the training data is used to select the tuning parameters $(\lambda_1,\lambda_2,\lambda_3)$. Note that for $p=0$, we explicitly set $\beta\equiv 0$ in the penalized optimization $\eqref{eq:approx A beta_penalized}$, which reduces to the structured ridge regression in \eqref{eq:optimization_nolasso1} and can be solved efficiently with a closed-form solution as detailed in Section \ref{subsec:coordinateD}.

\subsection{Function-on-function regression}\label{subsec:ffr}
In this subsection, we set $p=0$ and compare the proposed RKHS estimator with two popular methods in the literature for function-on-function regression.

The first competitor \citep{Ramsay2005} estimates the coefficient function $A(r,s)$ based on a penalized B-spline basis function expansion. We denote this approach as FDA and it is implemented via the R package \texttt{fda.usc}~(\texttt{fregre.basis.fr} function). The second competitor \citep{ivanescu2015penalized} is the penalized flexible functional regression~(PFFR) in  and is implemented via the R package \texttt{refund}~(\texttt{pffr} function).  Unlike our proposed RKHS method,  neither FDA or PFFR has optimal theoretical guarantees.

Both FDA and PFFR are based on penalized basis function expansions and require a hyper-parameter $N_b$, which is the number of basis. Intuitively, the choice of $N_b$ is related to the complexity of $A^*(r,s)$, which is unknown in practice. In the literature, it is recommended to set $N_b$ at a large number to guard against the underfitting of $A^*(r,s)$. On the other hand, a larger $N_b$ can potentially incur higher estimation variance and will also increase the computational cost. For a fair comparison, for both FDA and PFFR we set $N_b=20$, which is sufficient to accommodate the model complexity across all simulation settings except for Scenario B with $q=50$. See more discussions later. We use a standard 5-fold CV to select the tuning parameter $\lambda_1$ of RKHS and the roughness penalty of FDA from the range $10^{(-15:0.5:2)}$. PFFR uses a restricted maximum likelihood~(REML) approach to automatically select the roughness penalty.

Besides the penalized method based on basis function expansions~(i.e. FDA and PFFR), another popular method for function-on-function regression in the literature is based on functional PCA~(i.e. the Karhunen--Lo\'eve expansion, FPCA), see for example \cite{yao2005functional}. However, FPCA does not seem to perform competitively in our simulation and real data analysis, see similar observations in \cite{cai2012minimax} and \cite{sun2018optimal}. For completeness, we report the performance of FPCA in \Cref{subsec:realdata_cv} and \Cref{subsec:withFPCA}.

\textbf{Numerical results}: For each method, Table \ref{tab:FR} reports its average nRMISE~(nRMISE$_{\mathrm{avg}}$) across 500 experiments under all simulation settings. For each simulation setting, we further report R$_{\mathrm{avg}}$, which is the percentage improvement of excess risk given by RKHS defined as 

\begin{figure}[h]
	\begin{subfigure}{0.32\textwidth}
		\includegraphics[angle=270, width=1.2\textwidth]{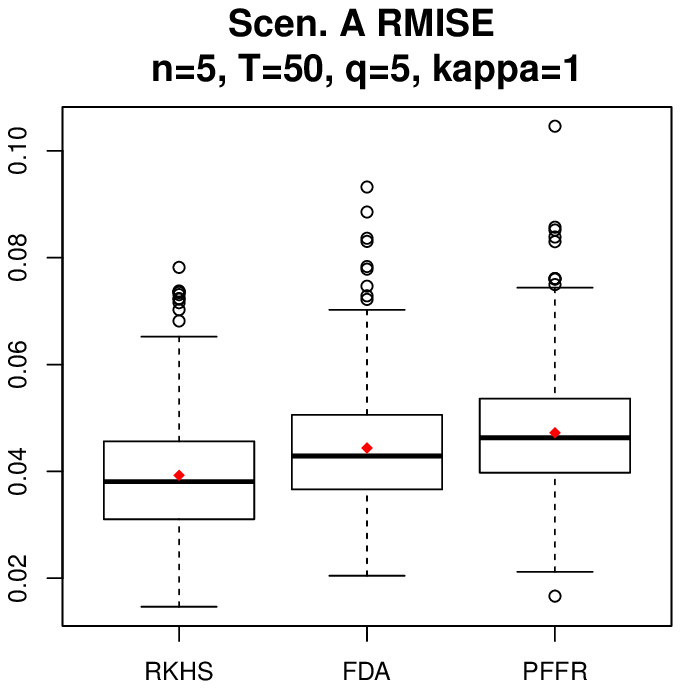}
		\vspace{-0.8cm}
	\end{subfigure}
	~
	\begin{subfigure}{0.32\textwidth}
		\includegraphics[angle=270, width=1.2\textwidth]{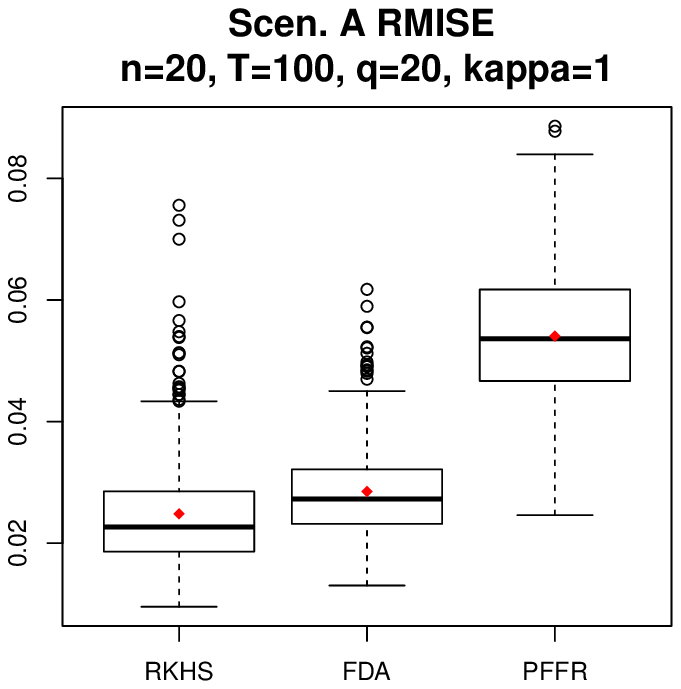}
		\vspace{-0.8cm}
	\end{subfigure}
	~
	\begin{subfigure}{0.32\textwidth}
		\includegraphics[angle=270, width=1.2\textwidth]{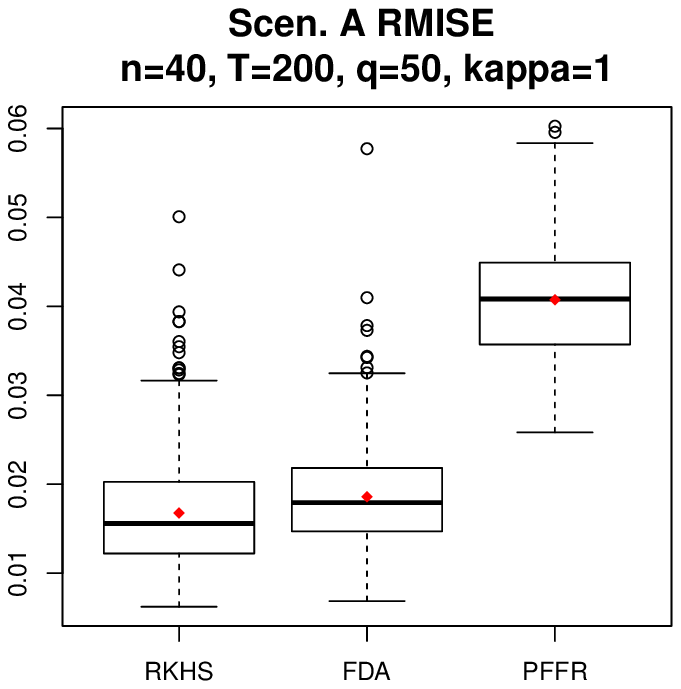}
		\vspace{-0.8cm}
	\end{subfigure}
    ~
    \begin{subfigure}{0.32\textwidth}
    	\includegraphics[angle=270, width=1.2\textwidth]{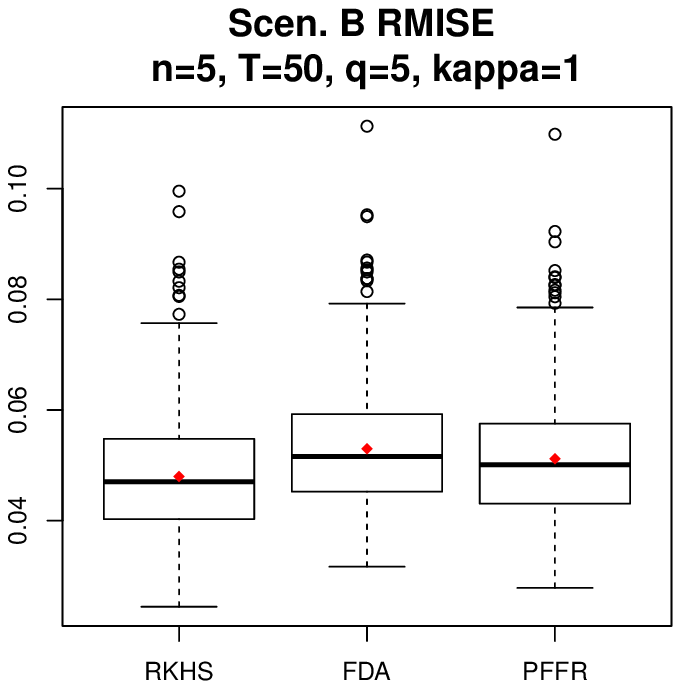}
		\vspace{-0.8cm}
    \end{subfigure}
    ~
    \begin{subfigure}{0.32\textwidth}
    	\includegraphics[angle=270, width=1.2\textwidth]{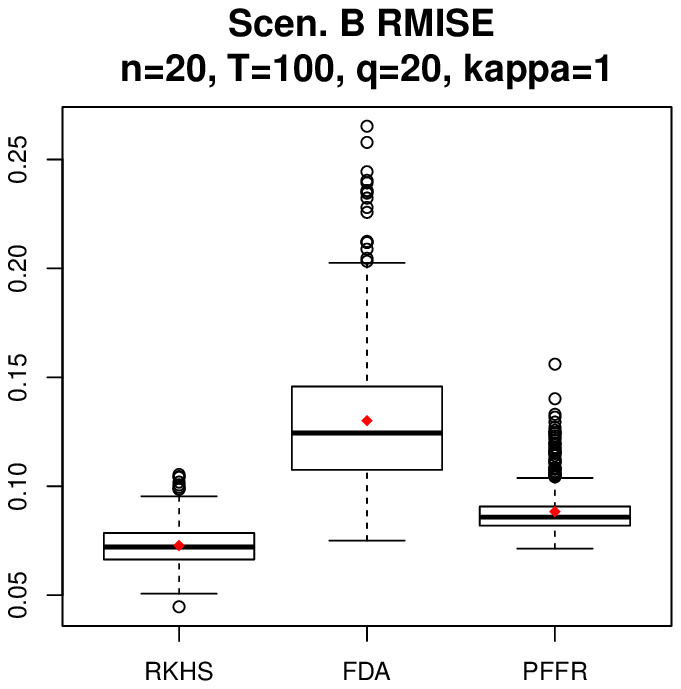}
		\vspace{-0.8cm}
    \end{subfigure}
    ~
    \begin{subfigure}{0.32\textwidth}
    	\includegraphics[angle=270, width=1.2\textwidth]{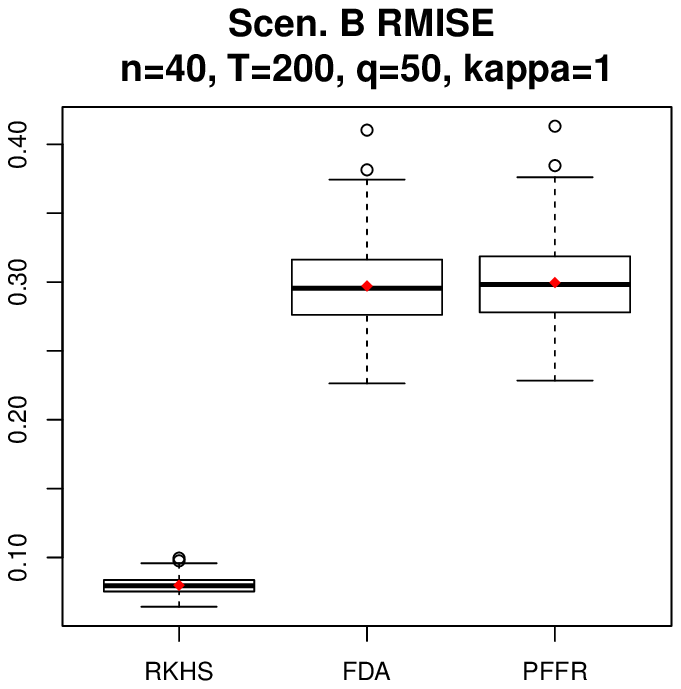}
		\vspace{-0.8cm}
    \end{subfigure}
	\caption{Boxplots of RMISE of RKHS, FDA and PFFR across 500 experiments under function-on-function regression with $\kappa=1$. Red points denote the average RMISE.}
	\label{fig:FR}
\end{figure}

\centerline{R$_{\mathrm{avg}}$=[$\min\{$nRMISE$_{\mathrm{avg}}$(FDA), nRMISE$_{\mathrm{avg}}$(PFFR)$\}$ / nRMISE$_{\mathrm{avg}}$(RKHS$) -1]\times 100\%$.}
\noindent In addition, we report R$_w$, which is the percentage of experiments (among 500 experiments) where RKHS gives the lowest RMISE. We further give the boxplot of RMISE in Figure \ref{fig:FR} under $\kappa=1$. The boxplots of RMISE under $\kappa\in\{0.5,2\}$ can be found in \Cref{subsec:additional_simu}.

As can be seen, RKHS in general delivers the best performance across almost all simulation settings. For all methods, the estimation performance improves with a larger SNR~(reflected in $\kappa$). Since the complexity of Scenario A is insensitive to $q$, the performance of all methods improve as $(n,T)$ increases under Scenario A. However, this is not the case for Scenario B, as a larger $q$ increases the estimation difficulty. Compared to Scenario A, the excess risk RMISE of Scenario B is larger, especially for $q \in \{20, 50\}$, due to the more complex nature of the operator $A^*.$ In general, the improvement of RKHS is more notable under high model complexity and SNR.

Note that for Scenario B with $q=50$, both FDA and PFFR incur significantly higher excess risks due to the underfitting bias caused by the insufficient basis dimension $N_b$, indicating the potential sensitivity of the penalized basis function approaches to the hyper-parameter. In comparison, RKHS automatically adapts to different levels of model complexity. In \Cref{subsec:withLargeBasis}, we further conduct the same simulation but with a larger number of basis dimension $N_b=50$ for FDA and $N_b=30$ for PFFR. For Scenario A, where the bivariate function $A^*(r,s)$ is a simple exponential function, FDA and PFFR give essentially the same performance, while for Scenario B with $q=50$, FDA and PFFR give much improved performance due to lower underfitting bias, though still having a notable performance gap compared to RKHS. In addition, note that a larger $N_b$ can significantly increase the computational cost of the penalized basis function approaches. We refer to \Cref{subsec:withLargeBasis} for more details.

In addition, \Cref{subsec:simu_me} collects the results of the same simulation but with the observed functional responses $Y_t(r)$ being additionally corrupted with measurement errors, i.e.\ the scenario discussed in Section \ref{subsec:extension_me}~(see Theorem \ref{thm-main-mixed noise}). It is seen that the performance of RKHS only worsens slightly with the additional measurement errors, providing numerical support for Theorem \ref{thm-main-mixed noise} that measurement errors do not affect the convergence rate of RKHS.


\begin{table}[h]
	\centering
	\begin{tabular}{rrrrrr|rrrrr}
		\hline\hline
		 & \multicolumn{5}{c|}{ Scenario A:  $n= 5 , T= 50, q=5$} & \multicolumn{5}{c}{ Scenario B:  $n= 5 , T= 50, q=5 $} \\
		$\kappa$ & RKHS & FDA & PFFR &R$_{\mathrm{avg}}$(\%) & R$_{w}$ (\%) & RKHS & FDA & PFFR &R$_{\mathrm{avg}}$(\%) & R$_{w}$ (\%)\\
		\hline
		$0.5$ & 15.98 & 17.97 & 19.80 & 12.46 & 64 & 20.86 & 22.98 & 22.28 & 6.79 & 63 \\ 
		$1$ & 9.14 & 10.31 & 10.98 & 12.81 & 68 & 10.42 & 11.50 & 11.12 & 6.63 & 71 \\ 
		$2$ & 4.99 & 6.18 & 5.70 & 14.37 & 72 & 5.21 & 5.75 & 5.56 & 6.68 & 74 \\ 	
		\hline
		& \multicolumn{5}{c|}{ Scenario A:  $n= 20 , T= 100, q=20 $} & \multicolumn{5}{c}{ Scenario B:  $n= 20 , T= 100, q=20 $} \\
		$\kappa$ & RKHS & FDA & PFFR &R$_{\mathrm{avg}}$(\%) & R$_{w}$ (\%) & RKHS & FDA & PFFR &R$_{\mathrm{avg}}$(\%) & R$_{w}$ (\%)\\
		\hline		
		$0.5$ & 10.61 & 12.15 & 24.84 & 14.57 & 80 & 37.41 & 45.91 & 36.51 & $-$2.41 & 47 \\ 
		$1$ & 5.89 & 6.75 & 12.84 & 14.66 & 81 & 18.39 & 32.41 & 22.32 & 21.37 & 92 \\ 
		$2$ & 3.60 & 4.26 & 6.89 & 18.45 & 84 & 9.19 & 27.71 & 12.56 & 36.58 & 100 \\ 
		\hline 
		& \multicolumn{5}{c|}{ Scenario A:  $n= 40 , T= 200, q=50 $} & \multicolumn{5}{c}{ Scenario B:  $n= 40 , T= 200, q=50 $} \\
		$\kappa$ & RKHS & FDA & PFFR &R$_{\mathrm{avg}}$(\%) & R$_{w}$ (\%) & RKHS & FDA & PFFR &R$_{\mathrm{avg}}$(\%) & R$_{w}$ (\%)\\
		\hline
		$0.5$ & 7.37 & 8.53 & 18.93 & 15.72 & 81 & 39.22 & 79.37 & 78.92 & 101.21 & 100 \\ 
		$1$ & 3.98 & 4.41 & 9.68 & 10.96 & 75 & 20.75 & 76.74 & 77.40 & 269.73 & 100 \\ 
		$2$ & 2.34 & 2.36 & 5.05 & 0.71 & 49 & 12.59 & 75.95 & 77.09 & 503.19 & 100 \\ 
		\hline\hline
	\end{tabular}
	\caption{Numerical performance of RKHS, FDA and PFFR under function-on-function regression. The reported nRMISE$_{\mathrm{avg}}$ is multiplied by 100 in scale. R$_{\mathrm{avg}}$ reflects the percent improvement of RKHS over the best performing competitor, and R$_{w}$ reflects the percentage of experiments in which RKHS achieves the lowest RMISE.} 
	\label{tab:FR}
\end{table}

\subsection{Functional regression with mixed predictors}\label{subsec:mfr}
In this subsection, we set $p=3$ and compare the performance of RKHS with PFFR, as FDA cannot handle mixed predictors.  For RKHS, we use the standard 5-fold CV to select the tuning parameter $(\lambda_1,\lambda_2,\lambda_3)$ from the range $10^{(-4:1:0)}\times10^{(-4:1:0)}\times 10^{(-1:1:3)}$ for Scenario A and from the range $10^{(-17:1:-13)}\times10^{(-17:1:-13)}\times 10^{(-1:1:3)}$ for Scenario B.  PFFR uses REML to automatically select the roughness penalty and we set $N_b=20$ for PFFR as before.

\begin{figure}[h]
	\begin{subfigure}{0.32\textwidth}
		\includegraphics[angle=270, width=1.2\textwidth]{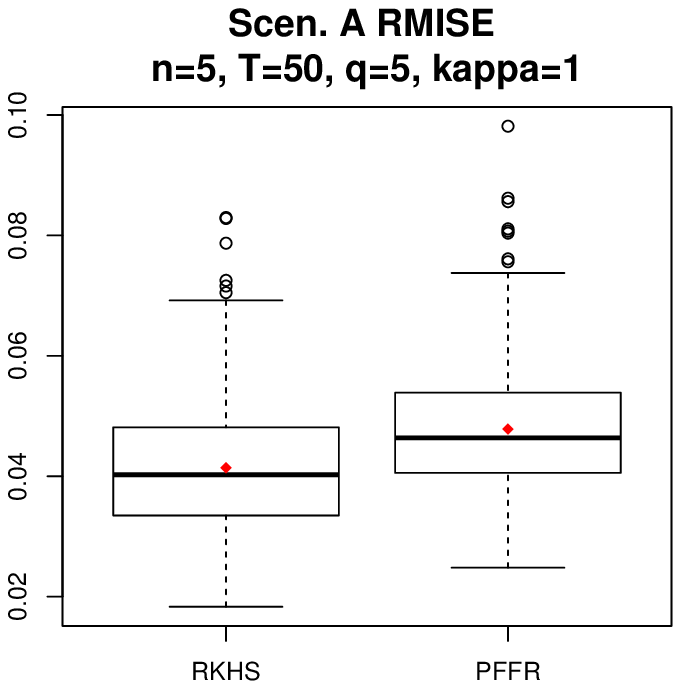}
		\vspace{-0.8cm}
	\end{subfigure}
	~
	\begin{subfigure}{0.32\textwidth}
		\includegraphics[angle=270, width=1.2\textwidth]{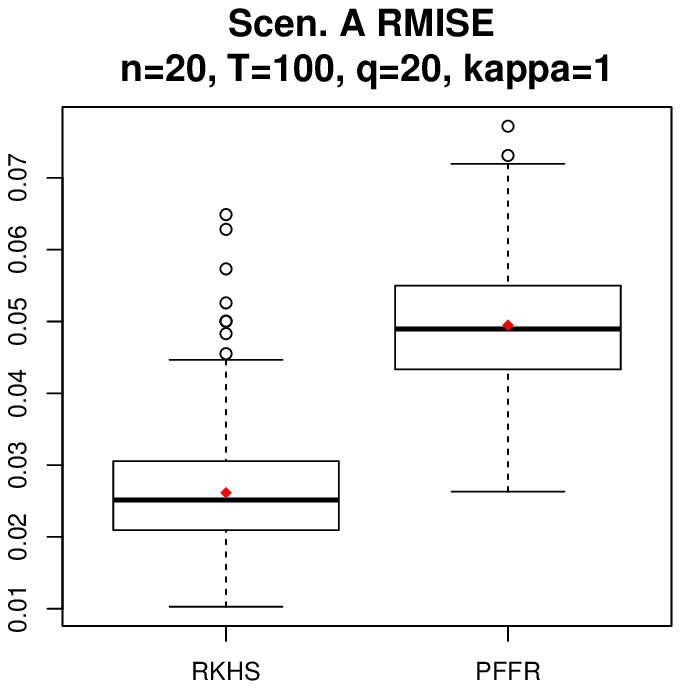}
		\vspace{-0.8cm}
	\end{subfigure}
	~
	\begin{subfigure}{0.32\textwidth}
		\includegraphics[angle=270, width=1.2\textwidth]{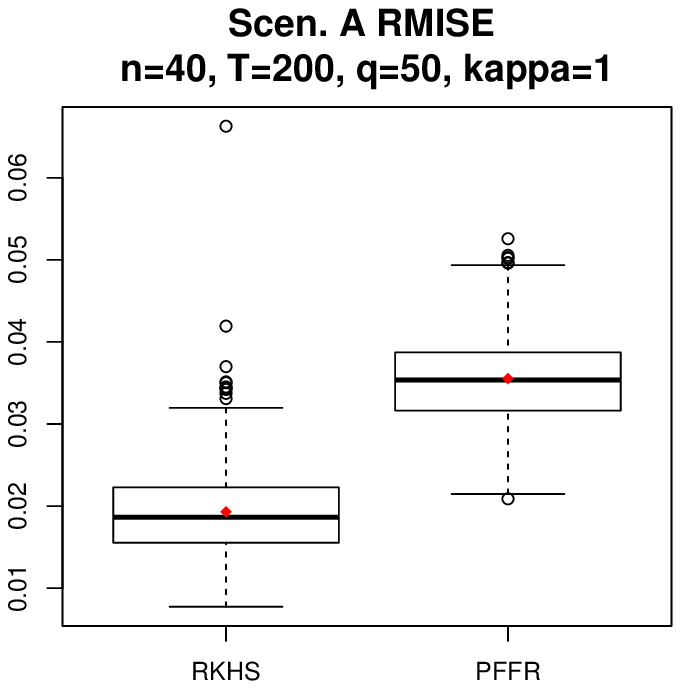}
		\vspace{-0.8cm}
	\end{subfigure}
    ~
    \begin{subfigure}{0.32\textwidth}
    	\includegraphics[angle=270, width=1.2\textwidth]{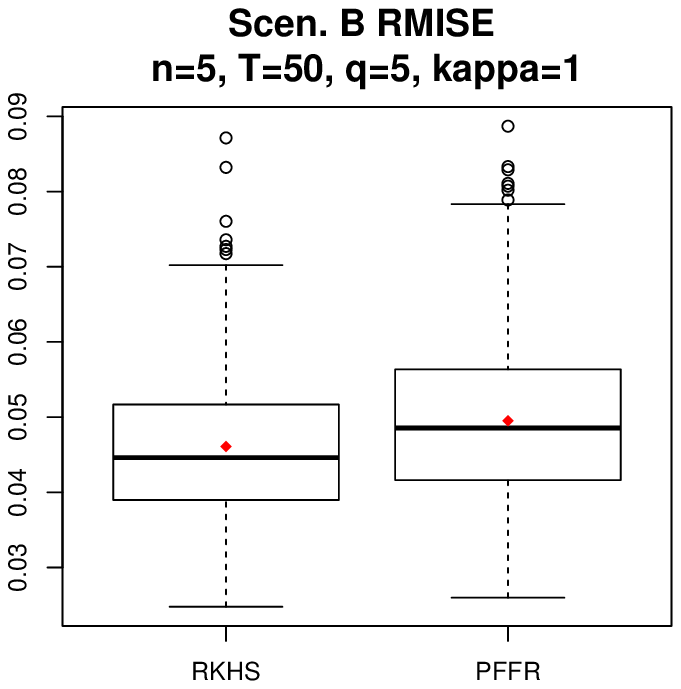}
		\vspace{-0.8cm}
    \end{subfigure}
    ~
    \begin{subfigure}{0.32\textwidth}
    	\includegraphics[angle=270, width=1.2\textwidth]{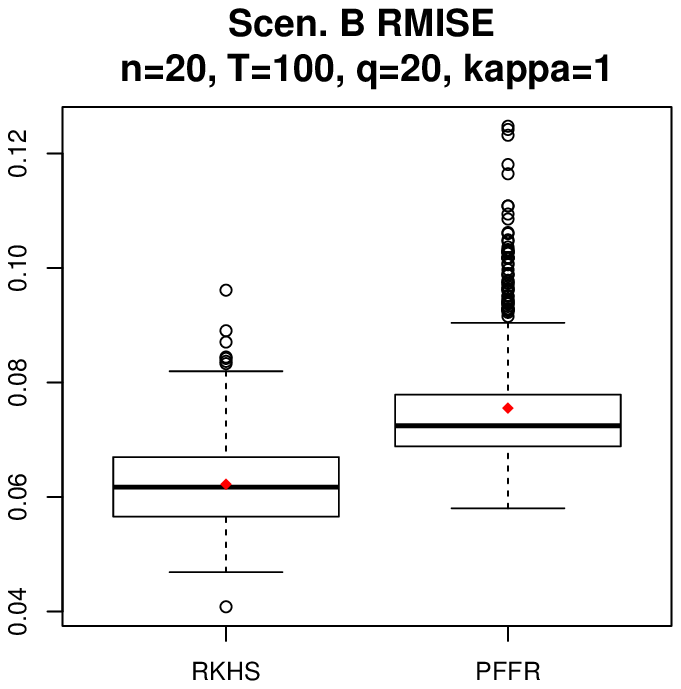}
		\vspace{-0.8cm}
    \end{subfigure}
    ~
    \begin{subfigure}{0.32\textwidth}
    	\includegraphics[angle=270, width=1.2\textwidth]{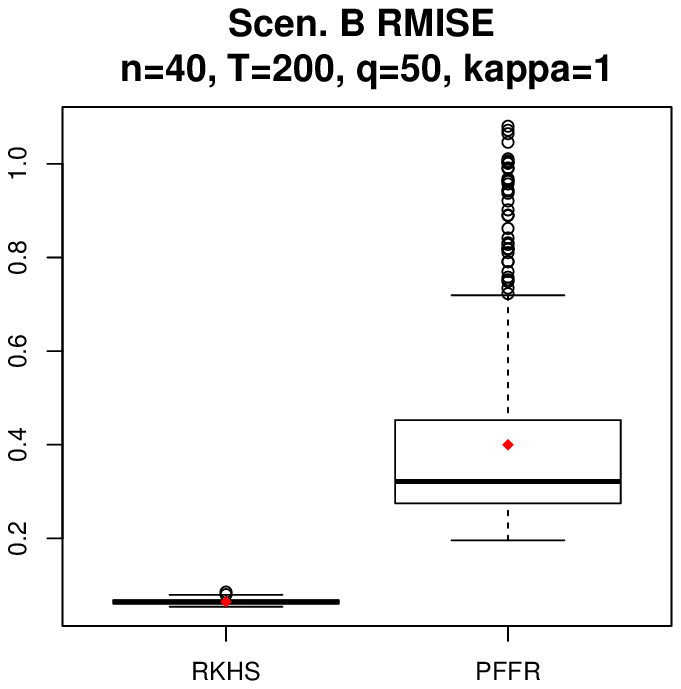}
		\vspace{-0.8cm}
    \end{subfigure}
	\caption{Boxplots of RMISE of RKHS and PFFR across 500 experiments under functional regression with mixed predictors with $\kappa=1$. Red points denote the average RMISE.}
	\label{fig:FR_mixed}
\end{figure}

\textbf{Numerical results}: Table \ref{tab:FR_mixed} reports the average nRMISE~(nRMISE$_{\mathrm{avg}}$) across 500 experiments for RKHS and PFFR under all simulation settings. Table \ref{tab:FR_mixed} further reports R$_{\mathrm{avg}}$, the percentage improvement of RKHS over PFFR, and R$_w$, the percentage of experiments (among 500 experiments) where RKHS returns lower RMISE. We further give the boxplot of RMISE in Figure \ref{fig:FR_mixed} under $\kappa=1$. The boxplots of RMISE under $\kappa=0.5,2$ can be found in \Cref{subsec:additional_simu}.  The result is consistent with the one for function-on-function regression, where RKHS delivers the best performance across almost all simulation settings with notable improvement.

\begin{table}[h]
	\centering
	\begin{tabular}{rrrrr|rrrr}
		\hline\hline
		& \multicolumn{4}{c|}{Scenario A:  $n= 5 , T= 50, q=5 $} & \multicolumn{4}{c}{Scenario B:  $n= 5 , T= 50, q=5 $} \\
		$\kappa$ & RKHS & PFFR & R$_{\mathrm{avg}}$(\%) & R$_{w}$ (\%) & RKHS & PFFR & R$_{\mathrm{avg}}$(\%) & R$_{w}$ (\%)\\
		\hline
		$0.5$ & 17.16 & 20.14 & 17.35 & 77 & 14.81 & 15.99 & 8.02 & 70 \\ 
		$1$ & 9.43 & 10.89 & 15.53 & 76 & 7.42 & 7.97 & 7.46 & 75 \\ 
		$2$ & 5.03 & 5.60 & 11.35 & 74 & 3.72 & 3.99 & 7.29 & 75 \\ 
		\hline
		& \multicolumn{4}{c|}{Scenario A:  $n= 20 , T= 100, q=20 $} & \multicolumn{4}{c}{Scenario B:  $n= 20 , T= 100, q=20 $} \\
		$\kappa$& RKHS & PFFR & R$_{\mathrm{avg}}$(\%) & R$_{w}$ (\%) & RKHS & PFFR & R$_{\mathrm{avg}}$(\%) & R$_{w}$ (\%)\\
		\hline
		$0.5$ & 11.25 & 21.95 & 95.11 & 100 & 23.95 & 25.74 & 7.48 & 63 \\ 
		$1$ & 5.95 & 11.25 & 89.31 & 100 & 11.83 & 14.44 & 22.03 & 90 \\ 
		$2$ & 3.38 & 5.94 & 75.77 & 100 & 5.92 & 7.85 & 32.47 & 99 \\ 
		\hline
		& \multicolumn{4}{c|}{Scenario A:  $n= 40 , T= 200, q=50 $} & \multicolumn{4}{c}{Scenario B:  $n= 40 , T= 200, q=50 $} \\
		$\kappa$& RKHS & PFFR & R$_{\mathrm{avg}}$(\%) & R$_{w}$ (\%) & RKHS & PFFR & R$_{\mathrm{avg}}$(\%) & R$_{w}$ (\%)\\
		\hline
		$0.5$ & 8.50 & 15.76 & 85.46 & 99 & 27.52 & 77.13 & 180.22 & 100 \\ 
		$1$ & 4.36 & 8.02 & 84.03 & 100 & 14.23 & 76.19 & 435.38 & 100 \\ 
		$2$ & 2.33 & 4.16 & 78.36 & 100 & 8.20 & 75.99 & 827.09 & 100 \\
		\hline\hline
	\end{tabular}
	\caption{Numerical performance of RKHS and PFFR under functional regression with mixed predictors. The reported nRMISE$_{\mathrm{avg}}$ is multiplied by 100 in scale. R$_{\mathrm{avg}}$ reflects the percent improvement of RKHS over PFFR and R$_{w}$ reflects the percentage of experiments in which RKHS achieves the lower RMISE.} 
	\label{tab:FR_mixed}
\end{table}

In the following, we further examine the performance of RKHS for high-dimensional scalar predictors. Specifically, keep the simulation setting identical as above, we increase the number of scalar predictors to $p=10,50,100,200$. As a reminder, the sparsity of $\beta^*=(\beta_1^*,\beta_2^*,\ldots,\beta_p^*)$ is 1 as only $\beta_1^*$ is a non-zero function. In other words, the additionally introduced scalar predictors $(Z_{t2},Z_{t3},\ldots, Z_{tp})$ are purely noise. Thus, compared to the case of $p=3$, the performance of RKHS and PFFR are expected to worsen with an increasing dimension $p.$ However, thanks to the group Lasso-type penalty, which induces sparsity of the estimated coefficient functions $\widehat{\beta}$, we expect the excess risk of RKHS to grow at the rate of $O\big(\log (p)\big)$, as suggested in Theorem \ref{thm-main-mixed}.
	
\begin{figure}[h]
	\hspace{-1cm}
	\centering
	\begin{subfigure}{0.32\textwidth}
		\includegraphics[angle=270, width=1.2\textwidth]{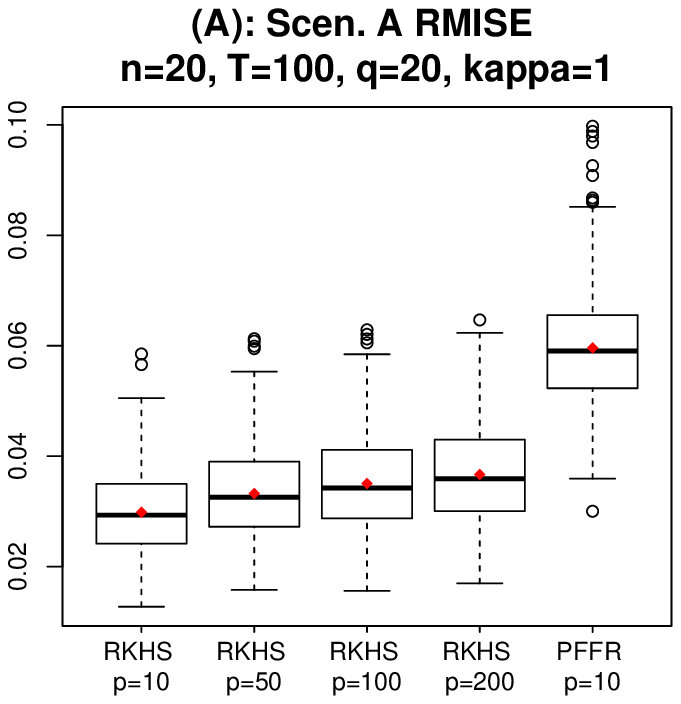}
		\vspace{-0.3cm}
	\end{subfigure}
	~
	\begin{subfigure}{0.32\textwidth}
		\includegraphics[angle=270, width=1.2\textwidth]{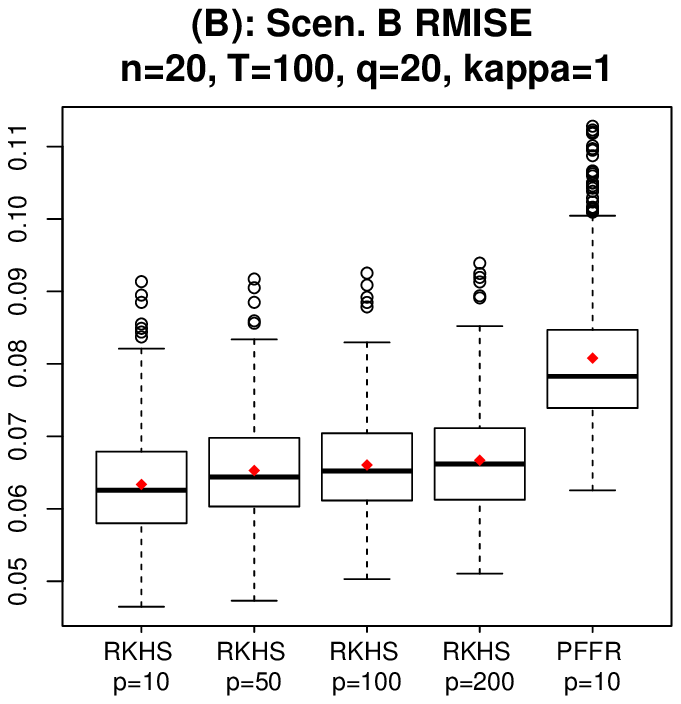}
		\vspace{-0.3cm}
	\end{subfigure}
	~
	\begin{subfigure}{0.32\textwidth}
		\includegraphics[angle=270, width=1.2\textwidth]{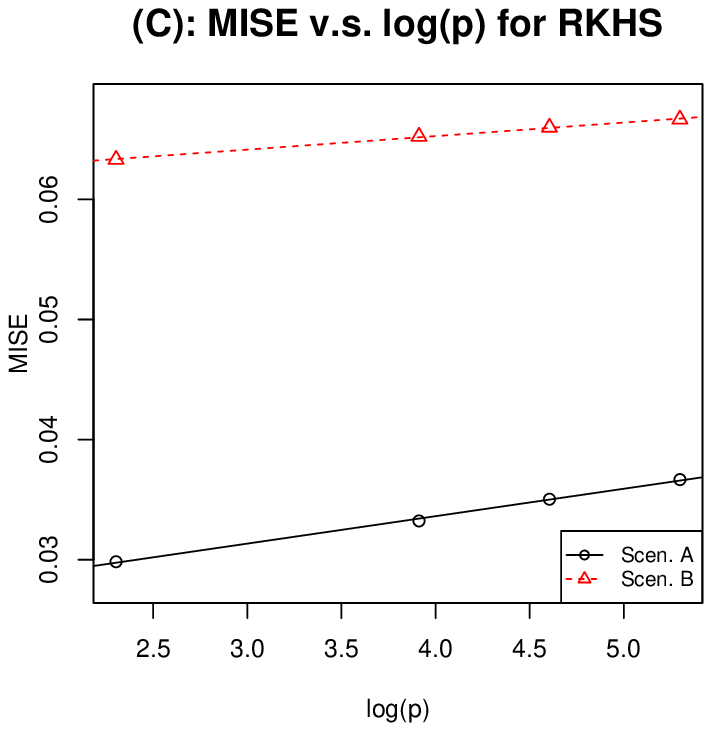}
		\vspace{-0.3cm}
	\end{subfigure}
	\caption{Boxplots of RMISE of RKHS and PFFR across 500 experiments under functional regression with mixed predictors with $\kappa=1$ for $p=10,50,100,200$. Red points denote the average RMISE. For plot (C), MISE = RMISE$^2$ and the two lines are fitted via OLS.}
	\label{fig:FR_mixed_HD}
\end{figure}

To conserve space, we present the result under the simulation setting with $(n,T,q)=(20,100,20)$ and $\kappa=1.$ Results under other settings are similar and thus omitted. Due to the lack of sparsity penalty, PFFR may not be suitable for the setting of high-dimensional predictors. For comparison, we only implement PFFR for $p=10$. Figure \ref{fig:FR_mixed_HD} (A) and (B) give the boxplot of RMISE for RKHS and PFFR across 500 experiments. As expected, the RMISE of RKHS increases as the dimension $p$ increases though at a rate slower than $p$. Indeed, RKHS at $p=200$ still gives better performance than PFFR at $p=10$. Figure \ref{fig:FR_mixed_HD} (C) gives the plot between average MISE\footnote{To match the result in Theorem \ref{thm-main-mixed}, we compute MISE, which is the squared RMISE with MISE = RMISE$^2$.} across 500 experiments and $\log(p)$, where the relationship is seen to be roughly linear. This confirms that the excess risk of RKHS increases in the order of $O\big(\log(p)\big)$ as suggested in Theorem \ref{thm-main-mixed}.

\subsection{Real data applications}\label{subsec:realdata}
In this section, we conduct real data analysis to further demonstrate the promising utility of the proposed RKHS-based functional regression in the context of crowdfunding. In recent years, crowdfunding has become a flexible and cost-effective financing channel for start-ups and entrepreneurs, which helps expedite product development and diffusion of innovations.

We consider a novel dataset collected from one of the largest crowdfunding websites, kickstarter.com, which provides an online platform for creators, e.g.~start-ups, to launch fundraising campaigns for developing a new product such as electronic devices and card games. The fundraising campaign takes place on a webpage set up by the creators, where information of the new product is provided, and typically has a 30-day duration with a goal amount $G$ preset by the creators. Over the course of the 30 days, backers can pledge funds to the campaign, resulting in a functional curve of pledged funds $\{P(r), r\in [0,30]\}$. At time $r$, the webpage displays real-time $P(r)$ along with other information of the campaign, such as its number of creators $Z_1$, number of product updates $Z_2$ and number of product FAQs $Z_3$.

A fundraising campaign succeeds if $P(30)\geq G$. Importantly, only creators of successful campaigns can be awarded the final raised funds $P(30)$ and the platform kickstarter.com charges a service fee~($5\%\times P(30)$) of successful campaigns. Thus, for both the platform and the creators, an accurate prediction of the future curve $\{P(r), r\in (s, 30]\}$ at an early time $s$ is valuable, as it not only reveals whether the campaign will succeed but more importantly suggests timing along $(s,30]$ for potential intervention by the creators and the platform to boost the fundraising campaign and achieve greater revenue.

The dataset consists of $T=454$ campaigns launched between Dec-01-2019 and Dec-31-2019. For each campaign $t=1, \ldots, T$, we observe its normalized curve $N_t(r)=P_t(r)/G_t$\footnote{Note that the normalized curve $N_t(r)$ is as sufficient as the original curve $P_t(r)$ for monitoring the fundraising process of each campaign.} at 60 evenly spaced time points over its 30-day duration and denote it as $\{N_t(r_i)\}_{i=1}^{60}$. See Figure \ref{fig:representative_plot} for normalized curves of six representative campaigns. At a time $s\in(0,30)$, to predict $\{N_t(r), r\in (s, 30]\}$ for campaign $t$, we employ functional regression, where we treat $\{N_t(r_i), r_i\in (s, 30]\}$ as the functional response $Y_t$, use $\{N_t(r_i), r_i\in [0, s]\}$ as the functional covariate $X_t$ and $(Z_1,Z_2,Z_3)$ as the vector covariate. We compare the performances of RKHS, FDA and PFFR. Note that FDA only allows for one functional covariate, thus for RKHS and PFFR, we implement both the function-on-function regression and the functional regression with mixed predictors~(denoted by RKHS$_{\text{mixed}}$ and PFFR$_{\text{mixed}}$). The implementation of each method is the same as that in Sections~\ref{subsec:ffr} and \ref{subsec:mfr}.

We vary the prediction time $s$ such that $s=7\text{th}, 8\text{th},\ldots, 20\text{th}$ day of a campaign. Note that we stop at $s=20$ as prediction made in the late stage of a campaign is not as useful as early forecasts. To assess the out-of-sample performance of each method, we use a 2-fold CV, where we partition the 454 campaigns into two equal-sized sets and use one set to train the functional regression and the other to test the prediction performance, and then switch the role of the two sets. For each campaign $t$ in the test set, given its prediction $\{\widehat{N}_t(r), r\in (s,30]\}$ generated at time $s$, we calculate its RMSE and MAE with respect to the true value $\{N_t(r_i), r_i\in (s, 30]\}$, where
	\[
		\text{RMSE}_{t,s}=\sqrt{\frac{1}{\# \{r_i\in (s,30]\} }\sum_{r_i\in (s,30]} (\widehat{N}_t(r_i)-N_t(r_i))^2}
	\]
	and
	\[
	\text{MAE}_{t,s}=\frac{1}{\# \{r_i\in (s,30]\} }\sum_{r_i\in (s,30]} \left|\widehat{N}_t(r_i)-N_t(r_i)\right|.
	\]

\begin{figure}[h]
	\hspace*{-1cm}         
	\centering                                                  
	\begin{subfigure}{0.4\textwidth}
		\includegraphics[angle=270, width=1.2\textwidth]{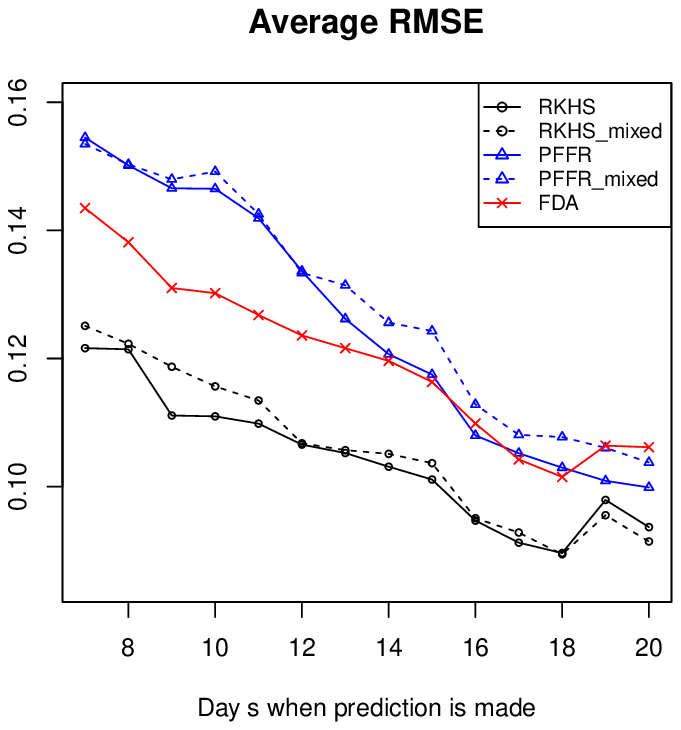}
		\vspace{-0.3cm}
	\end{subfigure}
	~                                                          
	\begin{subfigure}{0.4\textwidth}
		\includegraphics[angle=270, width=1.2\textwidth]{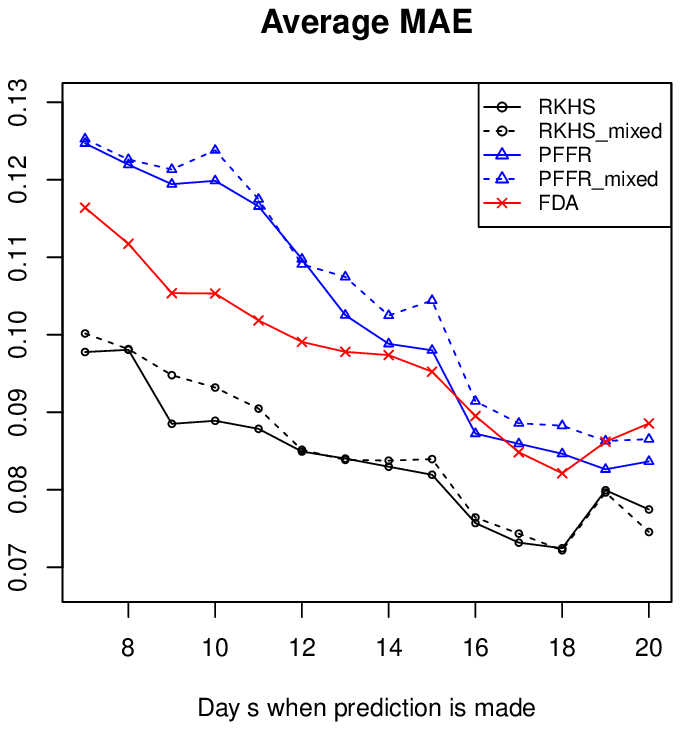}
		\vspace{-0.3cm}
	\end{subfigure}
	\caption{Average RMSE and MAE achieved by different functional regression methods. }
	\label{fig:realdata}
\end{figure}

Figure \ref{fig:realdata} visualizes RMSE$_s=T^{-1}\sum_{t=1}^T\text{RMSE}_{t,s}$ and MAE$_s=T^{-1}\sum_{t=1}^T\text{MAE}_{t,s}$  achieved by different functional regression methods across $s \in \{7,8,\ldots,20\}.$ In general, the two RKHS-based estimators consistently achieve the best prediction accuracy. As expected, the performance of all methods improve with $s$ approaching 20. Interestingly, the functional regression with mixed predictors does not seem to improve the prediction performance, which is especially evident for PFFR. Thanks to the group Lasso-type penalty, RKHS can perform variable selection on the vector covariate. Indeed, among the 28~(2 folds $\times$ 14 days) estimated functional regression models based on RKHS$_{\text{mixed}}$, 19 models select no vector covariate and thus reduce to the function-on-function regression, suggesting the potential irrelevance of the vector covariate.  We further provide a robustness check of the above analysis in \Cref{subsec:realdata_cv}, where we repeat the 2-fold CV procedure 100 times for RKHS, FDA and PFFR. It is seen that RKHS consistently provides the best performance, confirming the robustness of our findings. We refer to \Cref{subsec:realdata_cv}  for more details.

\begin{figure}[h]
	\hspace*{-0.9cm}                                                           
	\includegraphics[angle=270, width=1.1\textwidth]{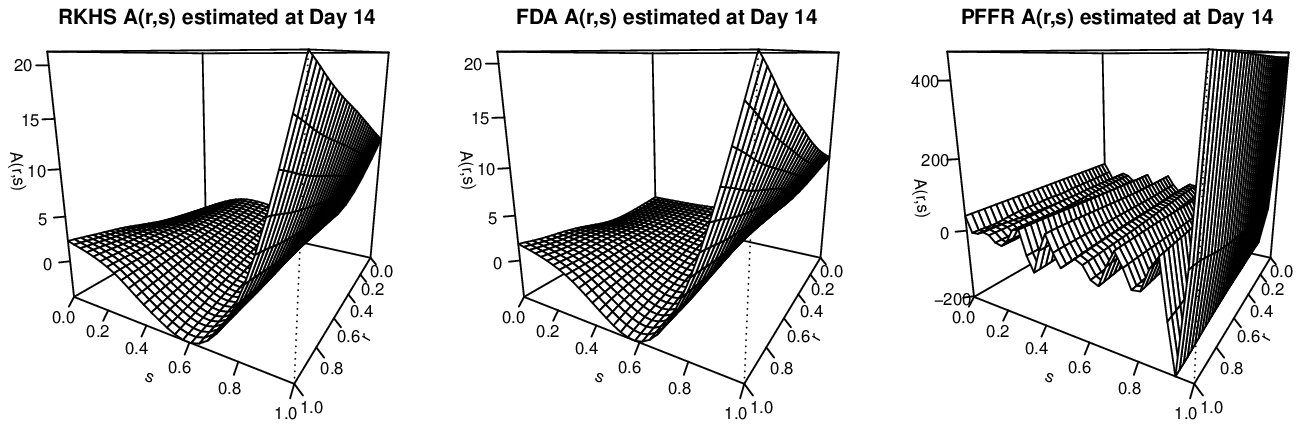}
	\vspace{-0.3cm}
	\caption{Estimated bivariate coefficient functions $\widehat{A}(r,s)$ by RKHS, FDA and PFFR at Day 14 ($(r,s)$ is rescaled to $[0,1]^2$).}
	\label{fig:A_est}
\end{figure}

For more intuition, Figure \ref{fig:representative_plot} plots the normalized fundraising curves $\{N_t(r_i)\}_{i=1}^{60}$ of six representative campaigns and further visualizes the functional predictions given by RKHS, FDA and PFFR at $s=14$th day. Note that FDA and RKHS provide more similar prediction while the prediction of PFFR seems to be more volatile. This is indeed consistent with the estimated bivariate coefficient functions visualized in Figure \ref{fig:A_est}, where $\widehat{A}(r,s)$ of RKHS and FDA resembles each other while PFFR seems to suffer from under-smoothing.

Figure \ref{fig:realdata_witharima} in \Cref{subsec:additional_simu} further gives the prediction performance of the classical time series model ARIMA, which is significantly worse than the predictions given by RKHS, FDA and PFFR, suggesting the advantage of functional regression for handling the current application.

%% file: conclusion.tex
In this paper, we study a functional linear regression model with functional responses and accommodating both functional and high-dimensional vector covariates. We provide a minimax lower bound on its excess risk. To match the lower bound, we propose an RKHS-based penalized least squares estimator, which is based on a generalization of the representer lemma and achieves the optimal upper bound on the excess risk. Our framework allows for partially observed functional variables and provides finite sample probability bounds. Furthermore, the result unveils an interesting phase transition between a high-dimensional parametric rate and a nonparametric rate in terms of the excess risk. A novel iterative coordinate descent based algorithm is proposed to efficiently solve the penalized least squares problem. Simulation studies and real data applications are further conducted, where the proposed estimator is seen to provide favorable performance compared to existing methods in the literature.

Throughout the paper, we assume knowledge of the kernel functions $\mathbb{K}$ and $\mathbb{K}_{\beta}$.  In practice, it would be ideal if one would be able to learn the kernel functions from data.  However, to our best knowledge, assuming knowing the true kernel functions is adopted in all RKHS-based functional data analysis literature.  Learning kernels is beyond the scope of our paper and we would like to pursue this direction in the future research.  We would also like to point out that our estimator $(\widehat{A}, \widehat{\beta})$ is a constrained estimator, with the constraints related with the true kernels.  A mis-specification might lead to an over/under penalization, which will affect the final prediction error rate.

Finally, we briefly discuss important distinctions between non-parametric regression and functional regression for interested readers. To make the comparison easier, we consider a simpler setting: the scalar response functional linear regression and the classical non-parametric regression.  

In non-parametric regression, we have
		\[
			y_i = f(x_i) +\epsilon_i, 
		\]
		where $ \{ x_i\}_{i=1}^n    $ are   assumed to be random or fixed grid points in $[0,1]$.  In scalar response functional linear regression   with fully observed functional data, we have 
		\begin{equation}\label{eq:classical functional regression}
			y_i = \int X_i(s) \beta(s) \dint s +\epsilon_i  ,
		\end{equation}
		where $ \{ X_i\}_{i=1}^n$ are  assumed to be i.i.d.~stochastic processes in $\mathcal L^2$. Note that we can treat $X_i$ and $\beta$ as infinite-dimensional vectors w.r.t.\ some $\mathcal L^2$ basis system, thus model \eqref{eq:classical functional regression} can be rewritten as  
	    \begin{equation*}
			y_i =    \widetilde{X}_i^\top \widetilde{\beta} +\epsilon_i ,
		\end{equation*}		
		where $ \widetilde X_i\in \mathbb R^\infty ,\widetilde \beta \in \mathbb R^\infty$. Under the standard assumption that  $\mathbb E(\|X_i\|_{\mathcal L^2}^2) <\infty  $,  
	the  eigenvalue sequence of the covariance matrix $\widetilde \Sigma_X =\mathbb E(\widetilde X_i\widetilde X_i^\top )   $  converges  to $0$. As such, model \eqref{eq:classical functional regression}  was described  as a high-dimensional or infinitely-dimensional ``ill-posed'' linear regression model in \cite{hall2007methodology}, which is fundamentally different from the standard non-parametric regression.
 
\section*{Acknowledgments}
We would like to thank the action editor, Dr. Pradeep Ravikumar, as well as the four anonymous reviewers for their thoughtful assessment and constructive comments which helped us to improve the quality and the presentation of our paper.   Z.Z. would like to acknowledge support from NSF DMS-2014053.  Y.Y. would like to acknowledge support from DMS-EPSRC EP/V013432/1 and EPSRC EP/W003716/1.  
R.W.  would like to acknowledge support from AFOSR FA9550‐18‐1‐0166,
DOE DE‐AC02‐06CH113575, NSF DMS‐1925101, ARO W911NF-17-1-0357, and
NGA HM0476-17-1-2003.

%% file: appendix.tex
\section{More discussions on bivariate functions, Assumptions~\ref{assume:mixed} and \ref{assume:joint smooth}}\label{sec-app-disc-assump}

\subsection{Bivariate functions and compact linear operators}\label{sec-bivfunc-linmap}
 
\begin{remark} \label{remark:bivariate functions}
Fro any  compact linear operator $A_2: \hk \to \hk$   denote 
\begin{align*}
	A_2[f,  g] = \langle A_2[g], f\rangle_\hk, \quad f, g \in \hk . 
\end{align*}
Note that $A_2[f, g]$ is well defined for any $f, g \in \hk $ due to the compactness of $A_2$. Define $a_{ij} = A_2[\psi_i, \psi_j] = \lb A_2[\psi_j], \psi_i\rb_{\mathcal{H}(\mathbb{K})}$, $i, j \in \mathbb{N}_*$.  We thus have for any $f, g \in \hk $, it holds that
\begin{align*}
	A_2[f, g] & = \lb A_2[g], f\rb_\hk = \sum_{i = 1}^{\infty} \langle f, \psi_i\rangle_{\hk } \langle A_2[g], \psi_i\rangle_{\hk }\nonumber\\
	& = \sum_{i = 1}^{\infty} \langle f, \psi_i\rangle_{\hk } \Bigg\langle A_2\left[\sum_{j = 1}^{\infty} \langle g, \psi_j\rangle_{\hk } \psi_j\right], \psi_i \Bigg\rangle_{\hk } \nonumber \\
	& = \sum_{i, j = 1}^{\infty} \langle f, \psi_i\rangle_{\hk } \langle g, \psi_j\rangle_{\hk } \lb A_2[\psi_j], \psi_i\rb_\hk = \sum_{i, j = 1}^{\infty} a_{ij} \langle f, \psi_i\rangle_{\hk } \langle g, \psi_j\rangle_{\hk }, 
\end{align*}
where the fourth identity follows from the compactness of $A_2.$ This justifies \Cref{eq-comp}.
 
 \end{remark}

\subsection{\Cref{assume:mixed}(d)}\label{sec-appendix-A.2}

In \Cref{assume:mixed}(d), we assume that
	\[
		\mathbb{E}( \lb X_1, f \rb_\lt Z^\top _1 v) \le \frac{3}{4}\sqrt { \mathbb{E}\{\langle X_1, f\rangle_{\mathcal{L}^2}^2\} (v^\top  \Sigma_Zv)},
	\]
	i.e.~the functional and vector covariates are allowed to be correlated up to $O(1)$.  In this subsection, we emphasize that the correlation cannot be equal to one.  We first note that 
	\begin{equation}\label{eq-holder-inequ}
		\mathbb{E}( \lb X_1, f \rb_\lt Z^\top _1 v) \le \sqrt { \mathbb{E}\{\langle X_1, f\rangle_{\mathcal{L}^2}^2\} (v^\top  \Sigma_Zv)}.
	\end{equation}

    For model identification, we require the inequality in \eqref{eq-holder-inequ} to be strict.  To see this, we suppose that there exist $g \in \hk$ and $u \in \mathbb{R}^p$ such that 
	\[
		\mathbb{E}( \lb X_1, g \rb_\lt Z^\top _1 u) =  \sqrt { \mathbb{E}\{\langle X_1, g\rangle_{\mathcal{L}^2}^2\} (u^\top  \Sigma_Zu)}.
	\]
	Since $\lb X_1, g\rb_\lt $ and  $Z^\top _1 u$ are both normal random variables, the above equality implies that the correlation between $\lb X_1, g\rb_\lt $ and $Z^\top _1 u$ is 1 and thus $\lb X_1, g\rb_\lt  = a Z^\top _1 u$ for some constant $a \in \mathbb{R}$.  This means that the  model defined  in  \eqref{eq-model-mixed-cov} is no longer identifiable,  because  the covariates are perfectly correlated.  Even in finite-dimensional regression problems, if the covariates are perfectly correlated, the solution to the linear system is ill-defined. 

\subsection{\Cref{assume:joint smooth}(b)}
Recall that \Cref{assume:joint smooth} is required for handling the case where the functional observations are on discrete sample points.  \Cref{assume:joint smooth}(a) formalizes the sampling scheme and \Cref{assume:joint smooth}(b) requires that the second moment of  the random variable $\|X\|_{\wot }$ is finite. 

In the following, via a concrete example, we show that the  second moment assumption in  \eqref{eq:X almost surely} holds under mild conditions. Specifically in \Cref{example:joint smooth}, we provide a simple sufficient condition on the eigen-decay of the covariance operator $\Sigma_X$, under which \eqref{eq:X almost surely} holds for $W^{\alpha,2}$ with $\alpha>1/2$.
 
\begin{example}\label{example:joint smooth}
Let $\alpha>1/2$ and  suppose   the Sobolev space $\wot$ is generated by the kernel 
$$\mathbb K_{\alpha }(r,s)  = \sum_{k=1}^\infty \omega_k \psi_k ^\alpha (r) \psi_k^\alpha (s)    , $$
where due to the property of the Sobolev space, we have $\omega_k \asymp  k^{-2\alpha } $ and  $\{ \psi_k^\alpha \}_{k=1}^\infty  $ can be taken as the Fourier basis such that $ \| \psi_k^\alpha \|_\wot ^2   \asymp k^{2\alpha}$. To better understand the second moment condition  
$$ \mathbb E (\| X_1\|_\wot ^2)<\infty ,$$
we proceed by following the same strategy as in \cite{cai2006prediction} and \cite{yuan2010reproducing} and  assume   that  $\mathbb K_\alpha$ and $\Sigma_X $ are perfectly aligned, i.e., they share the same set of eigen-functions. 

Let $z_k = \lb X_1, \psi_k^\alpha \rb_\lt  $. We have that $\{z_k\}_{k=1}^\infty $ is a collection of independent Gaussian random variables such that $z_k \sim N(0,\sigma_k^2) $, where $ \{ \sigma_k^2\}_{k=1}^\infty $ are the eigenvalues of $\Sigma_X$. Thus we have  
\begin{align*} \| X_1\|_\wot^2  =\sum_{k=1}^\infty z_k^2 \| \psi_k^\alpha \|_\wot ^2   \asymp \sum_{k=1}^\infty z_k ^2  k^{2\alpha},
\end{align*}
and
$$ \mathbb E( \| X_1\|_\wot^2 ) \asymp \sum_{k=1}^\infty \mathbb E (z_k ^2)   k^{2\alpha}  =   \sum_{k=1}^\infty \sigma_k^2 k^{2\alpha}.$$
To assure that $\mathbb E( \| X_1\|_\wot^2 ) < \infty $, it suffices to have that 
$$ \sigma_k^2 \asymp k^{-2\alpha -1 -\nu} $$ for any constant $\nu>0$.  
Note that when $\alpha$ is close to $1/2$, this means that the eigenvalues of $\Sigma_X$ decay at a polynomial rate slightly faster than 2.
\end{example}

\section{Proofs of the results in \Cref{sec-main-results}} \label{sec-app-proof-sec3}

\subsection{Proof of \Cref{lemma:representer 3} }  \label{appendix representer}

\begin{proof} [Proof of \Cref{lemma:representer 3}]  

For any linear subspaces $R, S \subset \hk$, let $\mathcal{P}_R$ and $\mathcal{P}_S$ be the projection mappings of the spaces $R$ and $S$, respectively with respect to $\|\cdot\|_\hk $.  For any $f, g \in \hk $ and any compact linear operator $A: \hk  \to \hk $, denote
	\[
		A\vert _{R \times S}[f, g] = A[\mathcal P_{R}f, \mathcal P_{S }g].
	\]

Let $(\widehat B, \widehat \alpha)$ be any solution to \eqref{eq:approx A beta}.  Let $Q = \text{span} \{\mk_\beta  (r_j,  \cdot)\}_{j = 1}^n \subset \hb $, $R = \text{span} \{\mathbb K (r_j,  \cdot)\}_{j = 1}^n \subset \hk $ and $S = \text{span} \{\mathbb K (s_j,  \cdot )\}_{j = 1}^n \subset \hk $.  Denote $\widehat \beta_l = \mathcal {P}_{Q  }  (\widehat  \alpha_ l )  $, $l \in \{1, \ldots, p\}$, and $ \widehat A [\cdot, \cdot] = \widehat B|_{R \times S}[\cdot, \cdot] = \widehat{B}[\mathcal{P}_R \cdot, \mathcal{P}_S \cdot]$.  

Let $S^\perp$ and $R^{\perp}$ be the orthogonal complements of $S$ and $R$ in $\h$, respectively.  Then for any compact linear operator $A$, we have the decomposition 
	\[
		A =   A\vert _{R  \times S  }  + A\vert _{R  \times S^\perp  } +  A\vert _{R^\perp \times S  } +  A\vert _{R^\perp  \times S ^\perp  }.
	\]
	Due to \Cref{lemma:constrained operator}, there exist $\{a_{ij}\}_{i, j = 1}^n \subset \mathbb{R}$ such that 
	\[
		A\vert  _{R   \times S  }  [f, g] = \sum_{i, j = 1}^n  a_{ij } \lb \mk (  r _i ,  \cdot   ), f \rb_\hk \lb \mk (s_j ,  \cdot   ),g \rb_\hk.
	\] 
\
\\
 {\bf Step 1.}  In this step, it is shown that for any compact linear operator $A$, its associated bivariate function $A(\cdot, \cdot)$ satisfies that $\{A(r_i, s_j)\}_{i, j=1}^n$ only depend on $A\vert _{R \times S}$.  The details of the bivariate function is explained in \Cref{sec-bivfunc-linmap}.
\
\\
Observe that 
	\[
		A\vert_{ R \times S^\perp  } (r_i , s _j)=     A\vert_{R \times S^\perp  }[ \mathbb K( r _i, \cdot),  \mathbb K(s_j, \cdot)]  = A[  \mathcal P_R  \mathbb K(r _i, \cdot) ,  \mathcal P_{S^\perp }  \mathbb K(s_j, \cdot) ] = 0.  
	\]
	Similar arguments also lead to that $A\vert_{R^\perp\times S  }[ \mathbb K(r _i, \cdot) , \mathbb K(s_j, \cdot) ]  = 0$ and $   A\vert_{R ^\perp\times S^\perp   }[ \mathbb K( r _i, \cdot) , \mathbb K(s_j, \cdot) ]  = 0 $, $i, j = 1, \ldots, n$.

\medskip
\noindent {\bf Step 2.}  By {\bf Step 1}, it holds that $\widehat A   (r _i,s_j) =\widehat B (r _i,s_j)$, $i, j = 1, \ldots, n$.  By \Cref{lemma:constrained operator}, we have that there exist $\{\widehat{a}_{ij}\}_{i, j = 1}^n \subset \mathbb{R}$ such that for any $f, g \in \hk$, 
	\[
		\widehat {    A }   [f ,g]   = \sum_{i, j = 1}^n \widehat a_{ ij } \lb \mk (r _i ,  \cdot   ), f \rb_\hk \lb \mk (s_j ,  \cdot ),g \rb_\hk .
	\]
	Therefore the associated bivariate function satisfies that, for all $r, s \in [0, 1]$,
	\begin{align*}
		\widehat{A}(r, s) & = \sum_{i, j = 1}^{\infty} \widehat{A}[\psi_i, \psi_j] \psi_i(r)\psi_j(s) = \sum_{i, j = 1}^{\infty} \sum_{k, l = 1}^n \widehat{a}_{kl} \psi_i(r_k) \psi_j(s_l) \psi_i(r)\psi_j(s) \\
		& = \sum_{k, l = 1}^n \widehat{a}_{kl} \left\{\sum_{i = 1}^\infty \psi_i(r)\psi_i(r_k)\right\} \left\{\sum_{i = 1}^\infty \psi_i(s)\psi_i(s_l)\right\} = \sum_{k, l = 1}^n \widehat{a}_{kl} \mathbb{K}(r, r_k) \mathbb{K}(s, s_l),
	\end{align*}
	where $K(r,s) = \sum_{i=1}^\infty \psi(s) \psi(r) $ is used in the last inequality. 

\medskip
\noindent {\bf Step 3.} For any $j \in \{1, \ldots, n\}$ and $l \in \{1, \ldots, p\}$, by the definition of  $\widehat \beta_l  $,
we have that  
	\[
		\widehat{\alpha}_l (r_j) = \lb \widehat \alpha _  l , \mathbb K_ \beta  (r_ j , \cdot) \rb_\h =\lb  \widehat \beta _  l , \mathbb K_\beta  (r_j,  \cdot) \rb_\h  = \widehat \beta_ l (r_ j ). 
	\]
 \
 \\
 Therefore 
	\begin{align*}  
		& \frac{1}{Tn} \sum_{t= 1}^T    \sum_{i=1}^ n  \wri   \left\{Y_{t} (r _i) -    \frac{1}{n }\sum_{j=1} ^n \wsj  \widehat A   (r _i, s_j)  X_{t  } (s_j)   - \lb  \widehat \beta  (r_i),  Z_t  \rb_p \right\}^2  + \lambda \sum_{ l =1}^p \|\widehat \beta_ l \|_n  \\
		=& \frac{1}{Tn} \sum_{t= 1}^T    \sum_{i=1}^ n   \wri   \left\{Y_{t} (r _i) -    \frac{1}{n }\sum_{j=1} ^n \wsj   \widehat B   (r _i, s_j)  X_{t  } (s_j)   - \lb  \widehat \alpha  (r_i),  Z_t  \rb_p \right\}^2 + \lambda \sum_{ l =1}^p \|\widehat \alpha_ l \|_n   .
	\end{align*}
\
\\
In addition, by \Cref{lemma:constrained operator} we have that 
	\begin{align*} 
		\| \widehat A\|_\fk \le \|\widehat B\|_\fk  \quad  \mbox{and} \quad  \| \widehat \beta _ l \|  _{\hb } \le \| \widehat \alpha  _l \|  _{\hb } \quad \text{ for all } l \in \{ 1,\ldots, p\} ,
	\end{align*} 
	which completes the proof.
\end{proof}

\begin{lemma}  \label{lemma:constrained operator}
Let $A: \hk  \to \hk $ be any Hilbert\textendash Schmidt operator.  Let $\{v_i, v'_i\}_{i = 1}^m \subset \hk $ be a collection of functions in $\hk $.  Let $R$ and $S$ be the subspaces of $\hk $ spanned by $\{v_i\}_{i = 1}^m$ and $\{v'_i\}_{i = 1}^m$, respectively.  Let $\mathcal{P}_R$ and $\mathcal{P}_S$ be the projection mappings of the spaces $R$ and $S$, respectively.  Then there exist $\{a_{ij}\}_{i, j=1}^m \subset \mathbb{R}$, which are not necessarily unique, such that for any $f, g \in \hk $,
	\[
		A\vert_{R \times S}[f, g] = A[\mathcal{P}_R f, \mathcal{P}_S g] = \sum_{i, j = 1}^m a_{ij} \langle f, v_i \rangle_{\hk } \langle g, v_i \rangle_{\hk },
	\]
	and that $\|A\vert_{R \times S}[f, g]\|_{\fk } \leq \|A\|_{\fk }$.
\end{lemma}

\begin{proof} 
Let $\{u_k\}_{k = 1}^K$ and $\{u'_l\}_{l =1}^L$ be orthogonal basis of subspaces $R$ and $S$ of $\mathcal{L}^2$, respectively, with $K, L \leq m$.  Since each $u_i$($u'_i$) can be written as a linear combination of $\{v_k\}_{k = 1}^m$($\{v'_k\}_{k = 1}^m$), it suffices to show that there exist $\{b_{ij}\}_{i = 1, j = 1}^{K, L} \subset \mathbb{R}$, such that for any $f, g \in \hk $,
	\[
		A\vert_{R\times S} [f, g] = \sum_{k = 1}^K \sum_{l = 1}^L b_{ij}\lb u_i, f\rb_{\hk  }\lb u_j', g\rb_{\hk  }.
	\]

Since $R$ and $S$ are linear subspaces of $\hk $, there exist $\{u_{k}\}_{k = K +1}^\infty, \{u'_l\}_{l = L + 1}^\infty \subset \hk $, such that $\{u_k\}_{k=1}^\infty$ and $\{u_l'\}_{l = 1}^\infty$ are two orthogonal basis of $\hk $, and that 
	\[
		A[f, g] = \sum_{i, j = 1}^\infty b_{ij}\lb u_i, f\rb_{\hk  }\lb u_j' , g\rb_{\hk },
	\]
	where $b_{ij} = A[u_i, u_j'] = \lb A[u_j'], u_i \rb_\hk $. Therefore
	\begin{align*} 
		A\vert_{R\times S} [f, g] = & \sum_{i, j=1}^\infty  b_{ij}\lb u_i, \mathcal P_R  f\rb_{\hk  }\lb u_j ', \mathcal P_S g\rb_{\hk  }  = \sum_{i,j=1}^\infty  b_{ij }\lb \mathcal P_R  u_i,  f\rb_{\hk  }\lb \mathcal P_S u_j' ,  g\rb_{\hk  } \\
		= & \sum_{k = 1}^K \sum_{l = 1}^L b_{kl}\lb u_k, f\rb_{\hk  }\lb u_l' , g\rb_{\hk  }.
	\end{align*}
	Moreover, we have that
	\[
		\|A\|_\fk ^2 = \sum_{k, l = 1}^{\infty} b_{kl}^2 \ge \sum_{k = 1}^K \sum_{l = 1}^L b_{kl}^2  = \|A\vert_{R \times S}\|_\fk ^2,
	\]
	which concludes the proof.
\end{proof}

\subsection{Proof of \Cref{thm-main-mixed}}
\label{subsection:proof of main}

\begin{proof}[Proof of \Cref{thm-main-mixed}] 

First note that $0 < \delta_T, \zeta_n < 1$, then $\delta_T^2 \le \delta_T $ and $\zeta_n^2 \le \zeta_n $.  These inequalities will be used repeatedly in the rest of this proof.  Recall the notation that $\Delta_\beta (r) = \widehat \beta(r) -\beta^*(r)$ and $\Delta_A(r,s) = \widehat A(r,s) - A^*(r,s)$, $r, s \in [0, 1]$.  Note that for any $ g\in \lt $, it holds that
	\[
		\mathbb{E}(\lb g, X^*\rb_\lt^2 ) =  \iint_{[0, 1]^2} g(r)\Sigma_X(r,s) g(s) \dint r \dint  s.
	\]
	The above expression will also be used repeatedly in the rest of the proof.  
	Observe that 
for any $r\in[0,1]$, it holds that 
$$ \| \Delta_A(r,\cdot) \|_{W^{\alpha, 2}} \le \| \Delta(r, \cdot) \|_{\hk} \le C_{\mathbb K}\| \Delta\|_{\fk} \le C_{\mathbb K} C_A  ,$$	
	where the first inequality follows from \Cref{assume:joint smooth} {\bf b}, which assumes that $ \hk \subset W^{\alpha,2 }$, the second inequality follows from \Cref{lemma:bound of A 1} and the last inequality follows from the assumption that $A^*, \widehat A \in \mathcal C_A$.  Throughout the proof, we assume without loss of generality that $ C_\mk =1. $
In addition, observe  that   by \Cref{lemma:Gaussian process sobolev norm}, it holds that 
	$$\mathbb P \big ( \| X_t \|_\wot  \ge 4C_X \sqrt {\log(T) } \big)\le T^{-6}. $$
 	By a union bound argument,  
	$$\mathbb P ( \mathcal D ) =  \mathbb P \big ( \| X_t \|_\wot \ge 4C_X \sqrt {\log(T)}  \ \text{for all } 1\le t \le T \big)\le T^{-5} .   $$ 
The rest of the proof is shown in the event of $\mathcal D $.
	\\
	\\
Let $\lambda =C_\lambda \sqrt {\frac{\log(p\vee T)}{T }} $ for sufficiently large constant $C_\lambda $.	From the     minimizer  property, we have that
	\begin{align*}
		& \frac{1}{Tn}\sum_{t=1}^T \sum_{i=1}^n \wri  \left\{Y_t(r_i) -  \frac{1}{n }   \sum_{j=1}^n \wsj   \widehat A(r_i, s_j) X_t(s_j)    -  \lb Z_t ,  \widehat \beta(r_i) \rb _p \right\}^2   + \lambda \sum_{j=1}^p  \|\widehat \beta_j   \|_n \\
		\le &  \frac{1}{Tn}\sum_{t=1}^T \sum_{i=1}^n  \wri \left\{Y_t(r_i) -      \frac{1}{n }   \sum_{j=1}^n  \wsj    A^* (r_i, s_j) X_t(s_j) - \lb Z_t , \beta ^* (r_i) \rb _p \right\}^2  + \lambda \sum_{j=1}^p \|\beta_j  ^*  \|_n  ,
	\end{align*}
	which implies that 
	\begin{align} 
 		& \frac{1}{Tn}\sum_{t=1}^T \sum_{i=1}^n \wri  \left\{\frac{1}{n }   \sum_{j=1}^n \wsj  \Delta _A (r_i, s_j) X_t(s_j)      \right\}^2 \label{eq:CAFR term 0} \\ 
		  + & \frac{1}{Tn}\sum_{t=1}^T \sum_{i=1}^n \wri \lb \Delta_\beta (r_i) , Z_t\rb^2 _p   \\
		  +& \frac{2}{Tn}\sum_{t=1}^T \sum_{i=1}^n \wri \left\{\frac{1}{n }   \sum_{j=1}^n \wsj   \Delta _A (r_i, s_j) X_t(s_j)\lb \Delta_\beta (r_i) , Z_t\rb _p \right\} \label{eq:CAFR term 2} \\
 		\le & \frac{2}{Tn} \sum_{t=1}^T \sum_{i=1}^n     \wri   \left( \frac{1}{n }   \sum_{j=1}^n \wsj   \Delta _A (r_i, s_j) X_t(s_j)  +     \lb \Delta_\beta (r_i) , Z_t\rb  _p  \right)\epsilon_t(r_i)   \label{eq:CAFR term 4} 
 		\\
 		+& \lambda \sum_{j=1}^p \|\beta_j  ^*  \|_n  - \lambda\sum_{j=1}^p  \| \widehat \beta_j     \|_n  \label{eq:CAFR term 5} .
	\end{align} 
	In the rest of the proof, for readability, we will assume that 
	$$\wri =1  \text{ for all } 1\le i \le n \quad \text{and} \quad  \wsj= 1\text{ for all }   1\le j \le n  .$$
	This assumption is equivalent to the case that $\{ r_i\}_{i=1}^n $ and $\{s_j\}_{j=1}^n $ are equally spaced grid. 
We note that in view of \Cref{lemma:discrete approximation}, the general case follows straightforwardly from the same argument and thus will be omitted for readability. 
\\
\\	
	 {\bf Step 1.}
Observe that 
	\begin{align}
		\eqref{eq:CAFR term 2} = & \frac{2}{Tn}\sum_{t=1}^T \sum_{i=1}^n \left\{ \frac{1}{n } \sum_{j=1}^n \Delta_A  (r_i, s_j) X_t(s_j) - \int_{[0, 1]} \Delta_A(r_i, s) X_t (s) \dint s\right\} \lb \Delta_\beta (r_i) , Z_t\rb_p  \label{eq:CAFR term 21}  \\
		& \hspace{2cm} + \frac{2}{Tn}\sum_{t=1}^T \sum_{i=1}^n \left\{\int_{[0, 1]} \Delta_A(r_i, s) X_t (s) \dint s  \right\}\lb \Delta_\beta (r_i) , Z_t\rb_p. \label{eq:CAFR term 22} 
	\end{align}
\
\\
As for the term \eqref{eq:CAFR term 21}, we have that for any $j \in \{1, \ldots, n\}$ and $t \in \{1, \ldots, T\}$, 
$$
	  \|\Delta_A(r_j, \cdot) X_t(\cdot)\|_{\wot } \leq   \|\Delta_A(r_j, \cdot)\|_\wot   \| X_t \|_ \wot \le 8C_AC_X \sqrt {\log(T) },
$$
where the first inequality follows from \Cref{lemma:basis properpties of sobolev} and  event $\mathcal D$.
	By   \Cref{lemma:discrete approximation},   it holds that  for all $1\le i\le n , 1\le t\le T   $,
	$$\bigg| \frac{1}{n } \sum_{j=1}^n \Delta_A  (r_i, s_j) X_t(s_j) - \int_{[0, 1]} \Delta_A(r_i, s) X_t (s) \dint s \bigg| 
	\le C_1   \zeta_n  \sqrt {\log(T) }   $$
	where $\zeta_n = n^{ -\alpha +1/2} $ and  $C_1$ depends on $C_A$ and $C_X$ only. 
	So,
	\begin{align*}
		\eqref{eq:CAFR term 21} \ge & - \frac{C _1 }{Tn}\sum_{t=1}^T \sum_{i=1}^n     \zeta_n  \sqrt {\log(T) }    \left|  \lb \Delta_\beta (r_i) , Z_t\rb _p   \right| 
		\\ 
		\ge&  - C _1      \zeta_n  \sqrt {\log(T) }   \sqrt {  \frac{1}{Tn}\sum_{t=1}^T \sum_{i=1}^n \lb \Delta_\beta (r_i) , Z_t\rb_p^2} \\
		\ge & - \frac{1}{640 Tn}\sum_{t=1}^T \sum_{i=1}^n \lb \Delta_\beta (r_i) , Z_t\rb  _p     ^2   -  640C_1^2  \zeta_n ^2    \log(T)  
		\\
		=
		 & - \frac{1}{640 Tn}\sum_{t=1}^T \sum_{i=1}^n \lb \Delta_\beta (r_i) , Z_t\rb  _p     ^2   -  C_1'    \zeta_n ^2    \log(T)      ,
	\end{align*}
	for some absolute constant $C '_1  > 0$.

In addition, for \Cref{eq:CAFR term 22},   note that  with probability at least $  1-(p\vee T)^{-4}$,
	\begin{align} \nonumber 
		&  \left\{\int_{[0, 1]} \Delta_A(r_i, s)  \frac{1}{T} \sum_{t=1}^T X_t(s)  Z_t^\top   \Delta_\beta(r_i) \dint s  \right\} 
- 			  \left\{\int_{[0, 1]} \Delta_A(r_i, s)  \E \{  X_t (s)  Z_t^\top  \}  \Delta_\beta(r_i)    \dint s  \right\}  
\\\nonumber 
= &  \int_{[0, 1]} \Delta_A(r_i, s)   \bigg( \frac{1}{T} \sum_{t=1}^T X_t(s)  Z_t^\top-  \E \{  X_t (s)  Z_t^\top  \}   \bigg)   \Delta_\beta(r_i) \dint s  
\\\nonumber 
		 \ge &  -\| \Delta_A (r_i,\cdot )  \|_\lt \bigg\|     \bigg( \frac{1}{T} \sum_{t=1}^T X_t(s)  Z_t^\top- \E \{  X_t (s)  Z_t^\top  \}   \bigg)   \Delta_\beta(r_i)   \bigg\| _\lt \\\nonumber 
		 \ge 
		 & -2 C_A \sum_{1\le j  \le p }\bigg\|       \frac{1}{T} \sum_{t=1}^T X_t(s)  Z_{t, j} ^\top- \E \{  X_t (s)  Z_{t, j  } ^\top  \}   \bigg\| _\lt   |\Delta_{\beta_ j } (r_i)|     
		 \\\nonumber  
		\ge &-2 C_A   \sum_{1\le  j  \le p } |\Delta_{\beta_ j }(r_i)  | \max_{1\le  j \le p }\bigg\|       \frac{1}{T} \sum_{t=1}^T X_t(s)  Z_{t,  j  } ^\top- \E \{  X_t (s)  Z_{t,  j } ^\top  \}   \bigg\| _\lt   
		\\ \label{eq:main step 1 term 2}
		\geq & -C _1'    \sum_{1\le  j  \le p } |\Delta_{\beta_j  }(r_i)  |  \sqrt { \frac{\log(p) }{T }} 
		\ge    -\frac{\lambda}{200} \sum_{1\le j  \le p }  |\Delta_{\beta_ j }(r_i)  |  ,
	\end{align}
where the second  inequality holds because $ \| \Delta_A (r_i,\cdot )  \|_\lt \le   \| \Delta_A (r_i,\cdot )  \|_ \hk \le 2C_A$, the fourth  inequality   follows from  \Cref{lemma:quadratic in functional space}, and the last inequality follows from the assumption that $ \lambda=C_\lambda \sqrt { \frac{\log(p\vee T )}{T}} $ with sufficiently large $C_\lambda$. In addition, by \Cref{assume:mixed}(c),
\begin{align*} 
 		 & \int_{[0, 1]} -\Delta_A(r_i, s)  \E \{  X_t (s)  Z_t^\top  \}  \Delta_\beta(r_i)    \dint s 
 		 \\
 		  \le &  \frac{3}{4} \sqrt { \Sigma_X [\Delta_A(r_i, \cdot) ,\Delta_A(r_i, \cdot )] } \sqrt { \Delta_\beta(r_i)^\top  \Sigma_Z\Delta_\beta(r_i) }
 		  \\
 		  \le & \frac{3}{8} \bigg (   \frac{6}{5}  \Sigma_X [\Delta_A(r_i, \cdot) ,\Delta_A(r_i, \cdot )]   + \frac{5}{6}  \Delta_\beta(r_i)^\top  \Sigma_Z\Delta_\beta(r_i) \bigg) 
 		  \\
 		=
 		  &  \frac{9}{20} \Sigma_X [\Delta_A(r_i, \cdot) ,\Delta_A(r_i, \cdot )]   + \frac{5}{16} \Delta_\beta(r_i)^\top  \Sigma_Z\Delta_\beta(r_i)
 		  \\
 		  =&  \frac{9}{20}  \iint_{[0,1]^2 } \Delta_A(r_i, r ) \Sigma_X(r,s )\Delta_A(r_i, s  )\dint r \dint s   + \frac{5}{16} \Delta_\beta(r_i)^\top  \Sigma_Z\Delta_\beta(r_i) .
	\end{align*} 
So  
\begin{align*} \nonumber 
		 			  &\eqref{eq:CAFR term 22}
		 			  \\
		 			   \ge &-  \frac{\lambda}{1 00} \frac{1}{n} \sum_{i=1}^n  \sum_{j=1 }^{  p }  |\Delta_{\beta_j }(r_i)  |  
		 			 -  \frac{1}{n} \sum_{i=1}^n \frac{9}{10}  \iint_{[0,1]^2}  \Delta_A(r_i, r ) \Sigma_X(r,s )\Delta_A(r_i, s  )\dint r \dint s    - \frac{1}{n} \sum_{i=1}^n  \frac{5}{8} \Delta_\beta(r_i)^\top  \Sigma_Z\Delta_\beta(r_i) 
		 			 \\
		 			 \ge &- \frac{\lambda}{1 00}   \sum_{j=1 }^{  p }   \|\Delta_{\beta_ j   } \|_n  -  \frac{1}{n} \sum_{i=1}^n \frac{9}{10}  \iint_{[0,1]^2}  \Delta_A(r_i, r ) \Sigma_X(r,s )\Delta_A(r_i, s  )\dint r \dint s    - \frac{1}{n} \sum_{i=1}^n  \frac{5}{8} \Delta_\beta(r_i)^\top  \Sigma_Z\Delta_\beta(r_i)  
		 			 \\
		 			=  & - \frac{\lambda}{1 00}  \big(  \sum_{  j\in S }  \|\Delta_{\beta_j   } \|_n  
		 			  +     \sum_{  j \in S^c  }  \| \widehat \beta_ j   \|_n \big)   
		 			  \\ -&   \frac{1}{n} \sum_{i=1}^n \frac{9}{10}  \iint_{[0,1]^2} \Delta_A(r_i, r ) \Sigma_X(r,s )\Delta_A(r_i, s  )\dint r \dint s    - \frac{1}{n} \sum_{i=1}^n  \frac{5}{8} \Delta_\beta(r_i)^\top  \Sigma_Z\Delta_\beta(r_i)  
	\end{align*} 
	Putting  calculations for \Cref{eq:CAFR term 21} and \Cref{eq:CAFR term 22} together,  we have that
	\begin{align} \nonumber 
		\eqref{eq:CAFR term 2} \ge & -\frac{1}{640 Tn}\sum_{t=1}^T \sum_{i=1}^n \lb \Delta_\beta (r_i) , Z_t\rb  _p     ^2   -  C_1'    \zeta_n ^2    \log(T)   - \frac{\lambda}{1 00}  \big(  \sum_{  j\in S }  \|\Delta_{\beta_j   } \|_n  
		 			  +     \sum_{ j \in S^c  }  \| \widehat \beta_j  \|_n \big)    
		\\
		      -  &  \frac{1}{n} \sum_{i=1}^n \frac{9}{10}  \iint_{[0,1]^2}  \Delta_A(r_i, r ) \Sigma_X(r,s )\Delta_A(r_i, s  )\dint r \dint s    - \frac{1}{n} \sum_{i=1}^n  \frac{5}{8} \Delta_\beta(r_i)^\top  \Sigma_Z\Delta_\beta(r_i)  .\label{eq:CAFR term 23} 
	\end{align}
	\
	\\
	 {\bf Step 2.}  It follows from   \eqref{eq: function regression rev 2}  that with probability at least 
$
		1 -  T^{-4}  
$,
	it holds that
	\begin{align}   
		\eqref{eq:CAFR term 0} 
		 \ge &\frac{ 159  }{ 160 n    }   \sum_{i=1}^n \iint_{[0, 1]^2} \Delta_A (r_i,  r  ) \Sigma (r,s)  \Delta_A   (r_i,  s )\dint  r \dint   s  - C _2  \left(  \log(T) \zeta_n^2 + \log(T) \delta_T \right),  \label{eq: function regression rev2}
	\end{align}
	where $C_2  > 0$ is an absolute constant.  

It follows from  \eqref{eq:FFR step 2 2} and \eqref{eq:sc form term 2 bound 2}, we have that with probability at least 
	$1 -2T^{-4}$,
	it holds that
	\begin{align*}
		 \eqref{eq:CAFR term 4} \le&  \frac{1 }{3 2 0 n} \sum_{i=1}^n \iint_{[0, 1]^2} \Delta_A (r_i,  r) \Sigma_X  (r,s)  \Delta_A(r_i,  s ) \dint r \dint s 
		  + C_2'   \Bigg\{   \zeta_n   +   \log(T)\delta_T \Bigg\} 
		  \\
		  + & \frac{  \lambda }{ 320 }       \sum_{j \in S }   \| \Delta_{\beta_j }\|_n    +   \frac{\lambda }{ 320 }     \sum_{j \in S ^c}   \|  \widehat \beta_j  \|_n  
	\end{align*}
	where $C_2' > 0$ is an absolute constant.

\medskip 
\noindent {\bf  Step 3.} Note that 
for \Cref{eq:CAFR term 5}, 
\begin{align*}\lambda\sum_{j=1}^p  \|\beta_j  ^*  \|_n  -\sum_{j=1}^p  \lambda \| \widehat \beta_j     \|_n  = &  \lambda \bigg\{ \sum_{j \in S }   \|\beta_j  ^*  \|_n  -  \sum_{j \in S }   \| \widehat \beta_j     \|_n -\sum_{j \in S^c  }   \| \widehat \beta_j     \|_n  \bigg\} 
\\
\le & \lambda \sum_{j \in S }   \| \Delta_{\beta_j}      \|_n  - \lambda \sum_{j \in S^c  }  \|  \widehat \beta_j   \|_n . 
\end{align*} 
\
\\
{\bf Step 4.} Putting all previous calculations  together leads to that
	\begin{align}
		& \frac{31}{320 n} \sum_{i=1}^n \iint_{[0, 1]^2} \Delta_A (r_i,  r  ) \Sigma_X (r,s)  \Delta_A  (r_i,  s ) \dint r \dint  s \nonumber \\ 
		  +&  \frac{639}{640Tn}\sum_{t=1}^T \sum_{i=1}^n \lb \Delta_\beta (r_i) , Z_t\rb_p^2  -  \frac{5}{8n }\sum_{i=1}^n  \Delta_\beta(r_i)^\top  \Sigma_Z\Delta_\beta(r_i) \label{eq:mixed eq2 2} \\
		\le & \big( 1 +  \frac{1}{100}  +  \frac{1}{320}   \big) \lambda   \sum_{  j \in S  }  \|  \Delta_{\beta_ j}  \| _{n }    -  \lambda  \big( 1- \frac{1}{100} - \frac{1}{320}   \big) \sum_{  j \in S ^c   }  \|   \widehat \beta_ j   \| _{n }  
		  +  C_3   \big \{    \log(T) \zeta_n     +\log(T)\delta_T \big \}    \nonumber,
		  \\
		  \le &  2\lambda   \sum_{  j \in S  }  \|  \Delta_{\beta_ j}  \| _{n }   +C_3   \big \{    \log(T) \zeta_n     +\log(T)\delta_T \big \}   
	\end{align}
	where  $C_3 > 0$ is  some absolute constant.
	\\
	\\
 For \Cref{eq:mixed eq2 2}, note that by constrain set $\mathcal C_\beta$, we have that 
$$ \sum_{ 1\le j  \le p  }    \|  \widehat \beta_ j   \|_n    \le   \sum_{ 1\le j  \le p  }    \| \widehat \beta_ j    \|_\infty       \le \sum_{ 1\le  j \le p  }    \|  \widehat \beta_ j  \|_\hb  
	\le C_\beta. 
$$
	Therefore 
	$$  \sum_{ 1\le j  \le p  }    \|  \Delta_{\beta_ j  }  \|_n   
 \le \sum_{ 1\le j  \le p  }     \|  \widehat \beta_ j   \|_n  + 
	 \sum_{ 1\le  j  \le p  }     \|    \beta_ j ^*    \|_n     \le 2C_\beta.   $$   From \Cref{theorem:restricted eigenvalue  condition},  it holds that with probability at least $1-\exp(-cT) $, 
	\begin{align*} 
		 \frac{639}{640T }\sum_{t=1}^T   \lb \Delta_\beta (r_i) , Z_t\rb _p^2  \ge &  \frac{639}{640  }  \cdot  \frac{2}{3} \Delta_\beta(r _i ) ^\top    \Sigma_Z \Delta_\beta(r _i )  - C_3' \frac{\log(p)}{T } \bigg(  \sum_{ 1\le j  \le p  }    \|  \Delta_{\beta_ j  }  \|_n          \bigg)  ^2 \\
		 \ge  &   \frac{213}{320 }    \Delta_\beta(r _i ) ^\top    \Sigma_Z \Delta_\beta(r _i )  -  4 C_3 'C_\beta^2  \frac{\log(p)}{T }  .
	\end{align*}
	 Substitute the above inequality into \Cref{eq:mixed eq2 2} gives 
\begin{align}
		& \frac{31}{320 n} \sum_{i=1}^n \iint_{[0, 1]^2} \Delta_A (r_i,  r  ) \Sigma_X (r,s)  \Delta_A  (r_i,  s ) \dint r \dint  s \nonumber \\ 
		  +&  \frac{13}{320 n} \sum_{i=1}^n  \Delta_\beta(r_i)^\top  \Sigma_Z\Delta_\beta(r_i)   \nonumber   \\
		\le & 2  \lambda   \sum_{  j \in S  }  \|  \Delta_{\beta_ j}  \| _{n }      \label{eq:mixed eq3} \\
		  +& C_3''   \big \{    \log(T) \zeta_n    +\log(T)\delta_T + \frac{\log(p)}{T } \big \}    \nonumber.
	\end{align} 
	\
	\\
 {\bf Step 5.}   Note that  by \Cref{assume:mixed}{\bf c},
 \begin{align*}  
		\left(  \lambda \sum_{  j   \in S }    \|  \Delta_{\beta_ j  }  \| _n   \right)^2 \le   \lambda ^2  \mathfrak {s} \sum_{  j =1}^p     \|  \Delta_{\beta_ j   }  \|_n    ^2   =  \lambda ^2   \mathfrak {s} \frac{1}{n }  \sum_{i=1}^n \sum_{  j =1}^p \{\Delta_{\beta_  j  } (r_i )\}^2  \le    \frac{ \lambda ^2   \mathfrak {s}  }{ c_z  }   \frac{1}{n }  \sum_{i=1}^n  \Delta_\beta(r _i )  ^\top  \Sigma_Z \Delta_\beta(r _i ),
	\end{align*}  
	where $ \mathfrak {s}= |S|$, the cardinality of  the support set $S$.
This gives 
	\begin{align} \label{eq:size of sparse 2}
		 \sum_{ j \in S }    \|  \Delta_{\beta_j }  \| _n      \le  \sqrt { \frac{  \mathfrak {s} }{c_z  }   \frac{1}{n }  \sum_{i=1}^n  \Delta_\beta(r _i )^\top   \Sigma_Z \Delta_\beta(r _i )} ,
	\end{align} 
 which implies that  
\begin{align*}\frac{1}{160} \lambda   \sum_{  j \in S  }  \|  \Delta_{\beta_ j}  \| _{n }   \le &  \frac{1}{320n} \sum_{i=1}^n  \Delta_\beta(r_i)^\top  \Sigma_Z\Delta_\beta(r_i)   + C_4  \mathfrak {s} \lambda ^ 2  
 \\
 \le &  \frac{1}{320n} \sum_{i=1}^n  \Delta_\beta(r_i)^\top  \Sigma_Z\Delta_\beta(r_i)   + C_4  C_\lambda^2  \frac{ \mathfrak {s} \log(p\vee T) }{T }  .
 \end{align*}
	 Substituting  the above inequality into   \eqref{eq:mixed eq3} gives 
	 \begin{align}
		& \frac{31}{320 n} \sum_{i=1}^n \iint_{[0, 1]^2} \Delta_A (r_i,  r  ) \Sigma_X (r,s)  \Delta_A  (r_i,  s ) \dint r \dint  s 
		  +   \frac{12}{320 n} \sum_{i=1}^n  \Delta_\beta(r_i)^\top  \Sigma_Z\Delta_\beta(r_i)  \nonumber \\
		\le &   
		    C_4'    \big \{    \log(T) \zeta_n     +\log(T)\delta_T + \frac{\log(p)}{T } +  \frac{  \mathfrak {s}\log(p\vee T) }{T }  \big \}  
		  \label{eq:main step 4 sample}
	\end{align}  
\
\\
{\bf Step 6.}
Note that by definition,  $\mathcal{E}^*(\widehat{A}, \widehat{\beta}) \ge 0 $ for any estimators $\widehat{A}$ and $\widehat{\beta}$, and that
    \begin{align} \nonumber 
		\mathcal E^* (\widehat A, \widehat \beta) \le & 4\int_{[0, 1]} \mathbb{E} \left\{\int_{[0, 1]} \Delta_A (r, s)   X^ *(s) \dint s\right\} ^2  \dint r + 4 \mathbb{E}\int_{[0, 1]} \left\{( Z_i^*)^\top  \Delta_\beta (r) \right\}^2 \dint r \\
		= & 4 \int_{[0, 1]} \bigg\{\iint_{_{[0,1]^2} }  \Delta_A  (r,s) \Sigma_X (s, u) \Delta_A (r, u) \dint s\dint u + \Delta_\beta ^\top  (r)  \Sigma_Z \Delta_\beta    (r)        \bigg\} \dint r, \label{eq:bound of excess risk theorem 1}
	\end{align}     
	where $\Sigma _X $ and $\Sigma _Z$ are the covariance of $X_t$ and $Z_t$,  as defined in \Cref{assume:mixed}. 
\\
\\
 It follows from \eqref{eq:main 1 step 5 lower bound 2} that  with probability  at least $1-T^{-4}   $,
\begin{align}\nonumber  
&\frac{1}{  n}  \sum_{ i=1}^ n  \iint_{[0, 1]^2}  \Delta_A (r_i,  r  ) ^\top \Sigma_X (r,s)  \Delta_A  (r_i,  s )\dint r \dint s
\\ \nonumber  
=&  		
		\mathbb{E}_{X^*}\left\{\frac{1}{n}\sum_{i=1}^n\left( \int_{[0, 1]}\Delta_A   (r_i,   s  )   X^*(s) \dint s \right)^2 \right\} 
		\\
		\ge &  c_5\mathbb{E}_{X^*}\left\{\int_{[0, 1]} \left  ( \int_{[0, 1]} \Delta_A  (s,r) X^*(r)  \dint r  \right )^2 \dint s  \right\}  -  C_5 \log(T) \zeta_n  ,   \label{eq:step 5 first}
\end{align} 
with $c_5, C_5 > 0$ being absolute constants. 
 In addition,  with absolute constants $c_5', C_5' > 0$, due to \eqref{eq:main 2 last step inequality 2}, we have that  
\begin{align}  \frac{1}{ n }    \sum_{i=1}^n  \Delta_\beta(r _i )  ^\top   \Sigma_Z \Delta_\beta(r _i ) \ge 
	c_5'        \int_{[0, 1]} \Delta_{\beta } (r )^\top \Sigma_Z \Delta_{\beta } (r )    \dint r  -C_5'  \zeta _n  \label{eq:step 5 second}
	\end{align}
	\
	\\
Finally, \eqref{eq:main step 4 sample},  \eqref{eq:bound of excess risk theorem 1},  \eqref{eq:step 5 first} and \eqref{eq:step 5 second} together imply that     with probability at least
$		1 - 6T^{-4}  ,$
	it holds that
	\begin{align*}  
		\mathcal{E}^*(\widehat{A}, \widehat{\beta}) \le   C _5''   \big \{    \log(T) \zeta_n       +\log(T)\delta_T + \frac{\log(p)}{T } + \frac{ \mathfrak {s}\log(p\vee T) }{T }    \big \}  .
	\end{align*} 
	where $C  $ is some  absolute constant.  This directly leads to the desired result.
\end{proof}

\subsection{Additional proofs related to  \Cref{thm-main-mixed} }\label{sec-proofs-special}
In this section, we present the technical results related to  \Cref{thm-main-mixed}.  
   Recall in the proof of  \Cref{thm-main-mixed}, we set  that $\Delta_\beta (r) = \widehat \beta(r) -\beta^*(r)$ and $\Delta_A(r,s) = \widehat A(r,s) - A^*(r,s)$, $r, s \in [0, 1]$.  Note that for any $ g\in \lt $, it holds that
	\[
		\mathbb{E}(\lb g, X^*\rb_\lt^2 ) =  \iint_{[0, 1]^2} g(r)\Sigma_X(r,s) g(s) \dint r \dint  s.
	\]
	We also assume that the following good event holds:
	$$   \| X_t \|_ \wot  \le 4C_X \sqrt {\log(T)}  \ \text{for all } 1\le t \le T   .   $$ 
	It was justified in the proof of \Cref{thm-main-mixed} that 
	$$ \mathbb P(\| X_t \|_\wot \le 4C_X \sqrt {\log(T)}  \ \text{for all } 1\le t \le T ) \ge 1- T^{-5}.$$

  \begin{lemma}Let $\zeta_n=   n^{-\alpha +1/2  } $.  
  Under the same conditions as in \Cref{thm-main-mixed}, with probability at least 
$
		1 -  T^{-4}  
$,  it holds that 
  \begin{align}	  \nonumber  &\frac{1}{Tn}\sum_{t=1}^T \sum_{i=1}^n \left\{\frac{1}{n }   \sum_{j=1}^n  \Delta _A (r_i, s_j) X_t(s_j)      \right\}^2  
  \\
   \ge &\frac{ 159  }{ 160 n    }   \sum_{i=1}^n \iint_{[0, 1]^2} \Delta_A (r_i,  r  ) \Sigma (r,s)  \Delta_A   (r_i,  s )\dint  r \dint   s  - C _2  \left(  \log(T) \zeta_n^2 + \log(T) \delta_T \right).  \label{eq: function regression rev 2}
  \end{align}
   \end{lemma} 
   
   \begin{proof}
   Let $\zeta_n=   n^{-\alpha +1/2  } $ and $\delta_T =T^{-2r/(2r+1 )}$.
Observe that for any $i \in \{1, \ldots, n\}$ and $t \in \{1, \ldots, T\}$, we have
	\begin{align*}
		     \| \Delta_A  (r_i, \cdot )  X_t(\cdot )   \|_ \wot  \le   \| \Delta_A  (r_i, \cdot )  \|_\wot \| X_t(\cdot )   \|_\wot  \le 8 C_AC_X \sqrt {\log (T)},
	\end{align*}
	where the first inequality follows   \Cref{lemma:basis properpties of sobolev}.
	By \Cref{lemma:discrete approximation}, we have that  
	\[
		\left |\frac{1}{n} \sum_{j =1}^n \Delta_A(r_i, s_j) X_t (s_j) - \int_{[0, 1]} \Delta_A (r_i, s) X_t(s)  \dint   s \right| \le C_4  \zeta_n \sqrt {\log(T) },
	\]
	where $C_4 > 0$ is an absolute constant.  Therefore, we have that  
	\begin{align} 
		\left\{\frac{1}{n} \sum_{ j =1}^n  \Delta_A (r_i, s_j) X_t (s_j) \right\}^2  & \ge \frac{319}{320}  \left\{\int_{[0, 1]} \Delta_A(r_i, s) X_t(s) d  s \right\}^2 -   C _ 4' \zeta_n^2 \log(T)  \nonumber \\
		& = \frac{319}{320} \lb\Delta_A(r_i, \cdot), X_t (\cdot)\rb_\lt^2 -  C _4'   \zeta_n^2 \log(T)  \label{eq:main 1 step 4 first 2} .
	\end{align} 
By  \Cref{lemma:functional restricted eigenvalues} and the fact that 
$$ 
\| \Delta_A(r_i, \cdot)\| _\hk \le \|A^* \|_\hk + \|\widehat A\|_\hk \le 2C_A,
$$
we have that  with probability at least $1-T^{-4} $,
	it holds that uniformly for all $ 1\le i \le n$, 
	\begin{align*} 
		& \left|\frac{1}{T}\sum_{t=1}^T \lb\Delta_A(r_i, \cdot), X_t(\cdot)\rb_\lt^2 - \iint_{[0, 1]^2}  \Delta_A (r_i,  r)\Sigma _X (r,s) \Delta _A  (r_i,  s ) \dint r \dint s   \right|	\\
		\leq &   \frac{1}{320} \iint_{[0, 1]^2} \Delta_A (r_i,  r  ) \Sigma_X (r,s)  \Delta_A   (r_i,  s ) \dint r \dint s      +C_4' \log(T) \delta_T 
	\end{align*}
	where $C_4' > 0$ is an absolute constant.  Thus  the above display implies that 
	\begin{align}\label{eq:main 1 step 4 second 2} 
		 & \frac{1}{T}\sum_{t=1}^T   \lb  \Delta_A  (r _i ,\cdot) , X_t\rb _\lt ^ 2 \ge  \frac{319}{320} \iint_{[0, 1]^2} \Delta_A (r_i,  r  ) \Sigma (r,s)  \Delta_A   (r_i,  s ) \dint   r \dint  s 
	  - C _4''  \log(T) \delta_T .
	\end{align} 
	Therefore,
	\begin{align}    \nonumber   
		 & \sum_{t=1}^T \sum_{i=1}^{n }  \left( \frac{1}{n } \sum_{j=1}^{n }   \Delta_A  (r_i, s_j) X_t (s_j) \right)^2 
		 \\  \geq &  \sum_{t=1}^T \sum_{i=1}^n \frac{319}{320}   \lb  \Delta_A  (r _i ,\cdot) , X_t\rb _\lt ^2 - C _4'  Tn \log(T) \zeta_n^2    \nonumber \\
		\ge & T  \bigg(\frac{ 319  }{ 320  }\bigg)^2   \sum_{i=1}^n \iint_{[0, 1]^2} \Delta_A (r_i,  r  ) \Sigma (r,s)  \Delta_A   (r_i,  s )\dint r \dint s- C _4''' Tn \left(  \log(T) \zeta_n^2 + \log(T) \delta_T \right),  \nonumber 
		\\
			\ge &  T\frac{ 159  }{ 160    }   \sum_{i=1}^n \iint_{[0, 1]^2} \Delta_A (r_i,  r  ) \Sigma (r,s)  \Delta_A   (r_i,  s )\dint r \dint s- C _4''' Tn \left(  \log(T) \zeta_n^2 + \log(T) \delta_T \right), \label{eq: function regression rev}
	\end{align}
	where the first inequality follows from \eqref{eq:main 1 step 4 first 2} and the second inequality follows from \eqref{eq:main 1 step 4 second 2}.

  \end{proof}

   \begin{lemma}Let $\zeta_n=   n^{-\alpha +1/2  } $ and $\delta_T =T^{-2r/(2r+1 )}$.  
  Under the same conditions as in \Cref{thm-main-mixed}, if in addition, $\delta_T \ge \frac{\log(T)}{T} $, then with probability at least
$
		1  -  T^{-4}
$, it holds that  
  \begin{align} \nonumber 
& \frac{ 2}{Tn }\sum_{t=1}^T \sum_{i=1}^n  \left\{\frac{1}{n} \sum_{j=1}^n  \Delta _A (r_i, s_j) X_t (s_j) \right\}  \epsilon_t(r _i)\\ 
  \le &     \frac{1 }{3 2 0 n} \sum_{i=1}^n \iint_{[0, 1]^2} \Delta_A (r_i,  r) \Sigma_X  (r,s)  \Delta_A(r_i,  s ) \dint r \dint s 
		  + C    \Bigg\{   \zeta_n    +\log(T)\delta_T \Bigg\}    \label{eq:FFR step 2 2} .
\end{align} 
In addition,  let $\{ \mathfrak{E} _{t,i} \}_{t=1,i=1}^{T,n}   $ be a collection of standard normal random variables independent of $\{ X_t\}_{t=1}^T $,$\{ r_i\}_{i=1}^n  $ and $\{ s_i\}_{i=1}^n  $.  Then with probability at least
$1  -  T^{-4}$, 
it holds that  
  \begin{align} \nonumber 
& \frac{ 2}{Tn }\sum_{t=1}^T \sum_{i=1}^n  \left\{\frac{1}{n} \sum_{j=1}^n  \Delta _A (r_i, s_j) X_t (s_j) \right\} \mathfrak{E} _{t,i}  \\ 
  \le &     \frac{1 }{3 2 0 n} \sum_{i=1}^n \iint_{[0, 1]^2} \Delta_A (r_i,  r) \Sigma_X  (r,s)  \Delta_A(r_i,  s ) \dint r \dint s 
		  + C    \Bigg\{   \zeta_n    +\log(T)\delta_T \Bigg\}    \label{eq:FFR step 2 3} .
\end{align}  
  \end{lemma}
  
  \begin{proof} Let $\zeta_n=   n^{-\alpha +1/2  } $. 
  For  \Cref{eq:FFR step 2 2}, note that 
	\begin{align}
	\nonumber 	& 2\sum_{t=1}^T \sum_{i=1}^n  \left\{\frac{1}{n} \sum_{j=1}^n  \Delta _A (r_i, s_j) X_t (s_j) \right\}  \epsilon_t(r _i)
		\\ = & 2\sum_{t=1}^T \sum_{i=1}^n \left\{\frac{1}{n} \sum_{j =1}^n \Delta_A (r_i, s_j) X_t (s_j) - \int_{[0, 1]} \Delta_A (r_i, s) X_t(s) \dint s\right\}\epsilon_t(r_i)   \label{eq:third term 1 2} \\
		   +&  2 \sum_{t=1}^T \sum_{i=1}^n \epsilon_t(r_i) \int_{[0, 1]}\Delta_A(r_i, s) X_t(s) \dint  s. \label{eq:third term 2 2}
	\end{align}
	Since  for all $r_i$
	$$ \| \Delta _A  (r_i, \cdot ) X_t (\cdot)   \|_\wot \le 8C_AC_X  \sqrt {\log(T) },$$
\Cref{lemma:sobolev covering 0} implies that with probability $1 -  T^{-4}  $, it holds that uniformly for all $1\le i \le n$,
 	\begin{align*}
		\frac{1}{T} \sum_{t=1}^T \left\{\frac{1}{n} \sum_{   j =1}^n  \Delta _A  (r_i, s_j) X_t (s_j)   - \int_{[0, 1]} \Delta _A  (r_i,s) X_t(s) \dint  s \right\}\epsilon_t(r_i) \le C _2  \zeta_n\log(T)   ,
	\end{align*}
	where $C _2 > 0$ is an absolute constant.  Therefore 
	\[
		\eqref{eq:third term 1 2} \le 2 C _ 2 Tn  \zeta_n\log(T)   .
	\]

To control the term \eqref{eq:third term 2 2}, we deploy \Cref{lemma:functional deviation bound}. Since
$$ 
\| \Delta_A(r_i, \cdot)\| _\hk \le \|A^* \|_\hk + \|\widehat A\|_\hk \le 2C_A,
$$ 
then with probability at least 
	$		1 -T^{-4}$, it holds that 
	that
	\begin{align*}
		  \eqref{eq:third term 2 2} 
		\leq &  C _2 'Tn   \Bigg\{ \frac{1}{n}\sum_{i=1}^n    \sqrt {\log(T)\delta_T \iint_{[0, 1]^2} \Delta_A (r_i,  r)\Sigma_X (r,s) \Delta_A (r_i, s)\dint r \dint s} + \log(T)  \delta_T   \Bigg\}
		\\ 
			\leq  & 2C _2'Tn   \Bigg\{ 320 C_2' \log(T)  \delta_T   + \frac{1}{320 C_2' n}\sum_{i=1}^n \iint_{[0, 1]^2} \Delta_A (r_i,  r)\Sigma_X (r,s) \Delta_A (r_i, s)\dint r \dint s  \Bigg\},  
	\end{align*}
	where $C_2' > 0$ is an absolute constant.
\\
\\
Therefore, we have that with probability at least
$
		1   -2 T^{-4}
$
	that
	\begin{align}   \nonumber 
		& 2\sum_{t=1}^T \sum_{i=1}^n  \left\{\frac{1}{n} \sum_{j=1}^n  \Delta _A (r_i, s_j) X_t (s_j) \right\}  \epsilon_t(r _i) 
		\\
		  \le & Tn   \frac{1 }{3 2 0 n} \sum_{i=1}^n \iint_{[0, 1]^2} \Delta_A (r_i,  r) \Sigma_X  (r,s)  \Delta_A(r_i,  s ) \dint r \dint s 
		  + C_2'' Tn \Bigg\{   \zeta_n  + \frac{\log^2 (T)}{T} +\log(T)\delta_T \Bigg\},  \nonumber  
	\end{align}
	where $C_3 > 0$ is an absolute constant.     \Cref{eq:FFR step 2 2} follows from the assumption that $ \delta_T \ge \log(T)/T .$
	\\
	\\
	 The argument of \Cref{eq:FFR step 2 3} is the same as that of  \Cref{eq:FFR step 2 2} and will be omitted for brevity. 
  \end{proof}

   \begin{lemma}Let $\zeta_n=   n^{-\alpha +1/2  } $ and $\delta_T =T^{-2r/(2r+1 )}$.   
  Under the same conditions as in \Cref{thm-main-mixed},  with probability at least $1-(T\vee p)^{-4}$, it holds that
\begin{align}   
		  \frac{2}{Tn} \sum_{t=1}^T \sum_{i=1}^n \lb \Delta_\beta (r_i) , Z_t\rb_p  \epsilon_t(r_i)
		  \le     \frac{\lambda }{ 320 }       \sum_{ l  \in S }   \| \Delta_{\beta_l  }\|_n    +   \frac{\lambda }{ 320 }     \sum_{ l  \in S ^c}   \|  \widehat \beta_ l   \|_n  \label{eq:sc form term 2 bound 2}.
	\end{align} 
	In addition,  let $\{ \mathfrak{E} _{t,i} \}_{t=1,i=1}^{T,n}   $ be a collection of standard normal random variables independent of $\{ X_t\}_{t=1}^T $, $\{ r_i\}_{i=1}^n  $ and $\{ s_i\}_{i=1}^n  $.  Then with probability at least
$1  -  T^{-4}$, 
it holds that  
  \begin{align}  \frac{2}{Tn} \sum_{t=1}^T \sum_{i=1}^n \lb \Delta_\beta (r_i) , Z_t\rb_p  \mathfrak{E} _{t,i}
		  \le     \frac{\lambda }{ 320 }       \sum_{ l  \in S }   \| \Delta_{\beta_l  }\|_n    +   \frac{\lambda }{ 320 }     \sum_{ l  \in S ^c}   \|  \widehat \beta_ l   \|_n    \label{eq:sc form term 2 bound 23}.
\end{align}  
  \end{lemma}
  
  \begin{proof} Let $\zeta_n=   n^{-\alpha +1/2  } $.  For \Cref{eq:sc form term 2 bound 2}, note that 
	 since $Z_{t,j} $ is centered Gaussian with variance bounded by $C_Z $ and  $\epsilon_t (r_i) $ is Gaussian with variance bounded by $C_\epsilon $, $Z_{t,j} \epsilon _t(r_i)$ is sub-exponential with parameter $C_\epsilon C_Z $. Therefore by    a union bound argument,  there exists  an absolute constant $C_1  > 0$ such that, for any $i \in \{1, \ldots, n\}$,
	\begin{align} \label{eq:main 2 step 1 first 2}
		\mathbb{P}\left\{\left\|\frac{1}{T} \sum_{t=1}^T Z_t \epsilon_t(r_i)\right\|_{\infty} \ge C  _1 \sqrt { \frac{\log(p \vee T ) }{T }} \right\} \le (Tp)^{-5}.
	\end{align} 
	Therefore, with probability at least $1-(T\vee p)^{-5}$, it holds that 
	\begin{align*}   \nonumber 
		& \frac{2}{Tn} \sum_{t=1}^T \sum_{i=1}^n \lb \Delta_\beta (r_i) , Z_t\rb_p  \epsilon_t(r_i)
		\\ \le & \frac{2}{n} \sum_{i=1}^n \left\{\left\|\frac{1}{T} \sum_{t=1}^T Z_t  \epsilon_t(r_i)\right\|_{\infty} \sum_{ l =1}^p |\Delta_{\beta_ l  } (r_i)| \right\} \nonumber \\ 
		\le & \frac{2C_1 }{n} \sqrt { \frac{\log(p\vee T ) }{T }}\sum_{i=1}^n   \sum_{ l =1}^p |  \Delta_{\beta_ l  } (r_i)|
		\\
		 \le   &   C_2 \sqrt { \frac{\log(p\vee T ) }{T }  } \sum_{ l =1}^p        \| \Delta_{\beta_ l  }\|_n   
		 \\
		 \le  &   \frac{\lambda }{ 320 }       \sum_{ l  \in S }   \| \Delta_{\beta_l  }\|_n    +   \frac{\lambda }{ 320 }     \sum_{ l  \in S ^c}   \|  \widehat \beta_ l   \|_n 
	\end{align*}
	where the second inequality follows from \eqref{eq:main 2 step 1 first 2}, the third inequality follows from H\"{o}lder's inequality,  and the last inequality follows from the assumption that 
	$\lambda =C_\lambda \sqrt { \frac{\log(p\vee T) }{T}}$ for sufficiently large constant $ C_\lambda$.
	\\
	\\
The argument of \Cref{eq:sc form term 2 bound 23} is the same as that of  \Cref{eq:sc form term 2 bound 2} and will be omitted for brevity.  
  \end{proof} 
 	
  \begin{lemma}Let $\zeta_n=   n^{-\alpha +1/2  } $.  
  Under the same conditions as in \Cref{thm-main-mixed}, with probability at least 
$
		1 -  T^{-4}  
$,  it holds that 
\begin{align} \nonumber   & \frac{1}{n} \sum_{i=1}^n \iint_{[0, 1]^2} \Delta_A(r_i,  r) \Sigma_X (r, s) \Delta_A  (r_i,  s ) \dint r \dint s 
\\
\ge   & c  \mathbb{E}_{X^*}\left\{\int_{[0, 1]} \left  ( \int_{[0, 1]} \Delta_A  (s,r) X^*(r)  \dint   r  \right )^2 \dint s  \right\}  -  C '\log(T)  \zeta_n  .
\label{eq:main 1 step 5 lower bound 2}
\end{align} 
\end{lemma} 
  
  \begin{proof}
   We have 
	\begin{align}  
		& \mathbb{E}_{X^*}\left\{\frac{1}{n} \sum_{i=1}^n  \left( \int_{[0, 1]} \Delta_A(r_i, s) X^*(s) \dint s \right)^2 \right\} =  \frac{1}{n} \sum_{i=1}^n \iint_{[0, 1]^2} \Delta_A(r_i,  r) \Sigma_X (r, s) \Delta_A  (r_i,  s ) \dint r \dint s \nonumber   
		  \label{eq:step 3 er 3 1}
	\end{align}  
	where $X^*$ is the predictor in the test set.  By \Cref{lemma:Gaussian process sobolev norm},   we have with probability at least 
	$ 1- T^{-4} $ that 
	 $\| X ^* \|_\wot  \le 4 C_X \sqrt {\log(T) } .$
So it holds that 
	\begin{align*}
		& \left \|\lb\Delta_A(r_i, \cdot ), X^*\rb_\lt\right\|_ \wot  = \left\|\int_{[0, 1]} \Delta_A(r_i,  s)  X^*(s) \dint  s \right \| _\wot   \\
		\le & \int_{[0, 1]}  \| \Delta _A  (\cdot,s )\|_\wot   | X^*  (s)  |  \dint s    \le 2 C_A  \|X^*\|_\lt \le 8 C_AC_X  \sqrt {\log(T) } .     
	\end{align*}
	Let  $f (r) =  \lb  \Delta_A   (r ,  \cdot ) ,X^*(\cdot)\rb_\lt   $. Then
	$$\| f^2\|_\wot \le 64C_A^2 C_X^2 \log(T)  . $$
	Applying \Cref{lemma:discrete approximation} to $f ^2   $, we have  that 
	\[
		\frac{   1}{   n}  \sum_{i=1}^n    \lb  \Delta_A   (r_i,  \cdot ) ,X^*\rb_\lt^2   
		 \ge   \int_{[0, 1]}   \lb \Delta _A (r , \cdot ) ,X^*\rb_\lt^2 \dint r     - C_5'  \log (T) \zeta_n ,
	\]
	and this implies that 
\begin{align*}  
		& \mathbb{E}_{X^*}\left\{\frac{1}{n}\sum_{i=1}^n\left( \int_{[0, 1]}\Delta_A   (r_i,   s  )   X^*(s) \dint s \right)^2 \right\} \\
		\ge &   \mathbb{E}_{X^*}\left\{\int_{[0, 1]} \left  ( \int_{[0, 1]} \Delta_A  (s,r) X^*(r)  \dint   r  \right )^2 \dint s  \right\}  -  C_5'\log(T)  \zeta_n  
\end{align*} 
  \end{proof}
  
 \begin{lemma}
 Let $\zeta_n=   n^{-\alpha +1/2  } $.   
 Under the same conditions as in \Cref{thm-main-mixed}, it holds that 
\begin{align}  \frac{1}{ n }    \sum_{i=1}^n  \Delta_\beta(r _i )  ^\top   \Sigma_Z \Delta_\beta(r _i ) \ge 
	c         \int_{[0, 1]} \Delta_{\beta } (r )^\top \Sigma_Z \Delta_{\beta } (r )    \dint r  -C  \zeta _n     \label{eq:main 2 last step inequality 2}.
	\end{align} 
 \end{lemma} 
  \begin{proof}
   Note that  since the minimal eigenvalue of $\Sigma_Z$ is lower bounded by $c_z$,
	$$ \frac{1}{ n }    \sum_{i=1}^n  \Delta_\beta(r _i )  ^\top   \Sigma_Z \Delta_\beta(r _i )  \ge 
	\frac{c_z }{ n }    \sum_{i=1}^n   \sum_{l =1}^p \Delta^2 _{\beta_ l } (r _i )  . $$
	By  \Cref{lemma:discrete approximation},  there exists  an absolute constant      $C_4$ such that 
	\begin{align}\label{eq:step 4 main 2} 
	 \frac{1  }{n}\sum_{i=1}^n  \Delta^2 _{\beta_ l } (r _i )    \ge        \int_{[0, 1]} \Delta_{\beta_ l }^2 (r )  \dint r    - C _4 
	  \zeta _n       \| \Delta _{\beta_ l } ^2 \|_ \wot , \quad \forall  l  = 1, \ldots, p  .  
	  \end{align}
 As a result, 
 	\begin{align*}
 	\frac{1}{ n }    \sum_{i=1}^n  \Delta_\beta(r _i )  ^\top   \Sigma_Z \Delta_\beta(r _i )  \ge 	\frac{c_z }{ n }    \sum_{i=1}^n   \sum_{ l =1}^p \Delta^2 _{\beta_ l } (r _i )  \ge & c  _4  c_z   \sum_{   l =1}^p   \int_{[0, 1]} \Delta_{\beta_ l }^2 (r )  \dint r   - C _4  c_z   \sum_{  l =1}^p  \zeta _n       \| \Delta _{\beta_ l } ^2 \|_\wot    \\
 		\ge & c_4'      \sum_{ l =1}^p   \int_{[0, 1]} \Delta_{\beta_l }^2 (r )  \dint r  -C_4'  \zeta _n     , 
 	\end{align*}
	where the second inequality follows from  \eqref{eq:step 4 main 2}, and the third  inequality follows from 
	\begin{align*}
		 \sum_{ l =1}^p     \| \Delta _{\beta_l } ^2 \|_\wot     \le  \sum_{ l =1}^p     \| \Delta _{\beta_l }\|_\wot  ^2 \le&  \sum_{ l =1}^p     \| \Delta _{\beta_l }\|_\hb ^2 
		  \le  2 \sum_{ l =1}^p \big(      \| \beta _l ^*\|_\hb^2   +  \|  \widehat \beta_ l   \|_\hb   ^2 \big) 
		  \\   
		  \le&   2\big \{   \sum_{ l =1}^p      \| \beta_ l ^*\|_\hb \big  \}  ^2   + 2 \big\{  \sum_{ l =1}^p   \|  \widehat \beta_ l   \|_\hb   \big\} ^2  \le  4C_\beta ^2 .  
	\end{align*} 
	Therefore 
	\begin{align} \nonumber \frac{1}{ n }    \sum_{i=1}^n  \Delta_\beta(r _i )  ^\top   \Sigma_Z \Delta_\beta(r _i ) \ge 
	& c_4'      \sum_{ l =1}^p   \int_{[0, 1]} \Delta_{\beta_ l }^2 (r )  \dint r  -C_4'  \zeta _n    
\\
\ge &  \frac{c_4'}{C_z }        \int_{[0, 1]} \Delta_{\beta } (r )^\top \Sigma_Z \Delta_{\beta } (r )    \dint r  -C_4'  \zeta _n   \label{eq:main 2 last step inequality}  
	\end{align}
	where  the last inequality follows from the fact that $w^\top \Sigma w \le C_z  \|w\|_2^2  $ for all $w \in \mathbb R^p  $.

  \end{proof}

 \subsection{Extensions}
 \label{subsection:proof of main 2}
 \begin{corollary}\label{cor-main}
Define the discretized excess risk as	
	\begin{align}
		\mathcal{E}_{\mathrm{dis}}^*(\widehat{A}, \widehat{\beta}) = & \mathbb{E}_{X^*, Z^*, Y^*}\left\{\frac{1}{n_2}\sum_{j=1}^{n_2} \left(Y^*(r_j) - \frac{1}{n_1}\sum_{i=1}^{n_1} \widehat{A}(r_j, s_i) X^*(s_i) - \langle Z^*, \widehat{\beta}(r_j) \rangle_p\right)^2\right\} \nonumber \\
		-& \mathbb{E}_{X^*, Z^*, Y^*}\left\{\frac{1}{n_2}\sum_{j=1}^{n_2} \left(Y^*(r_j) -  \frac{1}{n_1}\sum_{i=1}^{n_1}   A^*(r_j, s_i) X^*(s_i)\dint s - \langle Z^*, \beta^*(r_j) \rangle_p \right)^2\right\}. \label{eq-excess-add-term}
	\end{align}
Suppose that Assumptions~\ref{assume:regularity}, \ref{assume:mixed} and \ref{assume:joint smooth} hold. 
Let $(\widehat{A}, \widehat{\beta})$ be any solution to \eqref{eq:approx A beta}  with the tuning parameter $\lambda =C_\lambda \sqrt {\frac{\log(p\vee T)}{T}}$ for some sufficiently large constant $C_\lambda$. Define $n=\min\{n_1,n_2\}$. 
For $ T\gtrsim \log(n)$, there exists absolute constants $C    > 0$ such that with probability at least 
 $1 -  8T^{-4}   $, it holds that
\begin{align}\label{eq-thm-main-mixed-result_dis}
	& \mathcal{E}_{\mathrm{dis}}^*(\widehat{A}, \widehat{\beta})  \le  C   \log(T)  \big \{   \delta_T +   \mathfrak {s} \log(p \vee T )/T    +   \zeta_n        \big \}    
\end{align}
where $\zeta_n = n^{-\alpha +1/2}$   and  $\delta_T = T^{\frac{-2r}{2r+1 }}$.
 \end{corollary}  
\Cref{cor-main} is a direct consequence of \Cref{thm-main-mixed} and provides a formal theoretical guarantee for using $\mathcal{E}_{\mathrm{dis}}^*(\widehat{A}, \widehat{\beta})$ to evaluate the proposed estimators in practice.

 \subsection{Proof of \Cref{thm-main-mixed noise}}
 \begin{proof}[Proof of \Cref{thm-main-mixed noise}]
 The proof is almost identical to the proof of \Cref{thm-main-mixed}. As a result, we only point out the difference. 

 From the     minimizer  property, we have that
	\begin{align*}
		& \frac{1}{Tn}\sum_{t=1}^T \sum_{i=1}^n \wri  \left\{ y_{t,i}-  \frac{1}{n }   \sum_{j=1}^n \wsj   \widehat A(r_i, s_j) X_t(s_j)    -  \lb Z_t ,  \widehat \beta(r_i) \rb _p \right\}^2   + \lambda \sum_{j=1}^p  \|\widehat \beta_j   \|_n \\
		\le &  \frac{1}{Tn}\sum_{t=1}^T \sum_{i=1}^n  \wri \left\{ y_{t,i} -      \frac{1}{n }   \sum_{j=1}^n  \wsj    A^* (r_i, s_j) X_t(s_j) - \lb Z_t , \beta ^* (r_i) \rb _p \right\}^2  + \lambda \sum_{j=1}^p \|\beta_j  ^*  \|_n  ,
	\end{align*}
	which implies that 
	\begin{align} 
 		& \frac{1}{Tn}\sum_{t=1}^T \sum_{i=1}^n \wri  \left\{\frac{1}{n }   \sum_{j=1}^n \wsj  \Delta _A (r_i, s_j) X_t(s_j)      \right\}^2 \label{eq:CAFR term 0 noise} \\ 
		  + & \frac{1}{Tn}\sum_{t=1}^T \sum_{i=1}^n \wri \lb \Delta_\beta (r_i) , Z_t\rb^2 _p   \\
		  +& \frac{2}{Tn}\sum_{t=1}^T \sum_{i=1}^n \wri \left\{\frac{1}{n }   \sum_{j=1}^n \wsj   \Delta _A (r_i, s_j) X_t(s_j)\lb \Delta_\beta (r_i) , Z_t\rb _p \right\} \label{eq:CAFR term 2 noise} \\
 		\le & \frac{2}{Tn} \sum_{t=1}^T \sum_{i=1}^n     \wri   \left( \frac{1}{n }   \sum_{j=1}^n \wsj   \Delta _A (r_i, s_j) X_t(s_j)  +     \lb \Delta_\beta (r_i) , Z_t\rb  _p  \right)\epsilon_t(r_i)   \label{eq:CAFR term 4 noise} 
 		\\
 		+& \lambda \sum_{j=1}^p \|\beta_j  ^*  \|_n  - \lambda\sum_{j=1}^p  \| \widehat \beta_j     \|_n  \label{eq:CAFR term 5 noise} 
 		\\
 		+ & \frac{2}{Tn} \sum_{t=1}^T \sum_{i=1}^n     \wri   \left( \frac{1}{n }   \sum_{j=1}^n \wsj   \Delta _A (r_i, s_j) X_t(s_j)  +     \lb \Delta_\beta (r_i) , Z_t\rb  _p  \right)\mathfrak{E}_{t,i }  . \label{eq:CAFR term 6 noise} 
	\end{align}
	Note that \Cref{eq:CAFR term 0 noise} - \eqref{eq:CAFR term 5 noise} are     identical to \Cref{eq:CAFR term 0} - \eqref{eq:CAFR term 5}. So it suffices to analyzed  \Cref{eq:CAFR term 6 noise}.
	Without loss of generality,  we will assume that 
	$$\wri =1  \text{ for all } 1\le i \le n \quad \text{and} \quad  \wsj= 1\text{ for all }   1\le j \le n  .$$
	\
	\\
It follows from  \eqref{eq:FFR step 2 3} and \eqref{eq:sc form term 2 bound 23}, we have that with probability at least 
	$1 -2T^{-4}$,
	it holds that
	\begin{align*}
		 \eqref{eq:CAFR term 6 noise} \le&  \frac{1 }{3 2 0 n} \sum_{i=1}^n \iint_{[0, 1]^2} \Delta_A (r_i,  r) \Sigma_X  (r,s)  \Delta_A(r_i,  s ) \dint r \dint s 
		  + C_1   \Bigg\{   \zeta_n   +   \log(T)\delta_T \Bigg\} 
		  \\
		  + & \frac{  \lambda }{ 320 }       \sum_{j \in S }   \| \Delta_{\beta_j }\|_n    +   \frac{\lambda }{ 320 }     \sum_{j \in S ^c}   \|  \widehat \beta_j  \|_n  
	\end{align*}
	where $C_1 > 0$ is an absolute constant. The rest of the argument is identical to \Cref{thm-main-mixed} and therefore is omitted. 
 \end{proof}

\section{Deviation bounds}  \label{sec-app-func}

\input{appendix/functional-deviation}

\section{Proof of the lower bound result}\label{sec-app-lb}

\begin{proof}  [Proof of \Cref{lem-lb111}] Since the proof is about lower bounds, it suffices to assume that  $1 \leq  \mathfrak {s} < p/3$ and   $T \geq  \mathfrak {s}^2  \log(2p)/24$.
\
\\
	{\bf Step 1.}  	Consider the following two special cases of model \eqref{eq-model-mixed-cov}:
	\begin{align} \label{eq:model 2}
		Y(r) =& \int_{[0, 1]} A^*(r, s)X(s)\dint s   + \epsilon(r), \quad r\in[0,1]
		\\
 \label{eq:model 3}
		Y(r) =& \sum_{j = 1}^p Z_j\beta^*_j(r) + \epsilon(r), \quad r\in[0,1].
	\end{align}
   Denote 
	\[
		\mathcal{E}^*(\widehat{A}) =  \int_{[0, 1]} \left\{\iint_{[0, 1] \times [0, 1]} \Delta_A(r, s_1) \Sigma_X(s_1, s_2) \Delta_A (r, s_2)\dint s_1 \dint s_2\right\}\dint r 
	\]
	and
	\[
		\mathcal{E}^*(\widehat{\beta}) =  \int_{[0, 1]} \Delta_\beta^\top(r) \Sigma_Z \Delta_\beta (r) \dint r,
	\]
	as the excess risks of model \eqref{eq:model 2} and model \eqref{eq:model 3} respectively. Since both model \eqref{eq:model 2} and model \eqref{eq:model 3} are special cases of model \eqref{eq-model-mixed-cov}, standard minimax analysis shows that 
\begin{align*}  \inf_{ \widehat A , \widehat \beta}    \sup_{\substack{A ^*  \in  \mathcal C_A, \\ \beta ^*  \in \mathcal{C}_\beta }} \mathbb   E\{ \mathcal{E}^*(\widehat{A}, \widehat{\beta})\}  
\ge&  \max\bigg\{ \inf_{ \widehat A , \widehat \beta } \sup_{\substack{A ^*  \in  \mathcal C_A, \\ \beta^*  = 0}}\mathbb E\{\mathcal{E}^*(\widehat{A} )\}  , \ 
\inf_{ \widehat A , \widehat \beta } \sup_{\substack{A ^* =0\   \\ \beta^*  \in \mathcal{C}_\beta }}\mathbb E\{\mathcal{E}^*(\widehat{\beta } )\}   \bigg\} 
\\
= & \max\bigg\{  \inf_{ \widehat A   } \sup_{\substack{A ^*  \in  \mathcal C_A, \\ \beta^*  = 0}}  \mathbb E\{\mathcal{E}^*(\widehat{A} )\} , \ \inf_{ \widehat \beta    } \sup_{\substack{A ^* =0, \\ \beta^*  \in \mathcal{C}_\beta} }  \mathbb E\{\mathcal{E}^*(\widehat{\beta } )\}  \bigg\}.
\end{align*} 
	 
By the above arguments, the task of finding a lower bound on the excess risk $\mathcal{E}^*(\widehat{A} ,\widehat{\beta})$ can be separated into two tasks of finding lower bounds on 
	\begin{align}\label{eq-eeeeee4}
		\inf_{ \widehat A   } \sup_{\substack{A ^*  \in   \mathcal C_A, \\ \beta^*  = 0}}  \mathbb E\{\mathcal{E}^*(\widehat{A} )\} \quad \mbox{and} \quad \inf_{ \widehat \beta    } \sup_{\substack{A ^* =0, \\ \beta^*  \in \mathcal{C}_\beta} }  \mathbb E\{\mathcal{E}^*(\widehat{\beta } )\}.
	\end{align}
  Theorem~2  in   \cite{sun2018optimal} shows that under the same eigen-decay assumption  in condition \eqref{eq:eigen decay} of  \Cref{thm-main-mixed},  it holds that
	\begin{align}\label{eq:lower bound 1}
		\inf_{ \widehat A   } \sup_{\substack{A ^*  \in   \mathcal C_A, \\ \beta^*  = 0}}\mathbb E\{ \mathcal{E}^*(\widehat{A} )\} \ge c  T^{-\frac{2r}{2r + 1 } },
	\end{align}
where $c > 0$ is a sufficiently small constant.  
Therefore, to provide a lower bound on the excess risk $\mathcal{E}^*(\widehat{A} ,\widehat{\beta})$, the only remaining task is to provide a lower bound of $\inf_{ \widehat \beta    } \sup_{\substack{A ^* =0, \\ \beta^*  \in \mathcal{C}_\beta} }  \mathbb E\{\mathcal{E}^*(\widehat{\beta } )\}$ based on model \eqref{eq:model 3}.

\
\\
{\bf Step 2.} For $ 0<\delta < \sqrt { 1/ (2 \mathfrak {s})   }$ to be specified  later, let $\widetilde \Theta (\delta,  \mathfrak {s})    $  denote 
	\begin{align*}
		\mathcal B ( \mathfrak {s})    = \left \{ \{\beta_j \}_{j=1}^p \subset \mathcal{H}_1; \, \sum_{j=1}^p \mathbbm{1}\{  \beta_j\not= 0 \} \le  \mathfrak {s},\,  \max_{j = 1, \ldots, p} \|\beta_j \|_{\h_1} \le  \delta\sqrt {2/ \mathfrak {s}}  \right\}, \quad  \mathfrak {s} \in \{1, \ldots, p\}.
	\end{align*}
	Since $\delta \sqrt{2 \mathfrak {s}} \le 1  $, 
  $\mathcal B( \mathfrak {s}) \subset \mathcal C_\beta$ with $C_\beta = 1  $. 
Let $\phi \in \hb$ be the leading eigenfunction of $\hb $.  Since $\hb  $ is generated by a bounded kernel, we have that $\|\phi\|_{\infty} \leq 1$ and $\|\phi\|_{\mathcal{L}^2} = 1$.
For any $\theta \in \mathbb R^p $ such that $\|\theta\|_{\infty} \le 1$, denote $\beta^{\theta} = (\beta^{\theta}_j, \, j = 1, \ldots, p) \subset \hb $ be such that 
	\[
		\beta_j^{\theta}(\cdot) = \theta_j \phi(\cdot), \quad j = 1, \ldots, p.
	\]
	Provided that $\|\theta\|_0 \le  \mathfrak s $, we have that $\beta^{\theta} \in \mathcal B( \mathfrak {s}) \subset  \mathcal{C}_\beta$.  For $t \in \{1, \ldots, T\}$, let $\epsilon_t (s) = \varepsilon_t \phi(s)$, where $ \{ \varepsilon_t\}_{t=1}^T | \overset {\mbox{i.i.d.}}{\sim} \mathcal{N}(0,1)$.  We thus have $\varepsilon_t $ is a Gaussian process with bounded second moment. 

Denote $P_\theta$ as the joint distribution of $\{ Y_t ,Z_t \} _{t=1}^T$ with $\beta^* = \beta^{\theta}$.  Since $ P_\theta $ is supported on the subspace spanned by $\phi$, for any $\theta, \theta'$, the Kullback--Leibler divergence between $ P_{\theta}  $ and $ P_{\theta'}  $ is 
	\[
		\mathrm{KL}( P_\theta  |P_{\theta'}  ) =  \mathbb{E}_{  P_\theta } \left\{\log \left ( \frac{\dint P_{\theta} }{ \dint P_{\theta'} }  \right)  \right\} =  \mathbb{E}_{  Q_\theta } \left\{\log \left ( \frac{\dint Q_{\theta} }{\dint Q_{\theta'} }  \right)  \right\},
	\]
	where $Q_\theta$ is the joint distribution of $\{y_t, Z_t \}_{t=1}^T $ with 
	\[
		y_t = \lb Z_t, \theta\rb_p +\varepsilon_t, \quad t \in \{1, \ldots, T\}.
	\] 
	Since $\{Z_t\}_{t=1}^T$ and $\{ \varepsilon_t\}_{t=1}^T  $ are independent, we have that
	\[
		\log \left ( \frac{\dint Q_{\theta} }{\dint Q_{\theta'} }   \right) = \sum_{t=1}^T \left( y_t\lb Z_t, \theta \rb_p   - \frac{1}{2} \lb Z_t ,\theta\rb_p ^2  -  y_t\lb Z_t, \theta' \rb_p   + \frac{1}{2} \lb Z_t ,\theta' \rb_p ^2  \right)   .
	\]
	Conditioning on $\{Z_t\}_{t=1}^T$, we have that $y_t \sim \mathcal{N}( \lb Z_t, \theta\rb_p, 1 )$. Therefore,
	\begin{align*}
	\mathrm{KL} ( P_\theta  |P_{\theta'}  ) = & \mathbb{E}_{  Q_\theta } \left\{\log \left ( \frac{\dint Q_{\theta} }{\dint Q_{\theta'} }  \right)  \right\} = \mathbb{E}_{\{ Z _t\}_{t=1}^T}  \left\{ \mathbb{E}_{Q_\theta| \{ Z _t\}_{t=1}^T}   \log \left ( \frac{\dint Q_{\theta} }{\dint Q_{\theta'} }   \right) \right\} \\
	= & \mathbb{E}_{\{ Z _t\}_{t=1}^T}  \left(\sum_{t=1}^T   \lb Z_t, \theta-\theta'\rb_p ^2  \right)  =T (\theta-\theta')^\top  \Sigma_Z (\theta-\theta') = T \| \theta-\theta'\|_2^2. 
	\end{align*}
	For $ 0<\delta < \sqrt {1/ (2 \mathfrak{s})  }$ to be specified  later, let $\widetilde \Theta (\delta,  \mathfrak {s})    $ be defined as in \Cref{lemma:packing of l_q ball}. Then for any $\theta \in  \widetilde \Theta (\delta,  \mathfrak {s}) $, we have that
	\[
		\|\theta\|_{1} \le \delta\sqrt {2  \mathfrak {s} }\le 1 
	\]
	and thus $\beta(\theta) \in \mathcal \mathcal{C}_\beta$. 
For $ \theta   \not =  \theta'   $ in $\widetilde { \Theta } (\delta, \mathfrak {s})  $, it holds that
\begin{align*}
 \sum_{j=1} ^p  \| \beta_j (\theta) -\beta_j (\theta' ) \|_\lt ^2 = \sum_{j=1}^p (  \theta_j -\theta_j')^2 \| \phi\|_\lt  ^2  =  \| \theta-\theta'\|_2^2  \ge  \delta^2
 \end{align*} 
 and that 
 \begin{align*}
 \mathrm{KL} ( P_\theta  |P_{\theta'}  ) 
  =  T \| \theta-\theta'\|_2^2 \le 8 T \delta^2 .
\end{align*}
\
\\
Using Fano's lemma \citep[e.g.][]{Yu1997}, we have that 
	\[
		\mathbb{P} \left(   \inf_{\widehat  \beta } \sup_{ \beta^* \in \mathcal{B}(s) }  \sum_{j=1}^p     \|  \widehat \beta_j  -\beta^*_j\|_\lt ^2 \ge \delta ^2 /4  \right) \ge 1- \frac {  \frac{ 2}{   M(M-1)  }  \sum_{\substack{\theta, \theta' \in \widetilde { \Theta } (\delta,s) \\ \theta \neq \theta'}} \mathrm{KL} ( P_\theta  |P_{\theta'}  )  + \log(2)     } { \log(M)  },
	\]
	where 
	\[
		M =  | \widetilde \Theta (\delta,  \mathfrak {s})| \geq \exp\left\{\frac{s}{2} \log\left(\frac{ p - \mathfrak {s} }{ \mathfrak {s}/2}\right)\right\}.
	\]
	Note  that 
	\begin{align*}
		\frac {  \frac{ 2}{   M(M-1)  }  \sum_{\substack{\theta, \theta' \in \widetilde { \Theta } (\delta, \mathfrak {s}) \\ \theta \neq \theta'}} \mathrm{KL} ( P_\theta  |P_{\theta'}  )  + \log(2)     } { \log(M)  } \le       \frac {   8T\delta^2 + \log(2)     } {\frac{ \mathfrak {s}}{2} \log( \frac{p-  \mathfrak {s} }{ \mathfrak {s}/2 })    } .
	\end{align*}
	Provided that $ \delta ^2  = \frac{   \mathfrak {s}}{48T }   \log( \frac{ p - \mathfrak {s} }{ \mathfrak {s}/2 })  $, we have that
	\[
		\mathbb{P} \left\{ \inf_{\widehat  \beta } \sup_{ \beta^* \in \mathcal{B}( \mathfrak {s}) }  \sum_{j=1}^p     \|  \widehat \beta_j  -\beta^*_j\|_\lt ^2 \ge  \frac{ \mathfrak {s}}{200T }   \log\left( \frac{p - \mathfrak {s} }{ \mathfrak {s}/2 }\right)  \right\} \ge 1-  1/2   .
	\]
	Since $ T\ge  \frac{ \mathfrak {s}^2 \log(2p) }{24}$, it holds that 
	\[
		2 \mathfrak {s} \delta ^2   = \frac{  \mathfrak {s}^2}{24 T }   \log\left( \frac{ p - \mathfrak {s} }{  \mathfrak {s}/2 }\right)   \le   \frac{  \mathfrak {s}^2  \log(2 p) }{24 T }    \le 1,
	\]
	hence $\delta < \sqrt { 1/(2 \mathfrak {s})    }$ as desired.  Since $\Sigma_Z$ is positive definite , it holds that $v^{\top}\Sigma_Zv \geq \kappa \|v\|^2$ for any $v \in \mathbb{R}^p$. 
\
\\
Finally, we have that  
	\[
		\mathbb{P} \left\{\inf_{\widehat  \beta } \sup_{ \beta^* \in \mathcal{B}( \mathfrak {s}) } \sum_{j=1}^p \int_{[0, 1]} \Delta_\beta^\top(r) \Sigma_Z \Delta_\beta (r) \dint r \ge  \frac{  \mathfrak {s}}{200T }   \log\left( \frac{p - \mathfrak {s} }{ \mathfrak {s}/2 }\right)  \right\} \ge   1/2,
	\]
	which directly implies the desired result. 
\end{proof}

\begin{lemma} \label{lemma:packing of l_q ball}
For every $\delta>0$ and $ \mathfrak {s}< p/3$, denote 
	\[
		\Theta(\delta, s) = \left\{ \theta\in \mathbb \{ -  \delta \sqrt {2 /    \mathfrak {s}    }, 0,    \delta  \sqrt { 2/  \mathfrak {s}    } \}^p:\, \| \theta\|_0 = \mathfrak {s} \right\} .
	\]
	There exists a set $ \widetilde \Theta(\delta,  \mathfrak {s})   \subset \Theta(\delta,   \mathfrak {s})$, such that 
	\[
		| \widetilde \Theta (\delta,  \mathfrak {s})   | \ge \exp\left\{\frac{ \mathfrak {s}}{2} \log\left(\frac{ p- \mathfrak {s} }{ \mathfrak {s}/2 }\right)\right\}.
	\] 
	This implies that for any $\theta, \theta' \in \widetilde \Theta (\delta,  \mathfrak {s})$, $\theta \neq \theta'$, it holds that 
	\[
		\|\theta-\theta'  \|_2^2 \ge \delta ^2 \quad \text{and} \quad  \|\theta-\theta'  \|_2^2 \le 8 \delta ^2.
	\]
\end{lemma}
\begin{proof}\Cref{lemma:packing of l_q ball} is Lemma 4 in \cite{raskutti2011minimax}.
\end{proof}

\section{Additional technical results}  \label{sec-app-additional}

   \begin{theorem}
  \label{theorem:restricted eigenvalue  condition}
 Suppose  that $Z_i \sim \mathcal N(0, \Sigma_Z)$. Then  there exist  universal constants $c,C>0$ such that 
   $$\mathbb  P\left(  \frac{1}{T} \sum_{t =1}^ T  (Z_t ^\top v)^2 \ge \frac{2}{3} v^\top \Sigma_Z v      - C  \frac{ \log(p) }{T }  \|v\|_1^2 
  \quad  \text{for all }  v \in \mathbb R^p   \right)  \le \exp(-cT ) .$$
 \end{theorem}  
 \begin{proof} This is the well known Restricted Eigenvalue condition. The   proof can be found in Theorem 1 of \cite{raskutti2010restricted}. See     \cite{loh2012high} for better constants. 
\end{proof}

\begin{lemma}
Suppose $\{ w_ i \}_{i=1}^n $ are i.i.d. centered sub-Gaussian random variables with variance 1 (denoted as SG(1)).Then 
$$ \mathbb P \left ( \left|  \frac{1}{n} \sum_{i=1}^n w_i \right|  \ge \gamma  \right )\le 2\exp(-2\gamma^2n  ) .$$

\end{lemma}

\begin{lemma} \label{lemma:sub-exponential tail}
Suppose $\{ z_ i \}_{i=1}^n $ are i.i.d.  centered sub-Exponential random variables with parameter $1$.   Then 
$$ \mathbb P \left ( \left|  \frac{1}{n} \sum_{i=1}^n z _i \right|  \ge \gamma  \right )\le 2\exp(-2n \min \{  \gamma^2, \gamma\}  ) .$$

\end{lemma} 
 
 See \cite{vershynin2018high} for the proofs of these   two lemmas. 

\begin{lemma} \label{lemma:tail bound with log factor}
Suppose $r> 1/2 $.
For any $p\ge D$, it holds that
$$   \sum_{k=D+1 }^\infty   \frac{   \log(k)     }{ k^{ 2  r  }  }     \le C_r   \left(  \frac{          \log(p) }{  D ^{ 2r- 1   }  }    +    \frac{ 1     }{ p ^{   r- 1/2  } }  \right)  . $$
\end{lemma}
\begin{proof}
Let $\nu \in \mathbb R$  $  \nu  =r-1/2>0   $.
\begin{align*}
    \sum_{k=D+1 }^\infty   \frac{   \log(k)     }{ k^{ 2  r  } }      
\le 
&      \sum_{k=D+1 }^ p   \frac{    \log(k)     }{ k^{ 2  r  } }     +    \sum_{k=p+1 }^ \infty    \frac{    \log(k)     }{ k^{ 2  r  } }        
\\
\le &    \frac{      \log(p) }{  D ^{ 2r- 1   } }    + C_r  \sum_{k=p+1 }^ \infty    \frac{ k^\nu     }{ k^{ 2  r  } }        
\\
\le 
&  \frac{          \log(p) }{  D ^{ 2r- 1   } }    +   C_ r  \frac{ 1     }{ p ^{ 2  r- 1- \nu  } }     ,
\end{align*}
where the second inequality holds because there exists a constant $ C_r $ such that  $\log(k) \le C_\nu  k^\nu  $ for all $k\ge 1 $.
\end{proof}

\begin{theorem}[Hanson--Wright]
Let $X=(X_1, \ldots, X_n) \in \mathbb R^n $ be a random vector where the components  $X_i$ are independent SG($K$). Let $A $ be any $n\times n$ matrix. Then 
 for every $t\ge 0$, 
 $$ \mathbb  P \left( \left| X^\top AX - E ( X^\top AX ) \right| \ge  t \right) \le 2 \exp \left [ -c \min \left(  \frac{t^2}{ K^4\| A\|_F^2 }      , \frac{ t } {   K^2 \|A\|_{ \text{op} }  }  \right)   \right]  .$$ 
\end{theorem} 
 \begin{proof} 
 The proof of this theorem can be found in    \cite{vershynin2018high}.
\end{proof}  
 
\begin{lemma} \label{lemma:Gaussian process norm}
Suppose $X $ is a centered Gaussian process with  covaraince function $\Sigma_X  $. Suppose that 
$$ \mathbb E (\|X\|_\lt ^2 ) =C_X <\infty  .$$
Then there exist  absolute constants $c, C>0$ only depending on $C_X$  such that 
$$\mathbb P \bigg ( \| X\|_\lt^2 \ge C   \eta \bigg ) \le \exp(-c\min\{ \eta, \eta^2\}).$$
\end{lemma}
\begin{proof}
 Let 
$$ \Sigma _X (r,s)  = \sum_{k=1}^\infty {\sigma_k^2}   \psi^x_k (r)  \psi^x_k (s)  $$
be the eigen-expansion of $\Sigma_X $.  
Let $z_k = \lb X, \psi_k^x\rb  $. Observe that $\{z_k\}_{k=1}^\infty$ is a collection of independent Gaussian random variables such that 
$z_k \sim N(0, \sigma_k^2) $. Since $E(\|X\|_\lt^2 )   =C_X $,
we have that 
$$ \mathbb E (\|X\|_\lt^2 ) = \sum_{k=1}^\infty \mathbb E (z_k^2 ) =   \sum_{k=1}^\infty \sigma_k^2  =  C_X. $$
Observe that 
 $ \|X \|_\lt^2 = \sum_{k=1}^\infty z_k^2$ 
is a sub-Exponential   random variable with parameter $\sum_{k=1}^\infty \sigma_k^4 \le C_X ^2   $. 
So
by sub-Exponential tail bound, it holds that for any $\eta>0$, 
$$\mathbb P \bigg ( \| X\|_\lt^2 \ge C' \big\{ \sum_{k=1}^\infty \sigma_k^4 \big\} \eta \bigg ) \le \exp(-c\min\{ \eta, \eta^2\}),$$
where $C'$ and $c$ are absolute constants. 
\end{proof} 
\begin{lemma} \label{lemma:Gaussian process sobolev norm}
Suppose $X $ is a centered Gaussian process with  
$$\mathbb E (\|X\|_\wot^2 ) =C_X <\infty . $$  
Then there exist   absolute constants $c, C>0$ depending only on $ C_X$ such that 
$$\mathbb P \bigg ( \| X\|_\wot ^2 \ge C   \eta \bigg ) \le \exp(-c\min\{ \eta, \eta^2\}).$$
\end{lemma}
\begin{proof} Let $\mathbb K_{\alpha } $ be the generating kernel of  $ \wot  $ with eigen-expansion 
$$\mathbb K_{\alpha }(r,s)  = \sum_{k=1}^\infty \omega_k   \psi_k^\alpha (r) \psi_k^\alpha  (s). $$
Note that we have $\omega_k \asymp  k^{-2\alpha } $ and  
$ \| \psi_k\|_\wot ^2  = \omega_k^{-1}  \asymp k^{2\alpha}$ due to the property of Sobolev space. Let 
 $ b_k = \lb X, \psi_k^\alpha  \rb_\lt .  $ So $ \{ b_k\}_{k=1}^\infty$ is a collection of  centered correlated Gaussian random variables  
and $X= \sum_{k=1}^\infty  b_k \psi_k^\alpha .$
\\
\\
{\bf Step 1.} Since $X\in \wot$ almost surely, by \Cref{remark:inverse of K},   
$Y : = L_{\mathbb K ^{-1/2 }_{\alpha }   } (X)  \in \lt $ is well defined. 
 So  $X= L_{\mathbb K ^{ 1/2 }_{\alpha }   } (Y) $ and that 
$$ \mathbb E (\| L_{\mathbb K ^{ 1/2 }_{\alpha }   }  (Y)\|_\wot ^2 )=C_X. $$
From 
 \Cref{lemma:connect L_2 to H}, $\| L_{\mathbb K ^{ 1/2 }_{\alpha }   }  (Y)\|_\wot ^2 = \| Y\|_\lt^2.$ 
 So $  E (\|   Y \|_\lt ^2  )=C_X. $  For any generic function $f\in \lt$ with $\|f\|_\lt <\infty$, 
 then 
 $$ \lb Y, f\rb _\lt  =  \bigg \lb L_{\mathbb K ^{-1/2 }_{\alpha }   } \big ( \sum_{k=1}^\infty b_k \psi^\alpha_k  \big )  ,f \bigg \rb_\lt = \sum_{k=1}^\infty b_k \bigg \lb L_{\mathbb K ^{-1/2 }_{\alpha }   } \big (   \psi^\alpha_k  \big )  ,f \bigg \rb _\lt , $$
 where $L_{\mathbb K ^{-1/2 }_{\alpha }   } \big (   \psi^\alpha_k  \big ) \in \lt  $ because $  \psi^\alpha_k  \in \wot $.
 Therefore $\lb Y, f\rb _\lt  $ is a  centered Gaussian random variable with  
 $$Var( \lb Y, f\rb_\lt  )  = \mathbb E(\lb Y, f\rb_\lt ^2  ) \le \mathbb  E \{ \|Y \|_\lt^2 \|f\|_\lt^2 \}  =C_X \|f\|_\lt ^2  <\infty. $$
 \
\\
 { \bf Step 2.}  Let $\Sigma_Y $ be the covariance function of $Y$. Since $E(\|Y\|_\lt^2) <\infty $, by the Hilbert–Schmidt theorem, the  eigen-expansion  of $ \Sigma_Y$ can be written as 
 $$ \Sigma_Y (r,s) = \sum_{k=1}^\infty \delta_k ^2 \psi_k^y (r) \psi_k^y(s) .$$
 Then by {\bf Step 1}, $ \lb Y, \psi_k^y\rb_\lt $ is Gaussian and that 
 $$  \mathbb E \bigg\{ \lb Y, \psi_k^y\rb_\lt \lb Y, \psi_l^y\rb_\lt  \bigg\} = \Sigma_Y [ \psi_k^y , \psi_l ^y] = 0 $$
 whenever $l \not = k  $. So $\{ \lb Y, \psi_k^y\rb_\lt \} _{k=1}^\infty $ is a collection of independent Gaussian random variables. 
 By \Cref{lemma:Gaussian process norm},
 $$\mathbb P \bigg ( \| Y\|_\lt^2 \ge C   \eta \bigg ) \le \exp(-c\min\{ \eta, \eta^2\}).$$
The desired result follows by observing that 
$$ \| Y\|_\lt^2 =  \|L_{\mathbb K ^{-1/2 }_{\alpha }   } (X) \|_\lt^2 =  \|X\|_\wot^2. $$
\end{proof}
 
 \begin{lemma} \label{lemma:quadratic in functional space}
 Suppose that $\{X_t\}_{t=1}^T $ is a collection of independent  centered Gaussian process with covariance operator $\Sigma_X $.  Let 
$$ \Sigma _X (r,s)  = \sum_{k=1}^\infty {\sigma_k^2}   \psi^x_k (r)  \psi^x_k (s)  $$
be the eigen expansion of $\Sigma_X $.  Suppose 
 $\E(\| X_t \|_\lt^2) =  \sum_{k=1}^\infty \sigma^2 _k    < \infty       $ and that $\{ \sigma_k\}_{k=1}^\infty  $ decay to $0$ at  polynomial rate.  Let $\{ z_t \} _{t =1}^T  \overset{i.i.d.}{\sim}  \mathcal N(0,  1) $. 
 If in addition, $ T \ge \log(p) $, then with probability  at least $1-10p^{-4}$
 $$  \bigg\|  \frac{1}{ T }\sum_{ t =1}^T X_ t  z_t  -  \E (X_1 z_ 1)\bigg \|_\lt  \le  C\sqrt {\frac{\log(p) }{ T  }}. $$  
 \end{lemma}
\begin{proof}
Let $ w_{ t ,k  } = \lb X_t  , \psi^x_k\rb_\lt   $. Then 
$ \{ w_{t ,k}\}_{1\le t  \le T , 1\le k<\infty }$ is a collection of centered independent random variables with 
$Var(w_{t ,k})={\sigma_k^2} $.
Note that 
\begin{align}  \label{eq:l_2 deviation basis expansion}
     \bigg\|  \frac{1}{ T  }\sum_{t  =1}^T X_t  z_t   -  \mathbb E (X_1 z_ 1)\bigg \|_\lt^2 
    =
     \bigg\|  \frac{1}{ T  }\sum_{ t =1}^T   \sum_{k=1}^\infty  \big( w_{ t ,k } z_ t   - \mathbb E(w_{1   ,k } z_ 1    ) \big)  \psi_k^x \bigg \|_\lt^2  
    = 
     \sum_{k=1}^\infty \bigg( \frac{1}{T } \sum_{t  =1}^T  \bigg\{ w_{t  ,k}z_t   - \mathbb E(w_{1,k}z_1 )  \bigg\} \bigg)^2  ,
\end{align}
where the last inequality follows from the fact that $\{ \psi^x_k \}$ is a collection of basis functions in $\lt$.
\\
 Note that $w_{t , k} z_t $ is SE with parameter ${\sigma_k^2} $. 
 So 
 $$\mathbb  P\bigg( \bigg| \frac{1}{T  } \sum_{t =1}^T   w_{t  ,k}z_t  -\mathbb E(w_{1,k}z_1 ) \bigg|  \ge  \eta    \sigma _k \bigg) \le \exp(- n\{ \eta   ,\eta ^2\} ). $$
 Therefore 
 $$ \mathbb P\bigg( \bigg| \frac{1}{T  } \sum_{t  =1}^T    w_{t ,k}z_t  -\mathbb E(w_{1  ,k}z_ 1   ) \bigg|  \ge   C  \sigma_k   \bigg\{  \sqrt {\frac{\log(p) +\log(k) }{T }  } +\frac{\log(p) +\log(k) }{T }    \bigg\} \bigg) \le  \frac{1}{k^3p^4} . $$
 Since $\sum_{k=1}^\infty k^{-3} <3, $ with probability at least $ 3p^{-4} $,
 $$  \bigg| \frac{1}{T  } \sum_{ t =1}^T    w_{t ,k}z_t   -\mathbb E(w_{1,k}z_1 ) \bigg|  \le   C   \sigma_k \bigg\{  \sqrt {\frac{\log(p) +\log(k) }{T }  } +\frac{\log(p) +\log(k) }{T }   \bigg\}  \quad \text{ for all } 1\le k<\infty. $$
 Under this good event,   
 \Cref{eq:l_2 deviation basis expansion} implies that
\begin{align}  \label{eq:l_2 deviation basis expansion 2}
     \bigg\|  \frac{1}{ T  }\sum_{ t =1}^T  X_t z_t - \mathbb E (X_1 z_ 1)\bigg \|_\lt^2 
     \le 2 C{\sigma_k^2} \bigg\{  \frac{\log(p) +\log(k) }{T   }    + \bigg( \frac{\log(p) +\log(k)  }{  T } \bigg)  ^2   \bigg\} .
\end{align}
Since $ \sum_{k=1}^\infty  {\sigma_k^2}   <\infty$, $\sigma_k \asymp k^{-1-\nu} $ for some $ \nu>0$. Therefore  $ \sum_{k=1}^\infty  {\sigma_k^2}  \log^2(k)  <\infty$. Since in addition,  
 $ \frac{\log(p)}{T}<1,$ 
 \Cref{eq:l_2 deviation basis expansion 2} implies that 
 $$ \bigg\|  \frac{1}{ T  }\sum_{ t=1}^T  X_ t z_ t  -\mathbb  E (X_1 z_ 1)\bigg \|_\lt^2 \le C' \frac{\log(p)}{T  } .$$
\end{proof}

\begin{theorem} \label{lemma:low rank expansion}
Suppose $A:\hk  \to  \hk  $ is a compact linear operator.  There exist $ \{ \psi_k\}_{k=1}^\infty, \{ \omega_k\}_{k=1}^\infty \subset \hk $, two orthogonal basis in $\mathcal{L}^2$, such that the associated bivariate function $A$ can be written as 
	\[
		A (r, s) = \sum_{k=1}^\infty a_k \psi_k(s) \omega_k(r), \quad r, s \in [0, 1],
	\]
	where $\{a_k\}_{k = 1}^{\infty} \subset \mathbb{R}$.   
\end{theorem}
 \begin{proof} This is the well known spectral theory for compact operators on Hilbert space. See Chapter 5 of    \cite{brezis2010functional}  for a detailed proof.
\end{proof}    

\begin{lemma}   \label{lemma:bound of A 1}
Let $A: \hk  \to \hk $ be any Hilbert-Schmidt operator.  We have that
	\[
		\max \left \{\sup_{r\in [0,1]} \| A(r, \cdot  )   \|_\hk , \,\sup_{s\in [0,1]} \| A(\cdot ,s )   \|_\hk , \,\|A(\cdot, \cdot)\|_{\infty} \right\} \le C_\mk\|  A \|_{\fk },
	\]
	where $C_\mk$ is some constant only depending on $\mk$.
\end{lemma}

\begin{proof} 
Since $\hk $ is generated by a bounded kernel $\mk $,  we have that for any $f \in \hk$ and any $ r\in [0,1]$,
\begin{align*} f(r) =   \lb f, \mk(r, \cdot) \rb_\hk  \le \| f\|_\hk  \| K(r,\cdot)\|_\hk  =  \| f\|_\hk  \sqrt {  \lb \mk (\cdot , r) ,  \mk(\cdot ,r) \rb_\hk } 
\le    \| f\|_\hk  \sqrt { \mk (r,r) } .
\end{align*}
Since $\mk $ is a bounded kernel, letting $C_\mk : =\sup_{r}\mk(r, r)  $, we have that 
$$\| f\|_\infty \le \|f\|_\hk  C_\mk.$$
Without loss of generality,  it suffices to assume that $  C_\mk=1$. 
By \Cref{lemma:low rank expansion}, there exist $ \{ \psi_k\}_{k=1}^\infty, \{ \omega_k\}_{k=1}^\infty \subset \hk $, two orthogonal basis in $\mathcal{L}^2$, and $\{a_k\}_{k = 1}^{\infty} \subset \mathbb{R}$, such that
	\[
		A(r, s) = \sum_{k=1}^\infty  a_k \psi_k(r) \omega_k(s). 
	\]
	Since  $\|\psi_k\|_{\hk } = \|\omega_k\|_{\hk } = 1$, $k \in \mathbb{N}_*$, $\| A\|_{\fk }^2  = \sum_{k=1}^\infty a_k^2$.  
	\\
Observe that for any $f \in \hk $ such that $\|f\|_\h= 1$, it holds that $f = \sum_{k=1}^\infty b_k \omega_k$, where $\sum_{k=1}^\infty b_k^2 = 1$.  Therefore, for any $r \in [0, 1]$, we have that
	\begin{align*}  
		\|A(r, \cdot )\|_\hk = & \sup_{\|f\|_\hk =1  } \Big\lb \sum_{k=1}^\infty  a_k \psi_k(r) \omega _ k(\cdot) , f(\cdot)\Big\rb_\hk  = \sup_{ \sum_{k=1}^\infty b_k^2 = 1 } \sum_{k=1}^\infty  a_k \psi_k(r)  b_k \\
		\le & \sup_{ \sum_{k=1}^\infty b_k^2 = 1 } \sum_{k=1}^\infty  |a_k| |b_k| \le \sup_{ \sum_{k=1}^\infty b_k^2 = 1 }  \sqrt { \sum_{k=1}^\infty  a_k^2 } \sqrt { \sum_{k=1}^\infty  b _k^2} = \|A\|_{\fk }, 
	\end{align*}
	where the first inequality is due to that $\| \psi_k\|_{\infty} \le \| \psi_k\|_\hk = 1$, $k \in \mathbb{N}_*$.  Similar arguments show that $\sup_{s\in [0,1]}\|A(\cdot, s )\|_\hk \le \| A\|_{\fk} $.

Moreover, observe that for any fixed $r, s \in [0, 1]$, it holds that 
	\[
		A(r, s) = \lb A(r, \cdot), \mathbb K(s, \cdot)  \rb_\hk .  
	\]	
	Therefore it holds that
	\begin{align*}
		\sup_{s\in [0,1]}|A(r, s)| & \le   \sup_{s\in [0,1] } |  \lb A(r,\cdot), \mathbb K(s, \cdot)  \rb_\hk | \le \|A(r, \cdot)\|_\hk  \sup_{s\in [0,1] } \| \mathbb K(s, \cdot)  \|_\hk \le \| A(r,\cdot)\|_\hk, 
	\end{align*}
	where  $ \|\mathbb K(s, \cdot) \|_\hk^2 =\mathbb K(s,s ) \le 1$, $s \in [0, 1]$, is used in the last inequality.  Finally, we have that
	\[
		\|A(\cdot, \cdot)\|_{\infty} \le \sup_{r\in [0,1] }\|A(r, \cdot )\|_\hk  \le   \|A\|_{\fk },
	\]
	which concludes the proof.
\end{proof}

\begin{lemma}\label{lem-derivation-of-eq8}
It holds that 
	\begin{align*}
		 \mathcal{E}^*(\widehat{A}, \widehat{\beta}) =  \mathbb{E}_{X^*, Z^*, Y^*} \left\{  \int_{[0, 1]}  \left (\int_{[0, 1]}  \Delta_A  (r,s)   X^ *(s) \dint s   +    \langle Z^*,  \Delta_\beta (r)\rangle_p     \right)^2          \dint r\right\}.
	\end{align*}
\end{lemma}
    
\begin{proof}
We have that 
\begin{align*}
		 & \mathcal{E}^*(\widehat{A}, \widehat{\beta}) =  \mathbb{E}_{X^*, Z^*, Y^*}\left\{\int_{[0, 1]} \left(Y^*(r) - \int_{[0, 1]} \widehat{A}(r, s) X^*(s)\dint s - \langle Z^*, \widehat{\beta}(r) \rangle_p\right)^2\dint r\right\} \\
		& \hspace{0mm} - \mathbb{E}_{X^*, Z^*, Y^*}\left\{\int_{[0, 1]} \left(Y^*(r) - \int_{[0, 1]} A^*(r, s) X^*(s)\dint s - \langle Z^*, \beta^*(r) \rangle_p \right)^2\dint r\right\}  \\
		= & \mathbb{E}_{X^*, Z^*, Y^*} \bigg\{\int_{[0, 1]} \bigg(\int_{[0, 1]} A^*(r, s) X^*(s)\dint s + \langle Z^*, \beta^*(r) \rangle_p  + \epsilon^*_t(r) \\
		&  - \int_{[0, 1]} \widehat{A}(r, s) X^*(s)\dint s - \langle Z^*, \widehat{\beta}(r) \rangle_p\bigg)^2\dint r\bigg\} - \mathbb{E}_{X^*, Z^*, Y^*}\left\{\int_{[0, 1]} \{\epsilon^*_t(r)\}^2\dint r\right\}  \\
		= & \mathbb{E}_{X^*, Z^*, Y^*} \left\{\int_{[0, 1]} \left(\int_{[0, 1]}  \Delta_A  (r,s)   X^ *(s) \dint s   +    \langle Z^*,  \Delta_\beta (r)\rangle_p \right)^2 \dint r\right\} \\
		& + 2\mathbb{E}_{X^*, Z^*, Y^*} \left\{\int_{[0, 1]} \left(\int_{[0, 1]}  \Delta_A  (r,s)   X^ *(s) \dint s   +    \langle Z^*,  \Delta_\beta (r)\rangle_p\right) \epsilon^*_t(r) \dint r\right\} \\
		= & \mathbb{E}_{X^*, Z^*, Y^*} \left\{\int_{[0, 1]} \left(\int_{[0, 1]}  \Delta_A  (r,s)   X^ *(s) \dint s   +    \langle Z^*,  \Delta_\beta (r)\rangle_p \right)^2 \dint r\right\} \\
		& + 2 \int_{[0, 1]} \mathbb{E}_{X^*, Z^*, Y^*}\left(\int_{[0, 1]}  \Delta_A  (r,s)   X^ *(s) \dint s + \langle Z^*,  \Delta_\beta (r)\rangle_p\right) \mathbb{E}_{X^*, Z^*, Y^*}\left\{\epsilon^*_t(r) \right\}\dint r \\
		= & \mathbb{E}_{X^*, Z^*, Y^*} \left\{\int_{[0, 1]} \left(\int_{[0, 1]}  \Delta_A  (r,s)   X^ *(s) \dint s   +    \langle Z^*,  \Delta_\beta (r)\rangle_p \right)^2 \dint r\right\},
\end{align*}
where the second identity is due to the definition of $Y^*(r)$, the fourth identity is due to the independence between $\epsilon_t(\cdot)$ and the rest, and the final identity is due to the mean-zero assumption of $\epsilon^*_t(\cdot)$.	
\end{proof}

 \subsection{Properties of Soboblev Spaces}\label{sec:Holder_smooth}

\begin{definition} \label{def:holder smooth}
Let $f : [0,1]\to \mathbb R $ be any real function. Then $f$ is said to be H\"{o}lder continuous of order $\kappa\in (0,1]$ if 
$$\sup_{x, x'\in [0,1]} \frac{|f(x)-f(x')|}{| x-x'|^{\kappa}} < \infty  . $$

\end{definition}

\begin{theorem}\label{theorem:Morrey}
Suppose that $f : [0,1]\to \mathbb R $ is such that 
$f\in \wot $ with $\alpha >1/2$. Then $f$ is H\"{o}lder continuous with $\kappa = \alpha -1/2 $ and  there exists an absolute constant $C$ such that 
$$\sup_{x, x'\in [0,1]} \frac{|f(x)-f(x')|}{| x-x'|^{\alpha -1/2 }} <C \| f\|_{\wot}. $$

\end{theorem}
\begin{proof}
This is the well-known  Morrey's inequality. See \cite{evans2010partial} for a proof.  
\end{proof}

\
\\
In  the  following lemma, we show that the integral of a H\"{o}lder smooth function $f$ can be approximated by the discrete sum of $f$ evaluated at the sample points   $\{s_i \}_{i=1}^n $. 
\begin{lemma}\label{lemma:discrete approximation}
Suppose that $\{s_i \}_{i=1}^n $ is a collection of  sample points satisfying \Cref{assume:joint smooth}{\bf a}. Suppose in addition  that $f : [0,1]\to \mathbb R $ is H\"{o}lder continuous with 
$$\sup_{x, x'\in [0,1]} \frac{|f(x)-f(x')|}{| x-x'|^{\alpha -1/2 }} <C_f  . $$  
Let $w_s(i) = (s_{i}-s_{i-1} ) n .  $
Then there exists an absolute constant $C$ such that 
$$\bigg|  \int_{[0,1]} f(s)ds  - \frac{1}{n} \sum_{i=1}^n \wsi  f(s_i) \bigg| \le CC_f  n^{-\alpha +1/2.} $$
Consequently  if  $f \in \wot, $
then 
$$\bigg|  \int_{[0,1]} f(s)ds  - \frac{1}{n} \sum_{i=1}^n \wsi  f(s_i) \bigg| \le C '\|f\|_\wot   n^{-\alpha +1/2.} $$
\end{lemma}
\begin{proof}Let $s_0=0$. 
Observe that 
\begin{align*}
\bigg| \int_{[0,1]} f(s)ds  - \frac{1}{n} \sum_{i=1}^n \wsi   f(s_i)  
\bigg| 
= &  \bigg| \sum_{ i=1}^{n } \int_{s_{i-1}}^{s_{i}} f(s )\dint s  -    \frac{1}{n} \sum_{i=1}^n  w_s(i) f(s_i) \bigg| 
\\
= & \bigg| \sum_{ i=1}^{n } \int_{s_{i-1}}^{s_{i}} f(s ) \dint  s  -    \sum_{ i=1}^{n } \int_{s_{i-1}}^{s_{i}}  f(s_i )\dint  s  \bigg| 
\\
\le & \sum_{ i=1}^{n }    \int_{s_{i-1}}^{s_{i}}\bigg|  f(s ) -f(s_i) \bigg|  \dint s     
\\
\le & \sum_{ i=1}^{n } (s_{i}- s_{i-1}) C_f {C_d}^{\alpha -1/2}    n^{-\alpha +1/2 }
\\
= & C_f C_d ^{\alpha -1/2}  n^{-\alpha +1/2 }.
\end{align*}
\end{proof}

 \begin{lemma} \label{lemma:basis properpties of sobolev} Suppose that $f, g: \mathbb [0,1]\to \mathbb R$  are such that $f,g \in W^{\alpha ,2}  $ for some $\alpha >1/2$. Then there exist absolute constants 
 $ C_1$ and $C_2$ such that  
 	\begin{align} \label{eq:product regularity} 
		\|fg\|_{W^{ \alpha,2}}  \le C_1  \|f\|_ \wot  \|g\|_\wot  \quad \text{and} \quad  \|f\|_\infty  \le C_2\|f\|_\wot  .
	\end{align} 
	 
\end{lemma}
\begin{proof}
These are well know Sobolev inequalities. See \cite{evans2010partial} for   proofs.
\end{proof}

\input{moresimulation}

%% file: appendix/functional-deviation.tex
\subsection{Functional deviation bounds}
\begin{lemma}  \label{lemma:sobolev covering 0} 
Let $\{s_ j \}_{ j  = 1}^n \subset [0, 1]$ be independent uniform random variables.  Let $\{\varepsilon_t\}_{t=1}^T$ be a collection of independent standard  Gaussian random variables independent of $\{X_t\}_{t=1}^T$ and $\{s_i\}_{i = 1}^n$.  
Under \Cref{assume:joint smooth},  it holds that 
	\begin{align*}
		\mathbb P\Bigg(  \bigg|\frac{1}{n }\sum_{t=1}^T \bigg( \frac{1}{T}  \sum_{j=1}^n g(s_j) X_t(s_j) - \int_{[0,1]} g(s) X_t(s) \dint s    \bigg) \varepsilon_t   \bigg| \ge C  \log(T)  \zeta_n \ \text{for all}  \  \|g\|_\wot \le 1     \Bigg)\le T^{-5}
 	\end{align*}
	where $C ,c > 0$ are   absolute constants depending only on $C_X$ and $C_\varepsilon$ and $\zeta_n=n^{-\alpha+1/2 } $. 
 
\end{lemma}

\begin{proof} 
{\bf Step 1.}	Note that   by \Cref{lemma:Gaussian process sobolev norm}, it holds that 
	$$ \mathbb P\big ( \| X_t \|_\wot  \ge 4C_X \sqrt {\log(T) } \big)\le T^{-6}. $$  	By a  union bound argument,  
	$$ \mathbb P(\mathcal E) =\mathbb  P\big ( \| X_t \|_\wot \ge 4C_X \sqrt {\log(T)}  \ \text{for all } 1\le t \le T \big)\le T^{-5} .   $$
	Under the  event $\mathcal E$, for any $g $ such that $\|g\|_\wot$, it holds that 
	$$\|g X_t  \|_\wot \le \|g\|_\wot \|X_t\|_\wot \le 4C_X \sqrt { \log (T) }. $$
	So under  the event $\mathcal E$ , for all $1\le  t \le T$ and all $\|g\|_\wot  \le 1, $ it holds that from \Cref{lemma:discrete approximation} that
	$$    \bigg|  \frac{1}{ n }  \sum_{j=1}^n g(s_j) X_t(s_j) - \int_{[0,1]}    g(s) X_t(s) \dint s   \bigg |\le C_1 n^{-\alpha + 1/2 } \sqrt {\log(T)} \  \text{for all}  \  \|g\|_\wot \le 1 .  $$
	Since 
	$$\mathbb P( |\varepsilon_t| \ge C_2\sqrt {\log(T)}) \le T^{-5} $$
 the desired result immediately follows.
 \end{proof} 
 \
 \\
 \\
   Through out this section, denote   the eigen-expansion of   linear map $L_{ \mk ^{1/2} \Sigma_X  \mk ^{1/2} }  $ as 
$$ L_{  \mk ^{1/2} \Sigma_X  \mk ^{1/2} } ( \Phi_k  )  = \sum_{k=1}^\infty  \xi_k   \Phi_k .  $$
     In addition, for any bilinear operator $\Sigma$ and any $\lt$ functions $f,g$,  denote 
$$     \Sigma  [  f ,  g ] = \iint_{[0,1]^2} \Sigma(r,s) f(r) g(s)\dint r \dint s.  $$
\begin{lemma}    \label{lemma:functional deviation bound}
Suppose $\{ X_t\} _{t=1} ^T $ are independent and   identically distributed centered Gaussian random processes and that 
the  eigenvalues  $\{ \xi_k\}_{k=1}^\infty  $ of  the linear operator  $L_{ \mk ^{1/2} \Sigma_X  \mk ^{1/2} }  $ satisfy that 
\begin{align*} \xi_k \asymp   k^{-2r }
\end{align*} for some $r>1/2$.  Let $\{ \varepsilon_ t \} _{t=1} ^T   \overset{i.i.d. } {\sim}  \mathcal  N(0, 1)$ and be  independent of $ \{ X_t\} _{t=1} ^T $. Then 
with probability at least  $1 -T^{-4}$,
$$ \left|  \frac{1}{T} \sum _{t=1} ^T \lb X_t , \beta\rb_\lt \varepsilon_t\right|  \le  C  \bigg(  \sqrt {  \Sigma_X [\beta, \beta  ]     \log(T)   \delta_T   }      
 +   \log(T)   \delta_T  \bigg)      \quad \text{ for all } \| \beta\|_{\hk } \le 1  , $$ 
where $\delta_T = T  ^{  - \frac{2 r}{ 2r+ 1 } }         $ and $C$ is some absolute constant independent of $T$.
\end{lemma}
\begin{proof}

Note that for any deterministic  $ \alpha \in \hk$.
$$ \mathbb  E(  \lb X_t , \alpha \rb_\lt \varepsilon_t ) =0  .$$
In addition, since   $ \lb X_t , \alpha \rb_\lt  \sim  \mathcal N( 0,  \Sigma _X [\alpha , \alpha ])$   
 and $ \varepsilon_t\sim \mathcal N(0, 1),   $ 
where 
$$ \Sigma _X [\alpha, \alpha]  = \iint_{[0,1]^2 } \alpha(s)  \Sigma _X (s,t)\alpha(t) \dint s \dint t .$$
Thus $ \lb X_t , \alpha \rb_\lt \varepsilon_t $ is  a centered sub-exponential random variable with parameter  $    \Sigma _X [\alpha , \alpha ]    $. By \Cref{lemma:sub-exponential tail}, for $ \gamma<1$.
\begin{align} 
\label{eq:one beta covering 13}
 \mathbb P\left( \left|  \frac{1}{T} \sum _{t=1} ^T \lb X_t , \alpha \rb_\lt \varepsilon_t \right| \ge   \gamma  \sqrt {  \Sigma _X [\alpha, \alpha]}      \right) \le  \exp( -2 \gamma ^2 T).  
\end{align}
\
\\
Denote 
\begin{align*} &\mathcal F :  = \text{span}  \left\{    L_\kh  (\Phi_k )  \right\} _{k=1}^D  \subset \hk 
\text{ and } 
\\
& \mathcal F^\perp  :  = \text{span}    \left\{   L_\kh  (\Phi_k )  \right\} _{k=D+1}^\infty \subset \hk .
 \end{align*}
Denote $\mathcal P_{\mathcal F}$ to be the projection  operator from $ \hk $ to $\mathcal F$ with respect to the $\hk $ topology. 
 \
 \\
 By \Cref{lemma:connect L_2 to H}, $ L_\kh  (\Phi_k ) $ is a collection of $\hk$ basis. 
For any $\beta$ such that $\| \beta\|_\hk \le 1  $,   since 
 $\{  L_\kh  (\Phi_k ) \}_{k=1}^\infty$ form a $\hk$ basis, 
 $ \| \mathcal  P_\mathcal F (  \beta )  \|^2_\hk +\|\mathcal  P_{\mathcal F ^\perp}   ( \beta )  \|^2_\hk  = \| \beta \|^2_\hk \le 1    $. 
Note that 
\begin{align*}    \frac{1}{T} \sum _{t=1} ^T \lb X_t , \beta\rb_\lt\varepsilon_t 
=&\frac{1}{T} \sum _{t=1} ^T \lb X_t , \mathcal  P_{\mathcal F  } (\beta) \rb_\lt\varepsilon_t   + \frac{1}{T} \sum _{t=1} ^T \lb X_t , \mathcal P_{\mathcal F ^\perp} (\beta)  \rb_\lt\varepsilon_t .
 \end{align*}  
 \
  \\
{\bf Step 1.}  Let $D \le T$ to be chosen later. For any $ J\in \mathbb Z$, $  J \ge 1  $ and $ 2^{-J}\ge  T^{-6} $, consider the sets 
\begin{align*} \mathcal F_J  = &  \left\{   \alpha = \sum_{k=1}^D  f_k  L_\kh (\Phi_k )   :   2^{-J-1} \le   \sqrt {   \Sigma _X [\alpha , \alpha ]  }    \le 2 ^{-J} , 
\sum_{k=1}^D f_k^2=\| \alpha  \|_\hk ^2   \le 1  \right\}
\end{align*} 
 \\
 \\
So $\mathcal F_J$  can be identified as a  unit ball in $ \mathbb R^D$. This means that for every $\delta>0$, there exists a collection $\{  \alpha_m\}_{m=1}^M $ such that for any $ \alpha \in \mathcal F_J$,
$$ \|\alpha_m -\alpha \|_\hk   \le \delta  $$
and $M\le \left( \frac{2}{\delta} \right) ^D $. 
 \\
 \\
Therefore for given $m$,  $\alpha _m\in \mathcal F_J$, and so 
\begin{align*} 
 \mathbb P\left( \left|  \frac{1}{T} \sum _{t=1} ^T \lb X_t , \alpha_m \rb_\lt \varepsilon_t \right| \ge \gamma  2^{-J}     \right) \le  \exp( -2 \gamma ^2 T).  
\end{align*}
So for any $\alpha \in \mathcal F_J$,
\begin{align*}
&\frac{1}{T} \sum _{t=1} ^T \lb X_t ,  \alpha   \rb_\lt \varepsilon_t 
\\
\le
&
\frac{1}{T} \sum _{t=1} ^T \lb X_t , \alpha- \alpha   _m   \rb_\lt \varepsilon_t  +  \sup_{1\le m \le M }\frac{1}{T} \sum _{t=1} ^T \lb X_t , \alpha_m   \rb_\lt \varepsilon_t 
&  
\\
\le
& \| \alpha - \alpha  _m\|_\lt  \left\|  \frac{1}{T} \sum _{t=1} ^T  X_t   \varepsilon_t   \right \|_\lt  +  \sup_{1\le m \le M }\frac{1}{T} \sum _{t=1} ^T \lb X_t , \alpha_m   \rb_\lt \varepsilon_t   
\\
\le
&  \frac{ \delta}{\sqrt T} +  \sup_{1\le m \le M }\frac{1}{T} \sum _{t=1} ^T \lb X_t , \alpha _m   \rb_\lt \varepsilon_t    .
\end{align*} 
 Therefore 
 $$  \mathbb  P \left( \left|  \frac{1}{T} \sum _{t=1} ^T \lb X_t , \alpha  \rb_\lt \varepsilon_t \right| \ge \gamma  2  \sqrt {  \Sigma _X [\alpha, \alpha ] }    + \frac{ \delta}{\sqrt T }   \text{ for all } \alpha \in \mathcal F_J   \right)  
 \le \exp\left( D \log(2/\delta ) -\gamma^2T  \right) .$$
Letting  $\gamma = C \sqrt { \frac{D \log(T)  }{T}  }$ and $\delta=T ^{-9/2} $ gives
 $$  \mathbb  P \left( \left|  \frac{1}{T} \sum _{t=1} ^T \lb X_t , \alpha  \rb_\lt \varepsilon_t \right| \ge 
   2C \sqrt { \frac{   \Sigma _X [\alpha, \alpha ] D\log(T)  }{ T } }      + \frac{1}{T^ 5 }      \text{ for all } \alpha \in \mathcal F_J \right)  
 \le T^{-6} .$$
 \\
 \\
 Let 
 $$\mathcal E  = \left\{ \alpha = \sum_{k=1}^D  f_k  L_\kh  (\Phi_k )  :      \sqrt {  \Sigma _X [\alpha, \alpha] }  \le \frac{1}{T^ 5 } , 
\sum_{k=1}^D f_k^2 =\| \alpha\|_\hk ^2   \le 1  \right\}.  $$
 The similar argument shows that 
  $$  \mathbb  P \left( \left|  \frac{1}{T} \sum _{t=1} ^T \lb X_t , \alpha  \rb_\lt \varepsilon_t \right| \ge 2 C
    \sqrt { \frac{   \Sigma _X [\alpha, \alpha ] D\log(T)  }{ T } }    + \frac{1}{T^ 5 }   \text{ for all } \alpha \in \mathcal E \right)  
 \le T^{-6} .$$
Since  $ \Sigma _X [\alpha, \alpha ]\le T ^ {-6}  $ and $D\le T $, 
  $$  \mathbb  P \left(  \left|  \frac{1}{T} \sum _{t=1} ^T \lb X_t , \alpha  \rb_\lt \varepsilon_t \right| \ge 
    \frac{   \sqrt {   \log(T) }   } { T ^  3  }      + \frac{1}{T^ 5 }      \text{ for all } \alpha \in \mathcal E   \right)  
 \le T^{-6} .$$
 Since 
 $$\{ \alpha \in \mathcal F : \|\alpha\|_\hk \le 1 \}   = \bigcup_{2 ^{-J} \ge \frac{1}{T}} \mathcal  F_J \cup \mathcal E, $$   by a  union bound argument,
\begin{align}   \nonumber  &\mathbb  P \left( \left|  \frac{1}{T} \sum _{t=1} ^T \lb X_t , \alpha  \rb_\lt \varepsilon_t \right| \ge 
  2 C   \sqrt { \frac{   D \Sigma_X [\alpha, \alpha ]  \log(T)  }{ T    } }        +    \frac{ \sqrt {  \log(T)  } } { T  ^ 3  }        + \frac{1}{T^ 5 }  \text{ for all } \alpha \in \mathcal F  , \|\alpha\|_\hk \le 1   \right)  
  \\
 \le&  \log(T) T^{-6} \label{eq:bias part pealing} 
 \end{align}

\
 \\
{ \bf Step 2.}  For any $\beta $ such that $\| \beta\|_\hk \le 1 $, it holds that $\| \mathcal P_\mathcal F (\beta) \|_\hk \le 1    $. Therefore 
\begin{align*}  
    &\frac{1}{T } \sum _{t=1} ^T \lb X_t ,  \mathcal P _{\mathcal F} (\beta )  \rb_\lt\varepsilon_t
    \\
\le 
&  2 C  \sqrt { \frac{   D \Sigma_X [ \mathcal P_{\mathcal F } ( \beta), \mathcal P_{\mathcal F } ( \beta) ]  \log(T)  }{ T    } }      +    \frac{ \sqrt {  \log(T)  } } {T  ^ 3  }         + \frac{1}{T^ 5 } 
\\
\le& 2 C  \sqrt { \frac{   D \Sigma_X [  \beta ,  \beta  ]  \log(T)  }{ T    } }       +    \frac{ \sqrt {  \log(T)  } } {T  ^ 3  }     
   + \frac{1}{T^ 5 } ,
 \end{align*} 
 where the first   inequality follows from \eqref{eq:bias part pealing},
 and the second  inequality holds because 
\begin{align*}
\Sigma_X [\beta  ,      \beta  ]  =  & \Sigma_X [  \mathcal P _{\mathcal F} ( \beta),  \mathcal P _{\mathcal F} (\beta) ]  + 2   \Sigma_X [  \mathcal P _{\mathcal F ^\perp } ( \beta),  \mathcal P _{\mathcal F} (\beta) ]  + \Sigma_X [  \mathcal P _{\mathcal F ^\perp } ( \beta),  \mathcal P _{\mathcal F  ^\perp } (\beta) ]   
\\
= & \Sigma_X [  \mathcal P _{\mathcal F} ( \beta),  \mathcal P _{\mathcal F} (\beta) ]    + \Sigma_X [  \mathcal P _{\mathcal F ^\perp } ( \beta),  \mathcal P _{\mathcal F  ^\perp } (\beta) ]    
\\
\ge & \Sigma_X [  \mathcal P _{\mathcal F} ( \beta),  \mathcal P _{\mathcal F} (\beta) ] 
\end{align*} 
 where the second equality follows   from  \Cref{lemma:perpendicular}, which implies that  $\Sigma_X[  \mathcal P _{\mathcal F ^\perp } ( \beta),  \mathcal P _{\mathcal F} (\beta) ]  =0  $

 \
 \\
 {\bf Step 3.} 
Denote   the eigen-expansion of   linear map $L_{ \mk ^{1/2} \Sigma_X  \mk ^{1/2} }  $ as 
$$ L_{  \mk ^{1/2} \Sigma_X  \mk ^{1/2} } ( \Phi_k  )  = \sum_{k=1}^\infty \xi_k   \Phi_k .  $$
 Let $$w_{t ,k}=\frac{  \lb  L_{\kmh  }(X_t ) , \Phi_k\rb_\lt  }{\sqrt {\xi _k }}. $$ 
 Then 
 \begin{align*}
 & \sup_{\| \beta\|_\hk \le 1     } \frac{1}{T }  \sum _{t=1} ^T \lb X_t ,  \mathcal P _{\mathcal F ^\perp } (\beta )  \rb_\lt\varepsilon_t 
 \\
 \le &
   \sup_{ \beta_ 2 \in \mathcal F^\perp, \|\beta_2\|_\hk \le 1    }  \frac{1}{T }  \sum _{t=1} ^T \lb X_t , \beta_ 2 \rb_\lt\varepsilon_t 
   \\
    \le &
   \sup_{ \beta_ 2 \in \mathcal F^\perp, \|\beta_2\|_\hk \le 1    }  \frac{1}{T }  \sum _{t=1} ^T \lb  L_\kh (X_t) ,L_{\kmh}  (\beta_ 2)  \rb_\lt\varepsilon_t 
   \\
   = & \sup_{ \beta_ 2 \in \mathcal F^\perp , \| \beta_2 \|_{\hk } \le 1   }\frac{1}{T} \sum _{t=1} ^T  \sum_{k= 1 }^\infty   \lb L_\kh(X_t) ,\Phi_k \rb_\lt \lb  L_\kmh  (  \beta_2 )    , \Phi_k  \rb_\lt \varepsilon_t    
\\
  =& \sup_{ \beta_ 2 \in \mathcal F^\perp , \| \beta_2 \|_{\hk } \le 1   }\frac{1}{T} \sum _{t=1} ^T  \sum_{k=D+1 }^\infty   \lb L_\kh(X_t) ,\Phi_k \rb_\lt \lb L_\kmh(\beta_2 ) ,\Phi_k\rb_\lt \varepsilon_t   
  \\
  =
  & \sup_{  \sum_{k=D+ 1}^\infty f_k^2 \le 1    }\frac{1}{T} \sum _{t=1} ^T  \sum_{k=D+1 }^\infty  \sqrt {\xi_k }w_{t,k } 
  f_k\varepsilon_t   , \end{align*}
  where the  fourth inequality holds because $\beta_2 \in \mathcal F^\perp$ and the last inequality follows from \Cref{lemma:change of basis 1}.
 \
 \\
 Note that since  $w_{t ,k}  $ and $\varepsilon_t$   are both centered Gaussian with variance $1$,
 \begin{align*}
  & \mathbb  P \left (   \left  | \frac{1}{T}    \sum_{ i=1}^n w_{t ,k} \varepsilon_t \right |   \le  \sqrt { \frac{ 4  \log(k) +12  \log(T)     }{T} } +    \frac{ 4  \log(k) +12  \log(T)     }{T}     \text{ for all } k\ge 1 \right)
   \\
   \le 
   &  \sum_{k=1}^\infty  \mathbb  P \left (   \left  | \frac{1}{T}    \sum_{ i=1}^n w_{t ,k} \varepsilon_t \right |   \le  \sqrt { \frac{ 4 \log(k)  +  12 \log(T)   }{T} } +    \frac{ 4  \log(k) +12  \log(T)     }{T}       \right) 
   \\
   \le 
   &  \sum_{k=1}^\infty \frac{1}{T^6 k^2}\le \frac{1}{T^5}.
 \end{align*}  
  \
  \\
 So for any $D' \ge D$, 
\begin{align} \nonumber  
&  \sup_{  \sum_{k= D + 1}^\infty f_k^2 \le 1    }\frac{1}{T} \sum _{t=1} ^T  \sum_{k= D +1 }^\infty  \sqrt { \xi _k }w_{t,k } 
  f_k\varepsilon_t  
  \\ \nonumber   
  =   &  \sup_{  \sum_{k= D + 1}^\infty f_k^2 \le 1    }  \sum_{k= D +1 }^\infty  \sqrt {\xi _k }
  f_k \left(  \frac{1}{T}\sum _{t=1} ^T  w_{t,k }   \varepsilon_t   \right) 
  \\   \nonumber 
  \le 
  & \sup_{  \sum_{k=D+ 1}^\infty f_k^2 \le 1    }  \sum_{k= D +1 }^\infty  k^{- r  } 
  \left( \sqrt { \frac{ 4  \log(k) +12  \log(T)     }{T} }     +      \frac{ 4  \log(k) +12  \log(T)     }{T}   \right) 
  f_k 
  \\
 \label{eq:term 1 deviation 1 step 3}
  \le 
  & \sup_{  \sum_{k=D+ 1}^\infty f_k^2 \le 1    }  \bigg( 
    \sqrt { \sum_{k=D +1 }^\infty   \frac{ 4  \log(k) +12  \log(T)     }{ k^{ 2 r }T}    }
 +        \sqrt { \sum_{k=D +1 }^\infty   \frac{ 4  \log ^2 (k) +12  \log ^2 (T)     }{ k^{ 2 r }T ^2  }     } 
     \bigg) \sqrt { \sum_{k=D +1 }^\infty f_k^2  }  
  \end{align}
  Note that 
 \begin{align*} 
 &   \sup_{  \sum_{k=D + 1}^\infty f_k^2 \le 1    }  \sqrt { \sum_{k=D +1 }^\infty   \frac{ 4  \log(k) +12  \log(T)     }{ k^{ 2 r }T}    } \sqrt { \sum_{k=D +1 }^\infty f_k^2  }
\\
 \le  
&  \sqrt {   \sum_{k=D+1 }^\infty   \frac{ 4  \log(k)     }{ k^{ 2  r  }T}     } + \sqrt {   \sum_{k=D +1 }^\infty   \frac{  12  \log(T)     }{ k^{ 2  r   }T}     }  
\\
\le 
&\sqrt {  \frac{ 4        \log(p) }{  D ^{ 2r- 1   }T}    +    \frac{ 1     }{ (D')  ^{    r- 1/2   }T}        } + \sqrt { \frac{12  \log(T) }{ T   D ^{ 2 r  - 1}}}   .
  \end{align*} 
 where the first inequality follows from \Cref{lemma:tail bound with log factor}.  The second term in  \Cref{eq:term 1 deviation 1 step 3} can be bounded in a similar way. 
  \\
  \\
{\bf Step 4.}  Putting {\bf Step 2} and {\bf Step 3} together, it holds that with high probability, for all $\beta \in \hk  $ such that $\|\beta\|_\hk \le 1  $,
\begin{align*} 
   \frac{1}{T} \sum _{t=1} ^T \lb X_t , \beta\rb_\lt\varepsilon_t   
\le   & C' \bigg(  \sqrt { \frac{   D\Sigma_X [\alpha, \alpha ]\log(T)  }{ T    }  }    +    \frac{ \sqrt {  \log(T)  } } {T  ^ 3  }  + \frac{1}{T^5 }+    \sqrt {  \frac{ 4        \log(p) }{  D ^{ 2r- 1   }T}    +    \frac{ 1     }{ (D')  ^{   r- 1/2  }T}        } + \sqrt { \frac{12  \log(T) }{ T  D ^{ 2 r  - 1}}} \bigg) 
  \end{align*} 
 Set $ D'  = \max\{ T^{ \frac{10}{   r   - 1/2  }},   D \}  $ gives
   \begin{align*} 
  \frac{1}{T} \sum _{t=1} ^T \lb X_t , \beta\rb_\lt\varepsilon_t   
\le   & C'  \bigg(   \sqrt { \frac{   D\Sigma_X [\alpha, \alpha ]\log(T)  }{ T    }  }     +    \frac{ \sqrt {  \log(T)  } } {T  ^ 3  }      + \frac{1}{T^5 }       +   \sqrt { \frac{ \log( T ) }{T}   D^{-2r + 1 }    }   +    \frac{ 1         }{   T ^ 5 }         + 
 \sqrt { \frac{12  \log(T) }{ T  D ^{ 2 r  - 1}}}   \bigg) 
  \end{align*}     
  \
  \\
  Taking $D =     T   ^{ \frac{1}{2r + 1 }} $ gives
  \begin{align*} 
   \frac{1}{T} \sum _{t=1} ^T \lb X_t , \beta\rb_\lt\varepsilon_t   
\le   & C''\bigg(  \sqrt {  \Sigma_X [\alpha, \alpha ]  \frac{    \log(T)  }{ T  ^{  \frac{2 r}{ 2r+ 1 } }     }  }       
 + \frac{    \log(T)  }{ T  ^{  \frac{2 r}{ 2r+ 1 } } }      +    \frac{  1       }{   T ^2  }    \bigg)   .
  \end{align*}  
  
 \end{proof}

\begin{lemma} \label{lemma:functional restricted eigenvalues}
Suppose $\{ X_t\} _{t=1} ^T $ are independent  identically distributed centered Gaussian random processes and that 
the  eigenvalues of $\{\xi_k\}_{k=1}^\infty  $ of  the linear operator  $L_{ \mk ^{1/2} \Sigma_X  \mk ^{1/2} }  $ satisfy
\begin{align*} \xi_k \asymp   k^{-2r }
\end{align*} for some $r>1/2$. Then 
with probability  at least  $1 -T^{-4}$,  it holds that for any $0<\tau<1 $,
$$ \left|   \frac{1}{T} \sum _{t=1} ^T \lb X_t , \beta\rb_\lt ^2  -  \Sigma_X [\beta, \beta ] \right|  \le   
\tau  \Sigma_X [\beta, \beta  ]     +C_\tau   \log       (T)  \delta_T    \quad  \text{for all} \quad  \| \beta\|_{\hk } \le 1   ,$$ 
where 
$ \delta_T =  T^{ -\frac{2r}{2r+1 }}  $,  and $C_\tau$ is some  constant only depending  on $\tau$ and  independent of $T$.
\end{lemma}

\begin{proof}

Note that for any deterministic  $ \alpha \in \hk$,  $ \lb X_t , \alpha \rb_\lt $ is  $SG( \Sigma _X[\alpha, \alpha]  )  $. 
Thus by Hanson-Wright, it holds that for all $\gamma<1$,  
\begin{align} 
\label{eq:one beta covering 1}
  \mathbb  P\left( \left|  \frac{1}{T} \sum _{t=1} ^T \lb X_t , \alpha \rb_\lt   ^2  -\Sigma_X [\alpha ,\alpha ] \right| \ge \gamma   \Sigma_X [\alpha, \alpha]  \right) \le  \exp( -2 \gamma ^2 T).  
\end{align}
\
\\
For  $D\le T$  to be chosen later, denote 
$$ \mathcal F :  =  \text{span} \left\{      L_\kh  (\Phi_k )   
  \right\}  _{k=1}^D \quad  \text{and}\quad  \mathcal F^\perp  :  =  \text{span} \left\{     L_\kh  (\Phi_k )  \right\} _{k=D+1}^\infty .  $$ 
 \
 \\
 Denote $\mathcal P_{\mathcal F}$ to be the projection  operator from $ \hk $ to $\mathcal F$ with respect to the $\hk $ topology. 
 Them it holds that for all $\beta \in \hk$, 
 \begin{align*}
&  \left|   \frac{1}{T} \sum _{t=1} ^T \lb X_t , \beta\rb_\lt ^2  - \Sigma_X[\beta, \beta ] \right|   
\\
\le & \left|   \frac{1}{T} \sum _{t=1} ^T \lb X_t ,  \mathcal P_{\mathcal F } ( \beta) \rb_\lt ^2  - \Sigma_X[ \mathcal P_{\mathcal F } ( \beta) ,  \mathcal P_{\mathcal F } ( \beta)  ] \right|    
\\
+ & 2   \left|   \frac{1}{T} \sum _{t=1} ^T  \lb X_t ,  \mathcal P_{\mathcal F } ( \beta) \rb_\lt  \lb X_t ,  \mathcal P_{\mathcal F ^\perp  } ( \beta) \rb_\lt   -\Sigma_X[ \mathcal P_{\mathcal F } ( \beta) ,  \mathcal P_{\mathcal F ^\perp  } ( \beta)  ] \right|     
\\
+ &   \left|   \frac{1}{T} \sum _{t=1} ^T \lb X_t ,  \mathcal P_{\mathcal F ^\perp  } ( \beta) \rb_\lt ^2  - \Sigma_X [ \mathcal P_{\mathcal F^\perp  } ( \beta) ,  \mathcal P_{\mathcal F ^\perp } ( \beta)  ] \right|    
\\
  \end{align*}

 \
 \\  
{\bf Step 1.}   For any $ J\in \mathbb Z$, $  J\ge 1  $ and $ 2^{-J}\ge  T^{-5}$, consider the set 
$$ \mathcal F_J = \left\{ \alpha = \sum_{k=1}^D  f_k  L_\kh  (\Phi_k )  :   2^{-J-1} \le   C[\alpha, \alpha]    \le 2^{-J} , 
\sum_{k=1}^D f_k^2 =\| \alpha\|_\hk ^2   \le 1  \right\}.  $$ 
 \\
 \\
So $\mathcal F_J$ can be viewed as a subset of unit ball in $ \mathbb R^D$. This means that for every $\delta>0$, there exists a collection $\{ \alpha_m\}_{m=1}^M $ such that for any $\alpha\in \mathcal F_J$,
$$ \|\alpha_m -\alpha \|_\hk \le \delta  $$
and $M\le \left( \frac{2}{\delta} \right) ^D $. Therefore for given $m$,  $\alpha_m\in \mathcal F_J$,  
\begin{align} 
\label{eq:one beta covering 12}
\mathbb P\left( \left|  \frac{1}{T} \sum _{t=1} ^T \lb X_t , \alpha_m  \rb_\lt   ^2  -\Sigma_X[\alpha_m ,\alpha_m ] \right| \ge \gamma    \Sigma_X[\alpha_m, \alpha_m ]  \right) \le  \exp( -2 \gamma ^2 T).  
\end{align}
Denote 
 $ \widehat \Sigma_X(r,s) = \frac{1}{T} \sum_{ t=1}^  T X_t (r) X_t(s). $ 
We have that  for any $\alpha \in \mathcal F_J$,
 \begin{align*}
 &\left|  \frac{1}{T} \sum _{t=1} ^T \lb X_t , \alpha   \rb_\lt   ^2  -\Sigma_X[\alpha  ,\alpha  ] \right|  
\\
=
&
 \left |      \iint  \alpha (r) \left( \frac{1}{ T }   \sum_{ t=1}^  T X_t (r) X_t(s)   -\Sigma_X  (r,s) \right)  \alpha (s)\dint r  \dint s     \right| 
\\
 \le 
 & 
 \bigg |      \iint  \alpha   (r) \left( \frac{1}{T }  \sum_{ t=1}^ T X_t (r) X_t(s)   - \Sigma_X (r,s) \right) \alpha (s)\dint r  \dint s   
 \\
 & -   \iint \alpha _m  (r) \left( \frac{1}{ T }   \sum_{ t=1}^  T X_t (r) X_t(s)   - \Sigma_X (r,s) \right) \alpha _ m  (s)\dint r  \dint s      \bigg | 
 \\
 +
 &\left |       \iint  \alpha _m (r) \left( \frac{1}{ T }   \sum_{ t=1}^ T  X_t (r) X_t(s)   -\Sigma_X (r,s) \right) \alpha _ m  (s)\dint r  \dint s      \right| 
 \\
 \le & 
   2  \|\alpha  -\alpha _m \|_\lt  \|\alpha \|_\lt \sqrt {   \iint   \left (  \frac{1}{ T }   \sum_{ t=1}^  T  X_t (r) X_t(s)   -\Sigma_X  (r,s) \right) ^2 \dint r  \dint s  } 
  \\
  +   &
  \sup_{1\le m \le M }  \left  | \iint \alpha _ m  (r) \left( \frac{1}{ T }   \sum_{ t=1}^  T  X_t (r) X_t(s)   -\Sigma_X  (r,s) \right) \alpha _ m  (s)\dint r  \dint s      \right| 
  \\
= & 
   2  \|\alpha  -\alpha _m \|_\lt  \|\alpha \|_\lt \big  \|   \widehat \Sigma_X   -\Sigma_X   \big \|_F
  \\
  +   &
  \sup_{1\le m \le M }  \left  | \iint \alpha _ m  (r) \left( \frac{1}{ T }   \sum_{ t=1}^  T  X_t (r) X_t(s)   -\Sigma_X  (r,s) \right) \alpha _ m  (s)\dint r  \dint s      \right| 
  \\
  \le 
  & 2 C_X\delta  \sqrt {\frac{\log( T )}{ T }} +
  \sup_{1\le  m  \le  M }  \left  | \iint  \alpha _ m  (r) \left( \frac{1}{ T }   \sum_{ t=1}^  T X_t (r) X_t(s)   -\Sigma_X  (r,s) \right)  \alpha _m (s)\dint r  \dint s      \right| 
  \\
  =
  &2 C_X\delta  \sqrt {\frac{\log( T )}{ T }} +\sup_{1\le  m  \le  M }  \left|  \frac{1}{T} \sum _{t=1} ^T \lb X_t , \alpha _m    \rb_\lt   ^2  -\Sigma_X[\alpha_m   ,\alpha_m   ] \right|   
\end{align*}   
\
\\
 Therefore 
 $$ \mathbb P \left(   \left|  \frac{1}{T} \sum _{t=1} ^T \lb X_t , \alpha   \rb_\lt   ^2  - \Sigma _X [\alpha  ,\alpha  ] \right| \ge 2  \gamma    \Sigma _X [\alpha , \alpha  ] +  2 C_X\delta  \sqrt {\frac{\log( T )}{ T }}    \text{ for all }  \alpha \in \mathcal F_J     \right)  
 \le \exp\left( D \log(2/\delta ) -\gamma^2T  \right) .$$
 So $\gamma = C \sqrt { \frac{D \log(T)  }{T}  }$ and $\delta=  T^{-9/2} $ gives
 $$ \mathbb P \left(    \left|  \frac{1}{T} \sum _{t=1} ^T \lb X_t , \alpha   \rb_\lt   ^2  -\Sigma_X [\alpha  ,\alpha  ] \right|  \ge 
    C_1 \bigg( \sqrt { \frac{  D\log(T)  }{ T } }     \Sigma_X [\alpha, \alpha ] + \frac{  \sqrt { \log(T)}  }{T ^5  }    \bigg)  \text{ for all }  \alpha \in \mathcal F_J   \right)  
 \le T^{-6} . $$
\
 \\
 Let 
 $$\mathcal E  = \left\{ \alpha = \sum_{k=1}^D  f_k  L_\kh  (\Phi_k ) \in \hk  :     \Sigma _X [\alpha, \alpha]  \le \frac{1}{T^ 5 } , 
\sum_{k=1}^D f_k^2 =\| \alpha\|_\hk  ^2  \le 1  \right\}.  $$
 The similar argument shows that 
  $$\mathbb  P \left(    \left|  \frac{1}{T} \sum _{t=1} ^T \lb X_t , \alpha   \rb_\lt   ^2  - \Sigma _X [\alpha  ,\alpha  ] \right|   \ge 
    C_1\bigg(  \sqrt { \frac{  D\log(T)  }{ T } }  \Sigma _X [\alpha, \alpha ]     + \frac{ \sqrt { \log(T) }  }{T ^5}  \bigg)   \text{ for all }  \alpha \in \mathcal E  \right)  
 \le T^{-6} .$$
Since  $  \Sigma _X [\alpha, \alpha ]\le  T ^{-5} $ when $\alpha \in \mathcal E$,   and $D\le T$, 
  $$\mathbb  P \left(   \left|  \frac{1}{T} \sum _{t=1} ^T \lb X_t , \alpha   \rb_\lt   ^2  -  \Sigma _X [\alpha  ,\alpha  ] \right|  \ge 
   C_1'  \frac{   \sqrt {  \log(T)   } }{  T ^{ 5}    }    \text{ for all }  \alpha \in \mathcal E         \right)  
 \le T^{-6} .$$
 Since 
 $$\{ \alpha \in \mathcal F :    \|\alpha \|_\hk \le 1  \}= \bigcup_{2 ^{-J} \ge T^{-5 }} \mathcal  F_J \bigcup \mathcal E ,$$  by union bound, 
\begin{align} \label{eq:bias pealing in quadratic form} & \mathbb P \left(  \left|  \frac{1}{T} \sum _{t=1} ^T \lb X_t , \alpha   \rb_\lt   ^2  - \Sigma _X[\alpha  ,\alpha  ] \right|  \ge 
  C_1'\bigg(   \sqrt { \frac{   D  \log(T)  }{ T    } }  \Sigma _X [\alpha, \alpha ]      +   \frac{   \sqrt {   \log(T)  }  }{ T  ^{  5}   }    \bigg)   \text{ for all }  \alpha \in \mathcal F   , \|\alpha \|_\hk \le 1   \right)  \\ \nonumber 
   \le & \log(T) T^{-6} 
 \end{align} 
 \
 \\
 { \bf Step 2.}  For any $\beta $ such that $\| \beta\|_\hk \le 1 $, it holds that $\| \mathcal P_\mathcal F (\beta) \|_\hk \le 1    $. Therefore 
\begin{align*}  
 &  \left|  \frac{1}{T} \sum _{t=1} ^T \lb X_t , \mathcal P _ { \mathcal F}  (\beta  )    \rb_\lt   ^2  -  \Sigma _X [\mathcal P _ { \mathcal F}   (\beta  )    ,\mathcal P _ { \mathcal F}  (\beta  )      ] \right|   
   \\
\le 
&  C_1'\bigg( \sqrt { \frac{   D   \log(T)  }{ T    } }        \Sigma _X [ \mathcal P _{\mathcal F} (\beta ) ,   \mathcal P _{\mathcal F} (\beta )  ]      +  \frac{   \sqrt {   \log(T)  }  }{ T  ^{  5}   }   \bigg) 
\\
\le 
&C_1' \bigg(  \sqrt { \frac{   D  \log(T)  }{ T    } }     \Sigma _X[   \beta  ,      \beta  ]    +    \frac{   \sqrt {   \log(T)  }  }{ T  ^{  5}   }   \bigg)  
 \end{align*} 
 where the  first   inequality follows from \eqref{eq:bias pealing in quadratic form},
 and the last inequality holds because 
\begin{align*}
 \Sigma _X [\beta  ,      \beta  ]  =  &   \Sigma _X [  \mathcal P _{\mathcal F} ( \beta),  \mathcal P _{\mathcal F} (\beta) ]  + 2    \Sigma _X[  \mathcal P _{\mathcal F ^\perp } ( \beta),  \mathcal P _{\mathcal F} (\beta) ]  + \Sigma _X [  \mathcal P _{\mathcal F ^\perp } ( \beta),  \mathcal P _{\mathcal F  ^\perp } (\beta) ]   
\\
= &   \Sigma _X [  \mathcal P _{\mathcal F} ( \beta),  \mathcal P _{\mathcal F} (\beta) ]    +  \Sigma _X [  \mathcal P _{\mathcal F ^\perp } ( \beta),  \mathcal P _{\mathcal F  ^\perp } (\beta) ]    
\\
\ge &   \Sigma _X [  \mathcal P _{\mathcal F} ( \beta),  \mathcal P _{\mathcal F} (\beta) ] 
\end{align*} 
 where  by \Cref{lemma:perpendicular}, $ \Sigma _X [  \mathcal P _{\mathcal F ^\perp } ( \beta),  \mathcal P _{\mathcal F} (\beta) ]  =0  $.
 \
 \\
 \\
 {\bf Step 3.} Denote   the eigen-expansion of   linear map $L_{ \mk ^{1/2} \Sigma_X  \mk ^{1/2} }  $ as 
$$ L_{  \mk ^{1/2} \Sigma_X  \mk ^{1/2} } ( \Phi_k  )  = \sum_{k=1}^\infty \xi_k   \Phi_k .  $$
  Let $$w_{t ,k}=\frac{  \lb  L_{\kh }(X_t ) , \Phi_k\rb_\lt  }{\sqrt {\xi _k }}. $$ 
For any $\beta$, let $ f_k = \lb L_\kmh (\Phi_k) , \beta\rb _\lt. $  Note that if 
$\| \beta\|_\hk  ^2 \le 1$, then $\|\mathcal P _  { \mathcal F^ \perp}  (\beta) \| _\hk  \le 1 $. 
$$   \lb X_t ,  \mathcal P_{\mathcal F ^\perp   }  (\beta)  \rb _\lt  =\sum_{k=D+1 }^\infty   \lb L_\kh(X_t) ,\Phi_k \rb_\lt \lb L_\kmh(\beta) ,\Phi_k\rb_\lt     = \sum_{k=D+1 }^\infty  \sqrt {\xi_k }w_{t ,k} f_k .  $$ 
Note that if 
$\| \beta\|_\hk \le 1$, then $\|\mathcal P _  { \mathcal F^ \perp}  (\beta) \| _\hk  \le 1 $. 
Therefore
 \begin{align*}
 &  \sup_{\|\beta\|_\hk \le 1 } \left|   \frac{1}{T} \sum _{t=1} ^T \lb X_t ,  \mathcal P_{\mathcal F ^\perp  } ( \beta) \rb_\lt ^2  - \Sigma_X[ \mathcal P_{\mathcal F^\perp  } ( \beta) ,  \mathcal P_{\mathcal F ^\perp } ( \beta)  ] \right|     
  \\
\le 
  & \sup_{  \sum_{k  =D+ 1}^\infty f_k^2 \le 1    }  \left| 
   \sum_{k , l =D+1 }^\infty   \sqrt {\xi_k } f_k \sqrt {\xi_ l } f_ l
 \left(   \frac{1}{T} \sum _{t=1} ^T    w_{i, k } w_{i,  l }   - \frac{1}{T} \sum _{t=1} ^T  E(w_{i, k}w_{i,l } )  
  \right)    \right|,
 \end{align*}
 where the last inequality follows from \Cref{lemma:change of basis 1}.
  \
 \\
 Let $\eta =   \sqrt { \frac{ 4  \log(k) + 4\log(l)  +12  \log(T)     }{T} }  $. Note that since  $w_{t ,k}  w_{i,l } $ is SE(1), 
 \begin{align*}
  &  \mathbb P \left (   \left  | \frac{1}{T}    \sum_{ i=1}^n w_{t ,k} w_{i, l } - \frac{1}{T} \sum _{t=1} ^T  E(w_{i, k}w_{i,l } )      \right |   \le \eta +\eta^2    \text{ for all } k, l \ge 1 \right)
   \\
   \le 
   &  \sum_{k,l =1}^\infty P \left (   \left  | \frac{1}{T}    \sum_{ i=1}^n w_{t ,k}w_{i, l }   - \frac{1}{T} \sum _{t=1} ^T  E(w_{i, k}w_{i,l } )    \right |   \le  \eta +\eta^2        \right) 
   \\
   \le 
   &  \sum_{k,l =1}^\infty \frac{1}{T^6 k^2 l^2 }\le \frac{1}{T^5}.
 \end{align*}
  As a result
  {\small
\begin{align} \nonumber 
& \sup_{  \sum_{k  =D+ 1}^\infty f_k^2 \le 1    }   \left| 
   \sum_{k , l =D+1 }^\infty   \sqrt {\xi_k } f_k \sqrt {\xi_ l } f_ l
 \left(   \frac{1}{T} \sum _{t=1} ^T    w_{i, k } w_{i,  l }   -E(w_{i, k}w_{i,l } )  
  \right)    \right|  
\\ \nonumber 
  \le  &  \sup_{  \sum_{k  =D+ 1}^\infty f_k^2 \le 1    }   \left| 
   \sum_{k , l =D+1 }^\infty   \sqrt {\xi_k } f_k \sqrt {\xi_ l } f_ l
 \left(     \eta +\eta^2    
  \right)    \right|  
  \\ \label{eq:restricted eigenvalues tail bound step 3}
  \le 
  &   \sup_{  \sum_{k  =D+ 1}^\infty f_k^2 \le 1    }   \frac{1}{ \sqrt T}  \Bigg| 
   \sum_{k , l =D+1 }^\infty   \sqrt {\xi_k } f_k \sqrt {\xi_ l } f_ l
 \bigg(     \sqrt {  4  \log(k)  } +  \sqrt { 4\log(l)}   + \sqrt { 12  \log(T)     }   +   \frac{ 4  \log(k)   
 +    4\log(l)    + 12  \log(T)     }{\sqrt T }
  \bigg)    \Bigg|  
\end{align}  }
 Observe that for $D' \ge D$, 
  \begin{align*}
 &  \sup_{  \sum_{k  =D+ 1}^\infty f_k^2 \le 1    }   \left| 
   \sum_{k , l =D+1 }^\infty   \sqrt {\xi_k } f_k \sqrt {\xi_ l } f_ l
   \sqrt {     \log(k)  }  \right|
   \\
   \le &\sup_{  \sum_{k  =D+ 1}^\infty f_k^2 \le 1    }     \left(   \sum_{l=D+1}^\infty \sqrt {\xi_l }f_l  \right)  \left(   \sum_{k=D+1}^\infty \sqrt {\xi_k}f_k \sqrt {     \log(k)  }    \right) 
   \\
   \le 
   &  \sup_{  \sum_{k  =D+ 1}^\infty f_k^2 \le 1    }   \left(   \sqrt{ \sum_{l=D+1}^\infty l^{-2r}} \sqrt { \sum_{l=D+1}^\infty  f_l  ^2   } 
    \right)   \left(   \sqrt{ \sum_{ k=D+1}^\infty  k ^{-2r} \log(k) } \sqrt { \sum_{l=D+1}^\infty  f_k  ^2   }   \right) 
    \\
    \le  & \sqrt {  D^{-2r+1 }  } \sqrt {  \left(  \frac{          \log(p) }{  D ^{ 2r- 1   }  }    +    \frac{ 1     }{ (D')  ^{    r- 1/2   } }  \right) },
  \end{align*}
  where the last inequality follows from \Cref{lemma:tail bound with log factor}. 
  In addition, 
  \begin{align*}
 &  \sup_{  \sum_{k  =D+ 1}^\infty f_k^2 \le 1    }   \left| 
   \sum_{k , l =D+1 }^\infty   \sqrt {\xi_k } f_k \sqrt {\xi_ l } f_ l
   \sqrt {     \log(T)  }  \right|
   \\
   \le & \sqrt { \log(T)}  \sup_{  \sum_{k  =D+ 1}^\infty f_k^2 \le 1    }     \left(   \sum_{l=D+1}^\infty \sqrt {\xi_l }f_l  \right)  \left(   \sum_{k=D+1}^\infty \sqrt {\xi_k}f_k     \right) 
   \\
   \le 
   & \sqrt { \log(T)}  \sup_{  \sum_{k  =D+ 1}^\infty f_k^2 \le 1    }   \left(   \sqrt{ \sum_{l=D+1}^\infty l^{-2r}} \sqrt { \sum_{l=D+1}^\infty  f_l  ^2   } 
    \right)   \left(   \sqrt{ \sum_{ k=D+1}^\infty  k ^{-2r}   } \sqrt { \sum_{l=D+1}^\infty  f_k  ^2   }   \right) 
    \\
    \le  & \sqrt { \log(T)}  D^{-2r+1 }  .
  \end{align*}
 The rest terms in \Cref{eq:restricted eigenvalues tail bound step 3} can be handled in a similar way.  So  
  \begin{align}
 &\sup_{\|\beta\|_\hk \le 1 }   \left|   \frac{1}{T} \sum _{t=1} ^T \lb X_t ,  \mathcal P_{\mathcal F ^\perp  } ( \beta) \rb_\lt ^2  -  \Sigma _X [ \mathcal P_{\mathcal F^\perp  } ( \beta) ,  \mathcal P_{\mathcal F ^\perp } ( \beta)  ] \right|      
 \\
 \le&  C_2 \frac{1}{\sqrt T } \left(   \sqrt {  D^{-2r+1 }  \left(  \frac{          \log(p) }{  D ^{ 2r- 1   }  }    +    \frac{ 1     }{  (D')  ^{    r- 1/2   } }  \right) }  + \sqrt { \log(T)}  D^{-2r+1 }   \right) 
 \\
 \le& C_2 \bigg(  \sqrt {     \frac{        D^{-4r+2 }       \log( T ) }{   T }    +    \frac{ D^{-2r +1 }    }{  T ^{10}     }    } +    
   \sqrt {  \frac{\log( T ) }{ T  } }  D^{-2r+1 } \bigg)   \le C_2' \bigg(  D^{-2 r+1 }  \sqrt {  \frac{\log( T )   }{ T  } }  +  \sqrt { \frac{ D^{-2r +1 }    }{  T ^{10}   }    } \bigg),
  \end{align}
  where the second inequality follows by setting $ D' = \max\{ n^{ \frac{  9 }{   r   - 1/2  }},   D \}  $.
  \\
  \\
 {\bf Step 4.} So with high probability, it holds that for all $\beta$ such that $\|\beta\|_\hk \le 1$,
 \begin{align*}
 \left|   \frac{1}{T} \sum _{t=1} ^T \lb X_t ,  \mathcal P_{\mathcal F } ( \beta) \rb_\lt ^2  -  \Sigma _X [ \mathcal P_{\mathcal F } ( \beta) ,  \mathcal P_{\mathcal F } ( \beta)  ] \right|     \le &  
C_3 \bigg(\sqrt { \frac{   D  \log(T)  }{ T    } }     \Sigma _X[   \beta  ,      \beta  ]    \bigg)
 +   \frac{ \sqrt {  \log(T)  } } { T  ^5  }     
\\
     \left|   \frac{1}{T} \sum _{t=1} ^T \lb X_t ,  \mathcal P_{\mathcal F ^\perp  } ( \beta) \rb_\lt ^2  -  \Sigma _X [ \mathcal P_{\mathcal F^\perp  } ( \beta) ,  \mathcal P_{\mathcal F ^\perp } ( \beta)  ] \right|     & \le  C_3\bigg( D^{-2 r+1 }  \sqrt {  \frac{\log( T )   }{ T  } } 
      +  \sqrt { \frac{  D  ^ {-2r +1 }    }{  T ^{10}    }    } \bigg). 
 \end{align*}
 Taking $D =   \frac{cT}{\log(T) } $ for  sufficiently small constant $c$ gives, 
  \begin{align}\label{eq:eigen condition term 1}
 \left|   \frac{1}{T} \sum _{t=1} ^T \lb X_t ,  \mathcal P_{\mathcal F } ( \beta) \rb_\lt ^2  -  \Sigma _X [ \mathcal P_{\mathcal F } ( \beta) ,  \mathcal P_{\mathcal F } ( \beta)  ] \right|   
   \le &   \frac{\tau}{2}    \Sigma _X [   \beta  ,      \beta  ]    +      C_4    T^{-4}   
\\ \nonumber 
     \left|   \frac{1}{T} \sum _{t=1} ^T \lb X_t ,  \mathcal P_{\mathcal F ^\perp  } ( \beta) \rb_\lt ^2  -   \Sigma _X [ \mathcal P_{\mathcal F^\perp  } ( \beta) ,  \mathcal P_{\mathcal F ^\perp } ( \beta)  ] \right|      \le & C_4  \bigg( \left(  \frac{\log(T)}{T} \right)^{ 2r -1/2}   +   \frac{1}{T^4 } \bigg)  
     \\ \label{eq:eigen condition term 3}    \le  & C_4\bigg(   \left(  \frac{\log(T)}{T} \right)^{ \frac{2r}{2r+1 }}   +   \frac{1}{T^4 }  \bigg),
 \end{align} 
 where $r>1/2$ is used in the last inequality. 
 \\
 \\
 {\bf  Step 5.}  
Denote $$ \delta_T = \left ( \frac{\log(T)}{T} \right)^{ \frac{2r}{2r+1 }} . $$
Note that  if $k \not =l$, 
\begin{align*}
\mathbb E(w_{i, k}w_{i,l } )   = &  \mathbb E  \left( \frac{  \lb  L_{\kh }(X_t ) , \Phi_k\rb_\lt  }{\sqrt { \xi _k }} \frac{  \lb  L_{\kh }(X_t ) , \Phi_l \rb_\lt  }{\sqrt { \xi  _l }} \right) 
\\
= & \frac{1}{\sqrt {\xi_k \xi _l }}  E   \left(  \lb   X_t   , L_{\kh } ( \Phi_k) \rb_\lt  \lb   X_t   , L_{\kh } ( \Phi_l) \rb_\lt   \right) 
\\
 =  & \frac{1}{\sqrt { \xi_k \xi_l }}  \Sigma _X   [L_{\kh } ( \Phi_k)  , L_{\kh } ( \Phi_l)  ]
 \\
 =
 & \frac{1}{\sqrt {\xi_k\xi_l }}    \lb L_{\kh C \kh  } ( \Phi_k)  ,  \Phi_l   \rb _\lt =0.
\end{align*}
 Therefore 
    \begin{align} \nonumber 
 &  \sup_{\|\beta\|_\hk \le 1 }      \Sigma _X [ \mathcal P_{\mathcal F^\perp  } ( \beta) ,  \mathcal P_{\mathcal F ^\perp } ( \beta)  ]  
  \\
  \nonumber 
\le 
  & \sup_{  \sum_{k  =D+ 1}^\infty f_k^2 \le 1    }  \left| 
   \sum_{k , l =D+1 }^\infty   \sqrt {\xi_k } f_k \sqrt {\xi_ l } f_ l
  \frac{1}{T} \sum _{t=1} ^T  E(w_{i, k}w_{i,l } )  
  \right|  
\\   
\nonumber 
  = 
& \sup_{  \sum_{k  =D+ 1}^\infty f_k^2 \le 1    }  \left| 
   \sum_{k   =D+1 }^\infty   \sqrt {\xi_k } f_k \sqrt {\xi_ k } f_ k
  \frac{1}{T} \sum _{t=1} ^T  E(w_{i, k}w_{i, k } )  
  \right|   
  \\
  \nonumber 
  \le
  & 2  \sup_{  \sum_{k  =D+ 1}^\infty f_k^2 \le 1    }  \left| 
   \sum_{k   =D+1 }^\infty       f_ k^2  k^{-2r} 
  \right|   
  \\
  \nonumber 
  \le 
  & 2\sup_{  \sum_{k  =D+ 1}^\infty f_k^2 \le 1    }   
    \sqrt { \sum_{k   =D+1 }^\infty       f_ k^4}  \sqrt {  \sum_{k   =D+1 }^\infty  k^{-4 r}  }     
    \\ \label{eq:population bias second inequality}
    \le  & 2 D^{-2r+1/2 }  \le C_5   \left(  \frac{\log(T) }{T}\right)^{2r-1/2 } \le C_5 \left ( \frac{\log(T)}{T} \right)^{ \frac{2r}{2r+1 }}  \le C_5 \log(T) \delta_T .
 \end{align} 
  \
\\
Observe that 
$$ \Sigma_X [ \mathcal P_{\mathcal F } ( \beta) ,  \mathcal P_{\mathcal F ^\perp  } ( \beta)  ] =0 .$$
So  for any $\| \beta\|_\hk \le 1$, 
\begin{align}\nonumber 
 &  \left|   \frac{1}{T} \sum _{t=1} ^T  \lb X_t ,  \mathcal P_{\mathcal F } ( \beta) \rb_\lt  \lb X_t ,  \mathcal P_{\mathcal F ^\perp  } ( \beta) \rb_\lt   - \Sigma _X [ \mathcal P_{\mathcal F } ( \beta) ,  \mathcal P_{\mathcal F ^\perp  } ( \beta)  ] \right|     
 \\ \nonumber 
 \le 
 &  \left|  \sqrt {   \frac{1}{T} \sum _{t=1} ^T  \lb X_t ,  \mathcal P_{\mathcal F } ( \beta) \rb_\lt ^2} \sqrt {   \frac{1}{T} \sum _{t=1} ^T  \lb X_t ,  \mathcal P_{\mathcal F ^\perp  } ( \beta) \rb_\lt ^2 }     \right|      
 \\
 \le &  \nonumber     \left(   \Sigma _X [ \mathcal P_{\mathcal F } ( \beta) ,  \mathcal P_{\mathcal F } ( \beta)  ]    +
\frac{\tau}{2} \Sigma _X [  \beta  ,    \beta   ]    +C_4  T^{-4}   \right)^{1/2 } 
 \left( \Sigma _X [ \mathcal P_{\mathcal F^\perp  } ( \beta) ,  \mathcal P_{\mathcal F ^\perp } ( \beta)  ]   + C_4\bigg( \delta_T    +   T ^{-4} \bigg) \right) ^{1/2 }   
 \\
\le &   \nonumber    \left( 2 \Sigma _X [  \beta  ,   \beta   ]      + T^{-4}   \right)^{1/2 } 
 \left(  C_ 6 \delta_T    +   T ^{-4} \right) ^{1/2 }    
 \\ \label{eq:eigen condition term 2 rate}
 \le&  \frac{\tau}{2}\Sigma _X [  \beta  ,   \beta   ] + C_7 \delta_T , 
\end{align}
 where the second inequality follows from \eqref{eq:eigen condition term 1} and \eqref{eq:eigen condition term 3}, and the third inequality follows from \eqref{eq:population bias second inequality}, and the last inequality holds if $C_7$ is a sufficiently large constant.
The desired results follows from \Cref{eq:eigen condition term 1},  \Cref{eq:eigen condition term 3} and  \Cref{eq:eigen condition term 2 rate}.
 
 \end{proof}

\subsection{Additional lemmas for kernel alignment }

\begin{remark}\label{remark:inverse of K} Note that following the same argument as that in the proof of  Theorem 2 in  \cite{cai2012minimax},  the function $ L_{ \mk ^{1/2}} : \lt \to \hk  $ admits a well-defined inverse function $L_{ \mk ^{-1/2}} :  \hk   \to \lt  $ such that 
 $$L_{ \mk ^{-1/2}} (\phi_k) = \mu_k^{-1/2} \phi_k, $$
 where the eigen-expansion of  $\mathbb K (r,s)$ satisfies $ \mathbb K (r,s)= \sum_{k=1}^\infty \mu_k \phi_k(r) \phi_k (s)$.
 \end{remark}
 \
 \\
Recall that the eigen-expansion of the linear map $L_{ \mk ^{1/2} \Sigma_X  \mk ^{1/2} }  $ takes the form 
$$ \mk ^{1/2} \Sigma_X  \mk ^{1/2}  = \sum_{k=1}^\infty  \xi_k   \Phi_k .  $$

\begin{lemma} \label{lemma:connect L_2 to H}
For any $f,g \in \hk$, 
$$\lb f ,g\rb_\hk    = \lb L_\kmh (f) , L_\kmh (g) \rb _\lt .  $$
In addition, for any $\lt$ basis $\{\Phi_k \}_{k=1}^\infty$, we have 
$\{ L_\kh (\Phi_k)  \}_{k=1}^\infty$ is a collection of  $\hk$ basis. (However, $\{ L_\kh (\Phi_k)  \}_{k=1}^\infty$ is  not necessary a collection of  $\lt $ basis.) 
\end{lemma}
\begin{proof}
Let $ f =    \sum_{k=1}^\infty a_k\phi  _k \in \hk $ and $ g =    \sum_{k=1}^\infty b_k\phi  _k \in \hk $. 
Note that $ L_\kmh (f) =\sum_{k=1}^\infty \frac{a_k}{\sqrt {\mu_k} } \phi_k $. 
As a result 
$$\| L_\kmh (f) \|_\lt^2 =  \sum_{k=1}^\infty \frac{a_k^2}{\mu_k} = \| f\| _\hk .$$
So $L_\kmh (f) \in\lt$. 
Then 
$$\lb f ,g\rb_\hk = \sum_{k=1}^\infty \frac{a_kb_k}{\mu_k  } =  \lb L_\kmh (f) , L_\kmh (g) \rb _\lt.  $$
In addition, from the above equality, 
$$\lb L_\kh  ( \Phi_i )  ,L_\kh  (\Phi_j)  \rb_\hk  =  \lb L_\kmh  L_\kh  (\Phi_i ) , L_\kmh L_\kh   (\Phi_j ) \rb _\lt =\begin{cases}
1 \text{ if } i=j;
\\
0 \text{ if } i \not=j. 
\end{cases} $$

\end{proof}

\begin{lemma}\label{lemma:change of basis 1} For any $\beta \in \hk$, 
it holds that
$$ \sum_{k= 1}^\infty  \lb L_\kmh (\beta) ,\Phi_k\rb^2 _\lt = \| \beta\|_\hk^2. $$
If  in addition  $ \beta \in \mathcal F ^ \perp$, where 
$$\mathcal F^\perp  :  =  \text{span} \left\{     L_\kh  (\Phi_k )  \right\} _{k=D+1}^\infty,$$
and $D$ is any positive integer,
then
 $$  \sum_{k=D+1}^\infty \lb L_\kmh (\beta) ,\Phi_k\rb^2 _\lt =\| \beta\|_\hk^2 $$
\end{lemma}
\begin{proof}
Denote  $\beta= \sum_{k= 1}^\infty L_\kh (\Phi_k)b_k $. Then 
$   \| \beta\|_\hk^2 =\sum_{k= 1}^\infty b_k^2 . $  So  
$$ \sum_{k= 1}^\infty  \lb L_\kmh (\beta) ,\Phi_k\rb^2 _\lt  =\sum_{k= 1}^\infty b_k^2 =\| \beta\|_\hk^2. $$
For the  second part, it suffices to observe  that 
$\beta= \sum_{k= D+1}^\infty L_\kh (\Phi_k)b_k$.
\end{proof} 

Suppose $W(r)$ is any Gaussian process with covariance operator 
$$ \E(W(r)W(s) ) = \Sigma _W (r,s ), $$ 
and  the  eigen-expansion 
 of  $\Sigma_W$ satisfies 
 $$\Sigma_W (r,s)   = \sum_{k =1}^\infty \sigma_k^2  h _k (r) h _k(s)  .$$ Then it holds that 
$$  W (r ) = \sum_{k=1}^\infty \alpha _k   h _k(r ) ,$$
where $ \{\alpha _k\}_{k=1}^\infty$ are independent  Gaussian random variables  with mean 0 and variance $\sigma_k^2$. Denote 
$ f   ( r ) = \sum_{k=1}^\infty b  _k  h _k(r ) \in \lt $. 
Then 
$\lb W,f \rb _\lt  $ is Gaussian with mean 0 and variance 
$$ \sum_{k=1}^\infty b _k^2  \sigma_k^2 = \E( \lb W,f \rb _\lt  ^2 ) =   \iint_{[0,1]^2 } f(r )   \Sigma_W ( r, s  ) f(s ) \dint r \dint s  =\lb L_{\Sigma_W }(f) ,f\rb_\lt . $$

\begin{lemma}  \label{lemma:variance of transpose of Gaussian}Suppose $R$ is symmetric.  Let $W=L_R  (X)  $.
Then 
$$ \E   (    W  (s) W(r ) ) = ( R\Sigma_X R)(s,r ).$$ Thus  
if $X$ is a  Gaussian process, then  $L_R  (X)$ is also a  Gaussian   process with covariance  function $R\Sigma_X R$.
\end{lemma}
\begin{proof}
For the first part, it suffices to observe that 
$$ \E(W(s)W( r ) ) = \iint_{[0,1]^2 }   R (s, u)\E (X(u)X(v)  )  R(r ,v) \dint  u\dint  v  = \iint _{[0,1]^2 }  R (s, u) \Sigma_X (u,v)  R(r ,v) \dint  u\dint  v . $$
For the second part, it suffices to observe that 
$ \lb L_R(X), f\rb _\lt = \lb X, L_R(f) \rb_\lt $ is  Gaussian with standard deviation 
$$     \sqrt {  \E(\lb X, L_{R  } ( f)  \rb _\lt  ^2 ) }   =  \sqrt { \E(\lb  L_{R  }  (X),  f   \rb _\lt  ^2 ) }  =\sqrt {  \iint _{[0,1]^2 } f(s) (R\Sigma_X R)(s,r )f(r )\dint  s\dint  r } .
  $$
\end{proof}
\begin{lemma} \label{lemma:decomposition conjugate}
Let $ X$ be a centered Gaussian process. Suppose eigenvalues of $\{ \xi_k\}_{k=1}^\infty  $ of  the linear operator  $L_{ \mk ^{1/2} \Sigma_X  \mk ^{1/2} }  $ satisfy
\begin{align*}  \xi_k \asymp   k^{-2r }
\end{align*} for some $r>1/2$. Let  $w_k = \frac{  \lb  L_{\kh }(X) , \Phi_k\rb_\lt  }{\sqrt {\xi _k }}  $. Then 
$w_k  $ is  Gaussian   with variance 1.
\end{lemma}
\begin{proof} Since $   \lb  L_{\kh }(X) , \Phi_k\rb_\lt   $ is  Gaussian  with standard deviation 
$$ \sqrt {  \iint_{[0,1]^2 }  \Phi_k (s) L_{ \kh  \Sigma _X  \kh }  (s,r )   \Phi_k (r )     \dint s \dint r }  = \sqrt {\xi_k }.$$
  So $w_k$ is  $\mathcal N(0,1)$.
\end{proof}
\begin{lemma} Under the same conditions as in \Cref{lemma:decomposition conjugate},  \label{lemma:perpendicular}
it holds that 
$$ \Sigma _X [  L_\kh  (\Phi_k ) ,  L_\kh  (\Phi_l ) ] =\begin{cases}
\xi_k,  &\text{if } k=l ; 
\\
 0, &\text{ if } k \not = l;
\end{cases} \quad k, l = 1,2,3 ,\ldots$$
\end{lemma} 
\begin{proof}
It suffices to observe that 
$$ \Sigma _X [  L_\kh  (\Phi_k ) ,  L_\kh  (\Phi_l ) ]  = 
\lb L_{   \Sigma _X }  ( L_{\kh }\Phi_k ) ,L_{\kh }  \Phi_l \rb_\lt 
 = \lb L_{\kh  \Sigma _X \kh }\Phi_k , \Phi_l \rb_\lt  =  \lb \xi_k \Phi_k , \Phi_l \rb_\lt  .$$
\end{proof}

%% file: moresimulation.tex
\section{Technical details of the optimization}\label{sec-app-F}
\subsection{Detailed formulation of the convex optimization }\label{subsec:convex_formulation}
With the notation defined as in Section \ref{subsec:formulation_convex}, in the following, we provide the details for rewriting the penalized optimization problem in \eqref{eq:approx A beta_penalized} into \eqref{eq:optimization}, which is a convex function of $R$ and $B=\left[\mb_1,\mb_2,\cdots,\mb_p\right]$. Specifically, we examine the three components of \eqref{eq:approx A beta_penalized} one by one.
\begin{itemize}
	\item The squared loss can be rewritten as
	\begin{align*}
		&\sum_{t= 1}^T    \sum_{j=1}^{n_2} w_r(j)   \left (Y_{t} (r _j) -    \frac{1}{n_1}\sum_{i=1} ^{n_1}w_s(i) A  (r _j, s_i)  X_{t  } (s_i)   - \lb \beta (r_j),  Z_t  \rb_p \right  )^2\\
		=&\sum_{t=1}^T \left\|W_R\left(Y_t- \frac{1}{n_1}K_1^\top RK_2 W_S X_t - K_1^\top B Z_t\right)\right\|^2_2= \left\|W_R\left(Y-\frac{1}{n_1}K_1RK_2 W_SX-K_1B Z\right)\right\|_\mathrm{F}^2,
	\end{align*}
	where $\|\cdot\|_2$ and $\|\cdot\|_{\mathrm{F}}$ are the $\ell_2$-norm of a vector and the Frobenius norm of a matrix.
	
	\item As for the Frobenius norm penalty $\|A\|_{\mathrm F(\mk)}^2$, first, it is easy to see ${A}^\top(r,s)=k_1(s)^\top Rk_2(r)$, where $A^\top$ is the adjoint operator of $A$. Let $u(s)=k_2(s)^\top c$, for any $c \in \mathbb{R}^{n_1}$.  We have that
	\begin{align*}
		&{A}^\top {A}[u](s)=\lb A^\top(s,r), A[u](r) \rb_\h=\lb k_1(r)^\top Rk_2(s), A[u](r) \rb_\h\\
		=&\lb k_2(s)^\top R^\top k_1(r), \lb A(r,s), u(s) \rb_\h \rb_\h=k_2(s)^\top R^\top \lb  k_1(r), k_1(r)^\top \rb_\h R \lb k_2(s), u(s)\rb_\h\\
		=&k_2(s)^\top R^\top K_1 R K_2 c.
	\end{align*}
	Thus, the eigenvalues of ${A}^\top {A}$ are the same as those of $R^\top K_1 R K_2$ and $\|A\|_{\mathrm{F}(\mk)}^2 = \mathrm{tr}(R^\top K_1 RK_2)$.
	
	\item As for the $\h(\mk)$-norm penalty $\|\beta_l\|_{\h(\mk)}$ and the group Lasso-type penalty $\|\beta_l\|_{n_2}$, we have $\|\beta_l\|_{\h(\mk)}=\sqrt{\lb \beta_l(r), \beta_l(r) \rb_{\h(\mk)}}=\sqrt{\mb_l^\top K_1 \mb_l}$ and
	\begin{align*}
		&\|\beta_l\|_{n_2} = \sqrt{\frac{1}{n_2}\sum_{j=1}^{n_2}w_r(j)\beta_l^2(r_j)}=\sqrt{\frac{1}{n_2}\sum_{j=1}^{n_2}w_r(j)\mb_l^\top k_1(r_j)k_1(r_j)^\top \mb_l}\\
		=&\sqrt{\mb_l^\top\frac{1}{n_2}\sum_{j=1}^{n_2} w_r(j)k_1(r_j)k_1(r_j)^\top \mb_l}=\sqrt{\frac{1}{n_2}\mb_l^\top K_1 W_R^2 K_1\mb_l}, \text{ for } l=1, \ldots, p. 
	\end{align*}
\end{itemize}
Combining all three components together, it is easy to see that the optimization problem in \eqref{eq:approx A beta_penalized} can be written as \eqref{eq:optimization}.

\subsection{Detailed optimization of the structured ridge regression \eqref{eq:optimization_nolasso1}}\label{subsec:opt_alg}
By simple linear algebra, the first order condition forms the linear system
\begin{align*}
		(S_1^\top S_1+\lambda_1 I_{n_1 n_2} )\me = S_1^\top \widetilde{\my}.
\end{align*}
Denote
\begin{align*}
	S(\lambda_1)=S_1^\top S_1+\lambda_1 I_{n_1 n_2}.
\end{align*}
Note that $S(\lambda_1)$ is a positive definite matrix and the optimization has a unique solution. However, when $n_1n_2$ is large, the inverse of the matrix $S(\lambda_1)$ becomes computationally expensive. We thus further exploit the structure of $S(\lambda_1)$.

Denote $S_4=\frac{1}{n_1}K_2^{1/2}W_SX \in \mathbb{R}^{n_1\times T}$, we have $S_1^\top S_1=(S_4S_4^\top)\otimes (K_1^{1/2} W_R^2 K_1^{1/2}),$ and thus
\begin{align*}
	S(\lambda_1)=(S_4S_4^\top)\otimes (K_1^{1/2} W_R^2 K_1^{1/2})+\lambda_1 I_{n_1n_2}.
\end{align*}
Denote the singular value decomposition~(SVD) of $S_4S_4^\top=\frac{1}{n_1^2}K_2^{1/2}W_SXX^\top W_S K_2^{1/2}$ as $S_4S_4^\top=U_1D_1U_1^\top$ and the SVD of $K_1^{1/2} W_R^2 K_1^{1/2}$ as $K_1^{1/2} W_R^2 K_1^{1/2}=U_2D_2U_2^\top$, we have
\begin{align*}
	S(\lambda_1)=(S_4S_4^\top)\otimes (K_1^{1/2} W_R^2 K_1^{1/2})+\lambda_1 I_{n_1n_2}=(U_1\otimes U_2)(D_1\otimes D_2 +\lambda_1 I_{n_1n_2})(U_1^\top\otimes U_2^\top).
\end{align*}
Thus, we have
\begin{align*}
	S^{-1}(\lambda_1)&=(U_1\otimes U_2)(D_1\otimes D_2 +\lambda_1 I_{n_1n_2})^{-1}(U_1^\top\otimes U_2^\top)\\
	&=\sum_{i=1}^{n_1}\sum_{j=1}^{n_2} \frac{1}{\lambda_1+D_{1i}D_{2j}}(U_{1i}\otimes U_{2j})(U_{1i}^\top \otimes U_{2j}^\top)\\
	&=\sum_{i=1}^{n_1}\sum_{j=1}^{n_2} \frac{1}{\lambda_1+D_{1i}D_{2j}}(U_{1i}U_{1i}^\top)\otimes (U_{2j} U_{2j}^\top).
\end{align*}
Thus, we have
\begin{align*}
	\me = S^{-1}(\lambda_1)S_1^\top \widetilde{\my} =
	\sum_{i=1}^{n_1}\sum_{j=1}^{n_2} \frac{1}{\lambda_1+D_{1i}D_{2j}}(U_{1i}U_{1i}^\top S_4)\otimes (U_{2j} U_{2j}^\top K_1^{1/2}W_R)\widetilde{\my},
\end{align*}
which implies that
\begin{align*}
	E=&\sum_{i=1}^{n_1}\sum_{j=1}^{n_2} \frac{1}{\lambda_1+D_{1i}D_{2j}}(U_{2j} U_{2j}^\top K_1^{1/2}W_R)\widetilde{Y}(S_4^\top U_{1i}U_{1i}^\top)\\
	=&\frac{1}{n_1}\sum_{i=1}^{n_1}\sum_{j=1}^{n_2} \frac{1}{\lambda_1+D_{1i}D_{2j}}(U_{2j} U_{2j}^\top) (K_1^{1/2}W_R\widetilde{Y}X^\top W_S K_2^{1/2}) (U_{1i}U_{1i}^\top),
\end{align*}
and we update $R=K_1^{-1/2}EK_2^{-1/2}.$

\section{Additional numerical analysis}\label{subsec:additional_simu}


\begin{figure}[ht]
	\begin{subfigure}{0.32\textwidth}
		\includegraphics[angle=270, width=1.25\textwidth]{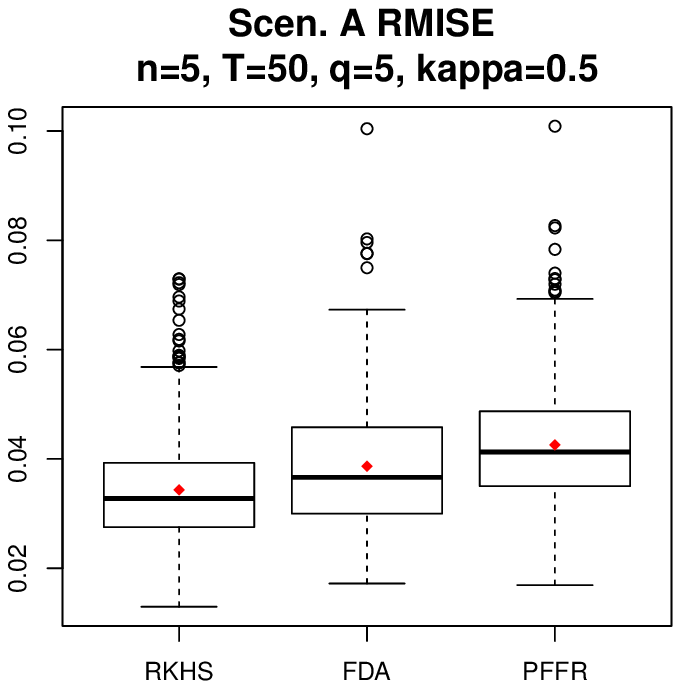}
		\vspace{-0.8cm}
	\end{subfigure}
	~
	\begin{subfigure}{0.32\textwidth}
		\includegraphics[angle=270, width=1.25\textwidth]{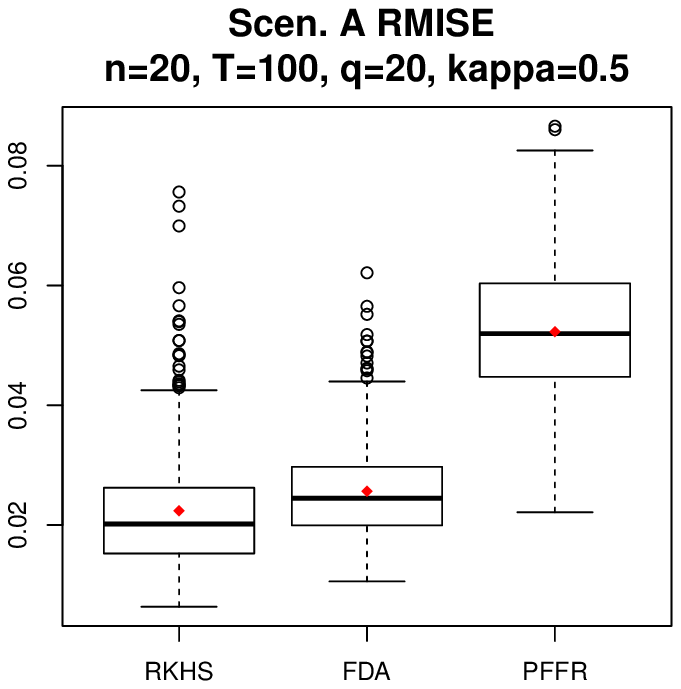}
		\vspace{-0.8cm}
	\end{subfigure}
	~
	\begin{subfigure}{0.32\textwidth}
		\includegraphics[angle=270, width=1.25\textwidth]{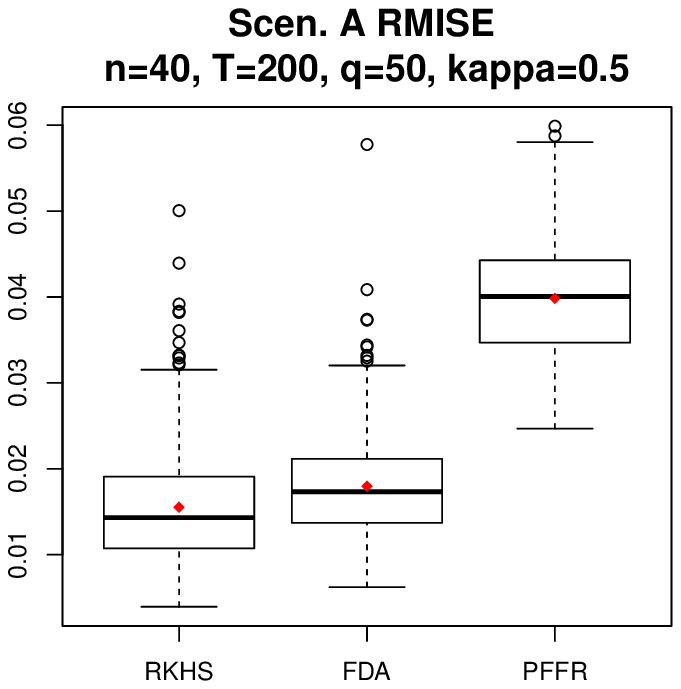}
		\vspace{-0.8cm}
	\end{subfigure}
	~
	\begin{subfigure}{0.32\textwidth}
		\includegraphics[angle=270, width=1.25\textwidth]{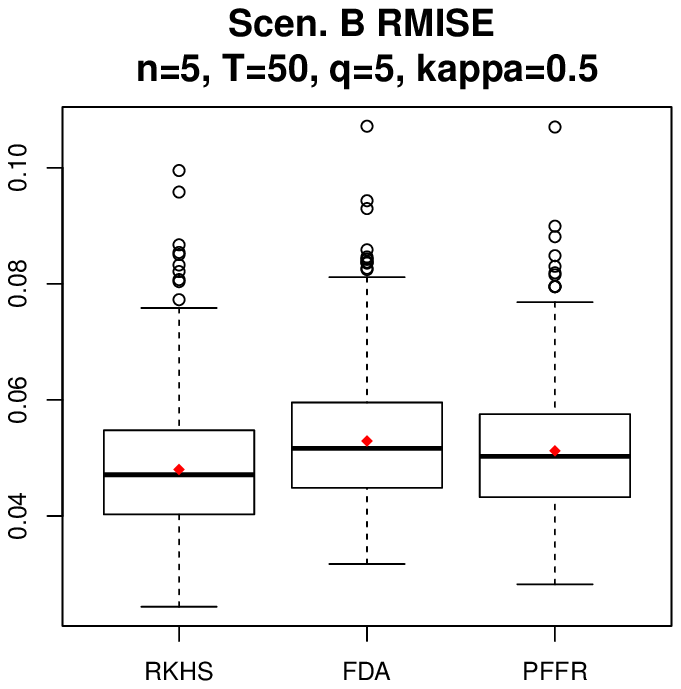}
		\vspace{-0.8cm}
	\end{subfigure}
	~
	\begin{subfigure}{0.32\textwidth}
		\includegraphics[angle=270, width=1.25\textwidth]{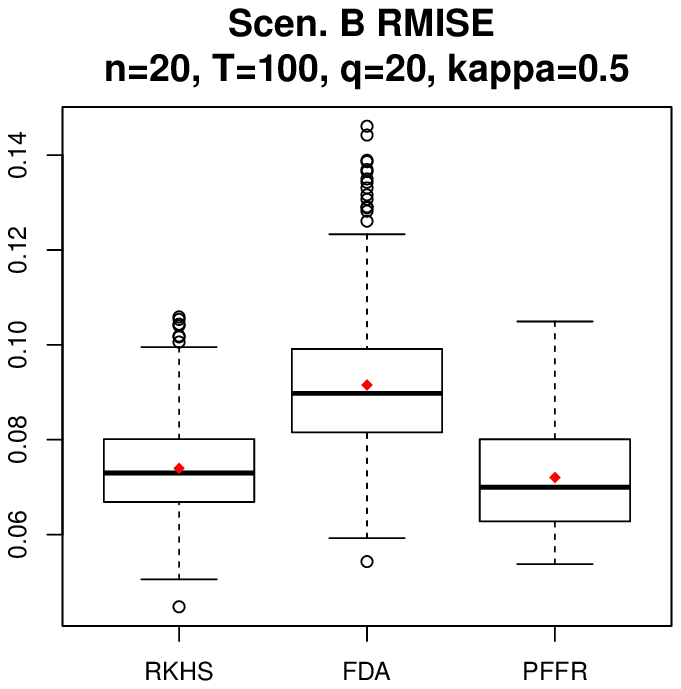}
		\vspace{-0.8cm}
	\end{subfigure}
	~
	\begin{subfigure}{0.32\textwidth}
		\includegraphics[angle=270, width=1.25\textwidth]{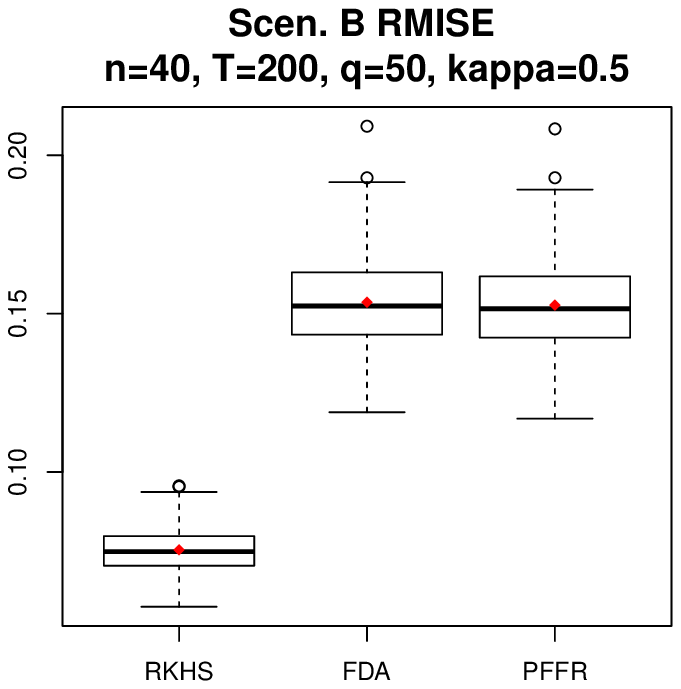}
		\vspace{-0.8cm}
	\end{subfigure}
	\caption{Boxplot of RMISE of RKHS, FDA and PFFR across 500 experiments under function-on-function regression with $\kappa=0.5$. Red points denote the average RMISE.}
	\label{fig:FR_k05}
\end{figure}

\begin{figure}[ht]
	\begin{subfigure}{0.32\textwidth}
		\includegraphics[angle=270, width=1.25\textwidth]{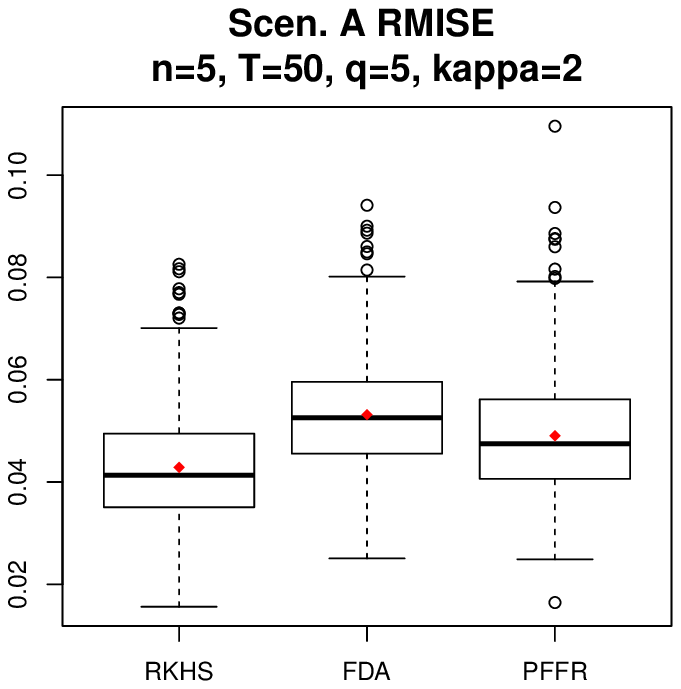}
		\vspace{-0.8cm}
	\end{subfigure}
	~
	\begin{subfigure}{0.32\textwidth}
		\includegraphics[angle=270, width=1.25\textwidth]{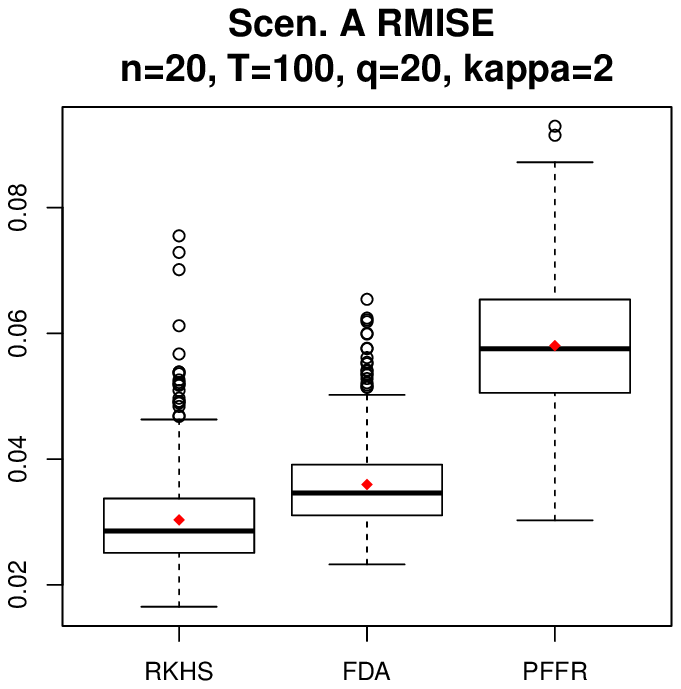}
		\vspace{-0.8cm}
	\end{subfigure}
	~
	\begin{subfigure}{0.32\textwidth}
		\includegraphics[angle=270, width=1.25\textwidth]{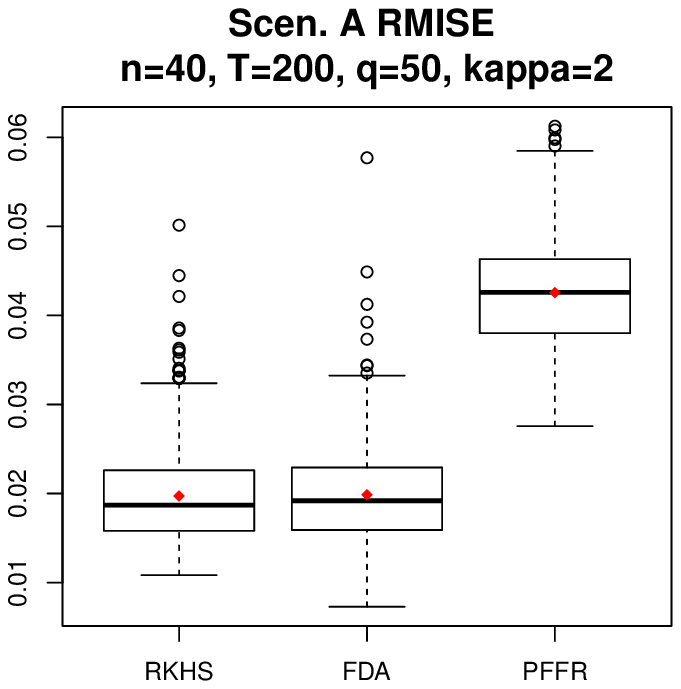}
		\vspace{-0.8cm}
	\end{subfigure}
	~
	\begin{subfigure}{0.32\textwidth}
		\includegraphics[angle=270, width=1.25\textwidth]{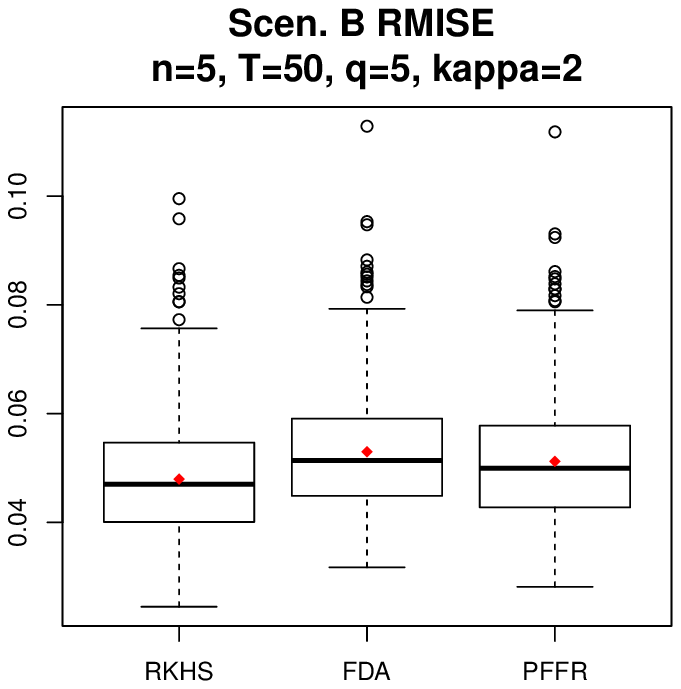}
		\vspace{-0.8cm}
	\end{subfigure}
	~
	\begin{subfigure}{0.32\textwidth}
		\includegraphics[angle=270, width=1.25\textwidth]{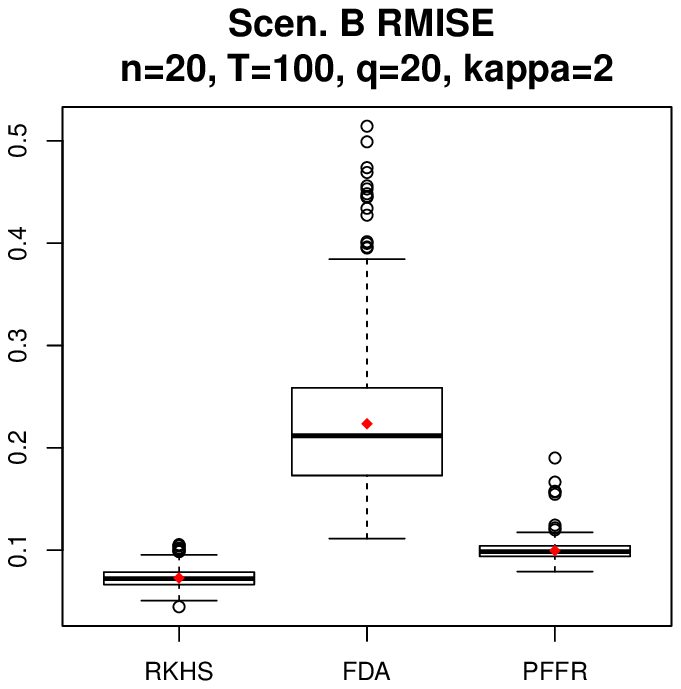}
		\vspace{-0.8cm}
	\end{subfigure}
	~
	\begin{subfigure}{0.32\textwidth}
		\includegraphics[angle=270, width=1.25\textwidth]{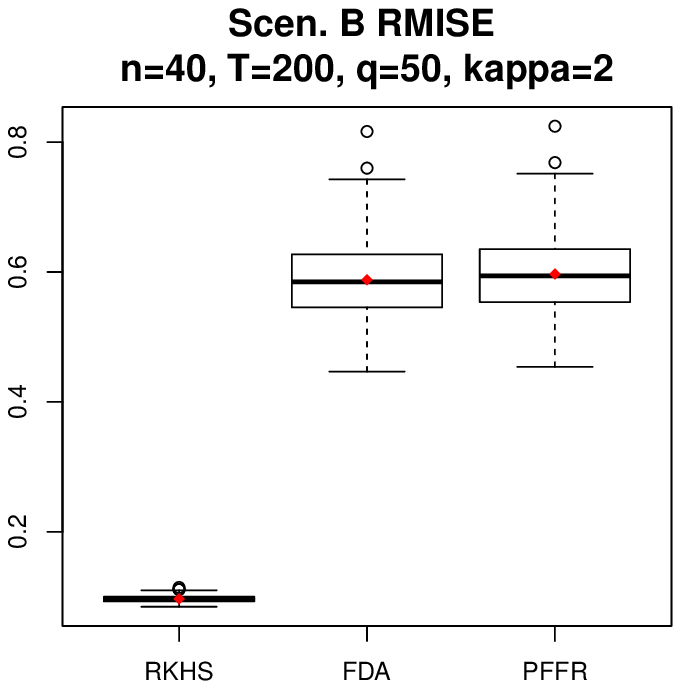}
		\vspace{-0.8cm}
	\end{subfigure}
	\caption{Boxplot of RMISE of RKHS, FDA and PFFR across 500 experiments under function-on-function regression with $\kappa=2$. Red points denote the average RMISE.}
	\label{fig:FR_k20}
\end{figure}

\begin{figure}[ht]
	\begin{subfigure}{0.32\textwidth}
		\includegraphics[angle=270, width=1.25\textwidth]{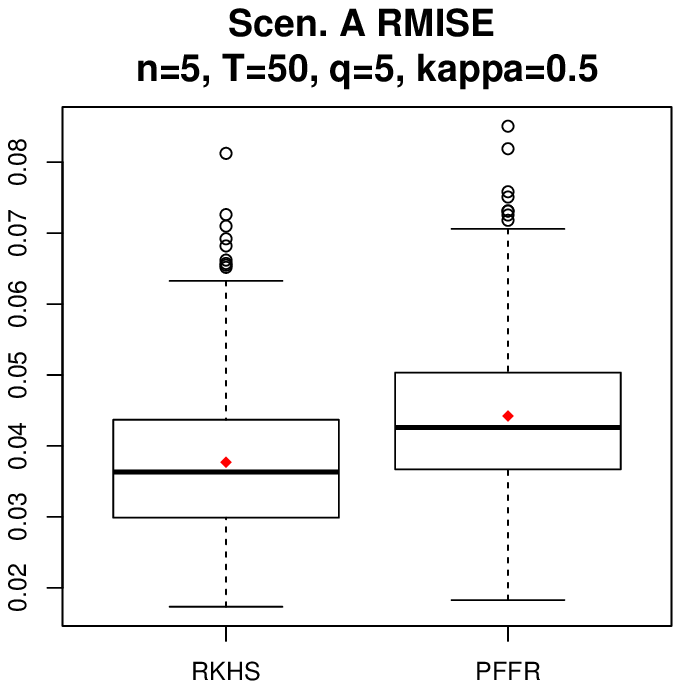}
		\vspace{-0.8cm}
	\end{subfigure}
	~
	\begin{subfigure}{0.32\textwidth}
		\includegraphics[angle=270, width=1.25\textwidth]{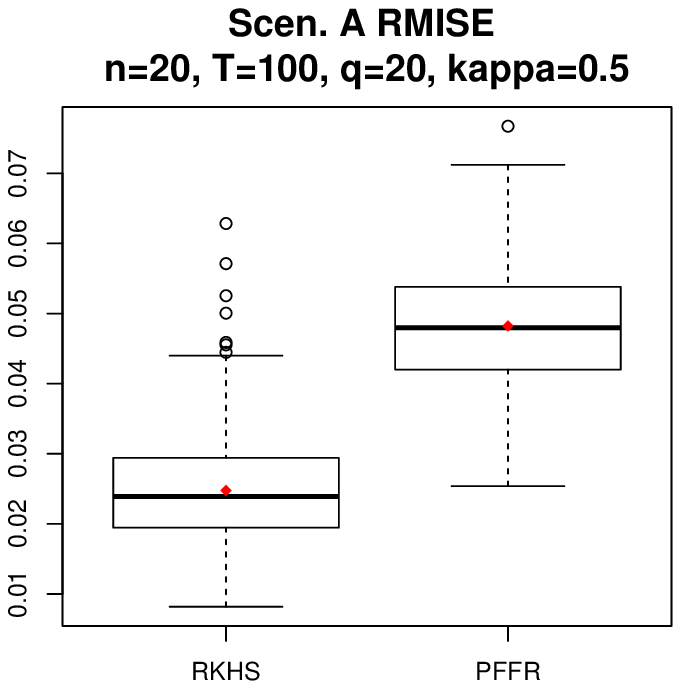}
		\vspace{-0.8cm}
	\end{subfigure}
	~
	\begin{subfigure}{0.32\textwidth}
		\includegraphics[angle=270, width=1.25\textwidth]{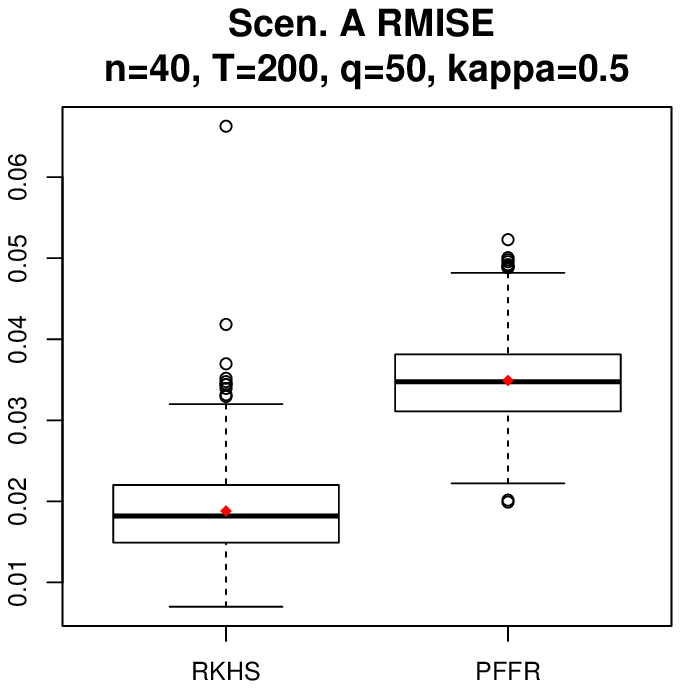}
		\vspace{-0.8cm}
	\end{subfigure}
	~
	\begin{subfigure}{0.32\textwidth}
		\includegraphics[angle=270, width=1.25\textwidth]{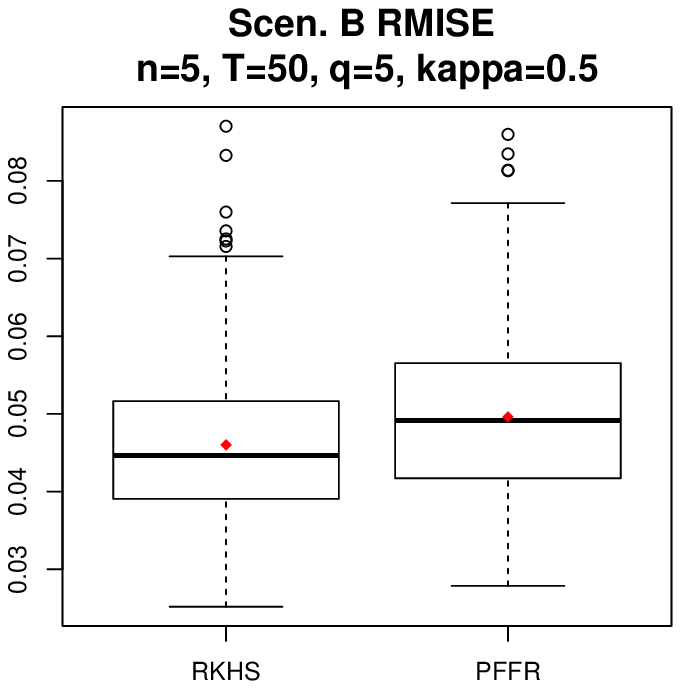}
		\vspace{-0.8cm}
	\end{subfigure}
	~
	\begin{subfigure}{0.32\textwidth}
		\includegraphics[angle=270, width=1.25\textwidth]{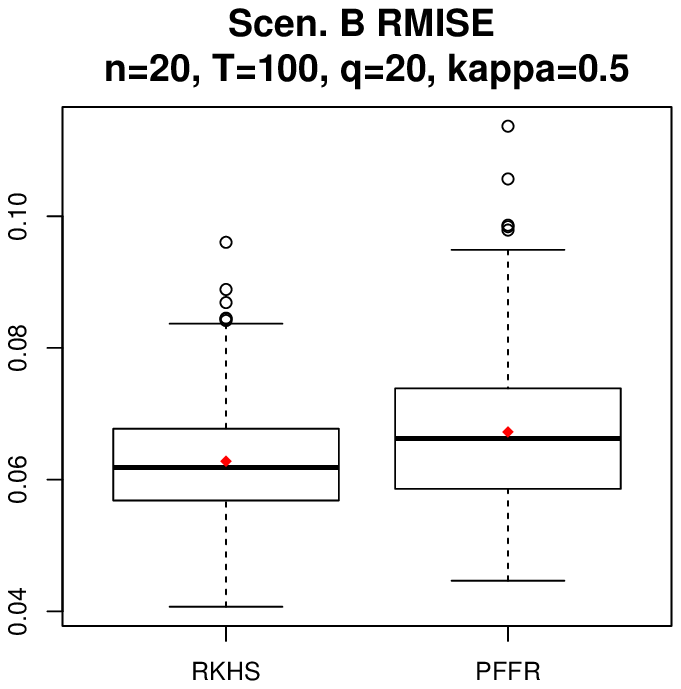}
		\vspace{-0.8cm}
	\end{subfigure}
	~
	\begin{subfigure}{0.32\textwidth}
		\includegraphics[angle=270, width=1.25\textwidth]{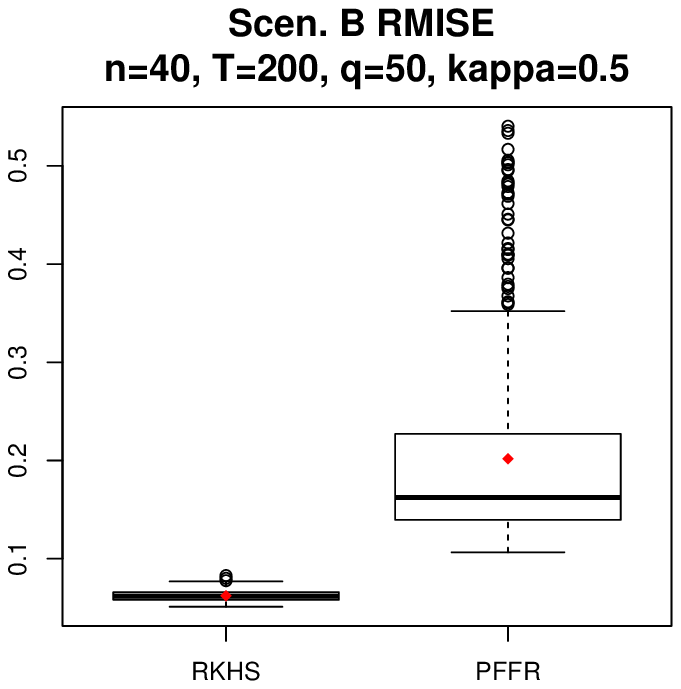}
		\vspace{-0.8cm}
	\end{subfigure}
	\caption{Boxplot of RMISE of RKHS and PFFR across 500 experiments under mixed functional regression with $\kappa=0.5$. Red points denote the average RMISE.}
	\label{fig:FR_mixed_k05}
\end{figure}

\begin{figure}[ht]
	\begin{subfigure}{0.32\textwidth}
		\includegraphics[angle=270, width=1.25\textwidth]{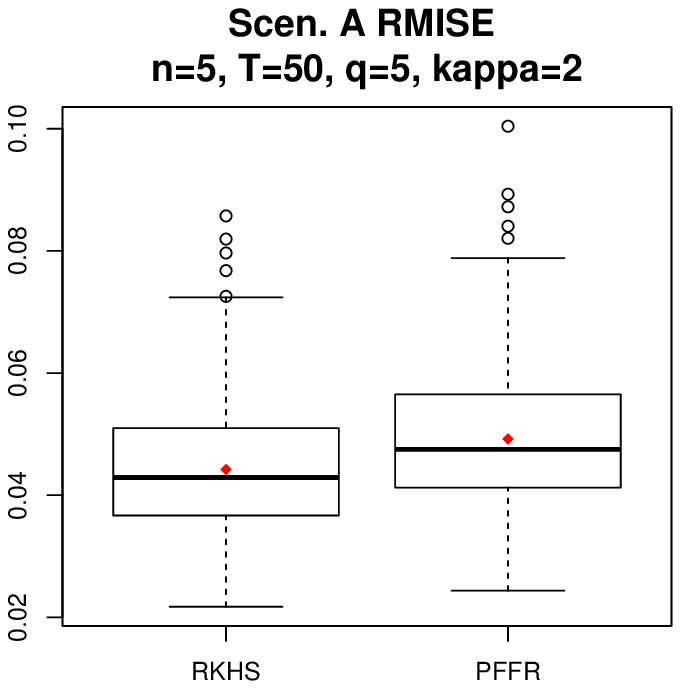}
		\vspace{-0.8cm}
	\end{subfigure}
	~
	\begin{subfigure}{0.32\textwidth}
		\includegraphics[angle=270, width=1.25\textwidth]{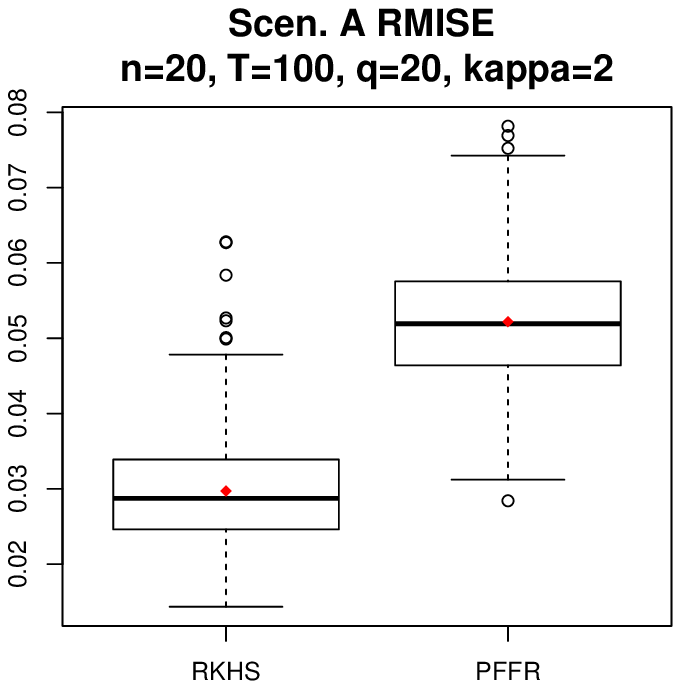}
		\vspace{-0.8cm}
	\end{subfigure}
	~
	\begin{subfigure}{0.32\textwidth}
		\includegraphics[angle=270, width=1.25\textwidth]{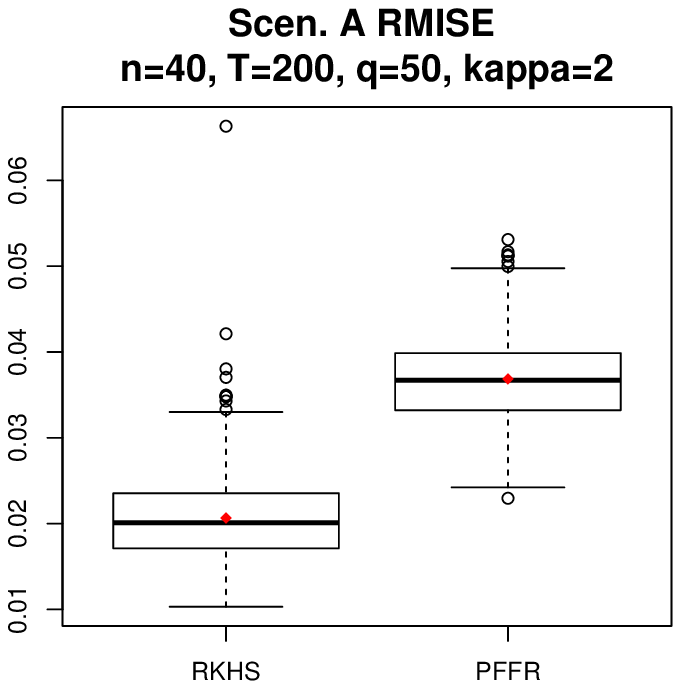}
		\vspace{-0.8cm}
	\end{subfigure}
	~
	\begin{subfigure}{0.32\textwidth}
		\includegraphics[angle=270, width=1.25\textwidth]{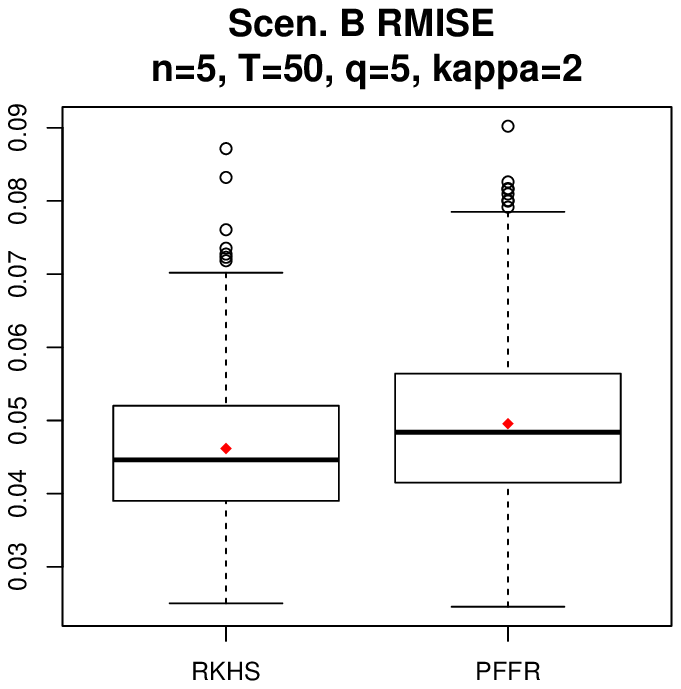}
		\vspace{-0.8cm}
	\end{subfigure}
	~
	\begin{subfigure}{0.32\textwidth}
		\includegraphics[angle=270, width=1.25\textwidth]{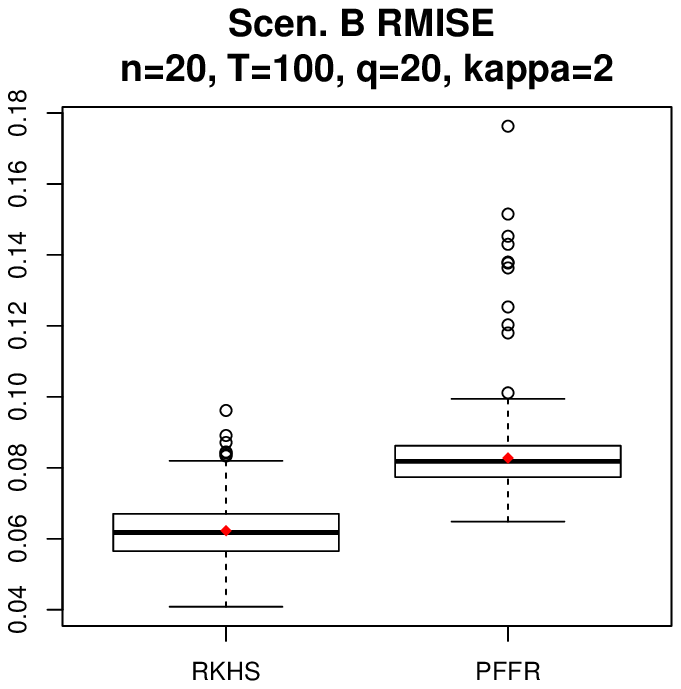}
		\vspace{-0.8cm}
	\end{subfigure}
	~
	\begin{subfigure}{0.32\textwidth}
		\includegraphics[angle=270, width=1.25\textwidth]{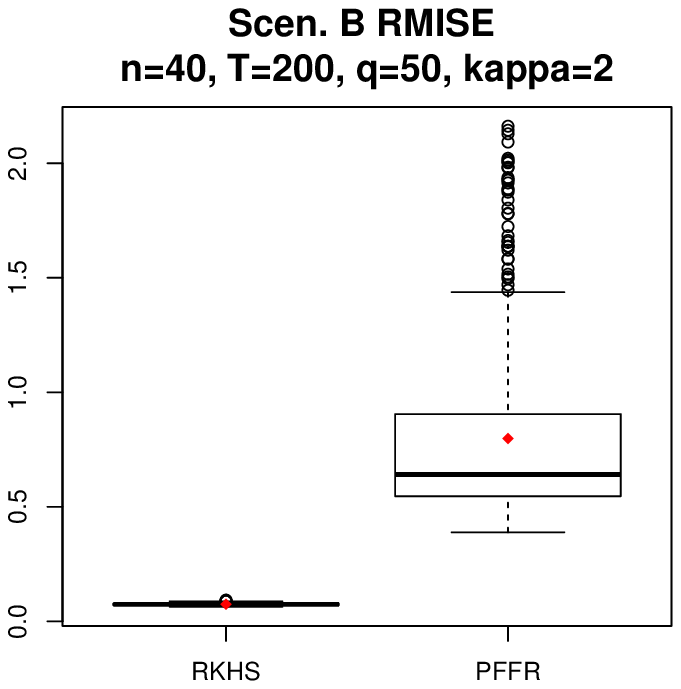}
		\vspace{-0.8cm}
	\end{subfigure}
	\caption{Boxplot of RMISE of RKHS and PFFR across 500 experiments under mixed functional regression with $\kappa=2$. Red points denote the average RMISE.}
	\label{fig:FR_mixed_k20}
\end{figure}

\begin{figure}[h]
	\hspace*{-1cm}
	\centering                                                    
	\begin{subfigure}{0.4\textwidth}
		\includegraphics[angle=270, width=1.2\textwidth]{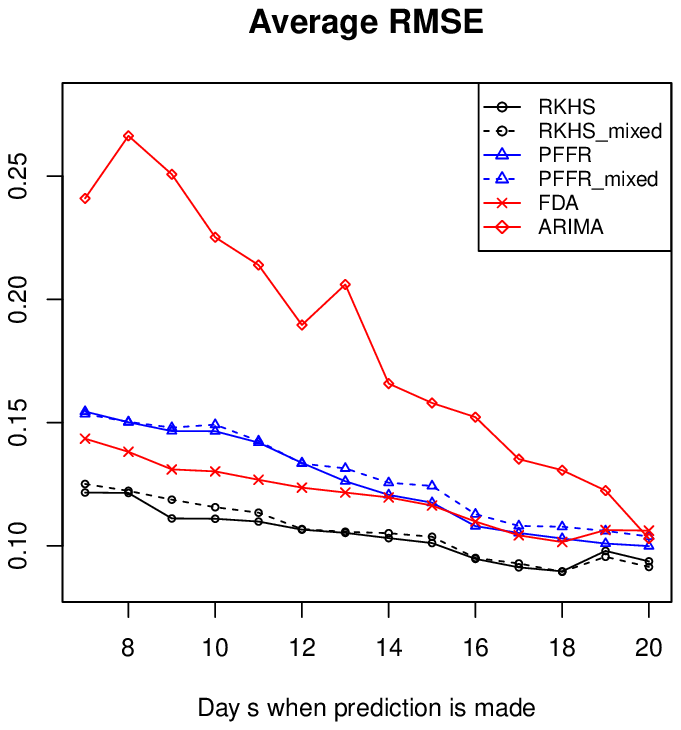}
		\vspace{-0.3cm}
	\end{subfigure}
	~
	\begin{subfigure}{0.4\textwidth}
		\includegraphics[angle=270, width=1.2\textwidth]{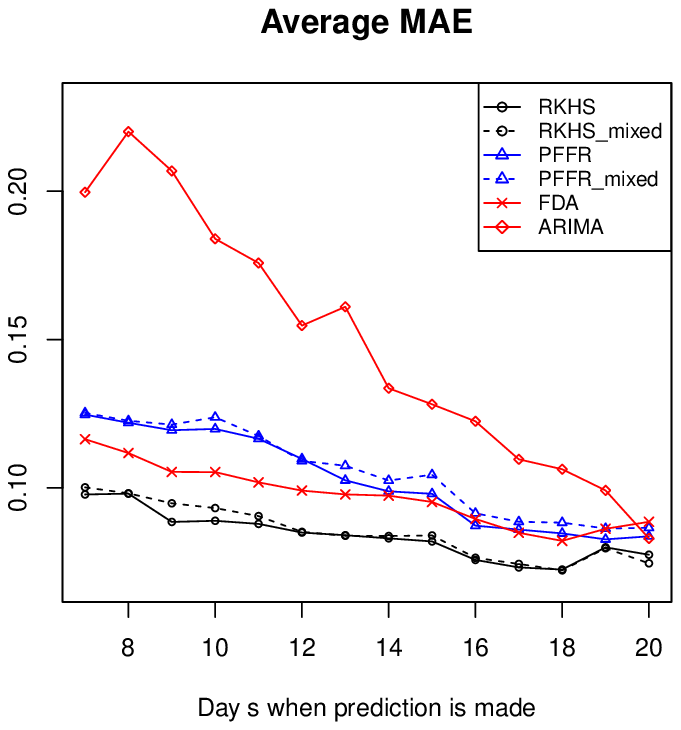}
		\vspace{-0.3cm}
	\end{subfigure}
	\caption{Average RMSE and MAE achieved by different functional regression methods. For campaign $t$ at prediction time $s$, ARIMA is trained on $\{N_t(r_i), r_i\in [0, s]\}$ for the prediction of $\{N_t(r_i), r_i\in (s, 30]\}$. Note that ARIMA can be used as $\{N_t(r_i)\}_{i=1}^{60}$ can be seen as an evenly spaced univariate time series.}
	\label{fig:realdata_witharima}
\end{figure}

\clearpage

\section{ Robustness check for real data analysis}\label{subsec:realdata_cv}
In this section, we further provide a robustness check for the real data analysis presented in \Cref{subsec:realdata}. Specifically, we repeat the two-fold CV procedure in Section \ref{subsec:realdata} independently for 100 times for RKHS, FDA and PFFR\footnote{For simplicity, we do not include RKHS$_{\text{mixed}}$ and PFFR$_{\text{mixed}}$ in the comparison, as the result in Section \ref{subsec:realdata}~(see Figure \ref{fig:realdata}) clearly indicates that the three scalar predictors do not seem to improve the prediction performance.}. In addition, we also implement FDA with 50 basis functions~(FDA50) and PFFR with 30 basis functions~(PFFR30). We further compare with the popular functional PCA based regression method~(FPCA) proposed in \cite{yao2005functional}.

Thus, for every method~(RKHS, FDA, PFFR, FPCA), on each prediction day $s\in \{7,8,9,\cdots,20\}$, we have 100 RMSE$_s$ and 100 MAE$_s$ obtained from the 100 times two-fold CV. We refer to Section \ref{subsec:realdata} for more detailed definition of RMSE and MAE.

\begin{figure}[h]
	\hspace*{-1cm}
	\centering                                                           
	\begin{subfigure}{0.4\textwidth}
		\includegraphics[angle=270, width=1.2\textwidth]{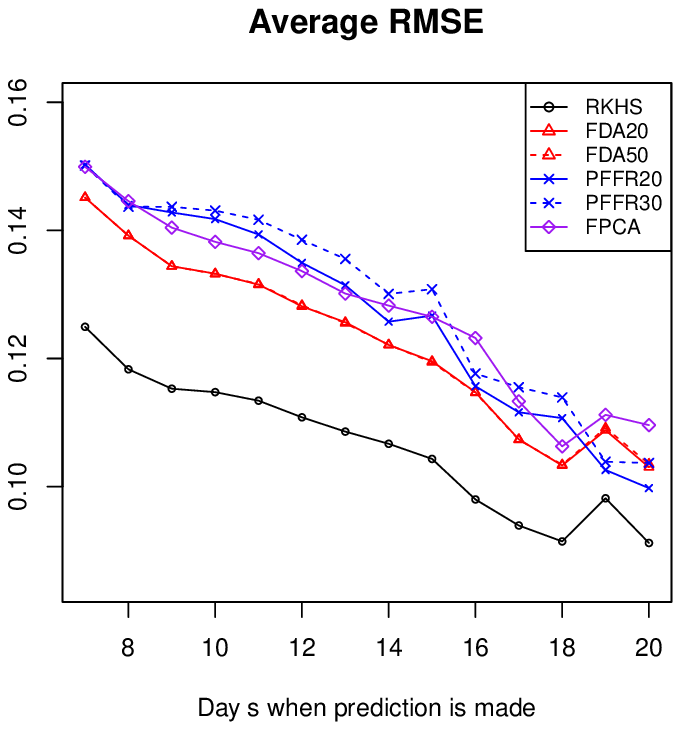}
		\vspace{-0.3cm}
	\end{subfigure}
	~
	\begin{subfigure}{0.4\textwidth}
		\includegraphics[angle=270, width=1.2\textwidth]{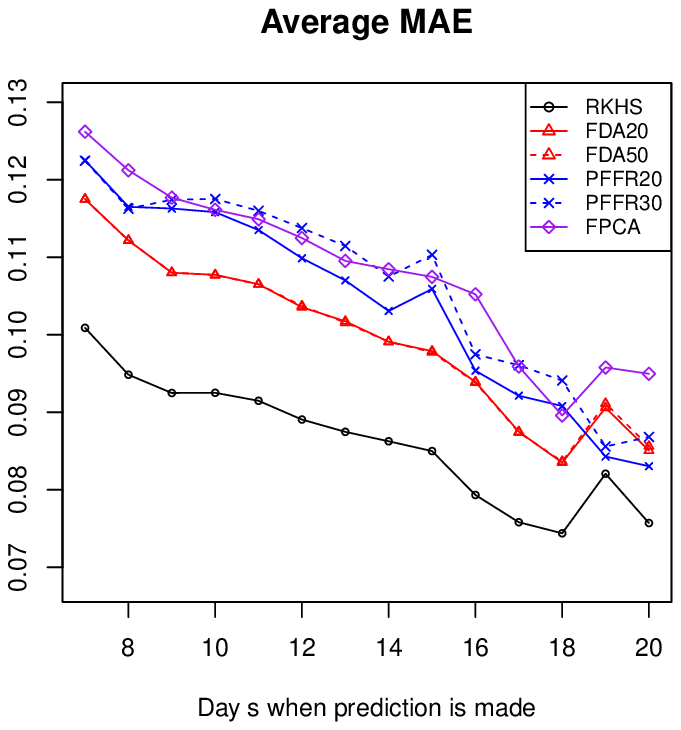}
		\vspace{-0.3cm}
	\end{subfigure}
	\caption{RMSE$_s$ and MAE$_s$ achieved by different functional regression methods averaged over 100 two-fold cross-validations for each $s\in \{7,8,9,\cdots,20\}$.}
	\label{fig:realdata_cv}
\end{figure}

Figure \ref{fig:realdata_cv} visualizes the average RMSE$_s$ and MAE$_s$ over the 100 times CV achieved by different functional regression methods across $s\in \{7,8,9,\cdots,20\}$. As can be seen clearly, the pattern exhibited in Figure \ref{fig:realdata_cv} is the same as the one in Figure \ref{fig:realdata} in Section \ref{subsec:realdata} and RKHS is the winner by a notable margin. Note that FDA20 and FDA50 give essentially the same performance, while PFFR30 indeed gives worse performance than PFFR20. This suggests that 20 basis dimension is sufficient for the current real data application.

Table \ref{tab:realdata_CV} further gives the number of two-fold CV (out of 100 times) where RKHS achieves the smallest RMSE$_s$ and MAE$_s$ for each $s\in \{7,8,9,\cdots,20\}$. As can be seen, RKHS is almost always the winner for every $s\in  \{7,8,9,\cdots,20\}$ except for $s=19$, where RKHS still provides the best performance more than two thirds of the time.

\begin{table}[h]
	\centering
	\begin{tabular}{lllllllllllllll}\hline\hline
		$s$    & 7   & 8   & 9   & 10  & 11  & 12  & 13  & 14  & 15  & 16  & 17  & 18  & 19 & 20 \\\hline
		RMSE$_s$ & 100 & 100 & 100 & 100 & 100 & 100 & 100 & 100 & 100 & 100 & 100 & 100 & 77 & 97\\
		MAE$_s$  & 100 & 100 & 100 & 100 & 100 & 100 & 100 & 100 & 100 & 100 & 100 & 100 & 66 & 96\\
		\hline\hline
	\end{tabular}
	\caption{Number of two-fold cross-validations (out of 100 times) where RKHS achieves the smallest RMSE$_s$ and MAE$_s$ for each $s\in \{7,8,9,\cdots,20\}$.} 
    \label{tab:realdata_CV}
\end{table}

\section{ Numerical results for observations with measurement error}\label{subsec:simu_me}
In this section, we conduct simulation experiments to examine the performance of different methods, including FDA, PFFR and the proposed RKHS, under the setting where the observed functional responses are additionally corrupted with measurement error.

The simulation analysis is done under the function-on-function regression setting as in \Cref{subsec:ffr}\footnote{Note that we use the function-on-function regression setting so that we can include FDA in the analysis as FDA cannot handle the general setting of functional regression with mixed covariates.}. Specifically, for the ease of comparison, we use exactly the same simulation setting that generates Table \ref{tab:FR} in Section \ref{subsec:ffr}. The only difference is that we assume we do not observe the functional response $Y_t(r_j)$ directly but instead a noisy version of it, denoted by $y_{tj}$, corrupted with measurement error such that
\begin{align*}
	y_{tj}=Y_t(r_j)+\mathfrak E_{tj},
\end{align*}
where $\mathfrak E_{tj}$ are independent Gaussian random variable with mean 0 and variance $\sigma^2_{\mathfrak E}$ for $t=1,\ldots, T$ and $j=1,\ldots, n_2.$ Note that for $\sigma^2_{\mathfrak E}=0$, we get back to the identical setting as in Section \ref{subsec:ffr}. It is easy to derive that the overall variability of the functional noise (i.e.\ $\mathbb{E}\int_{[0,1]}\epsilon_t(r) ^2\dint r$) under the simulation setting in Section \ref{subsec:ffr} is around $0.4^2/12\approx 0.115^2$, which can be calculated easily using the fact that $\epsilon_t(r)=\sum_{i=1}^{q} e_{ti}u_i(r)$. We refer to Sections \ref{subsec:simu_setting} and \ref{subsec:ffr} for more details.

Thus, in the following, we set the variance of the measurement error as $\sigma^2_{\mathfrak E}=0.12^2$ such that the variability of the functional noise and the measurement error are comparable. For each setting, we conduct 500 experiments. We again use RMISE to evaluate the excess risk of each estimator. Due to limited space, we only report the result for the case where the spectral norm $\kappa=1$. The results for $\kappa=0.5$ and $\kappa=2$ are similar and omitted. 

For each method, Table \ref{tab:FR_me} reports its average nRMISE across 500 experiments under all simulation settings. Note that Table \ref{tab:FR_me} is directly comparable with Table \ref{tab:FR} in Section \ref{subsec:ffr} (for $\kappa=1$), which reports the performance of each method for $\sigma_{\mathfrak E}^2=0.$ Thus, for ease of comparison, we copy the result reported in Table \ref{tab:FR} (for $\kappa=1$) to Table \ref{tab:FR_me} under the name $\sigma_{\mathfrak E}^2=0$.

Examining Table \ref{tab:FR_me}, it can be seen clearly that, compared to the case of no measurement error~(i.e.\ $\sigma^2_{\mathfrak E}=0$), the performance of all three estimators (RKHS, FDA and PFFR) only worsen slightly when the variability of the measurement error and functional noise are comparable (i.e.\ $\sigma^2_{\mathfrak E}=0.12^2$), suggesting that functional noise is more difficult to handle than the independent measurement error. This can also be viewed as numerical support for the result in Theorem \ref{thm-main-mixed noise}, i.e.\ the convergence rate of RKHS is not affected when the functional response is corrupted by independent mean zero measurement errors with a bounded variance.


\begin{table}[H]
	\centering
	\begin{tabular}{lrrrrr|rrrrr}
		\hline\hline
		& \multicolumn{5}{c|}{ Scenario A:  $n= 5 , T= 50, q=5$} & \multicolumn{5}{c}{ Scenario B:  $n= 5 , T= 50, q=5 $} \\
		$\sigma_{\mathfrak E}^2$ & RKHS & FDA & PFFR &R$_{avg}$(\%) & R$_{w}$ (\%) & RKHS & FDA & PFFR &R$_{avg}$(\%) & R$_{w}$ (\%)\\
		\hline
		$0$ & 9.14 & 10.31 & 10.98 & 12.81 & 68 & 10.42 & 11.50 & 11.12 & 6.63 & 71 \\ 
		$0.12^2$ & 10.14 & 11.40 & 12.20 & 12.45 & 64 & 12.61 & 13.87 & 13.23 & 4.96 & 64\\
		\hline
		& \multicolumn{5}{c|}{ Scenario A:  $n= 20 , T= 100, q=20 $} & \multicolumn{5}{c}{ Scenario B:  $n= 20 , T= 100, q=20 $} \\
		$\sigma^2_{\mathfrak E}$ & RKHS & FDA & PFFR &R$_{avg}$(\%) & R$_{w}$ (\%) & RKHS & FDA & PFFR &R$_{avg}$(\%) & R$_{w}$ (\%)\\
		\hline
		$0$ & 5.89 & 6.75 & 12.84 & 14.66 & 81 & 18.39 & 32.41 & 22.32 & 21.37 & 92 \\ 
		$0.12^2$ & 5.87 & 6.86 & 12.44 & 16.93 & 84 & 22.09 & 34.67 & 25.22 & 14.17 & 83 \\ 
		\hline 
		& \multicolumn{5}{c|}{ Scenario A:  $n= 40 , T= 200, q=50 $} & \multicolumn{5}{c}{ Scenario B:  $n= 40 , T= 200, q=50 $} \\
		$\sigma^2_{\mathfrak E}$ & RKHS & FDA & PFFR &R$_{avg}$(\%) & R$_{w}$ (\%) & RKHS & FDA & PFFR &R$_{avg}$(\%) & R$_{w}$ (\%)\\
		\hline
		$0$ & 3.98 & 4.41 & 9.68 & 10.96 & 75 & 20.75 & 76.74 & 77.40 & 269.73 & 100 \\ 
		$0.12^2$ & 4.02 & 4.47 & 9.75 & 11.32 & 77 & 23.71 & 76.95 & 77.52 & 224.49 & 100 \\ 
		\hline\hline
	\end{tabular}
	\caption{Numerical performance of RKHS, FDA and PFFR under function-on-function regression with no measurement error~($\sigma^2_{\mathfrak E}=0$) and with measurement error~($\sigma^2_{\mathfrak E}=0.12^2$). The reported nRMISE$_{avg}$ is multiplied by 100 in scale. R$_{avg}$ reflects the percent improvement of RKHS over the best performing competitor, and R$_{w}$ reflects the percentage of experiments in which RKHS achieves the lowest RMISE.} 
	\label{tab:FR_me}
\end{table}

\section{FDA and PFFR with larger number of basis functions}\label{subsec:withLargeBasis}
As discussed in Section \ref{subsec:ffr}, the performance of the penalized basis function approaches, such as FDA and PFFR, may be sensitive to the hyper-parameter $N_b$, which is the number of basis dimension used in the estimation.

In this section, we further examine the performance of FDA and PFFR with a larger $N_b.$ Specifically, we use the identical simulation setting as that in Section \ref{subsec:ffr}. We keep the implementation of RKHS the same and the only difference is that we implement FDA with $N_b=50$ and PFFR with $N_b=30$ instead of $N_b=20$. (For PFFR we can only do $N_b=30$ as the computational cost of PFFR with $N_b>30$ is practically forbidden for large-scale comparison. See later for more details.)

\begin{figure}[h]
	\begin{subfigure}{0.32\textwidth}
		\includegraphics[angle=270, width=1.2\textwidth]{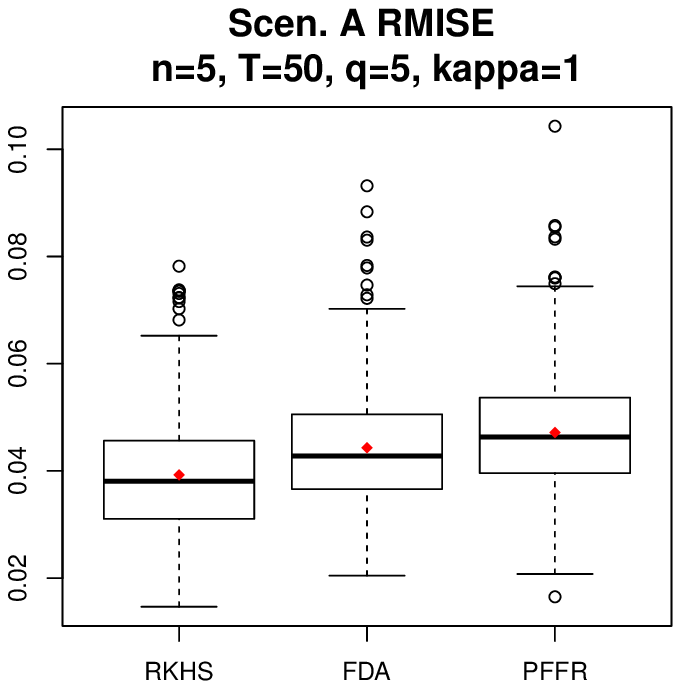}
		\vspace{-0.8cm}
	\end{subfigure}
	~
	\begin{subfigure}{0.32\textwidth}
		\includegraphics[angle=270, width=1.2\textwidth]{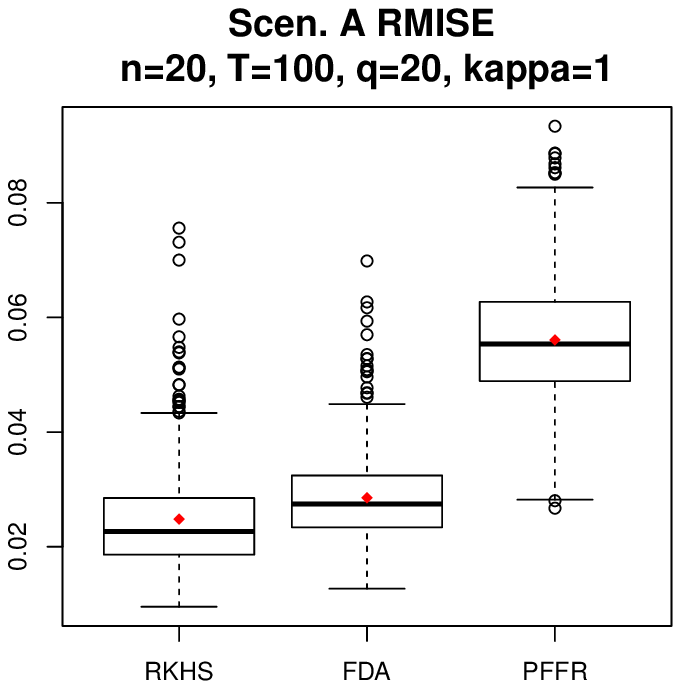}
		\vspace{-0.8cm}
	\end{subfigure}
	~
	\begin{subfigure}{0.32\textwidth}
		\includegraphics[angle=270, width=1.2\textwidth]{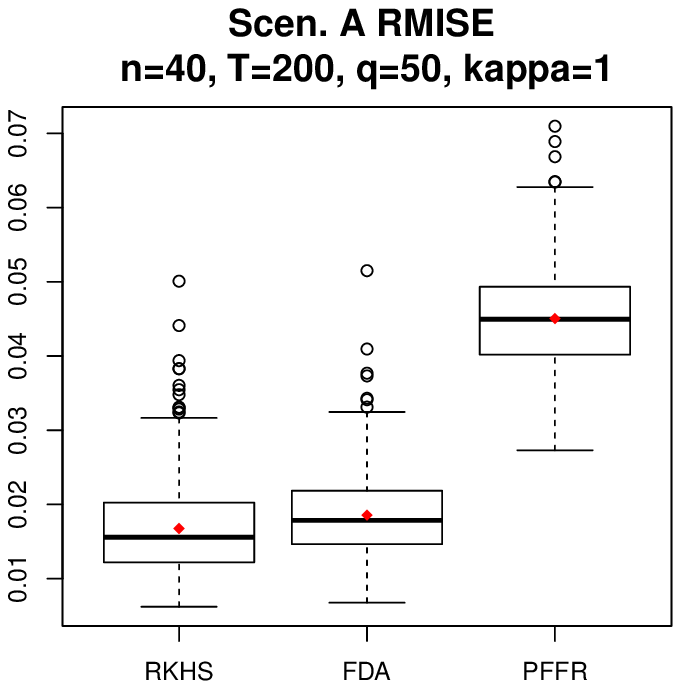}
		\vspace{-0.8cm}
	\end{subfigure}
	~
	\begin{subfigure}{0.32\textwidth}
		\includegraphics[angle=270, width=1.2\textwidth]{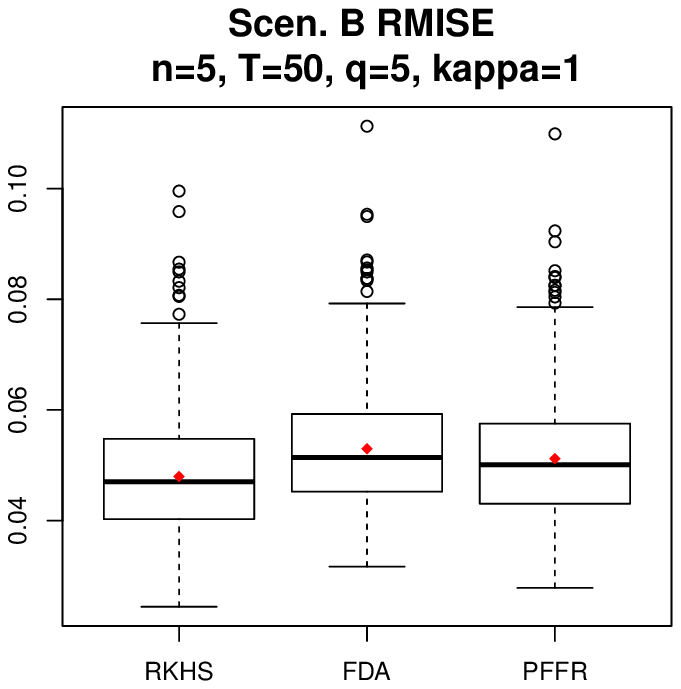}
		\vspace{-0.8cm}
	\end{subfigure}
	~
	\begin{subfigure}{0.32\textwidth}
		\includegraphics[angle=270, width=1.2\textwidth]{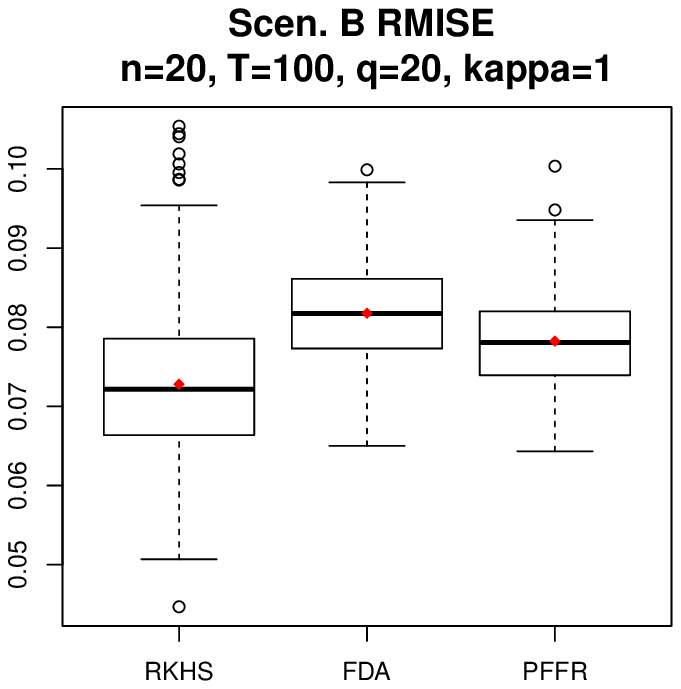}
		\vspace{-0.8cm}
	\end{subfigure}
	~
	\begin{subfigure}{0.32\textwidth}
		\includegraphics[angle=270, width=1.2\textwidth]{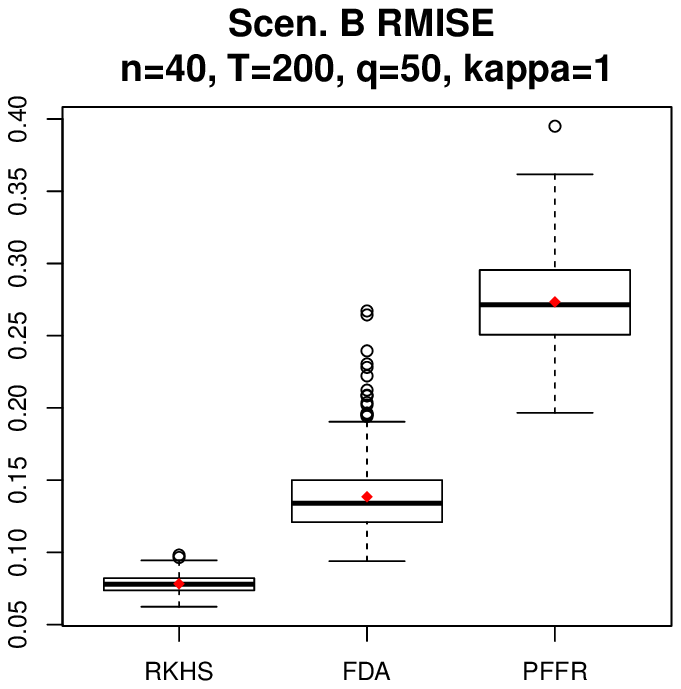}
		\vspace{-0.8cm}
	\end{subfigure}
	\caption{Boxplots of RMISE of RKHS, FDA50 and PFFR30 across 500 experiments under function-on-function regression with $\kappa=1$. Red points denote the average RMISE.}
	\label{fig:FR_LargeBasis}
\end{figure}

For each method, Table \ref{tab:FR_LargeBasis} reports the average nRMISE~(nRMISE$_{avg}$) across 500 experiments under all simulation settings. Table \ref{tab:FR_LargeBasis} is directly comparable with Table \ref{tab:FR} in Section \ref{subsec:ffr}. Cross-examining Table \ref{tab:FR_LargeBasis} and Table \ref{tab:FR}, we have the following observations. For Scenario A, where the bivariate function $A^*(r,s)$ is a simple exponential function, FDA50 and PFFR30 give essentially the same performance as $N_b=20$. The same applies to Scenario B with $q=5$ and with $q=20$ (for low SNR parameter $\kappa=0.5$). This indicates that $N_b=20$ is sufficient for simulation settings with low model complexity. On the other hand, for Scenario B with $q=20$~(for high SNR parameter $\kappa=1,2$) and $q=50$, compared to $N_b=20$, FDA50 and PFFR30 give much improved performance due to lower approximation bias, though still having a notable performance gap compared to RKHS. We further give the boxplot of RMISE for each method in Figure \ref{fig:FR_LargeBasis} under $\kappa=1$. The general pattern is the same as that exhibited in Figure \ref{fig:FR} for $N_b=20$, where RKHS provides the best overall performance.

\begin{table}[h]
	\centering
	\begin{tabular}{rrrrrr|rrrrr}
		\hline\hline
		& \multicolumn{5}{c|}{ Scenario A:  $n= 5 , T= 50, q=5$} & \multicolumn{5}{c}{ Scenario B:  $n= 5 , T= 50, q=5 $} \\
		$\kappa$ & RKHS & FDA50 & PFFR30 &R$_{avg}$(\%) & R$_{w}$ (\%) & RKHS & FDA50 & PFFR30 &R$_{avg}$(\%) & R$_{w}$ (\%)\\
		\hline
		$0.5$ & 15.98 & 17.98 & 19.77 & 12.41 & 64 & 20.86 & 22.97 & 22.30 & 6.88 & 64 \\ 
		$1$ & 9.14 & 10.32 & 10.97 & 12.74 & 68 & 10.42 & 11.51 & 11.12 & 6.66 & 71 \\ 
		$2$ & 4.99 & 6.17 & 5.69 & 14.22 & 72 & 5.21 & 5.75 & 5.55 & 6.69 & 74 \\ 
		\hline
		& \multicolumn{5}{c|}{ Scenario A:  $n= 20 , T= 100, q=20 $} & \multicolumn{5}{c}{ Scenario B:  $n= 20 , T= 100, q=20 $} \\
		$\kappa$ & RKHS & FDA50 & PFFR30 &R$_{avg}$(\%) & R$_{w}$ (\%) & RKHS & FDA50 & PFFR30 &R$_{avg}$(\%) & R$_{w}$ (\%)\\
		\hline		
		$0.5$ & 10.61 & 11.94 & 25.93 & 12.55 & 77 & 37.41 & 40.65 & 40.53 & 8.34 & 58 \\ 
		$1$ & 5.89 & 6.77 & 13.32 & 14.96 & 81 & 18.39 & 20.63 & 19.75 & 7.43 & 73 \\ 
		$2$ & 3.60 & 3.89 & 6.98 & 7.99 & 72 & 9.19 & 11.75 & 9.98 & 8.52 & 87 \\ 
		\hline 
		& \multicolumn{5}{c|}{ Scenario A:  $n= 40 , T= 200, q=50 $} & \multicolumn{5}{c}{ Scenario B:  $n= 40 , T= 200, q=50 $} \\
		$\kappa$ & RKHS & FDA50 & PFFR30 &R$_{avg}$(\%) & R$_{w}$ (\%) & RKHS & FDA50 & PFFR30 &R$_{avg}$(\%) & R$_{w}$ (\%)\\
		\hline
		$0.5$ & 7.37 & 8.54 & 21.00 & 15.86 & 82 & 39.22 & 48.27 & 72.42 & 23.40 & 94 \\ 
		$1$ & 3.98 & 4.40 & 10.70 & 10.75 & 75 & 20.75 & 35.82 & 70.61 & 76.05 & 100 \\ 
		$2$ & 2.34 & 2.34 & 5.57 & 0.11 & 47 & 12.59 & 29.77 & 70.20 & 151.92 & 100 \\ 
		\hline\hline
	\end{tabular}
	\caption{Numerical performance of RKHS, FDA and PFFR under function-on-function regression. The reported nRMISE$_{avg}$ is multiplied by 100 in scale. R$_{avg}$ reflects the percent improvement of RKHS over the best performing competitor, and R$_{w}$ reflects the percentage of experiments in which RKHS achieves the lowest RMISE.} 
	\label{tab:FR_LargeBasis}
\end{table}

\textbf{Computational Time}: Note that a larger $N_b$ can significantly increase the computational cost of the penalized basis function approaches. For example, for FDA with the same number of basis functions $N_b$ for the first and second argument of the bivariate coefficient function $A(r,s)$, the computation of FDA involves inversion of an $N_b^2\times N_b^2$-dimension matrix~(see equations (16.14) and (16.15) in \cite{Ramsay2005} for more details), which may not scale well with the number of basis $N_b$. PFFR uses restricted MLE~(REML) for its model estimation by recasting the functional regression model as a penalized additive model with mixed effects~(see \cite{ivanescu2015penalized} for more details), which seems to be slow in terms of computation and also does not scale well with $N_b$.

\begin{figure}[h]
	\hspace{-1cm}
	\centering
	\includegraphics[angle=270, width=1\textwidth]{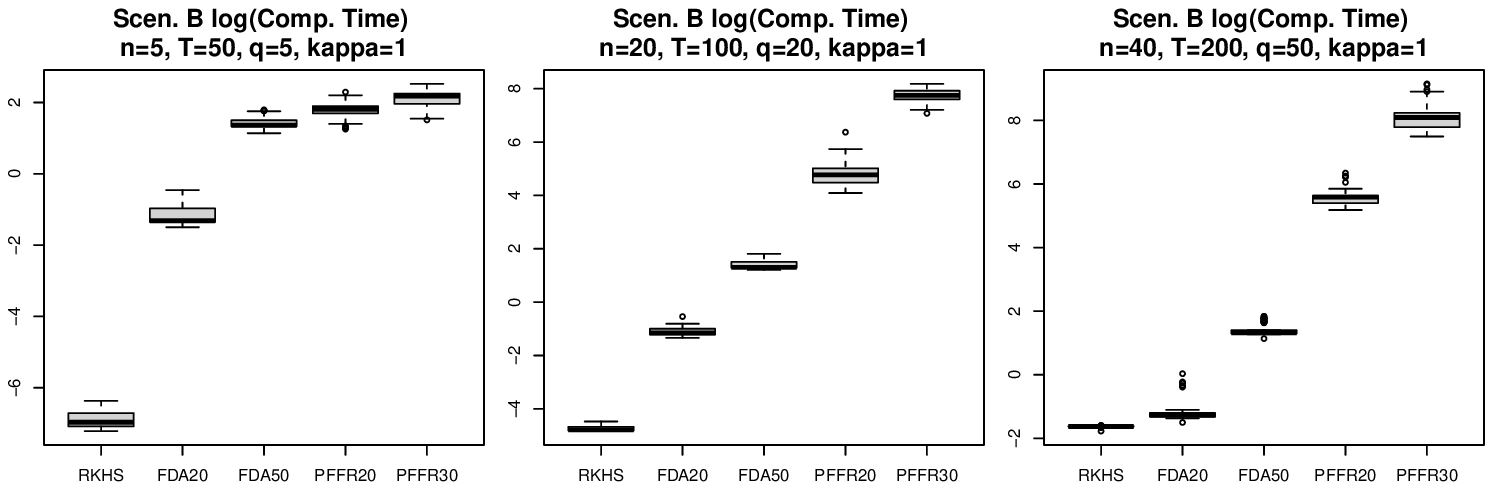}
	\vspace{-0.3cm}
	\caption{Boxplots of computational time (in log-scale) for RKHS, FDA20, FDA50 with a fixed tuning parameter $(\lambda=10^{-15})$ and for PFFR20 and PFFR30 across 500 experiments.}
	\label{fig:comp_time}
\end{figure}

For illustration, Figure \ref{fig:comp_time} gives the boxplot of computational time (in log-scale) for RKHS, FDA20, FDA50 with a fixed tuning parameter $(\lambda=10^{-15})$ and for PFFR20 and PFFR30 across 500 experiments. To conserve space, we present the case for Scenario B with $\kappa=1.$ Results under other settings are similar and thus omitted. Note that for FDA and RKHS, the computational time does not depend on $\lambda$ as we have closed-form solutions. For PFFR, it uses REML to automatically select the roughness penalty and does not require cross-validation. Thus, for fair comparison, if RKHS and FDA requires a 5-fold cross-validation to search for the best tuning parameters among 50 candidates, we need to add $\log(250)\approx 4.6$ to the log time for FDA and RKHS on Figure \ref{fig:comp_time}. 

As can be seen from Figure \ref{fig:comp_time}, FDA50 incurs much higher computational cost compared to FDA20. Similarly, PFFR30 has noticeably higher computational cost than PFFR20. For the current simulation settings, RKHS has the lowest computational cost. 

\section{Additional numerical comparison with FPCA}\label{subsec:withFPCA}
In this section, we further provide numerical comparison with the popular functional PCA based regression method~(FPCA) proposed in \cite{yao2005functional}, which estimates the functional linear model based on Karhunen--Lo\'eve expansion. FPCA is implemented via the \texttt{Matlab} package \texttt{PACE} (\texttt{FPCreg} function). We keep all tuning parameters as default values suggested in the \texttt{PACE} package\footnote{See http://www.stat.ucdavis.edu/PACE for more details.}.

\textbf{Function-on-function regression}: We first consider the function-on-function regression setting as in Section \ref{subsec:ffr} of the main text. Specifically, we use the identical simulation setting that generates Figure \ref{fig:FR} in Section \ref{subsec:ffr}. The only difference is that besides RKHS, FDA and PFFR, we further implement FPCA in the simulation.

\begin{figure}[h]
	\begin{subfigure}{0.32\textwidth}
		\includegraphics[angle=270, width=1.2\textwidth]{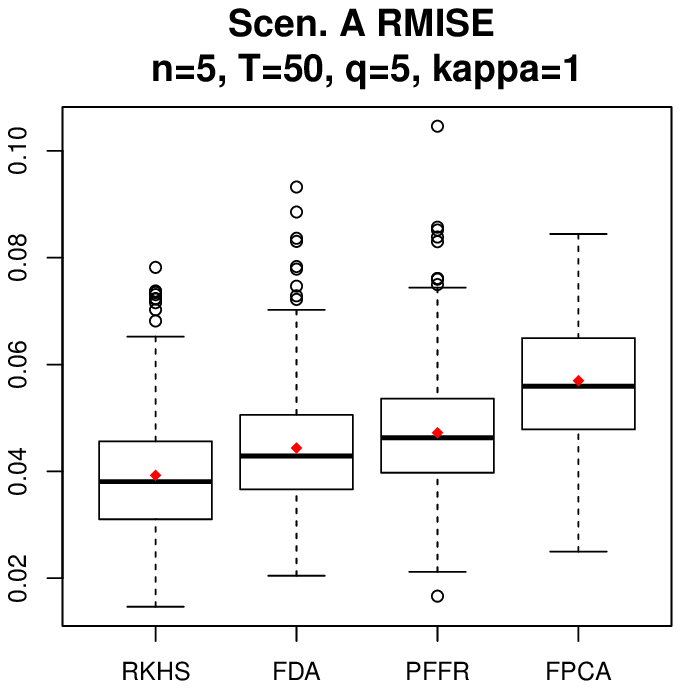}
		\vspace{-0.8cm}
	\end{subfigure}
	~
	\begin{subfigure}{0.32\textwidth}
		\includegraphics[angle=270, width=1.2\textwidth]{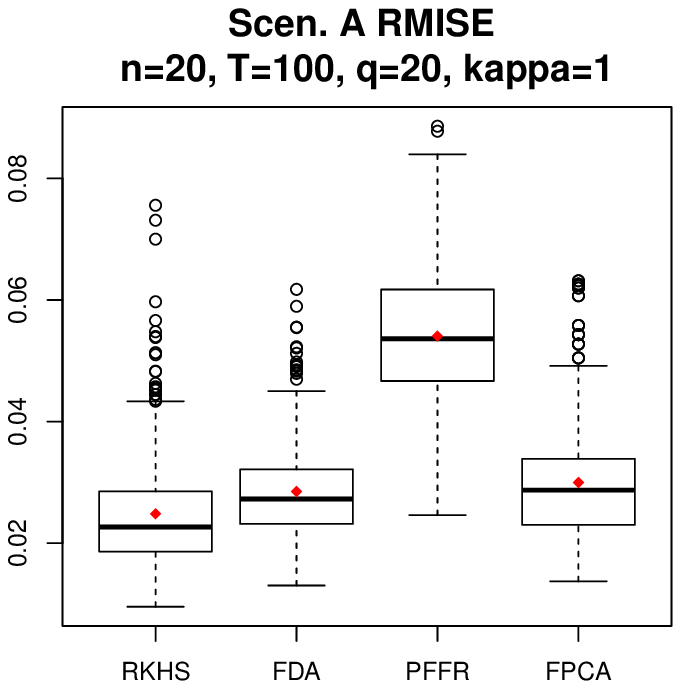}
		\vspace{-0.8cm}
	\end{subfigure}
	~
	\begin{subfigure}{0.32\textwidth}
		\includegraphics[angle=270, width=1.2\textwidth]{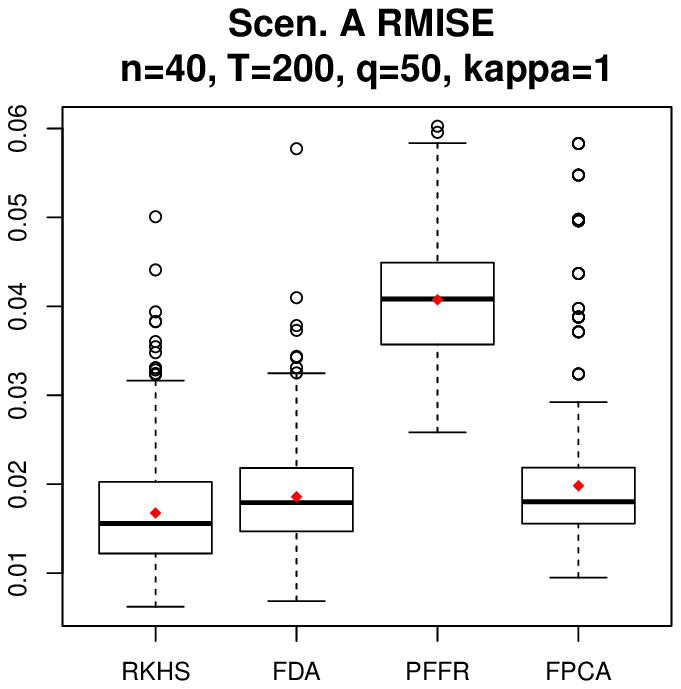}
		\vspace{-0.8cm}
	\end{subfigure}
	~
	\begin{subfigure}{0.32\textwidth}
		\includegraphics[angle=270, width=1.2\textwidth]{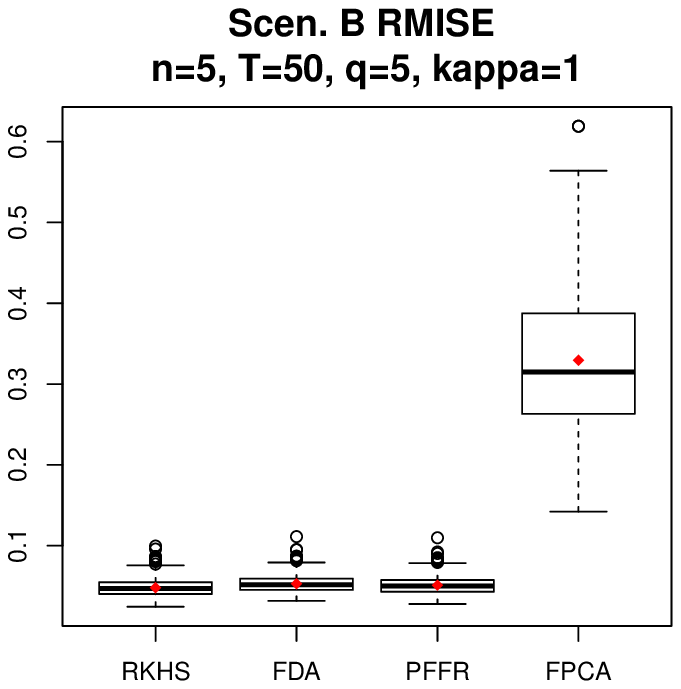}
		\vspace{-0.8cm}
	\end{subfigure}
	~
	\begin{subfigure}{0.32\textwidth}
		\includegraphics[angle=270, width=1.2\textwidth]{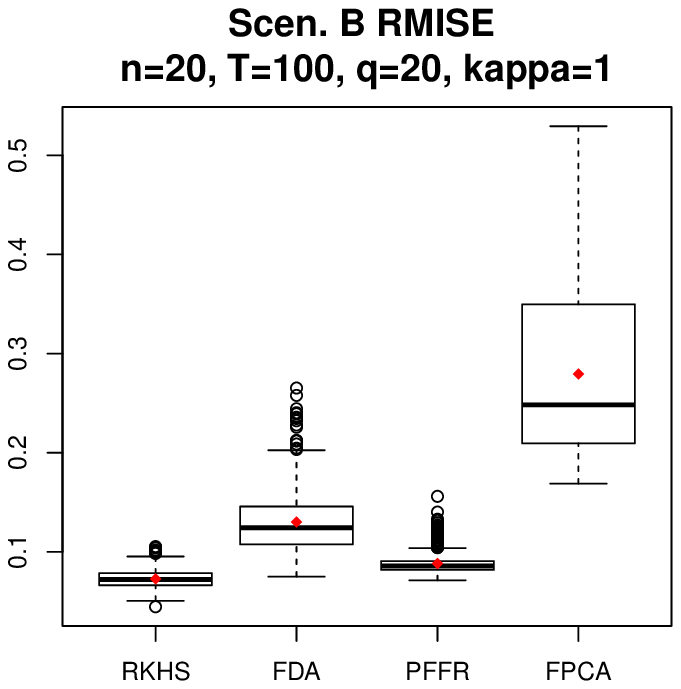}
		\vspace{-0.8cm}
	\end{subfigure}
	~
	\begin{subfigure}{0.32\textwidth}
		\includegraphics[angle=270, width=1.2\textwidth]{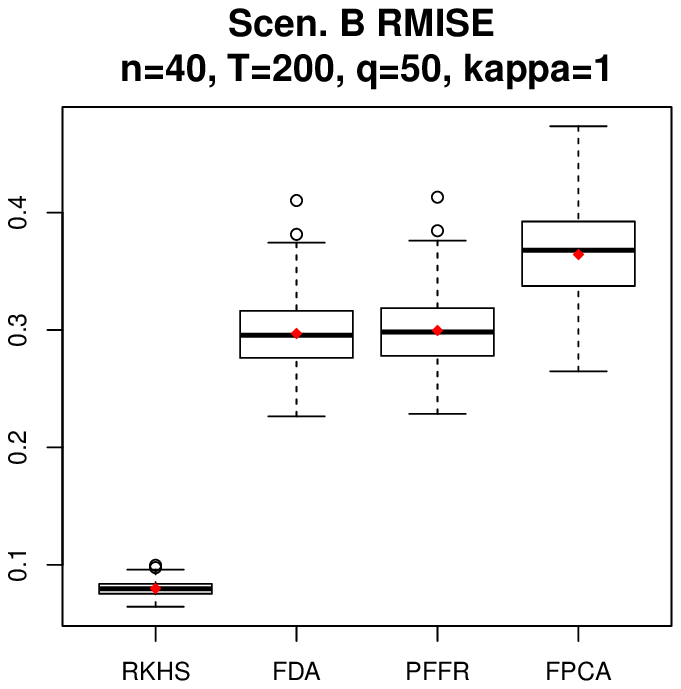}
		\vspace{-0.8cm}
	\end{subfigure}
	\caption{Boxplots of RMISE of RKHS, FDA, PFFR and FPCA across 500 experiments under function-on-function regression with $\kappa=1$. Red points denote the average RMISE.}
	\label{fig:FR_withFPCA}
\end{figure}

For each setting, we conduct 500 experiments. We again use RMISE to evaluate the excess risk of each estimator. Due to limited space, we only report the result for the case where the spectral norm $\kappa=1$. The results for $\kappa=0.5$ and $\kappa=2$ are similar and omitted. For each method (RKHS, FDA, PFFR and FPCA), Figure \ref{fig:FR_withFPCA} visualizes the boxplot of its RMISE across 500 experiments under each simulation setting. Note that the only different between Figure \ref{fig:FR} and Figure \ref{fig:FR_withFPCA} is that Figure \ref{fig:FR_withFPCA} further includes the boxplot of RMISE for FPCA.

As can be seen, for Scenario A, where the coefficient function $A(s,t)$ is the simple exponential function, FPCA provides reasonable performance. However, its performance is not satisfactory when the coefficient function $A(s,t)$ becomes more complex as in Scenario B. This phenomenon is well-known in the literature as the functional space spanned by the eigenfunctions associated with the leading principle components may not be able to well approximate the coefficient function $A(s,t)$, see \cite{cai2012minimax} and \cite{sun2018optimal} for similar observations.

\textbf{Functional regression with mixed predictors}: We then consider the functional regression with mixed covariates setting as in Section \ref{subsec:mfr} of the main text. Specifically, we use the identical simulation setting that generates Figure \ref{fig:FR_mixed} in Section \ref{subsec:mfr}. The only difference is that besides RKHS and PFFR, we further implement FPCA in the simulation.

Note that the original \texttt{FPCreg} function in the \texttt{Matlab} package \texttt{PACE} only supports function-on-function regression and does not support the general setting of functional regression with mixed covariates. Thus, we modify the original code in \texttt{FPCreg} to implement FPCA for functional regression with mixed covariates. The modification is straightforward, where we further include the observed scalar predictors $(Z_{t1},Z_{t2},Z_{t3})$ together with the estimated FPCs of the functional predictor $X_t$ to predict the FPCs of the functional response $Y_t.$

\begin{figure}[h]
	\begin{subfigure}{0.32\textwidth}
		\includegraphics[angle=270, width=1.2\textwidth]{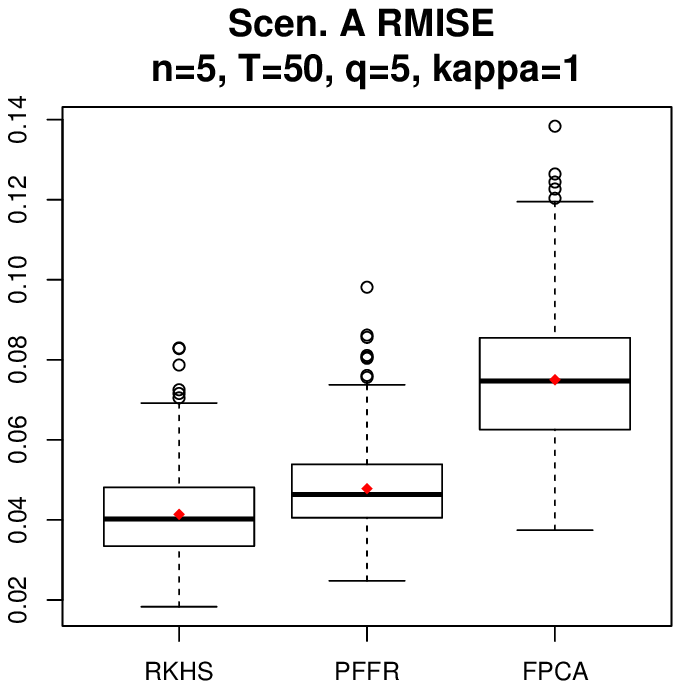}
		\vspace{-0.8cm}
	\end{subfigure}
	~
	\begin{subfigure}{0.32\textwidth}
		\includegraphics[angle=270, width=1.2\textwidth]{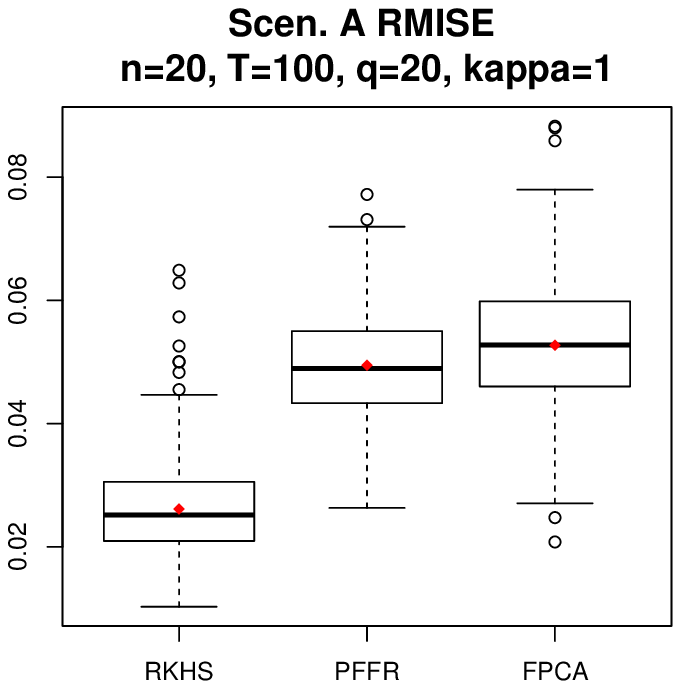}
		\vspace{-0.8cm}
	\end{subfigure}
	~
	\begin{subfigure}{0.32\textwidth}
		\includegraphics[angle=270, width=1.2\textwidth]{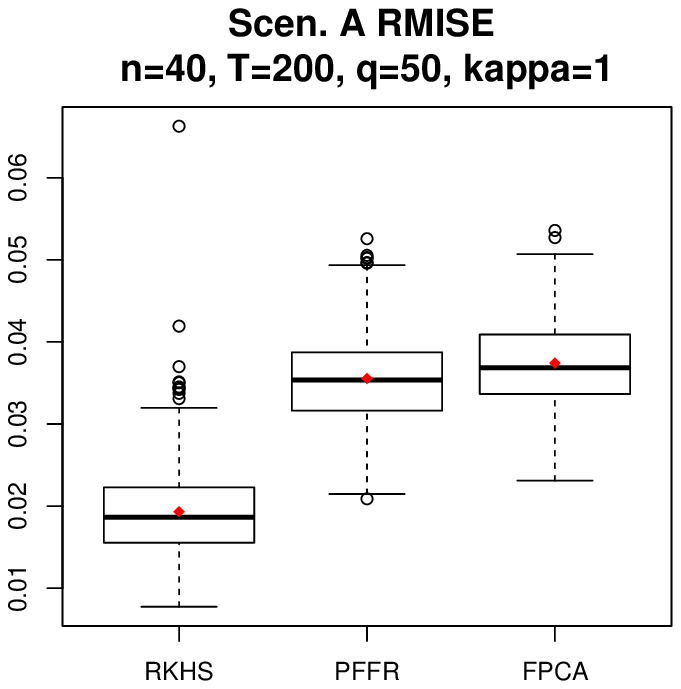}
		\vspace{-0.8cm}
	\end{subfigure}
	~
	\begin{subfigure}{0.32\textwidth}
		\includegraphics[angle=270, width=1.2\textwidth]{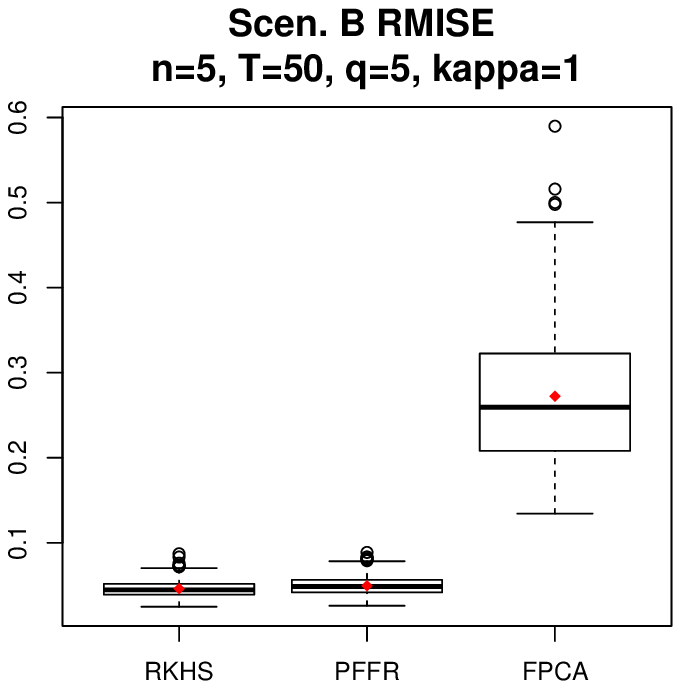}
		\vspace{-0.8cm}
	\end{subfigure}
	~
	\begin{subfigure}{0.32\textwidth}
		\includegraphics[angle=270, width=1.2\textwidth]{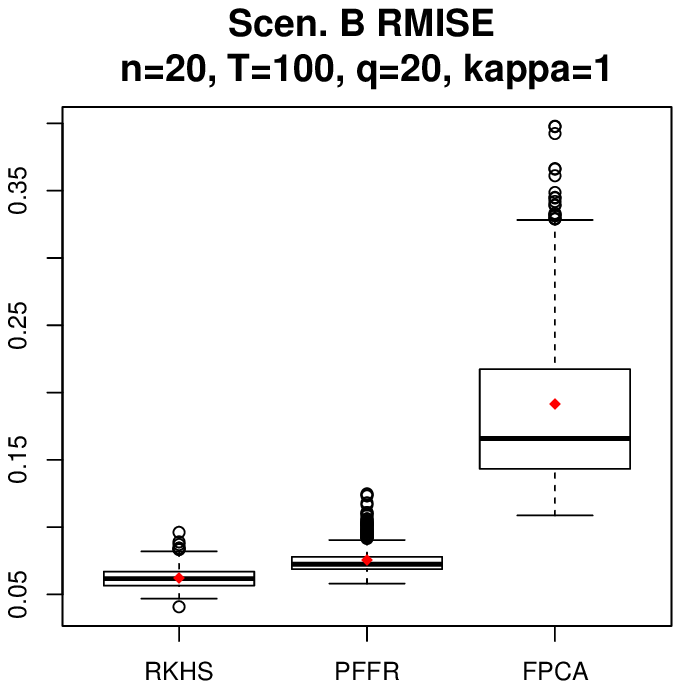}
		\vspace{-0.8cm}
	\end{subfigure}
	~
	\begin{subfigure}{0.32\textwidth}
		\includegraphics[angle=270, width=1.2\textwidth]{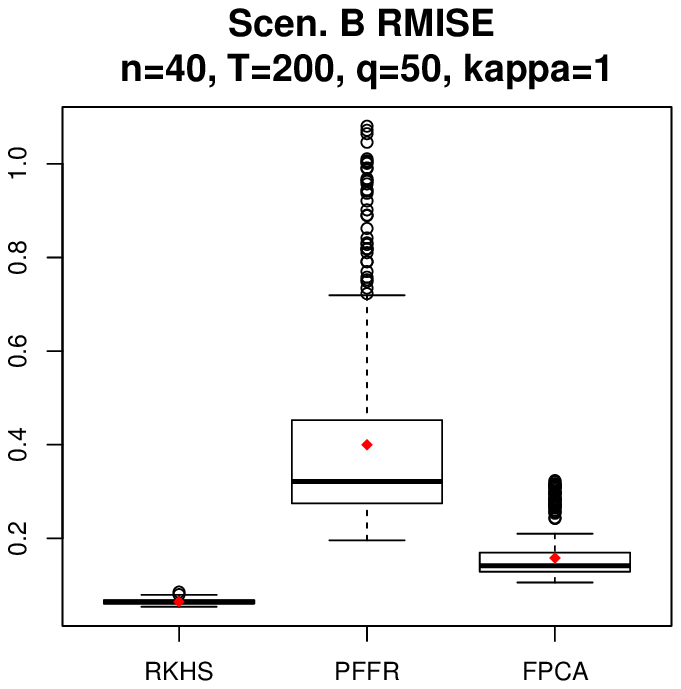}
		\vspace{-0.8cm}
	\end{subfigure}
	\caption{Boxplots of RMISE of RKHS, PFFR and FPCA across 500 experiments under functional regression with mixed predictors with $\kappa=1$. Red points denote the average RMISE.}
	\label{fig:FRmixed_withFPCA}
\end{figure}

For each setting, we conduct 500 experiments. We again use RMISE to evaluate the excess risk of each estimator. Due to limited space, we only report the result for the case where the spectral norm $\kappa=1$. The results for $\kappa=0.5$ and $\kappa=2$ are similar and omitted. For each method (RKHS, PFFR and FPCA), Figure \ref{fig:FRmixed_withFPCA} visualizes the boxplot of its RMISE across 500 experiments under each simulation setting. Note that the only different between Figure \ref{fig:FR_mixed} and Figure \ref{fig:FRmixed_withFPCA} is that Figure \ref{fig:FRmixed_withFPCA} further includes the boxplot of RMISE for FPCA.

As can be seen, the phenomenon is essentially the same as the one observed under the function-on-function regression setting. Specifically, for Scenario A, where the coefficient function $A(s,t)$ is the simple exponential function, FPCA provides reasonable performance. However, its performance is not satisfactory when the coefficient function $A(s,t)$ becomes more complex as in Scenario B. Again, this phenomenon is due to the fact that the functional space spanned by the eigenfunctions associated with the leading principle components may not be able to well approximate the coefficient function $A(s,t)$, see \cite{cai2012minimax} and \cite{sun2018optimal} for similar observations.

\section{ Illustration of phase transition}\label{subsec:phase_transition}
In this section, we provide an indirect numerical illustration for the phase transition phenomenon discovered in Theorem \ref{thm-main-mixed}. Recall Theorem \ref{thm-main-mixed} suggests that putting aside the discretization error $\zeta_n$, there is a phase transition between a nonparametric rate $\delta_T$ and a high-dimensional parametric rate $\mathfrak{s}\log(p)/T$.

Here, the nonparametric rate $\delta_T$ takes the form $\delta_T=T^{-2r/(2r+1)}$ with some $r>1/2$ and is associated with the estimation of the bivariate coefficient function $A^*$. In other words, $\delta_T=T^{-c}$ for some $c>0$ and $c$ is strictly smaller than 1. The high-dimensional parametric rate $\mathfrak{s} \log(p)/T$ is associated with the estimation of the high-dimensional sparse univariate coefficient functions $\beta^*.$ In the following, we numerically verify these two rates by considering two special cases of the proposed functional linear model with mixed predictors.

Specifically, we follow the same simulation setting as Scenario B in Section \ref{subsec:simu_setting} of the main text and simulate data from the functional linear regression model
\begin{align}\label{eq:FR_mixed1}
	Y_t(r)=\int_{[0,1]} A^*(r,s)X_t(s)\dint s + \sum_{j=1}^{p} \beta^*_j(r)Z_{tj} + \epsilon_t(r),~  r \in [0,1].
\end{align}
The functional predictor $X_t$, scalar predictor $Z_t$ and the functional error $\epsilon_t$ are generated using the same setting in Section \ref{subsec:simu_setting}. Following Scenario B, the coefficient functions $A^*$ and $\beta_1^*,\ldots,\beta_p^*$ are generated from a $q$-dimensional subspace spanned by basis functions $\{u_i(s)\}_{i=1}^{q}$. We refer to Section \ref{subsec:simu_setting} for more detailed description of Scenario B.

In the following, we fix the dimension $q=20$ and the SNR parameter $\kappa=2$. To minimize the impact of the discretization error $\zeta_n$, we set the discrete sample points $\{r_j\}_{j=1}^{n_2}$ for $Y_t$ and $\{s_i\}_{i=1}^{n_1}$ for $X_t$ to be evenly spaced dense grids on $[0,1]$ with a large number of grids $n=n_1=n_2=50$.

To evaluate the excess risk given by the RKHS estimator, we follow Section \ref{subsec:simu_setting} and use the RMISE defined in \eqref{eq:RMISE}. To match the result in Theorem \ref{thm-main-mixed}, we compute MISE, which is the squared RMISE with MISE = RMISE$^2$.

\begin{figure}[h]
	\hspace*{-1cm}
	\centering                                                           
	\includegraphics[angle=270, width=0.8\textwidth]{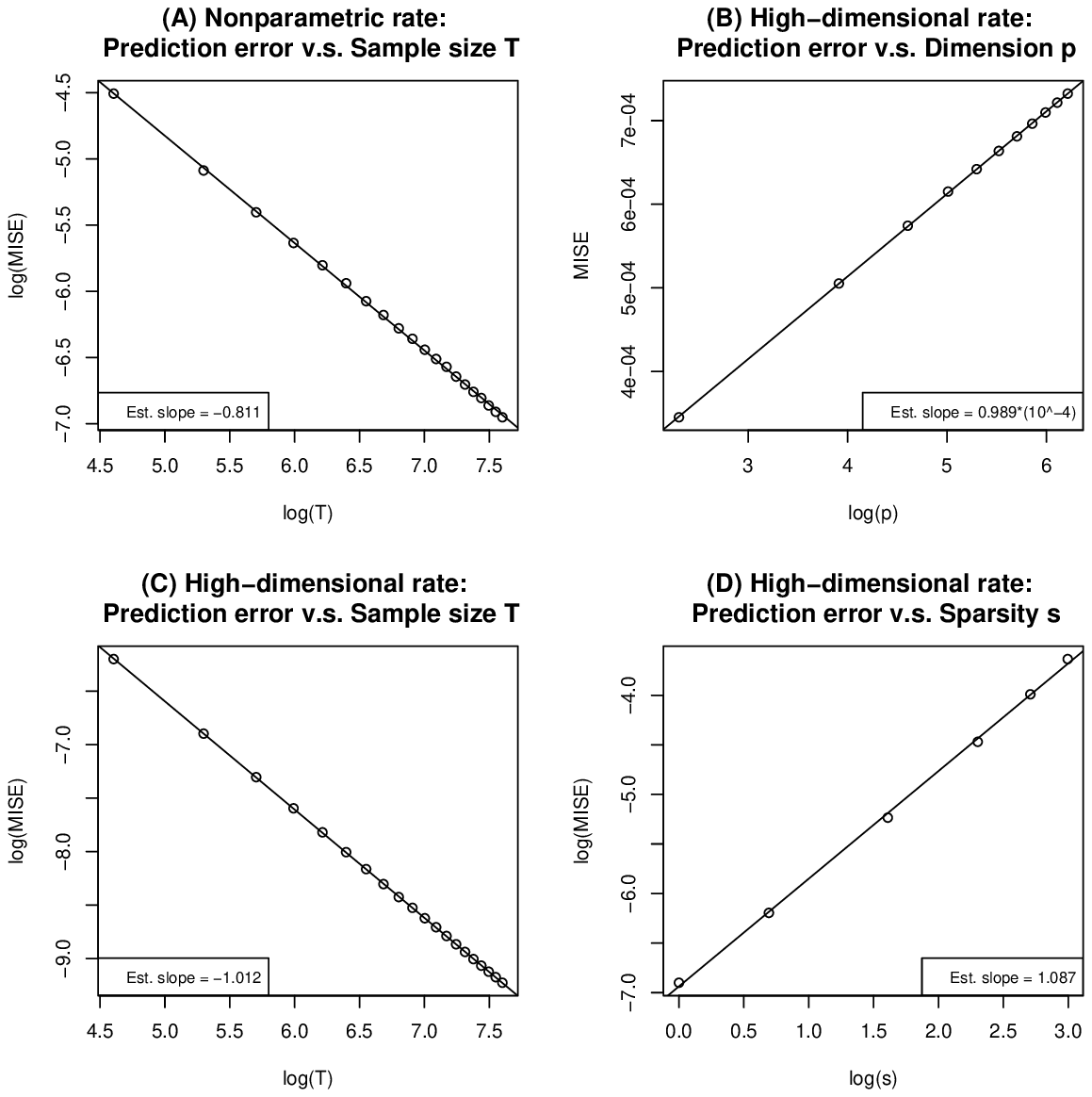}
	\vspace{-0.3cm}
	\caption{Numerical illustration of the nonparametric rate (A) and the high-dimensional parametric rate (B-D).}
	\label{fig:phase_transition}
\end{figure}

\textbf{Nonparametric rate $\delta_T$}: We first focus on the nonparametric rate $\delta_T$ and consider the special case of \eqref{eq:FR_mixed1} by setting $p=0$, i.e.\ the function-on-function regression. We vary the sample size $T=100,200,300,\ldots, 2000$ and for each $T$ we conduct 500 experiments. For each sample size $T$, we compute the average MISE across the 500 experiments and Figure \ref{fig:phase_transition}(A) gives the plot between $\log(\text{average MISE})$ and $\log(T)$, where the relationship is seen to be roughly linear. We further fit an ordinary least squared estimator (OLS) on the points, which gives an estimated slope of~$-0.811$, confirming that $\delta_T$ takes the nonparametric rate $O(T^{-c})$ with $c$ being a positive constant less than~1.

\textbf{High-dimensional parametric rate $\mathfrak{s} \log(p)/T$}: We now focus on the high-dimensional parametric rate and consider the special case of \eqref{eq:FR_mixed1} by setting $A^*=0$, i.e.\ we consider the functional linear regression with only scalar predictors $Z_{t1},Z_{t2},\cdots,Z_{tp}$.

We first fix $\mathfrak{s}=1$ and $T=200$, and vary the dimension $p=10,50,100,150,\ldots, 450,500.$ In other words, the coefficient functions $\beta_j^*\equiv 0$ for $j=2,3,\ldots, p$ and $\beta_1^*$ is generated as in Scenario~B of Section \ref{subsec:simu_setting} (see more detailed description above). For each dimension $p$, we conduct 500 experiments and compute the average MISE. Figure \ref{fig:phase_transition}(B) gives the plot between average MISE and $\log(p)$, where the relationship is seen to be roughly linear. This confirms that fixing the sparsity~$\mathfrak{s}$ and sample size $T$, the excess risk increases in the order of $O\big(\log(p)\big).$

We next fix $\mathfrak{s}=1$, $p=10$, and vary the sample size $T=100,200,300,\ldots,2000$. For each sample size $T$, we conduct 500 experiments and compute the average MISE. Figure \ref{fig:phase_transition}(C) gives the plot between $\log(\text{average MISE})$ and $\log(T)$, where the relationship is seen to be roughly linear. We further fit an OLS on the points, which gives an estimated slope of $-1.012$, confirming that fixing the sparsity $\mathfrak{s}$ and dimension $p$, the excess risk decreases in the order of roughly $O(1/T).$

Finally, we fix $T=200$, $p=100$, and vary the sparsity $\mathfrak{s}=1,2,5,10,15,20.$ In other words, the coefficient functions $\beta_j^*\equiv 0$ for $j=\mathfrak{s}+1,\mathfrak{s}+2,\ldots, p$ and $\beta_1^*=\beta_2^*=\cdots=\beta_{\mathfrak{s}}^*$ are generated as in Scenario B of Section \ref{subsec:simu_setting}. For each sparsity $\mathfrak{s}$, we conduct 500 experiments and compute the average MISE. Figure \ref{fig:phase_transition}(D) gives the plot between $\log(\text{average MISE})$ and $\log(\mathfrak{s})$, where the relationship is seen to be roughly linear. We further fit an OLS on the points, which gives an estimated slope of $1.087$, confirming that fixing the sample size $T$ and dimension $p$, the excess risk increases in the order of roughly $O(\mathfrak{s}).$

Thus, we numerically illustrate that the prediction error associated with the estimation of the high-dimensional sparse univariate coefficient functions $\beta^*$ is of order $O\big(\mathfrak s\log (p)/T\big)$. Putting these two cases together, we can conclude that there is a phase transition phenomenon between a nonparametric rate $O(T^{-c})$ and a high-dimensional parametric rate $O\big(\mathfrak s\log (p)/T\big)$.

